\documentclass[a4,10pt]{article}
\usepackage{amsmath,amscd,amssymb}

\newcommand{\nbiga}{\mathcal{A}}
\newcommand{\nbigb}{\mathcal{B}}
\newcommand{\nbigc}{\mathcal{C}}
\newcommand{\nbigd}{\mathcal{D}}
\newcommand{\nbige}{\mathcal{E}}
\newcommand{\nbigf}{\mathcal{F}}
\newcommand{\nbigg}{\mathcal{G}}
\newcommand{\nbigh}{\mathcal{H}}
\newcommand{\nbigi}{\mathcal{I}}
\newcommand{\nbigj}{\mathcal{J}}

\newcommand{\nbigl}{\mathcal{L}}
\newcommand{\nbigm}{\mathcal{M}}
\newcommand{\nbign}{\mathcal{N}}
\newcommand{\nbigo}{\mathcal{O}}
\newcommand{\nbigp}{\mathcal{P}}
\newcommand{\nbigq}{\mathcal{Q}}
\newcommand{\nbigr}{\mathcal{R}}
\newcommand{\nbigs}{\mathcal{S}}

\newcommand{\nbigu}{\mathcal{U}}
\newcommand{\nbigv}{\mathcal{V}}
\newcommand{\nbigw}{\mathcal{W}}

\newcommand{\nbigy}{\mathcal{Y}}

\newcommand{\proj}{\mathbb{P}}
\newcommand{\seisuu}{{\mathbb Z}}
\newcommand{\rnum}{{\mathbb Q}}

\newcommand{\cnum}{{\mathbb C}}
\newcommand{\real}{{\mathbb R}}
\newcommand{\hyperh}{\mathbb{H}}

\newcommand{\gbiga}{\mathfrak A}

\newcommand{\gbigc}{\mathfrak C}

\newcommand{\gbigf}{\mathfrak F}
\newcommand{\gbigg}{\mathfrak G}

\newcommand{\gminia}{\mathfrak a}
\newcommand{\gminib}{\mathfrak b}
\newcommand{\gminic}{\mathfrak c}

\newcommand{\gminio}{\mathfrak o}

\newcommand{\gminis}{\mathfrak s}

\newcommand{\vecxi}{{\boldsymbol \xi}}
\newcommand{\vecnu}{{\boldsymbol \nu}}
\newcommand{\vecmu}{{\boldsymbol \mu}}

\newcommand{\vece}{{\boldsymbol e}}

\newcommand{\vecv}{{\boldsymbol v}}
\newcommand{\vecu}{{\boldsymbol u}}
\newcommand{\vecw}{{\boldsymbol w}}

\newcommand{\veca}{{\boldsymbol a}}
\newcommand{\vecb}{{\boldsymbol b}}

\newcommand{\vecdelta}{{\boldsymbol \delta}}

\newcommand{\vecc}{{\boldsymbol c}}

\newcommand{\vech}{{\boldsymbol h}}

\newcommand{\vecm}{{\boldsymbol m}}

\newcommand{\vecN}{{\boldsymbol N}}

\newcommand{\vecx}{{\boldsymbol x}}

\newcommand{\vecn}{{\boldsymbol n}}

\newcommand{\vecV}{{\boldsymbol V}}

\newcommand{\veco}{{\boldsymbol o}}

\newcommand{\llarr}{\longleftarrow}

\newcommand{\lrarr}{\longrightarrow}

\newcommand{\pf}{{\bf Proof}\hspace{.1in}}
\newcommand{\qed}{\mbox{\rule{1.2mm}{3mm}}}

\def\Hom{\mathop{\rm Hom}\nolimits}

\def\End{\mathop{\rm End}\nolimits}

\def\Cok{\mathop{\rm Cok}\nolimits}

\def\Image{\mathop{\rm Im}\nolimits}

\def\Re{\mathop{\rm Re}\nolimits}

\def\Gr{\mathop{\rm Gr}\nolimits}

\def\GL{\mathop{\rm GL}\nolimits}

\def\Tot{\mathop{\rm Tot}\nolimits}

\def\rank{\mathop{\rm rank}\nolimits}

\def\Ker{\mathop{\rm Ker}\nolimits}

\def\Gr{\mathop{\rm Gr}\nolimits}
\def\Sym{\mathop{\rm Sym}\nolimits}

\def\ad{\mathop{\rm ad}\nolimits}
\def\Res{\mathop{\rm Res}\nolimits}

\def\ord{\mathop{\rm ord}\nolimits}

\def\ch{\mathop{ch}\nolimits}
\def\parchern{\mathop{\rm par\textrm{-}c}\nolimits}

\def\Tr{\mathop{\rm Tr}\nolimits}

\def\dvol{\mathop{\rm dvol}\nolimits}

\def\id{\mathop{\rm id}\nolimits}

\def\gcd{\mathop{\rm g.c.d.}\nolimits}

\def\ch{\mathop{\rm ch}\nolimits}

\def\Irr{\mathop{\rm Irr}\nolimits}

\newcommand{\del}{\partial}
\newcommand{\delbar}{\overline{\del}}

\newcommand{\nhom}{{\mathcal Hom}}

\newcommand{\harmonicbundle}{(E,\delbar_E,\theta,h)}

\newcommand{\barz}{\overline{z}}
\newcommand{\zbar}{\barz}
\newcommand{\zetabar}{\overline{\zeta}}

\newcommand{\Hbar}{\overline{H}}

\newcommand{\Poin}{{\bf p}}

\newcommand{\Par}{{\mathcal Par}}
\newcommand{\Sp}{{\mathcal Sp}}

\newcommand{\lefttop}[1]{{}^{#1}\!}

\newcommand{\vecnbign}{{\boldsymbol{\mathcal N}}}

\newcommand{\openclosed}[2]{]#1,#2]}

\newcommand{\Gtilde}{\widetilde{G}}

\newcommand{\vecwtilde}{\widetilde{\vecw}}

\newcommand{\wbar}{\overline{w}}

\newcommand{\kappatilde}{\widetilde{\kappa}}

\newcommand{\ftilde}{\widetilde{f}}

\newcommand{\Ftilde}{\widetilde{F}}

\newcommand{\stilde}{\widetilde{s}}

\newcommand{\nbigvtilde}{\widetilde{\nbigv}}

\newcommand{\Dtilde}{\widetilde{D}}

\def\ord{\mathop{\rm ord}\nolimits}

\def\Gal{\mathop{\rm Gal}\nolimits}

\def\Ad{\mathop{\rm Ad}\nolimits}
\def\ad{\mathop{\rm ad}\nolimits}
\def\RFM{\mathop{\rm RFM}\nolimits}
\def\Pic{\mathop{\rm Pic}\nolimits}
\def\VS{\mathop{\rm VS}\nolimits}
\def\VB{\mathop{\rm VB}\nolimits}

\def\Nahm{\mathop{\rm Nahm}\nolimits}

\newcommand{\Ubar}{\overline{U}}

\newcommand{\Ptilde}{\widetilde{P}}

\newcommand{\ubar}{\overline{u}}

\newcommand{\gtilde}{\widetilde{g}}

\newcommand{\Ybar}{\overline{Y}}
\newcommand{\wtilde}{\widetilde{w}}

\newcommand{\vecnbigi}{{\boldsymbol \nbigi}}

\newcommand{\nbigmlambda}{\nbigm^{\lambda}}

\newcommand{\taubar}{\overline{\tau}}

\newcommand{\nutilde}{\widetilde{\nu}}

\newcommand{\Xbar}{\overline{X}}

\newcommand{\nbigctilde}{\widetilde{\nbigc}}

\newcommand{\nbigmbar}{\overline{\nbigm}}

\newcommand{\Bbar}{\overline{B}}

\newcommand{\nubar}{\overline{\nu}}
\newcommand{\Sptilde}{\widetilde{\Sp}}

\newcommand{\chibar}{\overline{\chi}}
\newcommand{\Gammabar}{\overline{\Gamma}}

\newcommand{\gministilde}{\widetilde{\gminis}}

\newcommand{\vecIrr}{{\boldsymbol \Irr}}
\newcommand{\poincare}{\nbigp oin}
\newcommand{\vecnutilde}{\widetilde{\vecnu}}
\newcommand{\vecgminio}{{\boldsymbol\gminio}}
\newcommand{\poincaretilde}{\widetilde{\poincare}}
\newcommand{\nbigytilde}{\widetilde{\nbigy}}
\newcommand{\Spbar}{\overline{\Sp}}

\newtheorem{thm}{Theorem}[section]
\newtheorem{cor}[thm]{Corollary}

\newtheorem{rem}[thm]{Remark}
\newtheorem{lem}[thm]{Lemma}
\newtheorem{prop}[thm]{Proposition}

\newtheorem{assumption}[thm]{Assumption}

\begin{document}

\title{Asymptotic behaviour and the Nahm transform of \\
doubly periodic instantons
 with square integrable curvature}

\author{Takuro Mochizuki}
\date{}
\maketitle

\begin{abstract}
We study the asymptotic behaviour of 
doubly periodic instantons with square-integrable curvature.
Then, we establish the equivalence given by the Nahm transform
between the doubly periodic instantons with square integrable curvature
and the wild harmonic bundles on the dual torus.
\end{abstract}

\section{Introduction}
Let $T:=\cnum/L$,
where $L$ is a lattice of $\cnum$.
The product $T\times\cnum$
is equipped with the standard metric
$dz\,d\zbar+dw\,d\wbar$,
where $(z,w)$ is the standard local coordinate
of $T\times\cnum$.
In this paper, we shall study 
any $L^2$-instanton $(E,\nabla,h)$ on $T\times\cnum$,
i.e., the curvature $F(\nabla)$ satisfies the equation
$\Lambda F(h)=0$,
and it is $L^2$.

There is another natural decay condition around $\infty$.
That is the quadratic curvature decay,
i.e., $\bigl|F(\nabla)\bigr|=O(|w|^{-2})$
with respect to $h$ and the Euclidean metric
$dz\,d\zbar+dw\,d\wbar$.
M. Jardim \cite{Jardim1} studied the Nahm transform of 
some kind of harmonic bundles 
with tame singularity on the dual torus $T^{\lor}$
to produce instantons on $T\times\cnum$
satisfying the quadratic curvature decay.
O. Biquard and Jardim \cite{Biquard-Jardim}
studied the asymptotic behaviour of such instantons
with rank $2$.
Based on the results,
the inverse transform was constructed in \cite{Jardim2},
i.e., the Nahm transform of such instantons on $T\times\cnum$
to produce some type of harmonic bundles
with tame singularity on $T^{\lor}$.
See also \cite{Jardim3} and \cite{Jardim4}.

It is our purpose in this paper to generalize their results.
Namely, we will study the asymptotic behaviour of $L^2$-instantons,
and establish the equivalence between
the $L^2$-instantons on $T\times\cnum$
and harmonic bundles with wild singularity on $T^{\lor}$.
We shall also introduce algebraic counterparts of the transforms.
They are useful to describe the induced transformations
of the singular data, for example.

\subsection{Asymptotic behaviour of $L^2$-instanton}

\subsubsection{The dimensional reduction due to Hitchin}

Briefly speaking,
it is our goal in the study of the asymptotic behaviour of
$L^2$-instantons,
to show that they behave like wild harmonic bundles around $\infty$.
For the explanation,
let us recall the dimensional reduction due to N. Hitchin.
Let $U$ be any open subset of $\cnum$.
Let $(V,\delbar_V)$ be a holomorphic vector bundle on $U$
with a Higgs field $\theta$.
Let $h$ be a hermitian metric of $V$.
We have the Chern connection 
$\nabla_{V,h}=\delbar_V+\del_{V,h}$.
We have the adjoint $\theta^{\dagger}$ of $\theta$
with respect to $h$.
The tuple $(V,\delbar_V,h,\theta)$
is called a harmonic bundle,
if the Hitchin equation
$F(\nabla_{V,h})+[\theta,\theta^{\dagger}]=0$
is satisfied.

Let $p:T\times U\lrarr U$ be the projection.
We have the expression $\theta=f\,dw$
and $\theta^{\dagger}=f^{\dagger}\,d\wbar$,
where $f$ is a holomorphic endomorphism of $V$,
and $f^{\dagger}$ is the adjoint of $f$.
We set
$(E,h_E):=p^{\ast}(V,h)$.
Let $\nabla_E$ be the unitary connection given by
$\nabla_E=p^{\ast}(\nabla_{V,h})+f\,d\zbar-f^{\dagger}dz$.
Then,
$(E,h_E,\nabla_E)$ is an instanton,
if and only if
$(V,\delbar_V,h,\theta)$ is a harmonic bundle.
Indeed, 
we have the equivalence between
harmonic bundles on $U$,
and $T$-equivariant instantons on $T\times U$,
which is due to Hitchin.

\subsubsection{Examples and remarks}

We set $U:=\bigl\{w\in\cnum\,\bigr|\,|w|>R\bigr\}$.
We shall make $R$ larger without mention.
Let $\gminia$ be any holomorphic function on $U$.
We have the harmonic bundle $\nbigl(\gminia)$
obtained as the tuple of
the trivial line bundle $\nbigo_U\,e$,
the trivial metric $h(e,e)=1$
and the Higgs field $d\gminia$.
By using the dimensional reduction,
we have the associated instanton.
Its curvature is
$\del_w^2\gminia\,dw\,d\zbar
 +\overline{\del_w^2\gminia}dz\,d\wbar$,
which is $L^2$ if and only if
it has quadratic decay.

We can obtain more examples
by considering ramifications along $\infty$.
We set
$U_{\eta}:=\bigl\{\eta\in\cnum\,\big|\,
 |\eta|>R^{1/2}
 \bigr\}$.
We consider a harmonic bundle
$\nbigl(\gminia)$,
where $\gminia$ is a holomorphic function on $U_{\eta}$.
Let $\varphi:U_{\eta}\lrarr U$ be given by
$\varphi(\eta)=\eta^2$.
We obtain a harmonic bundle
$\varphi_{\ast}\nbigl(\gminia)$
of rank $2$
on $U$ obtained as the push-forward.
It is easy to check that 
the associated instanton is $L^2$
if and only if 
$\eta^{-2}\gminia(\eta)$ is holomorphic at $\infty$.
In that case,
the curvature $F$ satisfies the decay condition $O(|w|^{-3/2})$.
If $\gminia=\alpha\eta$ for $\alpha\neq 0$,
we have 
$0<C_1<|F|\,|w|^{3/2}<C_2$
for some constants $C_i$.

More generally, for any positive integer $p$,
we set $U^{\langle p\rangle}:=
 \bigl\{
 w_p\in\cnum\,\big|\,
 |w_p|>R^{1/p}
 \bigr\}$.
For a covering 
$\varphi_p:U^{\langle p\rangle}
 \lrarr U$ given by $\varphi_p(w_p)=w_p^p$
and for a holomorphic function $\gminia$ on $U^{\langle p\rangle}$,
we obtain a harmonic bundle
$\varphi_{p\ast}\nbigl(\gminia)$ 
of rank $p$ on $U$.
The associated instanton is $L^2$
if and only if 
$\varphi_p^{\ast}(w)^{-1}\gminia$ is holomorphic at $\infty$.
If $\gminia$ is a polynomial of $w_p$,
then it is described as a condition on the slope
$\deg_{w_p}(\gminia)/p\leq 1$.
In that case,
the curvature $F$ satisfies $O(|w|^{-1-1/p})$.
It is easy to construct an example satisfying
$0<C_1<|F|\,|w|^{1+1/p}<C_2$
for some $C_i>0$.

\vspace{.1in}

Let $\harmonicbundle$ be
a wild harmonic bundle on $U$,
i.e.,
the (possibly multi-valued) eigenvalues of $\theta$
are meromorphic at $\infty$.
As the above examples suggest,
the $L^2$-condition and
the quadratic decay condition
can be described in terms of 
the eigenvalues of the Higgs field.
If we take an appropriate covering
$\varphi_p:
 U^{\langle p\rangle}\lrarr U$,
we have a holomorphic decomposition
\begin{equation}
\label{eq;13.1.30.1}
 \varphi_p^{\ast}(E,\theta)
=\bigoplus_{\gminia\in w_p\cnum[w_p]}
 (E_{\gminia},\theta_{\gminia})
\end{equation}
such that
$\theta_{\gminia}-d\gminia$ are tame,
i.e.,
for the expression
$\theta_{\gminia}-d\gminia=f_{\gminia}\,dw_p/w_p$,
the eigenvalues of $f_{\gminia}$ are bounded.
We set $\Irr(\theta):=
 \bigl\{
 \gminia\,\big|\,
 E_{\gminia}\neq 0
\bigr\}$.
Then, 
by using the results in the asymptotic behaviour of
wild harmonic bundles \cite{Mochizuki-wild},
it is not difficult to show that 
the instanton associated to $\harmonicbundle$ is $L^2$,
if and only if $\deg_{w_p}(\gminia)/p\leq 1$
for any $\gminia\in\Irr(\theta)$.
It satisfies the quadratic decay condition,
if and only if the harmonic bundle is unramified, in addition,
i.e., it has a decomposition as in (\ref{eq;13.1.30.1})
on $U$, if $R$ is sufficiently large.

The condition can also be described in terms
of the spectral variety of $\theta$.
We have the expression $\theta=f\,dw$.
Let $\Sp(f)\subset\cnum_{\zeta}\times U$ 
denote the support of the cokernel of
$\nbigo_{\cnum_{\zeta}\times U}
\stackrel{\zeta-f}{\lrarr}
 \nbigo_{\cnum_{\zeta}\times U}$.
It induces a subvariety
$\Phi(\Sp(f))$,
where $\Phi:\cnum_{\zeta}\times U\lrarr T^{\lor}\times U$
denotes the projection.
Then, the instanton associated to $\harmonicbundle$
is $L^2$,
if and only if
the closure of $\Phi\bigl(\Sp(f)\bigr)$
in $T\times \Ubar$
is a complex subvariety,
where $\Ubar=U\cup\{\infty\}$.

\subsubsection{Brief description of the asymptotic behaviour of
$L^2$-instantons}
\label{subsection;13.10.16.10}

Let $(E,\nabla,h)$ be an $L^2$-instanton on $T\times U$.
Let $(E,\delbar_E)$ denote 
the underlying holomorphic vector bundle
on $T\times U$.
By using a theorem of Uhlenbeck,
we obtain $F(\nabla)=o(1)$.
It implies that the restrictions $(E,\delbar_E)_{|T\times\{w\}}$
are semistable of degree $0$
if $|w|$ is sufficiently large.
Hence, the relative Fourier-Mukai transform
of $(E,\delbar_E)$ gives 
an $\nbigo_{T^{\lor}\times U}$-module
whose support $\Sp(E)$
is relatively $0$-dimensional over $U$.
The first important issue in the study
is the following.
\begin{thm}[Theorem \ref{thm;12.6.22.101}]
$\Sp(E)$ is extended to 
a complex analytic subvariety $\overline{\Sp(E)}$
in $T\times \Ubar$.
\end{thm}
We use an effective control of
the spectrum of semistable bundles of degree $0$,
in terms of the eigenvalues of the monodromy transformations
of unitary connections with the small curvature
(Corollary \ref{cor;12.6.24.2}).
If we fix an embedding
$\Sym^{\rank E}(T^{\lor})\subset\proj^N$,
the spectrum induces
a holomorphic map from $U$ to $\proj^N$,
which we regard as a harmonic map.
We will observe that
the energy of the harmonic map
is dominated by the curvature.
Then, we obtain the desired extendability
of the spectral curve
from the regularity theorem of J. Sacks and K. Uhlenbeck
\cite{Sacks-Uhlenbeck}
for harmonic maps with finite energy.

Let $\pi:T\times U\lrarr U$ denote the projection.
We fix a lift of $\overline{\Sp(E)}$
to $\Sptilde(E)\subset \cnum\times \Ubar$.
Then,
we obtain a holomorphic vector bundle $V$
on $U$ with an endomorphism $g$,
with a $C^{\infty}$-isomorphism
$\pi^{\ast}V\simeq E$
such that
(i) $\pi^{\ast}\delbar_V+g\,d\zbar
=\delbar_E$,
(ii) $\Sp(g)=\Sptilde(E)$.
By the identification $E=\pi^{\ast}V$,
we obtain a $T$-action on $E$.

We set $\Sptilde_{\infty}(E):=
 (\cnum\times\{\infty\})
\cap
 \Sptilde(E)$.
We have the decomposition
$(V,g)=\bigoplus_{\alpha\in\Sptilde_{\infty}(E)}
 (V_{\alpha},g_{\alpha})$
such that
the eigenvalues of $g_{\alpha}(w)$
goes to $\alpha$ when $w\to \infty$.
We have the corresponding decomposition
$E=\bigoplus_{\alpha\in\Sptilde_{\infty}(E)}
 E_{\alpha}$.

The hermitian metric $h$ of $E$
is decomposed into the sum
$h=\sum h_{\alpha,\beta}$,
where $h_{\alpha,\beta}$
are the sesqui-linear pairings of
$E_{\alpha}$ and $E_{\beta}$.
By using the Fourier expansion,
we decompose $h_{\alpha,\beta}$ 
into the $T$-invariant part
and the complement.
Let $h^{\circ}$ denote the $T$-invariant part of
$\sum h_{\alpha,\alpha}$.
We shall prove that the complement 
$h^{\bot}:=h-h^{\circ}$
and its derivatives
have exponential decay.
\begin{thm}[Theorem \ref{thm;12.8.7.1}]
\label{thm;13.10.15.50}
For any polynomial $P(t_1,t_2,t_3,t_4)$
of non-commutative variables,
there exists $C>0$
such that 
\begin{equation}
 \label{eq;13.10.18.1}
 P(\nabla_z,\nabla_{\zbar},\nabla_w,\nabla_{\wbar})
 h^{\bot}=O\bigl(\exp(-C|w|)\bigr). 
\end{equation}
\end{thm}

We have a hermitian metric $h_V$ of $V$
induced by $h^{\circ}$.
As a result, $(V,\delbar_V,h_V,g\,dw)$ satisfies
the Hitchin equation up to an exponentially small term
(Proposition \ref{prop;13.1.11.21}).
Such a tuple $(V,\delbar_V,h_V,g\,dw)$
can be studied as in the case of
wild harmonic bundles \cite{Mochizuki-wild}
with minor modifications.
(See \S\ref{subsection;13.1.11.10}.)
Thus, we will arrive at a satisfactory stage of 
understanding of 
the asymptotic behaviour of $L^2$-instantons.
We state some significant consequences.
\begin{thm}[Theorem \ref{thm;13.1.30.2}]
\label{thm;13.10.18.2}
There exists $\rho>0$
such that the following holds:
\begin{equation}
 \label{eq;13.10.15.51}
 F(h)=O\Bigl(
 \frac{dz\,d\zbar}{|w|^{2}(-\log|w|)^2}
 \Bigr)
+O\Bigl(
 \frac{dw\,d\wbar}{|w|^{2}(-\log|w|)^2}
 \Bigr)
+O\Bigl(
 \frac{dw\,d\zbar}{|w|^{1+\rho}}
 \Bigr)
+O\Bigl(
 \frac{dz\,d\wbar}{|w|^{1+\rho}}
 \Bigr)
\end{equation}
In particular,
$F(\nabla)=O(|w|^{-1-\rho})$ for some $\rho>0$
with respect to $h$ and the Euclidean metric
$dw\,d\wbar+dz\,d\zbar$.
\end{thm}

The estimate (\ref{eq;13.10.15.51})
implies that $(E,\delbar_E,h)$ is acceptable,
i.e., $F(\nabla)$ is bounded with respect to
$h$ and the Poincar\'e like metric
$|w|^{-2}(\log|w|^2)^{-2}dw\,d\wbar
+dz\,d\zbar$  on $T\times U$.
By applying a general result
in \cite{Mochizuki-wild},
we obtain the following prolongation result.
\begin{cor}[Corollary \ref{cor;13.10.15.40}]
The holomorphic bundle
$(E,\delbar_E)$ is naturally extended 
to a filtered bundle $\nbigp_{\ast}E$
on $(T\times \Ubar,T\times\{\infty\})$.
\end{cor}

Here, the filtered bundle $\nbigp_{\ast}E$
on $(T\times\Ubar,T\times\{\infty\})$
is an increasing sequence 
$(\nbigp_aE\,|\,a\in\real)$
of locally free $\nbigo_{T\times\Ubar}$-modules
such that
(i) $\nbigp_a(E)_{|T\times U}=E$,
(ii) $\nbigp_a(E)/\nbigp_{<a}(E)$
 are locally free $\nbigo_{T\times\{\infty\}}$-modules,
 where $\nbigp_{<a}E=\sum_{b<a}\nbigp_bE$,
(iii) $\nbigp_a(E)=\nbigp_{a+\epsilon}(E)$ for some $\epsilon>0$,
(iv) $\nbigp_{a+1}(E)=
 \nbigp_{a}(E)\otimes\nbigo_{T\times \Ubar}
 \bigl(T\times\{\infty\}\bigr)$.
The sheaves $\nbigp_aE$ are obtained 
as the space of the holomorphic sections of $E$ 
whose norms with respect to $h$
have growth order $O(|w|^{a+\epsilon})$ for any $\epsilon>0$.

The filtered bundle is useful in the study of
the instanton.
For example,
it turns out that
$\frac{1}{8\pi^2}
 \int_{T\times\cnum}
 \Tr(F(h)^2)$
is equal to 
$\int_{T\times\proj^1}
 c_2(\nbigp_aE)$
for any $a\in\real$,
where $c_2$ denotes the second Chern class.
(Proposition \ref{prop;13.10.15.110}.)
In particular,
the number
 $\frac{1}{8\pi^2}
 \int_{T\times\cnum}
 \Tr(F(h)^2)$
is an integer.
(See \cite{Wehrheim}
for this kind of integrality
in a more general situation,
which was informed by the referee.)

We can also use this filtered bundle
to characterize the metric,
i.e., a uniqueness part of 
the so-called Kobayashi-Hitchin correspondence.
The stability condition for this type of
filtered bundles is defined,
as in \cite{Biquard-Jardim}.
(See \S\ref{subsection;13.10.15.200}.
Note that it is not a standard (slope-)stability condition
for filtered bundles.)

\begin{prop}[Proposition \ref{prop;13.1.24.3}, 
 Proposition \ref{prop;13.1.29.20}]
The associated filtered bundle
$\nbigp_{\ast}E$ is poly-stable of degree $0$.
The metric $h$ is uniquely determined 
as a Hermitian-Einstein metric of $(E,\delbar_E)$
adapted to $\nbigp_{\ast}E$,
up to obvious ambiguity.
\end{prop}

We observe that,
we need only a weaker assumption on the curvature decay,
if we assume the prolongation of the spectral curve.
\begin{thm}[Theorem \ref{thm;13.2.26.3}]
\label{thm;13.10.15.52}
Suppose that $F(\nabla)\to 0$ when $|w|\to\infty$,
and that the spectral curve $\Sp(E)$ is extended
to a complex subvariety of $T\times \Ubar$.
Then, $(E,\nabla,h)$ is an $L^2$-instanton.
In particular,
we may apply the results on
the asymptotic behaviour of $L^2$-instantons.
\end{thm}
More precisely, we directly prove
the claims of Theorem \ref{thm;13.10.15.50}
and Theorem \ref{thm;13.10.18.2}
under the assumption,
without considering $L^2$-condition.

\subsubsection{Some remarks}

In \cite{Biquard-Jardim},
Jardim and Biquard showed that 
an instanton of rank $2$ with quadratic decay 
is an exponentially small perturbation of 
a tuple $(V,\delbar_V,gdw,h_V)$
which satisfies the Hitchin equation up to an exponentially small term.
Our result could be regarded as a generalization of theirs.
But, the methods are rather different.
To obtain a decomposition into the $T$-invariant part
and the complement,
they started with the construction of a global frame
satisfying some nice property, which is an analogue of
the Coulomb gauge of Uhlenbeck.
Their method seems to require a stronger decay condition than $L^2$,
for example the quadratic decay condition as they imposed.
We use a more natural decomposition
induced by a standard method of the Fourier-Mukai transform
in complex geometry,
which allows us to consider $L^2$-instantons,
once we deal with the issue of the prolongation of the spectral curve. 
(See also \cite{Charbonneau} for the $L^2$-property
and the quadratic decay property of doubly periodic instantons.)

\vspace{.1in}
As mentioned above,
we shall establish that
an $L^2$-instanton is an exponentially small perturbation of
$(V,\delbar_V,h_V,\theta_V)$
which satisfies the Hitchin equation 
up to exponentially small term.
Interestingly to the author,
we can obtain a more refined result.
Namely,
we can naturally construct a harmonic metric $h'_V$
of $(V,\delbar_{V},g\,dw)$
defined on a neighbourhood of $\infty$,
from the $L^2$-instanton.
It is an analogue of the reductions from
wild harmonic bundles to tame harmonic bundles
studied in {\rm\cite{Mochizuki-wild}}.
We consider a kind of meromorphic prolongation of
the holomorphic vector bundle
on the twistor space associated to $T\times\cnum$,
then we encounter 
a kind of infinite dimensional Stokes phenomena.
By taking the graduation with respect to the Stokes structure,
we obtain a wild harmonic bundle.
Relatedly, in this paper, we consider only 
the product holomorphic structure of
$T\times\cnum$.
From the viewpoint of twistor theory,
the holomorphic vector bundle
with respect to the other holomorphic structures
should also be studied.
The prolongation of the twistor family of the holomorphic structure
is related with the issue in the previous paragraph.
The author hopes to return 
to this deeper aspect of the study elsewhere.

\vspace{.1in}
Although we do not use it explicitly,
we prefer to regard an instanton on
$T\times U$
as an infinite dimensional harmonic bundle
on $U$,
which is suggested by the reduction of Hitchin.
This heuristic is useful in our study of the asymptotic
behaviour of $L^2$-instantons.
From the view point,
Theorem {\rm \ref{thm;13.10.15.50}}
and Theorem {\rm\ref{thm;13.10.15.52}}
can be naturally regarded
as a variant of Simpson's main estimate
{\rm\cite{Simpson90}}.
(See also
{\rm\cite{mochi2}}
and {\rm\cite{Mochizuki-wild}}.)

\subsection{Nahm transform for 
wild harmonic bundles and $L^2$-instantons}

\subsubsection{Nahm transforms and algebraic Nahm transforms}

As an application of the study of the asymptotic behaviour,
we shall establish the equivalence between
$L^2$-instantons on $T\times\cnum$
and wild harmonic bundles on $T^{\lor}$
given by the Nahm transforms,
which is a differential geometric variant of 
the Fourier-Mukai transform.
(See \cite{Bartotti-Bruzzo-Ruiperez}  and \cite{Jardim4}
for the long history of a various versions of the Nahm transforms.)

Once we understand the asymptotic behaviour
of an $L^2$-instanton $(E,\nabla,h)$,
we can prove the desired property of
the associated cohomology groups
and harmonic sections.
Hence, the standard $L^2$-method allows us 
to construct the Nahm transform
$\Nahm(E,\nabla,h)$
which is a wild harmonic bundle
on $\bigl(T^{\lor},\Sp_{\infty}(E)\bigr)$.
(See \S\ref{subsection;13.2.3.10}.)
Conversely,
we may construct the Nahm transform of
any wild harmonic bundle
$(\nbige,\delbar_{\nbige},\theta,h_{\nbige})$
on $(T^{\lor},D)$
to $L^2$-instantons 
$\Nahm(\nbige,\delbar_{\nbige},\theta,h_{\nbige})$
on $T\times\cnum$,
by using the result on wild harmonic bundles on curves
(\cite{Mochizuki-wild}, \cite{sabbah2} and 
 \cite{Zucker}),
although we need some estimates to establish
the $L^2$-property (see \S\ref{subsection;13.1.27.10}).

\vspace{.1in}

To study their more detailed properties,
we introduce the algebraic Nahm transforms
for filtered Higgs bundles on $(T^{\lor},D)$
and filtered bundles on $(T\times\cnum,T\times\{\infty\})$,
which do not necessarily come from wild harmonic bundles
or $L^2$-instantons.
The constructions are 
based on the Higgs interpretation of the Nahm transforms.
It could be regarded as a filtered version of
the Fourier transform for Higgs bundles studied in
\cite{Bonsdorff1},
although we restrict ourselves to the case
that the base space is an elliptic curve.

\vspace{.1in}
As mentioned in \S\ref{subsection;13.10.16.10},
we obtain the filtered bundle $\nbigp_{\ast}E$
on $(T\times\proj^1,T\times\{\infty\})$
associated to any $L^2$-instanton $(E,\nabla,h)$,
and the metric $h$ is determined by $\nbigp_{\ast}E$
essentially uniquely.
We have the good filtered Higgs bundle
$(\nbigp_{\ast}\nbige,\theta)$ on $(T^{\lor},D)$
associated to any wild harmonic bundle
$(\nbige,\delbar_{\nbige},\theta,h_{\nbige})$,
and the metric $h_{\nbige}$ is determined by 
$(\nbigp_{\ast}\nbige,\theta)$
essentially uniquely
\cite{biquard-boalch}, \cite{Mochizuki-wild}.
So, it is significant to describe
the induced transformation between 
the underlying filtered bundles on $T\times\proj^1$
and the underlying good filtered Higgs bundles on $(T^{\lor},D)$,
that is given by the algebraic Nahm transform.
In particular, it is useful to understand
how the singular data are transformed.
It is also useful to prove that the Nahm transforms
are mutually inverse.

\subsubsection{Algebraic Nahm transform for
filtered Higgs bundles}

Let us briefly explain how the algebraic Nahm transform
is constructed for filtered Higgs bundles
$(\nbigp_{\ast}\nbige,\theta)$
on $(T^{\lor},D)$.
(The details will be given in \S\ref{section;13.1.30.10}
after the preliminary in \S\ref{section;13.10.29.1}.)
We should impose several conditions
to the filtered Higgs bundles.

\paragraph{Goodness and admissibility}

One of the conditions is the compatibility of 
the filtered bundle $\nbigp_{\ast}\nbige$
and the Higgs field $\theta$
at each $P\in D$.
Suppose that the filtered Higgs bundle
$(\nbigp_{\ast}\nbige,\theta)$ comes from
a good wild harmonic bundle.
Let $U_P$ be a small neighbourhood of $P$
with a coordinate $\zeta_P$ with
$\zeta_P(P)=0$.
If we take a ramified covering
$\varphi_p:U_P'\lrarr U_P$
given by $\varphi_p(u)=u^p=\zeta_P$,
for an appropriate $p$,
then we have the following decomposition:
\[
 \varphi_p^{\ast}(\nbigp_{\ast}\nbige,\theta)
=\bigoplus_{\gminia\in u^{-1}\cnum[u^{-1}]}
 (\nbigp_{\ast}\nbige'_{\gminia},\theta'_{\gminia})
\]
Here, $\theta'_{1,\gminia}-d\gminia$
are logarithmic
in the sense that
$(\theta'_{1,\gminia}-\gminia)\nbigp_{a}\nbige'_{\gminia}
\subset
 \nbigp_{a}\nbige'_{\gminia}\,du/u$.
Such a filtered Higgs bundle is called good.
This kind of filtered Higgs bundles
are closely related with
wild harmonic bundles
and $L^2$-instantons.

But, it seems natural to consider
an algebraic Nahm transform 
for a wider class of filtered Higgs bundles.
For $(p,m)\in\seisuu_{>0}\times\seisuu_{\geq 0}$
with $\gcd(p,m)=1$,
we say that a filtered Higgs bundle has type $(p,m)$
at $P$,
if we have that $u^{m}\cdot\varphi_p^{\ast}\theta$
gives a morphism of filtered bundles
$\varphi_p^{\ast}\nbigp_{\ast}\nbige
\lrarr
 \varphi_p^{\ast}\nbigp_{\ast}\nbige\,du/u$
on $U_P'$,
and that  it is an isomorphism 
in the case $(p,m)\neq (1,0)$.
We say that $(\nbigp_{\ast}\nbige,\theta)$ is admissible
at $P$,
if it is a direct sum 
$\bigoplus (\nbigp_{\ast}\nbige_{P}^{(p,m)},\theta_{P}^{(p,m)})$
of the filtered bundles of type $(p,m)$,
after $U_P$ is shrank appropriately.
We say that its slope is smaller (resp. strictly smaller) 
than $\alpha$,
if $\nbige^{(p,m)}_{P}=0$ for $m/p>\alpha$
(resp. $m/p\geq \alpha$).

Each $(\nbigp_{\ast}\nbige_{P}^{(p,m)},\theta_{P}^{(p,m)})$
has a refined decomposition,
as explained in \S\ref{subsection;13.1.21.200}.
In particular,
$(\nbigp_{\ast}\nbige_{P}^{(1,0)},\theta_{P}^{(1,0)})$
has a decomposition
\[
 (\nbigp_{\ast}\nbige_{P}^{(1,0)},\theta_{P}^{(1,0)})
=\bigoplus_{\alpha\in\cnum}
  (\nbigp_{\ast}\nbige^{(1,0)}_{P,\alpha},\theta_{P,\alpha}^{(1,0)}).
\]
Here, 
for the expression
$\theta^{(1,0)}_{P,\alpha}=
 f^{(1,0)}_{\alpha}\,d\zeta_P/\zeta_P$,
the eigenvalues of $f^{(1,0)}_{\alpha}$
goes to $\alpha$ when $\zeta_P\to 0$.
On $U_P$, we set
\[
 \nbigc^0(\nbigp_{\ast}\nbige,\theta)_P
=\bigoplus_{(p,m)\neq (1,0)}\nbigp_{-1/2-m/p}\nbige^{(p,m)}_{P}
\oplus
 \bigoplus_{\alpha\neq 0}\nbigp_{-1/2}\nbige^{(1,0)}_{P,\alpha}
\oplus
 \nbigp_{0}\nbige^{(1,0)}_{P,0}
\]
\[
  \nbigc^1(\nbigp_{\ast}\nbige,\theta)_P
=\Bigl(
 \bigoplus_{(p,m)\neq (1,0)}\nbigp_{1/2}\nbige^{(p,m)}_{P}
\oplus
 \bigoplus_{\alpha\neq 0}\nbigp_{1/2}\nbige^{(1,0)}_{P,\alpha}
\Bigr)
\otimes\Omega^1_{T^{\lor}}
\oplus
\Bigl(
 \nbigp_{<1}\nbige^{(1,0)}_{P,0}
 \otimes\Omega^1
+\theta^{(1,0)}_{P,0}
 \nbigp_{0}\nbige^{(1,0)}_{P,0}
\Bigr)
\]
The Higgs field $\theta$ gives a morphism
$\nbigc^0(\nbigp_{\ast}\nbige,\theta)_P
\lrarr
 \nbigc^1(\nbigp_{\ast}\nbige,\theta)_P$.
Thus, we obtain a complex
$\nbigc^{\bullet}(\nbigp_{\ast}\nbige,\theta)$
on $U_P$,
which is an extension of
$\nbige\stackrel{\theta}\lrarr \nbige\otimes\Omega^1$
on $U_P\setminus P$.

We say that $(\nbigp_{\ast}\nbige,\theta)$ on $(T^{\lor},D)$
is admissible,
if its restriction to a neighbourhood of each $P\in D$
is admissible.
By considering the extension at each $P\in D$,
we obtain a complex
$\nbigc^{\bullet}(\nbigp_{\ast}\nbige,\theta)$
on $T^{\lor}$,
as an extension of
$\nbige\lrarr\nbige\otimes \Omega^1$
on $T^{\lor}\setminus D$.

\paragraph{Vanishing of some cohomology groups}
For each $w\in\cnum$
and a holomorphic line bundle $L$ of degree $0$ 
on $T^{\lor}$,
we obtain a complex
$\nbigc^{\bullet}_{w,L}(\nbigp_{\ast}\nbige,\theta):=
 \nbigc^i(\nbigp_{\ast}\nbige\otimes L,\theta+w\,d\zeta)$.
To consider algebraic Nahm transform for
$(\nbigp_{\ast}\nbige,\theta)$,
it is natural to impose the following vanishing:
\begin{description}
\item[(A0)]
$\hyperh^i\bigl(
 T^{\lor},\nbigc^{\bullet}_{w,L}
 \bigr)=0$
unless $i=1$,
for any $w\in\cnum$
and a holomorphic line bundle $L$
of degree $0$ on $T^{\lor}$.
\end{description}

For $I\subset\{1,2,3\}$,
let $p_I$ denote the projection of
$T^{\lor}\times T\times\proj^1$
onto the product of the $i$-th components
$(i\in I)$.
Let $\poincare$ denote the Poincar\'e bundle
on $T^{\lor}\times T$.
We consider the following complex on 
$T^{\lor}\times T\times\proj^1$:
\[
\begin{CD}
 \nbigctilde^0:=
 p_1^{\ast}\nbigc^0\otimes
 p_{12}^{\ast}\poincare\otimes
 p_3^{\ast}\nbigo_{\proj^1}(-1)
@>{\theta+w\,d\zeta}>>
 \nbigctilde^1:=
 p_1^{\ast}\nbigc^1\otimes
 p_{12}^{\ast}\poincare
\end{CD}
\]
Then, it turns out that
$\vecN(\nbigp_{\ast}\nbige,\theta):=
 R^1p_{23\ast}\nbigctilde^{\bullet}$
is a locally free $\nbigo_{T\times\proj^1}$-module
on $T\times\proj^1$.
In particular,
we obtain a locally free
$\nbigo_{T\times\proj^1}(\ast(T\times\{\infty\}))$-module
\[
 \Nahm(\nbigp_{\ast}\nbige,\theta):=
 \vecN(\nbigp_{\ast}\nbige,\theta)
 \otimes\nbigo_{T\times\proj^1}
 \bigl(\ast (T\times\{\infty\})\bigr).
\]

\paragraph{Filtered bundles
on $(T\times U,T\times\{\infty\})$}

The algebraic Nahm transform
of $(\nbigp_{\ast}\nbige,\theta)$
is defined to be a filtered bundle
over the meromorphic bundle
$\Nahm(\nbigp_{\ast}\nbige,\theta)$,
i.e.,
an increasing sequence of 
$\nbigo_{T\times\proj^1}$-submodules
of $\Nahm(\nbigp_{\ast}\nbige,\theta)$.
For the construction of such a filtration,
it would be convenient to have a description of
filtered bundles $\nbigp_{\ast}E$ on
$(T\times U,T\times\{\infty\})$
satisfying the following condition,
where $U$ denotes a neighbourhood of $\infty$:
\begin{description}
\item[(A1)]
 $\Gr^{\nbigp}_c(E)$
 are semistable bundles of degree $0$
 on $T$ for any $c\in\real$.
\end{description}
If $U$ is shrank,
$\nbigp_c(E)_{|T\times\{w\}}$
are semistable of degree $0$
for any $w\in U$
and for any $c\in\real$.
We may assume it from the beginning.
We set
$\Sp_{\infty}(E):=
 \Sp(\nbigp_cE_{|T\times\{\infty\}})
\subset T^{\lor}$,
which is independent of $c\in\real$.
We fix a lift
$\Sptilde_{\infty}(E)\subset\cnum$
of $\Sp_{\infty}(E)$,
i.e.,
$\Sptilde_{\infty}(E)$ is mapped
bijectively to $\Sp_{\infty}(E)$
by the projection
$\cnum\lrarr T^{\lor}$.
Then, we have a filtered bundle
$\nbigp_{\ast}V$ on $(U,\infty)$
with an endomorphism $g$ 
such that $\Sp(g_{|\infty})=\Sptilde_{\infty}(E)$
corresponding to $\nbigp_{\ast}E$.
Namely, 
we have a $C^{\infty}$-isomorphism
$\nbigp_{\ast}E\simeq
 \pi^{\ast}\nbigp_{\ast}V$,
under which
$\delbar_{\nbigp_{\ast}E}
=\pi^{\ast}(\delbar_{\nbigp_{\ast}V})
+g\,d\zbar$,
where 
$\pi:T\times U\lrarr U$
denotes the projection.
The filtered bundle with an endomorphism
$(\nbigp_{\ast}V,g)$,
or equivalently,
the filtered Higgs bundle
$(\nbigp_{\ast}V,g\,dw)$
completely determines 
$\nbigp_{\ast}E$.
We have the decomposition
$(\nbigp_{\ast}V,g)=
 \bigoplus_{\alpha\in\Sptilde_{\infty}(E)}
 (\nbigp_{\ast}V_{\alpha},g_{\alpha})$
with
$\Sp(g_{\alpha|\infty})=\{\alpha\}$.
The filtered bundle $\nbigp_{\ast}E$ 
satisfying {\bf(A1)} is called admissible,
if the following holds:
\begin{description}
\item[(A2)]
 The filtered Higgs bundles
 $\bigl(\nbigp_{\ast}V_{\alpha},(g_{\alpha}-\alpha)\,dw\bigr)$
 are admissible for any $\alpha\in\Sptilde_{\infty}(E)$.
\end{description}
The slope of 
$\bigl(\nbigp_{\ast}V_{\alpha},(g_{\alpha}-\alpha)\,dw\bigr)$
is strictly smaller than $1$
by construction.

\paragraph{Local algebraic Nahm transform
and algebraic Nahm transform}

The local algebraic Nahm transform
$\nbign^{0,\infty}$ is a transform from
the germs of admissible filtered Higgs bundles
to the germs of admissible filtered Higgs bundles
whose slopes are strictly smaller than $1$.
It is an analogue of the local Fourier transform
$\gbigf^{0,\infty}$
for meromorphic flat bundles on $\proj^1$
in \cite{Bloch-Esnault} and \cite{Garcia-Lopez}.
(More precisely, it is an analogue of local Fourier
transform of minimal extension
of meromorphic flat bundles.)
It gives a procedure to produce
an admissible filtered bundle 
$\nbigp_{\ast}E_P$
on $(T\times U,T\times\{\infty\})$
such that $\Sp_{\infty}(E_P)=\{P\}$,
from an admissible filtered Higgs bundle
$(\nbigp_{\ast}\nbige,\theta)_{|U_P}$
on $(U_P,P)$.
From 
the local Nahm transform
$\bigoplus_{P\in D}\nbigp_{\ast}E_P$
and 
the meromorphic bundle
$\Nahm(\nbigp_{\ast}\nbige,\theta)$,
we obtain a filtered bundle on
$(T\times\proj^1,T\times\{\infty\})$,
denoted by
$\Nahm_{\ast}(\nbigp_{\ast}\nbige,\theta)$,
that is the algebraic Nahm transform
for admissible filtered Higgs bundles.

\subsubsection{Algebraic Nahm transform for
 admissible filtered bundles}

Let $\nbigp_{\ast}E$ be an admissible
filtered bundle on $(T\times\proj^1,T\times\{\infty\})$.
To define the algebraic Nahm transform
of $\nbigp_{\ast}E$,
we impose the following vanishing condition.
\begin{description}
\item[(A3)]
 $H^0\bigl(T\times\proj^1,\nbigp_0E\otimes L\bigr)=0$
 and 
 $H^2\bigl(T\times\proj^1,\nbigp_{<-1}E\otimes L\bigr)=0$
 for any holomorphic line bundle $L$ of degree $0$
on $T$.
\end{description}
We set $D:=\Sp_{\infty}(E)$.
It is easy to observe that
the condition {\bf (A3)} implies that
$H^i(T\times\proj^1,\nbigp_cE\otimes L^{\lor})=0$
$(i\neq 1)$
for any $c\in\real$, unless $L\in D$.
For any $I\subset\{1,2,3\}$,
let $p_{I}$ denote the projection
of $T^{\lor}\times T\times \proj^1$
onto the product of the $i$-th components
$(i\in I)$.
We define
\[
 \Nahm(\nbigp_{\ast}E):=
 R^1p_{1\ast}\bigl(
 p_{12}^{\ast}\poincare^{\lor}
 \otimes
 p_{23}^{\ast}\nbigp_0E
 \bigr)
 (\ast D)
\simeq
 R^1p_{1\ast}\bigl(
 p_{12}^{\ast}\poincare^{\lor}
 \otimes
 p_{23}^{\ast}\nbigp_{-1}E
 \bigr)
 (\ast D).
\]
It is a locally free
$\nbigo_{T^{\lor}}(\ast D)$-module.
The multiplication of $-w\,d\zeta$
gives a Higgs field $\theta$
of $\Nahm(\nbigp_{\ast}E)$.
Thus, we obtain a meromorphic Higgs bundle
on $(T^{\lor},D)$.
The algebraic Nahm transform of $\nbigp_{\ast}E$
is defined as a filtered bundle
over $\Nahm(\nbigp_{\ast}E)$ with $\theta$,
i.e.,
an increasing sequence of $\nbigo_{T^{\lor}}$-submodules
of $\Nahm(\nbigp_{\ast}E)$.

\vspace{.1in}

Let $U$ be a small neighbourhood of $\infty$
in $\proj^1$.
On $T\times U$,
we have a decomposition
$\nbigp_{\ast}E_{|T\times U}=
 \bigoplus_{P\in D}
 \nbigp_{\ast}E_P$
with 
$\Sp_{\infty}(E_P)=\{P\}$.
We have the corresponding
filtered bundles
$(\nbigp_{\ast}V_P,g_P)$.
We have the decomposition
$\bigl(\nbigp_{\ast}V_P,g_P-\Ptilde\bigr)
=\bigoplus
 (\nbigp_{\ast}V_P^{(p,m)},g_P^{(p,m)})$,
where 
$(\nbigp_{\ast}V_P^{(p,m)},g_P^{(p,m)}\,dw)$
has slope $(p,m)$ with $m/p<1$,
and $\Ptilde\in\Dtilde$ is a lift of $P$.
Moreover, we have the decomposition
\[
\bigl(
 \nbigp_{\ast}V_P^{(1,0)},g_P^{(1,0)}
\bigr)
=\bigoplus_{\alpha\in\cnum}
 \bigl(
 \nbigp_{\ast}V_{P,\alpha}^{(1,0)},g_{P,\alpha}^{(1,0)}
\bigr)
\]
It turns out that we have a decomposition
of the meromorphic bundle
\[
 \Nahm(\nbigp_{\ast}E)_{|U_P}
=
 \Nahm(\nbigp_{\ast}E)^{(1,0)}_{P,0}
\oplus
 \bigoplus_{\alpha\neq \cnum}
 \Nahm(\nbigp_{\ast}E)_{P,\alpha}^{(1,0)}
\oplus
\bigoplus_{(p,m)\neq (1,0)}
 \Nahm(\nbigp_{\ast}E)_P^{(p,m)},
\]
and $\Nahm(\nbigp_{\ast}E)_{P,\alpha}^{(1,0)}$ 
$(\alpha\neq 0)$
and 
$\Nahm(\nbigp_{\ast}E)_P^{(p,m)}$
($(p,m)\neq (1,0)$)
are determined by 
$(\nbigp_{\ast}V^{(1,0)}_{P,\alpha},g^{(1,0)}_{P,\alpha})$
and 
$(\nbigp_{\ast}V^{(p,m)}_{P,\alpha},g^{(p,m)}_{P})$.

\vspace{.1in}
We have the local algebraic Nahm transform
$\nbign^{\infty,0}$,
which is a transform of
admissible Higgs bundles
$(\nbigp_{\ast}V,\theta)$
such that 
(i) the slopes are strictly smaller than $0$,
(ii) $\nbigp_{\ast}V^{(1,0)}_{0}=0$.
It is an analogue of the local Fourier transform
$\gbigf^{\infty,0}$
in \cite{Bloch-Esnault} and \cite{Garcia-Lopez}.
It is an inverse of $\nbign^{0,\infty}$
except for the part 
$(p,m)=(1,0)$ and $\alpha=0$.
We may introduce filtrations of 
$\Nahm(\nbigp_{\ast}E)_{P,\alpha}^{(1,0)}$ 
$(\alpha\neq 0)$
and 
$\Nahm(\nbigp_{\ast}E)_P^{(p,m)}$
($(p,m)\neq (1,0)$)
by using the local algebraic Nahm transform
$\nbign^{\infty,0}$.
As for the part $(p,m)=(1,0)$ and $\alpha=0$,
we have an injection
$\nbigp_{0}V^{(1,0)}_{P,0|\infty}
\subset
 R^1p_{1\ast}\bigl(
 p_{12}^{\ast}\poincare^{\lor}
 \otimes
 p_{23}^{\ast}\nbigp_{-1}E
 \bigr)^{(1,0)}_{P,0|P}$,
by which we can introduce
a filtration on
$\Nahm(\nbigp_{\ast}E)^{(1,0)}_{P,0}$.
Thus, we obtain a filtered Higgs bundle
denoted by
$\Nahm_{\ast}(\nbigp_{\ast}E)$.
We obtain the following correspondence.
\begin{thm}[Proposition 
\ref{prop;13.10.18.2},
Proposition \ref{prop;13.10.15.11},
Proposition \ref{prop;13.1.22.2}]\mbox{{}}\\
The Nahm transforms $\Nahm_{\ast}$
give an equivalence of the following objects,
and they are mutually inverse:
\begin{itemize}
\item
Admissible filtered Higgs bundles
on $(T^{\lor},D)$
satisfying the condition {\bf (A0)}.
\item
Admissible filtered bundles
$\nbigp_{\ast}E$
on $(T\times\proj^1,T\times\{\infty\})$
with $\Sp_{\infty}(E)=D$
satisfying the condition {\bf (A3)}.
\end{itemize}
\end{thm}
We can also that
the Nahm transforms preserve the parabolic degrees
(Proposition \ref{prop;13.1.22.3}).

As already mentioned,
the filtered Higgs bundles associated to
wild harmonic bundles satisfy
a stronger condition called goodness.
Similarly, it turns out that
the filtered bundles associated to $L^2$-instantons
are also good,
in the sense that the corresponding filtered Higgs bundles
are good.
We can observe that
the algebraic Nahm transforms
preserve the goodness conditions.

\begin{thm}[Theorem
\ref{thm;13.2.3.2}]
The Nahm transforms $\Nahm_{\ast}$
gives an equivalence of the following objects:
\begin{itemize}
\item
Good filtered Higgs bundles
on $(T^{\lor},D)$ satisfying the condition {\bf (A0)}.
\item
Good filtered bundles
$\nbigp_{\ast}E$
on $(T\times\proj^1,T\times\{\infty\})$
with $\Sp_{\infty}(E)=D$
satisfying the condition {\bf (A3)}.
\end{itemize}
\end{thm}

\subsubsection{Application of the algebraic Nahm transform}

We have the following compatibility of 
the Nahm transform and the algebraic Nahm transform.

\begin{thm}[Theorem \ref{thm;13.1.23.20},
Theorem \ref{thm;13.1.23.100}]
\mbox{{}}\label{thm;13.10.17.1}
\begin{itemize}
\item
 Let $(E,\nabla,h)$ be an $L^2$-instanton
 on $T\times\cnum$.
 Let $\nbigp_{\ast}E$ be the associated
 filtered bundle
 on $(T\times\proj^1,T\times\{\infty\})$.
 Then, the underlying filtered Higgs bundle
 of the wild harmonic bundle 
 $\Nahm(E,\nabla,h)$ on $(T^{\lor},\Sp_{\infty}(E))$
 is naturally isomorphic to
 the algebraic Nahm transform
 $\Nahm(\nbigp_{\ast}E)$.
\item
 Let $(\nbige_,\delbar_{\nbige},\theta,h_{\nbige})$
 be a wild harmonic bundle on $(T^{\lor},D)$.
 Let $(\nbigp_{\ast}\nbige,\theta)$ be
 the associated good filtered Higgs bundle
 on $(T^{\lor},D)$.
 Then, the underlying filtered bundle of
 the $L^2$-instanton
 $\Nahm(\nbige,\delbar_{\nbige},h_{\nbige},\theta)$
 is naturally isomorphic to
 the algebraic Nahm transform
 $\Nahm(\nbigp_{\ast}\nbige,\theta)$.
\hfill\qed
\end{itemize}
\end{thm}

As an application, we obtain
the inversion property of the Nahm transforms.

\begin{cor}[Corollary \ref{cor;13.1.30.31}]
For an $L^2$-instanton
$(E,\nabla,h)$ on $T\times\cnum$,
we have an isomorphism
\[
 \Nahm\bigl(
 \Nahm(E,\nabla,h)
 \bigr)
\simeq
 (E,\nabla,h).
\]
For a wild harmonic bundle 
$(\nbige,\delbar_{\nbige},\theta,h_{\nbige})$
on $(T^{\lor},D)$,
we have an isomorphism
\[
 \Nahm\bigl(
 \Nahm(\nbige,\delbar_{\nbige},\theta,h_{\nbige})
 \bigr)
\simeq (\nbige,\delbar_{\nbige},\theta,h_{\nbige})
\]
\end{cor}
Indeed,
it follows from Theorem \ref{thm;13.10.17.1}
and the uniqueness of the Hermitian-Einstein metric
(resp. the harmonic metric)
 adapted to the filtered bundle
(resp. filtered Higgs bundle).

As another application of the compatibility,
we can easily compute the characteristic classes
of the bundles obtained as the algebraic Nahm transform,
which allows us to obtain the rank and the second Chern class
of the bundle obtained as the Nahm transform.
The local algebraic Nahm transform also leads us 
a rather complete understanding
of the transformation of singularity data
by the Nahm transform.

\subsubsection{Some remarks}

Recall that 
the hyperkahler manifold
$T\times\cnum$ has a twistor deformation.
Namely, for any complex number $\lambda$,
we have the moduli space $\nbigm^{\lambda}$
of the line bundles of degree $0$ with
flat $\lambda$-connection on $T^{\lor}$.
We have $\nbigm^0=T\times\cnum$.
The spaces $\nbigm^{\lambda}$
can also be regarded as the deformation
associated to the hyperkahler structure of
$T\times\cnum$.
An instanton on $T\times\cnum$ naturally induces
a holomorphic vector bundle.
If the instanton is $L^2$,
the holomorphic bundle with the metric
induces a filtered bundle on 
$(\nbigmbar^{\lambda},T_{\infty}^{\lambda})$,
where
$\nbigmbar^{\lambda}$ is a natural compactification
of $\nbigmlambda$,
and $T_{\infty}^{\lambda}\simeq T$
is the infinity. 
A wild harmonic bundle 
has the underlying good filtered
$\lambda$-flat bundle 
for each complex number $\lambda$.
It is also natural to study the transformation
of the underlying filtered bundles 
on $(\nbigmbar^{\lambda},T^{\lambda}_{\infty})$
and the underlying filtered $\lambda$-flat bundles.
It should be a filtered enhancement of
the generalized Fourier-Mukai transform
for elliptic curves
due to G. Laumon and M. Rothstein.
We would like to study this interesting aspect
elsewhere.

\vspace{.1in}

If we consider a counterpart of the algebraic Nahm transform
for the other non-product holomorphic structure of $T\times\cnum$
underlying the hyperKahler structure,
it is essentially a filtered version of
the generalized Fourier-Mukai transform
in {\rm\cite{Laumon-Fourier-Mukai}} and {\rm\cite{Rothstein}}.
Interestingly to the author,
we have an analogue of the stationary phase formula
even in this case.
The details will be given elsewhere.

\vspace{.1in}

In this paper,
we consider transforms between
filtered bundles on $T\times\proj^1$
and filtered Higgs bundles on $T^{\lor}$.
We may introduce similar transforms
for filtered Higgs bundles on $\proj^1$
with an additional work on the local Nahm transform
$\nbign^{\infty,\infty}$,
which is an analogue of the local Fourier transform
$\gbigf^{\infty,\infty}$.
It should be the Higgs counterpart of the Nahm transforms
between wild harmonic bundles on $\proj^1$,
which is given by the procedure for wild pure twistor $D$-modules
established in {\rm\cite{Mochizuki-wild}}.

Relatedly,
Sz. Szab\'{o} {\rm\cite{Szabo}} studied the Nahm transform
for some interesting type of harmonic bundles on $\proj^1$.
He also studied the transform of the underlying parabolic Higgs bundles,
which looks closely related with 
the regular version of ours in {\rm\S\ref{subsection;13.2.14.25}}.
K. Aker and Szab\'{o} {\rm\cite{Aker-Szabo}}
introduced a transformation of 
more general parabolic Higgs bundles on $\proj^1$,
which they call the algebraic Nahm transform.
Their method to define the transform is different from ours,
and the precise relation is not clear at this moment.

\subsection{Acknowledgement}

I thank the referee for his effort to read this paper
and for his valuable suggestions to improve it.
It is one of my motivations of this study 
to understand and generalize
some of the interesting results in \cite{Biquard-Jardim}.
I thank Masaki Tsukamoto
for his inspiring enthusiasm to instantons,
who attracted my attention to 
doubly periodic $L^2$-instantons.
I thank Max Lipyanskiy who kindly pointed out
an error in an earlier version,
in using a theorem of Uhlenbeck
on the curvature of instantons.
The argument is simplified after the issue is fixed.
I am grateful to Szilard Szab\'{o} for many discussions,
who first introduced me the Nahm transform.
I thank Hiraku Nakajima for his lecture on the Nahm transform.
I thank Claude Sabbah for his kindness and many discussions.
I am grateful to Christian Schnell for some motivating discussions
on the Fourier-Mukai-Laumon-Rothstein transform.
I thank Philip Boalch for his interest and encouragement
to this study.
I thank Akira Ishii and Yoshifumi Tsuchimoto
for their encouragement.
This study is supported by 
Grant-in-Aid for Scientific Research (C) 22540078.

\section{Preliminaries on filtered objects}
\label{section;13.10.29.1}

\subsection{Semistable bundles of degree $0$
 on elliptic curves}
\label{subsection;13.10.15.30}

\subsubsection{Elliptic curve and Fourier-Mukai transform}
\label{subsection;12.6.23.1}

For a variable $z$,
let $\cnum_z$ denote a complex line 
with the standard coordinate $z$.
For a $\cnum$-vector space $V$ and a $C^{\infty}$-manifold $X$,
let $\underline{V}_X$ denote the product bundle
$V\times X$ over $X$.
If $X$ is a complex manifold,
the natural holomorphic structure of
$\underline{V}_X$ is denoted just by $\delbar$.

We have a real bilinear map
$\cnum_z\times\cnum_{\zeta}\lrarr\real$
given by 
$(z,\zeta)\longmapsto
 \Image(z\zetabar)$.
Let 
$\tau=\tau_1+\sqrt{-1}\tau_2$ $(\tau_i\in\real)$
be a complex number such that $\tau_2\neq 0$.
Let $L:=\seisuu+\seisuu\tau\subset \cnum_z$.
In this paper, the dual lattice $L^{\lor}$ means the following:
\[
 L^{\lor}:=
 \bigl\{
 \zeta\in \cnum_{\zeta}\,\big|\,
 \Image(\chi\zetabar)\in \pi\seisuu
\,\,\,(\forall\chi\in L)
 \bigr\}
=\Bigl\{
 \frac{\pi}{\tau_2}
 (n+m\tau)\,\Big|\,
 n,m\in\seisuu
 \Bigr\}
\]
We obtain the elliptic curves
$T:=\cnum_z/L$
and 
$T^{\lor}:=\cnum_{\zeta}/L^{\lor}$.

For any $\nu\in L^{\lor}$,
we have $\rho_{\nu}\in C^{\infty}(T)$
given by
$\rho_{\nu}(z):=\exp\bigl(
 2\sqrt{-1}\Image(\nu\zbar)
 \bigr)
=\exp\bigl(
 \nu\zbar-\nubar z
 \bigr)$.
We have 
$\delbar_z\rho_{\nu}
=\rho_{\nu}\, \nu\,d\zbar$
and
$\del_z\rho_{\nu}=-\rho_{\nu}\,\nubar\,dz$.

We can naturally regard $T^{\lor}$
as the moduli space $\Pic_0(T)$
of holomorphic line bundles of degree $0$ on $T$.
Indeed,
$\zeta$ gives a holomorphic bundle
$\nbigl_{\zeta}=(\underline{\cnum}_T,\delbar+\zeta\,d\zbar)$.
It induces an isomorphism $T^{\lor}\simeq\Pic_0(T)$.
We have the isomorphism
$\Phi:\nbigl_{\zeta}\simeq\nbigl_{\zeta+\nu}$
given by $\Phi(f)=\rho_{-\nu}\cdot f$.

We have the unitary flat connection associated to
$\nbigl_{\zeta}$
with the trivial metric,
which is $d-\zetabar\,dz+\zeta\,d\zbar$.
The monodromy along the segment from $0$ to 
$\chi\in L$ is
$\exp\Bigl(
 \int_0^1
(-\zetabar\chi+\zeta\chibar)\,dt
 \Bigr)
=\exp\bigl(
 2\sqrt{-1}\Image(\zeta\chibar)
 \bigr)$.

We recall a differential geometric construction of the Poincar\'e bundle
on $T\times T^{\lor}$,
following \cite{Donaldson-Kronheimer}.
On $T\times \cnum_{\zeta}$,
we have the holomorphic line bundle
$\poincaretilde=\bigl(
 \underline{\cnum}_{T\times\cnum_{\zeta}},
 \delbar+\zeta\,d\zbar
 \bigr)$.
The $L^{\lor}$-action on $T\times \cnum_{\zeta}$
is naturally lifted to the action on $\poincaretilde$
given by
$\nu(z,\zeta,v)
=\bigl(z,\zeta+\nu,\rho_{-\nu}(z)\,v\bigr)$.
Thus, a holomorphic line bundle is induced on
$T\times T^{\lor}$,
which is the Poincar\'e bundle
denoted by $\poincare$.
The dual bundle $\poincare^{\lor}$
is induced by
$\poincaretilde^{\lor}=\bigl(
 \underline{\cnum}_{T\times\cnum_{\zeta}},
 \delbar-\zeta\,d\zbar
 \bigr)$
with the action
$\nu(z,\zeta,v)
=\bigl(z,\zeta+\nu,\rho_{\nu}(z)\,v\bigr)$.

\vspace{.1in}

Let $S$ be any complex analytic space.
For $I\subset\{1,2,3\}$,
let $p_I$ denote the projection of 
$T\times T^{\lor}\times S$
onto the product of the $i$-th components
$(i\in I)$.
For an object $\nbigm\in D^b\bigl(\nbigo_{T\times S}\bigr)$,
we set
\[
 \RFM_{\pm}(\nbigm):=Rp_{23\ast}\bigl(
 p_{13}^{\ast}(\nbigm)\otimes p_{12}^{\ast}\poincare^{\pm 1}\bigr)[1]
\in D^b(\nbigo_{T^{\lor}\times S}).
\]
For an object $\nbign\in D^b(\nbigo_{T^{\lor}\times S})$,
we set 
\[
 \widehat{\RFM}_{\pm}(\nbign):=
 Rp_{13\ast}\bigl(
 p_{23}^{\ast}(\nbign)\otimes p_{12}^{\ast}\poincare^{\pm 1}
 \bigr)
\in D^b(\nbigo_{T\times S}).
\]
Recall that we have a natural isomorphism
$\widehat{\RFM}_+\circ\RFM_-(\nbigm)\simeq\nbigm$.

\subsubsection{Semistable bundle of degree $0$}
\label{subsection;13.1.28.100}

For a holomorphic vector bundle $(E,\delbar_E)$ on $T$,
we have the degree $\deg(E):=\int_Tc_1(E)$
and the slope $\mu(E):=\deg(E)/\rank(E)$.
A holomorphic vector bundle $E$ on $T$ is called semistable,
if $\mu(F)\leq \mu(E)$ holds 
for any non-trivial holomorphic subbundle $F\subset E$.
In the following,
we shall not distinguish a holomorphic vector bundle
and the associated sheaf of holomorphic sections.

Let $E$ be a semistable bundle of degree $0$ on $T$.
It is well known that
the support $\Sp(E)$ of $\RFM_-(E)$ are finite points.
Indeed, $E$ is obtained as an extension of
the line bundles $\nbigl_{\zeta}$ $(\zeta\in\Sp(E))$.
We call $\Sp(E)$ the spectrum of $E$.
We have the decomposition
$E=\bigoplus_{\alpha\in \Sp(E)}
 E_{\alpha}$,
where the support of
$\RFM_-(E_{\alpha})$ is $\{\alpha\}$.
It is called the spectral decomposition of $E$.
We call a subset $\Sptilde(E)\subset\cnum$
is a lift of $\Sp(E)$,
if the projection $\Phi:\cnum\lrarr T^{\lor}$
induces the bijection $\Sptilde(E)\simeq\Sp(E)$.
If we fix a lift, 
an $\nbigo_{\cnum}$-module $\nbigm(E)$
is determined (up to canonical isomorphisms)
by the conditions
(i) the support of $\nbigm(E)$ is $\Sptilde(E)$,
(ii) $\Phi_{\ast}\nbigm(E)\simeq \RFM_-(E)$.
Such $\nbigm(E)$ is called a lift of $\RFM_-(E)$.
The multiplication of $\zeta$ on $\nbigm(E)$
induces an endomorphisms of $\RFM_-(E)$ and $E$.
The endomorphism of $E$ is denoted by $f_{\zeta}$.

\vspace{.1in}

Let $S$ be any complex analytic space.
Let $E$ be a holomorphic vector bundle on $T\times S$.
It is called semistable of degree $0$ relatively to $S$,
if $E_{|T\times\{s\}}$ are semistable of degree $0$
for any $s\in S$.
The support of $\RFM_-(E)$ is relatively $0$-dimensional over $S$.
It is denoted by $\Sp(E)$, and called the spectrum of $E$.
If we have a hypersurface 
$\Sptilde(E)\subset \cnum_{\zeta}\times S$
such that
the projection $\Phi:\cnum_{\zeta}\times S\lrarr T^{\lor}\times S$
induces $\Sptilde(E)\simeq \Sp(E)$,
we call $\Sptilde(E)$ a lift of $\Sp(E)$.
If we have a lift of $\Sp(E)$,
we obtain a lift $\nbigm(E)$ of $\RFM_-(E)$
as in the case that $S$ is a point.
We also obtain an endomorphism $f_{\zeta}$ of $E$
induced by the multiplication of $\zeta$
on $\nbigm(E)$.

\subsubsection{Equivalence of categories}
\label{subsection;13.1.21.10}

For a vector space $V$,
let $\underline{V}$ denote the product bundle
$T\times V$ over $T$,
and let $\delbar_0$ denote the natural holomorphic structure
of $\underline{V}$.
For any $f\in\End(V)$,
we have the associated holomorphic vector bundle
$\gbigg(V,f):=
(\underline{V},\delbar_0+f\,d\zbar)$.
We have a natural isomorphism
$\gbigg(V,f)\simeq \gbigg(V,f+\nu\id_V)$
for each $\nu\in L^{\lor}$,
induced by the multiplication of $\rho_{-\nu}$.
Let $\Sp(f)$ denote the set of the eigenvalues of $f$.

\begin{lem}
\label{lem;13.1.28.1}
$\gbigg(V,f)$ is semistable of degree $0$,
and we have
$\Sp(\gbigg(V,f))
=\Phi\bigl( \Sp(f) \bigr)$
in $T^{\lor}$,
where $\Phi:\cnum\lrarr T^{\lor}$
denotes the projection.
\end{lem}
\pf
We have only to consider the case $f$ has a unique 
eigenvalue $\alpha$.
In that case, $\gbigg(V,f)$ is an extension of
the line bundle $\nbigl_{\alpha}$.
Then, the claim is clear.
\hfill\qed

\vspace{.1in}

Let $VS^{\ast}$ denote the category of finite dimensional
$\cnum$-vector spaces with an endomorphism,
i.e., an object in $\VS^{\ast}$
is a finite dimensional 
vector space $V$ with an endomorphism $f$,
and a morphism $(V,f)\lrarr (W,g)$ in $\VS^{\ast}$
is a linear map $\varphi:V\lrarr W$
such that $g\circ\varphi-\varphi\circ f=0$.
For a given subset $\gministilde\subset\cnum$,
let $\VS^{\ast}(\gministilde)\subset \VS^{\ast}$ 
denote the full subcategory of
$(V,f)$ such that $\Sp(f)\subset\gministilde$.

Let $\VB^{ss}_0(T)$ denote the category of
semistable bundles of degree $0$ on $T$,
i.e., an object in $\VB^{ss}_0(T)$
is a semistable vector bundle of degree $0$ on $T$,
and a morphism $V_1\lrarr V_2$ in $\VB^{ss}_0(T)$
is a morphism of coherent sheaves.
For a given subset $\gminis\subset T^{\lor}$,
let $\VB^{ss}_0(T,\gminis)\subset\VB^{ss}_0(T)$ 
denote the full subcategory
of semistable bundles of degree $0$
whose spectrum are contained in $\gminis$.

We have the functor
$\gbigg:\VS^{\ast}\lrarr \VB^{ss}_0(T)$
given by the above construction.
If $\gministilde$ is mapped to $\gminis$
by the projection $\Phi:\cnum_{\zeta}\lrarr T^{\lor}$,
it induces a functor
$\gbigg:
\VS^{\ast}(\gministilde)\lrarr \VB^{ss}_0(T,\gminis)$.

\begin{prop}
\label{prop;13.1.8.22}
If $\Phi:\cnum\lrarr T^{\lor}$ 
induces a bijection $\gministilde\simeq\gminis$,
then $\gbigg$ gives an equivalence of the categories
$\VS^{\ast}(\gministilde)\simeq \VB^{ss}_0(T,\gminis)$.
\end{prop}
\pf
Let us show that it is fully faithful.
We set $E_f:=\gbigg(V,f)$.
We will not distinguish $E_f$
and the associated sheaf of holomorphic sections.
Suppose that $f$ has a unique eigenvalue $\alpha$
such that $\alpha\not\equiv 0$ modulo $L^{\lor}$.
Because $E_f$ is obtained as an extension of
the holomorphic line bundle $\nbigl_{\alpha}$,
we have $H^0(T,E_f)=H^1(T,E_f)=0$.
In particular, we obtain the following.
\begin{lem}
\label{lem;13.1.8.20}
Assume that $f_i\in \End(V)$ has a unique eigenvalue $\alpha_i$
for $i=1,2$.
If $\alpha_1\not\equiv\alpha_2$ modulo $L^{\lor}$,
any morphism
$E_{f_1}\lrarr E_{f_2}$ is $0$.
\hfill\qed
\end{lem}

Suppose that $f$ is nilpotent.
We have the natural inclusion $V\lrarr C^{\infty}(T,E_f)$
as constant functions.
We have a linear map
$V\lrarr C^{\infty}(T,E_f\otimes\Omega^{0,1})$
given by $s\longmapsto s\,d\zbar$.
They induce a chain map $\iota$ from 
$\nbigc_1=(V\stackrel{f}{\lrarr} V)$ to 
the Dolbeault complex 
$C^{\infty}(T,E_f\otimes\Omega_T^{0,\bullet})$
of $E_f$.
\begin{lem}
$\iota$ is a quasi-isomorphism.
\end{lem}
\pf
Let $W$ be the monodromy weight filtration of $f$.
It induces filtrations of $\nbigc_1$
and $C^{\infty}(T,E_f\otimes\Omega_T^{0,\ast})$,
and $\iota$ gives a morphism of filtered chain complex.
It induces a quasi-isomorphism of 
the associated graded complexes.
Hence, $\iota$ is a quasi isomorphism.
\hfill\qed

\vspace{.1in}
We obtain the following lemma 
as an immediate consequence.

\begin{lem}
\label{lem;13.1.8.21}
Assume that $f_i\in \End(V)$ are nilpotent $(i=1,2)$.
Then, holomorphic morphisms $E_{f_1}\lrarr E_{f_2}$
naturally correspond to holomorphic morphisms 
$\phi:E_0\lrarr E_0$ such that
$f_2\circ\phi-\phi\circ f_1=0$.

In particular,
if $f$ is nilpotent,
holomorphic sections of
$\End(E_f)$ bijectively corresponds
to holomorphic sections $g$ of $\End(E_0)$
such that $[f,g]=0$.
\hfill\qed
\end{lem}

The fully faithfulness of the functor
$\gbigg$ follows from Lemma \ref{lem;13.1.8.20}
and Lemma \ref{lem;13.1.8.21}.
Let us show the essential surjectivity of $\gbigg$.
Let $E\in \VB_0^{ss}(T,\gminis)$.
We have the $\nbigo_{\cnum_{\zeta}}$-module $\nbigm(E)$
and the endomorphism $f_{\zeta}$ of $E$
as in \S\ref{subsection;13.1.28.100}.
We have a natural isomorphism
$\widehat{\RFM}_+\circ \RFM_-(E)\simeq E$.
The functor
$\widehat{\RFM}_+$ is induced by
the holomorphic line bundle on $T\times T^{\lor}$,
obtained as the descent of
$\poincaretilde=
 (\underline{\cnum},\delbar_0+\zeta\,d\zbar)$.
Let $p$ and $q$ denote the projections
$T\times\cnum_{\zeta}\lrarr T$
and $T\times\cnum_{\zeta}\lrarr\cnum_{\zeta}$.
We have
$E\simeq
 p_{\ast}(q^{\ast}(\nbigm)\otimes\poincaretilde)$,
and the latter is naturally isomorphic to
$\bigl(
 \underline{H^0(\cnum_{\zeta},\nbigm(E))},
 \delbar_0+\,f_{\zeta}\,d\zbar
 \bigr)$.
Thus, we obtain the essential surjectivity of
$\gbigg$.
The proof of Proposition \ref{prop;13.1.8.22}
is finished.
\hfill\qed

\vspace{.1in}

As appeared in the proof of Proposition \ref{prop;13.1.8.22},
we have another equivalent construction of $\gbigg$.
Let $\nbign'(V,f)$ denote the cokernel of
the endomorphism 
$\zeta\id-f$
on $V\otimes\nbigo_{\cnum_{\zeta}}$.
It naturally induces an $\nbigo_{T^{\lor}}$-module
$\nbign(V,f)$.
We obtain $\widehat{\RFM}_+\bigl(\nbign(V,f)\bigr)$,
which is naturally isomorphic to
$\gbigg(V,f)$.
We obtain a quasi-inverse of $\gbigg$ as follows.
Let $E$ be a semistable bundle of degree $0$ on $T$.
We obtain a vector space $H^0(T^{\lor},\RFM_-(E))$.
If we fix a lift of $\Sp(E)$ to 
$\Sptilde(E)\subset \cnum$,
then the multiplication of $\zeta$ induces
an endomorphism $g_{\zeta}$ of $H^0(T^{\lor},\RFM_-(E))$.
The construction of 
$\bigl(H^0(T^{\lor},\RFM_-(E)),g_{\zeta}\bigr)$
from $E$ gives a quasi-inverse of $\gbigg$.

\vspace{.1in}
Let $(E,\delbar_{E})$ be a semistable bundle 
of degree $0$ on $T$.
Let $\gministilde$ be a subset of $\cnum$
which is injectively mapped to $T^{\lor}=\cnum/L^{\lor}$.
\begin{cor}
We have a unique decomposition 
$\delbar_{E}=\delbar_{E,0}+f\,d\zbar$
with the following property:
\begin{itemize}
\item
$(E,\delbar_{E,0})$ is holomorphically trivial,
i.e., isomorphic to
a direct sum of $\nbigo_T$.
\item
 $f$ is a holomorphic endomorphism of 
 $(E,\delbar_{E,0})$.
 We impose the condition that
 $\Sp\bigl(H^0(f)\bigr)\subset \gministilde$,
 where $H^0(f)$ is the induced endomorphism of
 the space of the global sections of $(E,\delbar_{E,0})$.
\end{itemize}
\end{cor}
\pf
The existence of such a decomposition
follows from the essential surjectivity of
$\gbigg$.
Let us show the uniqueness.
By considering the spectral decomposition,
we have only to consider the case $\gministilde=\{0\}$.
Suppose that $\delbar_E=\delbar'_{E,0}+g\,d\zbar$
is another decomposition with the desired property.
Because $f$ is holomorphic with respect to $\delbar_E$,
we have $\delbar'_{E,0}f=0$ and $[f,g]=0$
by Lemma \ref{lem;13.1.8.21}.
We put $h=f-g$, which is also nilpotent.
The identity induces an isomorphism
$(E,\delbar_{E,0}+h)
\simeq
 (E,\delbar'_{E,0})$.
Because $\gbigg$ is fully faithful,
we obtain $h=0$.
\hfill\qed

\paragraph{The family version}
We have a family version of the equivalence.
Let $S$ be any complex manifold.
Let $\pi_S:T\times S\lrarr S$ denote the projection.
Let $\VB^{\ast}(S)$ denote the category of
pairs $(V,f)$ of coherent locally free $\nbigo_S$-module $V$
and its endomorphism $f$.
A morphism $(V,f)\lrarr (V',f')$ in $\VB^{\ast}(S)$
is a morphism of $\nbigo_S$-modules
$g:V\lrarr V'$ such that $f'\circ g=g\circ f$.
Such $(V,f)$ naturally induces an 
$\nbigo_{\cnum_{\zeta}\times S}$-module
$\nbigm(V,f)$.
The support is denoted by $\Sp(f)$.
When we are given a divisor 
$\gministilde\subset \cnum_{\zeta}\times S$
which is finite and relatively $0$-dimensional over $S$,
then $\VB^{\ast}(S,\gministilde)$
denote the full subcategory of 
$(V,f)\in \VB^{\ast}(S)$
such that $\Sp(f)\subset\gministilde$.

Let $\VB^{ss}_0(T\times S/S)$
denote the full subcategory of 
$\nbigo_{T\times S}$-modules,
whose objects are semistable of degree $0$
relative to $S$.
When we are given a divisor 
$\gminis\subset T^{\lor}\times S$
which is relatively $0$-dimensional over $S$,
then let
$\VB^{ss}_0(T\times S/S,\gminis)$
denote the full subcategory of $E\in\VB^{ss}_0(T\times S/S)$
such that $\Sp(E)\subset\gminis$.

Let $V$ be a holomorphic vector bundle on $S$
with a holomorphic endomorphism $f$.
The $C^{\infty}$-vector bundle $\pi_S^{-1}V$
is equipped with a naturally induced holomorphic structure
obtained as the pull back,
denoted by $\delbar_0$.
We obtain a holomorphic vector bundle
$\gbigg(V,f):=(\pi_S^{-1}V,\delbar_0+f\,d\zbar)$.
By Lemma \ref{lem;13.1.28.1},
$\gbigg$ gives a functor
$\VB^{\ast}(S)\lrarr
 \VB^{ss}_0(T\times S/S)$.
If we are given $\gminis\subset T^{\lor}\times S$
and its lift $\gministilde\subset\cnum_{\zeta}\times S$,
it gives an equivalence of the categories
$\VB^{\ast}(S,\gminis)\lrarr
 \VB^{ss}_0(T\times S/S,\gministilde)$.

We have another equivalent description of $\gbigg$.
Let $(V,f)\in\VB^{\ast}(S)$.
We have the naturally induced 
$\nbigo_{\cnum_{\zeta}\times S}$-module
$\nbigm(V,f)$,
which induces an $\nbigo_{T^{\lor}\times S}$-module
$\nbign(V,f)$.
We have a natural isomorphism
$\gbigg(V,f)\simeq
 \RFM_-\bigl(\nbign(V,f)\bigr)$.

If we are given $\gminis\subset T^{\lor}\times S$
with a lift $\gministilde\subset\cnum_{\zeta}\times S$,
for an object $E\in \VB^{ss}_0(T\times S/S,\gminis)$,
we obtain an $\nbigo_{\cnum_{\zeta}\times S}$-module
$\nbigm(E)$ such that
(i) the support of $\nbigm(E)$ is contained in $\gministilde$,
(ii) $\Phi_{\ast}\nbigm(E)\simeq\RFM_-(E)$.
The multiplication of $\zeta$ induces
an endomorphism of $\RFM_-(E)$,
and hence an endomorphism of $\pi_{S\ast}(\RFM_-(E))$,
denoted by $g_{\zeta}$,
where $\pi_S:T^{\lor}\times S\lrarr S$.
The construction of 
$\bigl(\pi_{S\ast}\RFM_-(E),g_{\zeta}\bigr)$
from $E$ gives a quasi-inverse of $\gbigg$.

\subsubsection{Differential geometric criterion}

We recall a differential geometric criterion 
in terms of the curvature
for a metrized holomorphic vector bundle
to be semistable of degree $0$.
Let $(E,\delbar_E)$ be a holomorphic vector bundle
on $T$ with a hermitian metric $h$.
Let $F(h)$ denote the curvature.
We use the standard metric $dz\,d\zbar$ of $T$.

\begin{lem}
\label{lem;12.6.23.101}
There exists a constant $\epsilon>0$,
depending only on $T$ and $\rank E$,
with the following property:
\begin{itemize}
\item
 If $|F(h)|_h\leq \epsilon$,
 then $(E,\delbar_E,h)$ is semistable of degree $0$.
\end{itemize}
\end{lem}
\pf
The number $\deg(E)=\int \Tr F(h)$ is an integer.
We have $\int \bigl|\Tr F(h)\bigr|\leq 
|T|\rank E\epsilon$,
where $|T|$ is the volume of $T$.
Hence, we have $\int \Tr F(E_w)=0$,
if $\epsilon$ is sufficiently small.
For any subbundle $E'\subset E$,
by using the decreasing property of the curvature of subbundles,
we also obtain $\deg(E')<1$ and hence $\deg(E')\leq 0$.
\hfill\qed

\subsection{Filtered bundles}
\label{subsection;13.2.26.3}

\subsubsection{Filtered sheaves}
\label{subsection;13.10.28.2}

Let us recall the notion of filtered sheaves
and filtered bundles.
Let $X$ be a complex manifold
with a smooth hypersurface $D$.
(We restrict ourselves to the case that
$D$ is smooth,
because we are interested in only the case
in this paper.)
Let $\nbige$ be a coherent $\nbigo_X(\ast D)$-module.
Let $D=\coprod_{i\in\Lambda}D_i$
be the decomposition into the connected components.
A filtered sheaf $\nbigp_{\ast}\nbige$
over $\nbige$ is
a sequence of coherent $\nbigo_X$-submodules
$\nbigp_{\veca}\nbige\subset\nbige$ 
indexed by $\real^{\Lambda}$ 
satisfying the following.
\begin{itemize}
\item
$\nbigp_{\veca}\nbige_{|X\setminus D}=\nbige_{|X\setminus D}$.
We have
$\nbigp_{\veca}\nbige\subset\nbigp_{\veca'}\nbige$,
if $a_i\leq a_i'$ $(i\in\Lambda)$,
where $\veca=(a_i\,|\,i\in\Lambda)$
and $\veca'=(a_i'\,|\,i\in\Lambda)$.
\item
On a small neighbourhood $U$ of $D_i$ $(i\in\Lambda)$,
$\nbigp_{\veca}\nbige_{|U}$ depends only on $a_i$,
which we denote by $\nbigp_{a_i}(\nbige_{|U})$,
or 
$\lefttop{i}\nbigp_{a_i}(\nbige_{|U})$.
when we emphasize $i$.
\item
For each $i$ and $c\in\real$,
there exists $\epsilon>0$
such that
$\lefttop{i}\nbigp_{c}(\nbige_{|U})
=\lefttop{i}\nbigp_{c+\epsilon}(\nbige_{|U})$.

\item
We have 
$\nbigp_{\veca+\vecn}\nbige
=\nbigp_{\veca}\nbige(\sum n_iD_i)$,
where $\vecn=(n_i)\in\seisuu^{\Lambda}$.
\end{itemize}
The tuple
$(\nbige,\{\nbigp_{\veca}\nbige\,|\,\veca\in\real^{\Lambda}\})$
is denoted by $\nbigp_{\ast}\nbige$.
The filtration $\{\nbigp_{\veca}\nbige\,|\,\veca\in\real^{\Lambda}\}$
is also denoted by $\nbigp_{\ast}\nbige$.
We say that $\nbige$ is the $\nbigo_X(\ast D)$-module
underlying $\nbigp_{\ast}\nbige$.

For a small neighbourhood $U$ of $D_i$,
we set 
$\lefttop{i}\nbigp_{<a}(\nbige_{|U}):=
 \sum_{b<a}\nbigp_b(\nbige_{|U})$.
We also put
$\lefttop{i}\nbigp_a(\nbige)_{|D_i}:=
 \nbigp_a(\nbige_{|U})_{|D_i}$,
and
$\lefttop{i}\Gr^{\nbigp}_a(\nbige):=
 \lefttop{i}\nbigp_{a}(\nbige_{|U})\big/\lefttop{i}\nbigp_{<a}(\nbige_{|U})$,
which are coherent $\nbigo_{D_i}$-modules.
We set
\[
 \Par(\nbigp_{\veca}\nbige,i):=
 \bigl\{
 b\in\openclosed{a_i-1}{a_i}
 \,\big|\,\lefttop{i}\Gr^{\nbigp}_b(\nbige)\neq 0
 \bigr\},
\quad
 \Par(\nbigp_{\ast}\nbige,i):=
 \bigcup_{\veca\in\real^{\Lambda}}
 \Par(\nbigp_a\nbige,i).
\]

A morphism of filtered sheaves
$\nbigp_{\ast}\nbige_1\lrarr\nbigp_{\ast}\nbige_2$
is a morphism of $\nbigo_X$-modules
which is compatible with the filtrations.
A subobject $\nbigp_{\ast}\nbige_1\subset\nbigp_{\ast}\nbige$
is a subsheaf $\nbige_1\subset\nbige$
satisfying 
$\nbigp_{\veca}(\nbige_1)\subset
 \nbigp_{\veca}(\nbige)$
for any $\veca\in\real^{\Lambda}$.
It is called strict,
if $\nbigp_{\veca}(\nbige_1)=
 \nbige_1\cap\nbigp_{\veca}(\nbige)$
for any $\veca\in\real^{\Lambda}$.

\subsubsection{Filtered bundles and basic operations}
\label{subsection;13.10.27.1}

A filtered sheaf $\nbigp_{\ast}\nbige$ is called
a filtered bundle,
if (i) $\nbigp_{\veca}\nbige$ are locally free $\nbigo_X$-modules,
(ii) $\lefttop{i}\Gr^{\nbigp}_a(\nbige)$
 are locally free $\nbigo_{D_i}$-modules
for any $i\in\Lambda$
and $a\in\real$.
In that case,
for any $b\in\openclosed{a-1}{a}$,
we set
\[
 F_b\bigl(
 \lefttop{i}\nbigp_a(\nbige)_{|D_i}
 \bigr):=
 \Image\Bigl(
 \lefttop{i}\nbigp_b(\nbige)_{|D_i}
\lrarr 
 \lefttop{i}\nbigp_a(\nbige)_{|D_i}
 \Bigr).
\]
It is called the parabolic filtration.

The direct sum of filtered bundles
$\nbigp_{\ast}\nbige_i$ $(i=1,2)$
is defined as 
the locally free $\nbigo_X(\ast D)$-module
$\nbige_1\oplus\nbige_2$
with the $\nbigo_X$-submodules
$\nbigp_{\veca}(\nbige_1\oplus\nbige_2)
=\nbigp_{\veca}\nbige_1\oplus\nbigp_{\veca}\nbige_2$
$(\veca\in\real^{\Lambda})$.
The tensor product of filtered bundles
$\nbigp_{\ast}\nbige_i$ $(i=1,2)$
is defined as
the $\nbigo_X(\ast D)$-module
$\nbige_1\otimes\nbige_2$
with the $\nbigo_X$-submodules
$\nbigp_{\veca}(\nbige_1\otimes\nbige_2)
=\sum_{\vecb+\vecc\leq \veca}
 \nbigp_{\vecb}(\nbige_1)\otimes
 \nbigp_{\vecc}(\nbige_2)$.
The inner homomorphism 
is defined as
the $\nbigo_X(\ast D)$-module
$\nhom_{\nbigo_X(\ast D)}(\nbige_1,\nbige_2)$
with the $\nbigo_X$-submodules
$\nbigp_{\veca}\nhom(\nbige_1,\nbige_2)=
 \bigl\{
 f\in\nhom(\nbige_1,\nbige_2)\,\big|\,
 f(\nbigp_{\vecb}\nbige_1)\subset
 \nbigp_{\vecb+\veca}\nbige_2 \bigr\}$.

The most typical example is
the $\nbigo_X(\ast D)$-module $\nbigo_X(\ast D)$
with the $\nbigo_X$-submodules
$\nbigp_{\veca}(\nbigo_X(\ast D)):=
 \nbigo(\sum[a_i]D_i)$,
where $[a]:=\max\{n\in\seisuu\,|\,n\leq a\}$.
The filtered bundle is denoted just by
$\nbigo_X(\ast D)$.
For any filtered bundle $\nbigp_{\ast}\nbige$,
the dual
$\nbigp_{\ast}(\nbige^{\lor})$
is defined as 
$\nhom\bigl(\nbigp_{\ast}\nbige,\,
 \nbigo_{X}(\ast D)\bigr)$.
We have a natural isomorphism
$\nbigp_{\veca}(\nbige^{\lor})
\simeq
 \nbigp_{<-\veca+\vecdelta}(\nbige)^{\lor}$,
where $\vecdelta=(1,\ldots,1)$.

Let $\varphi:(X',D')\lrarr (X,D)$ be a ramified covering
with $D'=\coprod_{i\in\Lambda}D'_i$
and $D=\coprod_{i\in\Lambda}D_i$.
Let $e_i$ be the degree of the ramification.
Let $\nbigp_{\ast}\nbige$ be a filtered bundle over $\nbige$.
The pull back of a filtered bundle is defined as
the $\nbigo_{X'}(\ast D')$-module
$\varphi^{\ast}(\nbige)$
with the $\nbigo_{X'}$-submodules
$\nbigp_{\veca}\varphi^{\ast}\nbige
=\sum_{\vece\vecb+\vecn\leq\veca}
 \varphi^{\ast}(\nbigp_{\vecb}\nbige)
\otimes
 \nbigo_{X'}(\sum n_iD_i')$,
where $\vece\vecb=(e_ib_i\,|\,i\in\Lambda)$.
The filtered bundle is denoted by
$\varphi^{\ast}(\nbigp_{\ast}\nbige)$.

Let $\nbigp_{\ast}\nbige'$ be a filtered bundle
on $(X',D')$.
Then, we obtain a locally free
$\nbigo_{X}(\ast D)$-module
$\varphi_{\ast}\nbige'$
with the $\nbigo_X$-submodules
$\nbigp_{\vecc}(\varphi_{\ast}\nbige')$
such that
$\nbigp_{\vecc}(\varphi_{\ast}\nbige')_{|U}
=\varphi_{\ast}(\nbigp_{e_ic_i}\nbige'_{|\varphi^{-1}(U_i)})$.
The filtered bundle is denoted by
$\varphi_{\ast}(\nbigp_{\ast}\nbige)$.
Suppose that $\varphi:(X',D')\lrarr (X,D)$
is a Galois covering with the Galois group
$\Gal(\varphi)$,
and that $\nbigp_{\ast}\nbige'$
be a $\Gal(\varphi)$-equivariant filtered bundle,
then $\varphi_{\ast}(\nbigp_{\ast}\nbige)$
is equipped with an induced
$\Gal(\varphi)$-action.
The $\Gal(\varphi)$-invariant part is called
the descent of  $\nbigp_{\ast}\nbige'$
with respect to $\varphi$.

\subsubsection{The parabolic first Chern class}
\label{subsection;13.10.15.1}

Let $\nbigp_{\ast}\nbige$ be a filtered sheaf on $(X,D)$.
Suppose that $\nbige$ is torsion-free.
The parabolic first Chern class of 
$\nbigp_{\ast}\nbige$ is defined as
\[
 \parchern_1(\nbigp_{\ast}\nbige)=
 c_1(\nbigp_{\veca}\nbige)
-\sum_{i\in\Lambda}
 \sum_{b\in\Par(\nbigp_{\veca}\nbige,i)}
 b \dim\lefttop{i}\Gr^{\nbigp}_b(\nbige)\,[D_i].
\]
Here, $[D_i]$ is the cohomology class of $D_i$.
It is independent of the choice of $\veca$.

\vspace{.1in}

Let $U_i$ be a small neighbourhood of $D_i$.
Suppose that we are given a decomposition
$\nbigp_{\ast}\nbige_{|U_i}=
 \bigoplus_{k\in I(i)}\nbigp_{\ast}\nbige_{i,k}$
for each $i\in\Lambda$.
Let $\nbigu$ be a locally free $\nbigo_X$-submodule
of $\nbige$
such that 
$\nbigu_{|U_i}
=\bigoplus_{k\in I(i)}
 \nbigp_{a(i,k)}\nbige_{i,k}$,
where $a(i,k)\in\real$.
It is easy to check the following equality:
\[
 \parchern_1(\nbigp_{\ast}\nbige)
=c_1(\nbigu)
-\sum_{i\in\Lambda}
 \sum_{k\in I(i)}
 \delta(\nbigp_{\ast}\nbige_{i,k},a(i,k)),
\]
\[
 \delta(\nbigp_{\ast}\nbige_{i,k},a(i,k)):=
 \sum_{b\in\Par(\nbigp_{a(i,k)}\nbige_{i,k})}
 b\,\rank\Gr^{\nbigp}_b(\nbige_{i,k})
[D_i]
\]

\subsubsection{Compatible frames}

For simplicity,
we consider the case
that $X$ is a neighbourhood of $0$ in $\cnum$,
and $D=\{0\}$.
Let $\nbigp_{\ast}\nbige$ be a filtered bundle
on $(X,D)$.
For any section $f$ of $\nbige$,
we set $\deg^{\nbigp}(f):=
 \min\bigl\{
 a\in\real\,\big|\,
 f\in\nbigp_a\nbige
 \bigr\}$.
Let $\vecv=(v_1,\ldots,v_r)$
be a frame of $\nbigp_a\nbige$.
We say that it is compatible with the parabolic structure,
if for any $b\in\Par(\nbigp_a\nbige)$,
the set $\bigl\{v_i\,\big|\,\deg(v_i)=b\bigr\}$
induces a base of
$\Gr^{\nbigp}_b(\nbige)$.

Let $\varphi:(X',D')\lrarr (X,D)$
be a ramified covering given by $\varphi(u)=u^p$.
Let $\nbigp_{\ast}\nbige$ be a filtered bundle on $(X,D)$.
Let $\vecv$ be a compatible frame of $\nbigp_a\nbige$.
Let $c_i:=\deg^{\nbigp}(v_i)$.
We set $n_i:=\max\{n\in\seisuu\,|\,n+pc_i\leq pa\}$,
and 
$w_i:=u^{-n}\varphi^{\ast}v_i$.
Then, $\vecw=(w_1,\ldots,w_r)$
is a compatible frame of
$\varphi^{\ast}(\nbigp_{\ast}\nbige)$,
such that
$\deg^{\nbigp}(w_i)=n_i+pc_i$.

Let $\nbigp_{\ast}\nbige'$ be a filtered bundle on $(X',D')$.
Let $\vecv'$ be a compatible frame of $\nbigp_a\nbige'$.
Let $c_i:=\deg^{\nbigp}(v_i')$.
For $0\leq j<p$,
we set $w'_{ij}:=u^{j}v_i'$.
They naturally induce sections of
$\nbigp_{a/p}(\varphi_{\ast}\nbige)$,
denoted by $\wtilde'_{ij}$.
Then,
$\vecwtilde':=
 \bigl(\wtilde'_{ij}\,\big|\,
 1\leq i\leq\rank\nbige,\,
 0\leq j<p
 \bigr)$
gives a compatible
frame of $\nbigp_{a/p}(\varphi_{\ast}\nbige)$
such that
$\deg^{\nbigp}(\wtilde_{ij}')=(c_i-j)/p$.

\subsubsection{Adapted metric}
\label{subsection;13.10.28.3}

Let us return to the setting in \S\ref{subsection;13.10.28.2}.
Let $V$ be a holomorphic vector bundle 
on $X\setminus D$
with a hermitian metric $h$.
Recall that,
for any $\veca\in\real^{\Lambda}$,
we obtain a natural $\nbigo_{X}$-module
$\nbigp_{\veca}^{(h)}V$
on $X$ as follows.
Let $U$ be any open subset of $X$.
For any  $P\in U$,
we take a holomorphic coordinate neighbourhood
$(X_P,z_1,\ldots,z_n)$ around $P$
such that (i) $X_P$ is relatively compact in $U$,
(i) $X_P\cap D=X_P\cap D_i$ for some $i\in\Lambda$,
(ii) $X_P\cap D=\{z_1=0\}$.
Then, let $\nbigp^{\vech}_{\veca}(U)$ denote the space of
holomorphic sections $f$ of $V_{|U\setminus D}$
such that
$|f_{|X_P}|_h=O(|z_1|^{-a_i-\epsilon})$ for any $\epsilon>0$
and for any $P\in U$.
In general, $\nbigp_{\veca}V$ are not $\nbigo_X$-coherent.

Suppose that we are given a filtered bundle $\nbigp_{\ast}V$
on $(X,D)$,
and that $V:=\nbigp_aV_{|X\setminus D}$
is equipped with a hermitian metric $h$
such that
$\nbigp^{h}_{\veca}V=\nbigp_{\ast}V$.
In that case,
$h$ is called adapted to $\nbigp_{\ast}V$.

\subsection{Filtered Higgs bundles}

Let us recall the notion of filtered Higgs bundles on curves.
Let $X$ be a complex curve with a discrete subset $D$.
Let $\nbigp_{\ast}\nbige$ be a filtered sheaf
on $(X,D)$.
Let $\theta$ be a Higgs field of $\nbige$,
i.e.,
$\theta$ is an $\nbigo_X$-homomorphism
$\nbige\lrarr\nbige\otimes\Omega_X^1$.
Then, $(\nbigp_{\ast}\nbige,\theta)$ is called 
a filtered Higgs bundle.

We shall consider two conditions
on the compatibility of  $\theta$
and the filtration $\nbigp_{\ast}\nbige$.
One is the admissibility,
and the other is goodness.
The latter is what we are really interested in,
because it is closely related with wild harmonic bundles
and $L^2$-instantons.
The former is easier to handle,
and more natural when we consider
algebraic Nahm transforms.
We shall explain the easier one first.

The conditions are given locally around
each point of $D$.
So, we shall explain them
in the case
$X:=\bigl\{
 z\in\cnum\,\big|\,|z|<\rho_0
 \bigr\}$
and $D:=\{0\}$.

\subsubsection{Admissible filtered Higgs bundles}
\label{subsection;13.1.21.200}

For each positive integer $p$,
let $\varphi_p:X^{\langle p\rangle}=\{|z_p|<\rho_0^{1/p}\}\lrarr X$ 
be given by 
$\varphi_p(z_p)=z_p^p$.
Let $\nbigp_{\ast}V$ be a filtered bundle
on $(X,D)$ with a Higgs field $\theta$.
Let $m\in\seisuu_{\geq 0}$ and $p\in\seisuu_{>0}$
such that $\gcd(p,m)=1$.
We say that $(\nbigp_{\ast}V,\theta)$
has slope $(p,m)$,
if the following holds:
\begin{itemize}
\item
Let $(\nbigp_{\ast}V^{\langle p\rangle},\theta^{\langle p\rangle})$
be a filtered Higgs bundle
obtained as the pull back of $(\nbigp_{\ast}V,\theta)$
by $\varphi_p$.
Then, we have
$z_p^{m}\theta^{\langle p\rangle}(\nbigp_cV^{\langle p\rangle})
 \subset\nbigp_cV^{\langle p\rangle}dz_p/z_p$
for any $c\in\real$.
\item
Let $\Res(z_p^m\theta^{\langle p\rangle})$
denote the endomorphism of
$\Gr^{\nbigp}_c(V^{\langle p\rangle})$
obtained as the residue of $z_p^m\theta^{\langle p\rangle}$.
If $(p,m)\neq (1,0)$,
we impose that 
$\Res(z_p^m\theta^{\langle p\rangle})$ is invertible
for any $c$.
\end{itemize}
Although $\Res(z_p^m\theta^{\langle p\rangle})$ may depend
on the choice of a coordinate,
the above condition is independent.
Let $\nbigi(\theta)$ denote the set of the eigenvalues
of $\Res(z_p^m\theta^{\langle p\rangle})$.
We have $\Gal(\varphi_p)$-action
on $(\nbigp_{\ast}V^{\langle p\rangle},\theta^{\langle p\rangle})$ 
and $\nbigi(\theta)$.
The quotient set $\nbigi(\theta)/\Gal(\varphi)$
is denoted by $\vecnbigi(\theta)$.
We have the orbit decomposition
$\nbigi(\theta)=
 \coprod_{\veco\in\vecnbigi(\theta)}\veco$.
We say that $(\nbigp_{\ast}V,\theta)$ has
type $(p,m,\veco)$,
if moreover $\vecnbigi(\theta)=\{\veco\}$.

If $m\neq 0$,
then $\veco$ is an element of 
$\nbigj(p,m):=\cnum^{\ast}/\Gal(\varphi_p)$,
where the action is given by
$(t,\alpha)\longmapsto t^{m}\alpha$.
If $m=0$,
then $\veco$ is an element of $\nbigj(1,0):=\cnum$.
When $(\nbigp_{\ast}V,\theta)$ has slope $(p,m)$,
it has a decomposition 
$(\nbigp_{\ast}V,\theta)
=\bigoplus_{\veco\in\nbigj(p,m)}
 (\nbigp_{\ast}V_{\veco},\theta_{\veco})$
after $X$ is shrank appropriately,
such that
$(\nbigp_{\ast}V_{\veco},\theta_{\veco})$
has type $(p,m,\veco)$.

We say that
$(\nbigp_{\ast}V,\theta)$ is admissible,
if it has a decomposition 
$(\nbigp_{\ast}V,\theta)
=\bigoplus_{(p,m)}(\nbigp_{\ast}V^{(p,m)},\theta^{(p,m)})$
after $X$ is shrank appropriately,
such that
$(\nbigp_{\ast}V^{(p,m)},\theta^{(p,m)})$
has slope $(p,m)$.
In this paper,
we call it the slope decomposition.
We may also have a refined decomposition
$(\nbigp_{\ast}V,\theta)
=\bigoplus_{(p,m,\veco)}
 (\nbigp_{\ast}V_{\veco}^{(p,m)},\theta_{\veco}^{(p,m)})$
such that 
$(\nbigp_{\ast}V_{\veco}^{(p,m)},\theta_{\veco}^{(p,m)})$
has type $(p,m,\veco)$.
In this paper, we call it the type decomposition.
For $\alpha\in\rnum_{\geq 0}$,
we say that the slope of
$(\nbigp_{\ast}V,\theta)$ is smaller 
(resp. strictly smaller) than $\alpha$,
if $\nbigp_{\ast}V^{(p,m)}=0$ for $p/m>\alpha$
(resp. $p/m\geq \alpha$) in the slope decomposition.

Suppose that $(\nbigp_{\ast}V,\theta)$ has type $(p,m,\veco)$.
After $X$ is shrank appropriately,
we have the decomposition
$(\nbigp_{\ast}V^{\langle p\rangle},\theta^{\langle p\rangle})
=\bigoplus_{\alpha\in\veco}
 (\nbigp_{\ast}V^{\langle p\rangle}_{\alpha},
 \theta^{\langle p\rangle}_{\alpha})$
such that 
$\Res(z_p^m\theta^{\langle p\rangle}_{\alpha})$
has a unique eigenvalue $\alpha$.
For any $\alpha\in\veco$,
we have a natural isomorphism
$\varphi_{p\ast}(\nbigp_{\ast}V^{\langle p\rangle}_{\alpha},
 \theta^{\langle p\rangle}_{\alpha})
\simeq
 (\nbigp_{\ast}V,\theta)$.

The following lemma is easy to see.
\begin{lem}
Let $(\nbigp_{\ast}V,\theta)$ be an admissible 
filtered Higgs bundle on $(X,D)$.
Let $\nbigp_{\ast}V'$ be a strict filtered Higgs subbundle,
i.e., it is a strict filtered subbundle such that
$\theta(V')\subset V'\otimes\Omega_X^1$.
The restriction of $\theta$ to $V'$ is denoted by
$\theta'$.
Then, $(\nbigp_{\ast}V',\theta')$ is admissible.
\hfill\qed
\end{lem}

\subsubsection{Good filtered Higgs bundles}
\label{subsection;13.1.23.20}

We have a stronger condition.
Let $X$ and $D$ be as in \S\ref{subsection;13.1.21.200}.
We say that a filtered Higgs bundle 
$(\nbigp_{\ast}V,\theta)$ on $(X,D)$
is good,
if there exists a ramified covering
$\varphi_p:(X^{\langle p\rangle},D^{\langle p\rangle})
 \lrarr (X,D)$
given by $\varphi_p(z_p)=z_p^p$
with a decomposition
\begin{equation}
\label{eq;13.1.22.20}
 \varphi_p^{\ast}(\nbigp_{\ast}V,\theta)
=\bigoplus_{\gminia\in z_p^{-1}\cnum[z_p^{-1}]}
 (\nbigp_{\ast}V^{\langle p\rangle}_{\gminia},
 \theta^{\langle p\rangle}_{\gminia}),
\end{equation}
such that
$\theta^{\langle p\rangle}_{\gminia}
-d\gminia\id_{V^{\langle p\rangle}_{\gminia}}$
is logarithmic
in the sense that
it gives a morphism 
$\nbigp_{\ast}V^{\langle p\rangle}_{\gminia}
\lrarr
 \nbigp_{\ast}V^{\langle p\rangle}_{\gminia}\,
 dz_p/z_p$.
Let $\Irr(\varphi_p^{\ast}\theta)$ denote the set of
$\gminia$
such that $V^{\langle p\rangle}_{\gminia}\neq 0$.
The Galois group 
$\Gal(\varphi_p)$ naturally acts on
$\varphi_p^{\ast}(\nbigp_{\ast}V,\theta)$
and $\Irr(\varphi_p^{\ast}\theta)$.
The quotient set $\Irr(\varphi_p^{\ast}\theta)/\Gal(\varphi_p)$
is denoted by $\vecIrr(\varphi_p^{\ast}\theta)$.
We have the orbit decomposition
$\Irr(\varphi_p^{\ast}\theta)
=\coprod_{\vecgminio\in\vecIrr(\varphi_p^{\ast}\theta)}\vecgminio$.
We set
$(\nbigp_{\ast}V^{\langle p\rangle}_{\vecgminio},
 \theta^{\langle p\rangle}_{\vecgminio}):=
 \bigoplus_{\gminia\in\vecgminio}
 (\nbigp_{\ast}V^{\langle p\rangle}_{\gminia},
 \theta^{\langle p\rangle}_{\gminia})$.
We obtain a $\Gal(\varphi_p)$-equivariant decomposition
$\varphi_p^{\ast}(\nbigp_{\ast}V,\theta)
=\bigoplus_{\vecgminio\in\vecIrr(\varphi_p^{\ast}\theta)}
 (\nbigp_{\ast}V^{\langle p\rangle}_{\vecgminio},
 \theta^{\langle p\rangle}_{\vecgminio})$.
By the descent,
we obtain a decomposition
\begin{equation}
\label{eq;13.1.22.21}
(\nbigp_{\ast}V,\theta)
=\bigoplus_{\vecgminio\in\vecIrr(\varphi^{\ast}\theta)}
 (\nbigp_{\ast}V_{\vecgminio},\theta_{\vecgminio}).
\end{equation}

If we have a factorization
$\varphi_p=\varphi_{p_1}\circ\varphi_{p_2}$
such that
$\varphi_{p_2}^{\ast}(\nbigp_{\ast}V,\theta)$
has a decomposition as above,
$\varphi_{p_1}$ gives a bijection
$\Irr(\varphi_{p_2}^{\ast}\theta)
\simeq
 \Irr(\varphi_{p}^{\ast}\theta)$.
It induces a bijection of the quotient sets
by the Galois groups.
By the identification,
we denote them by
$\Irr(\theta)$ and $\vecIrr(\theta)$.
The decomposition (\ref{eq;13.1.22.21}) is independent
of the choice of $\varphi_p$.

For each $\vecgminio\in\vecIrr(\theta)$,
there exists a minimum $p_{\vecgminio}$
among the numbers $p$ such that
$\varphi_p^{\ast}(\nbigp_{\ast}V_{\vecgminio},\theta_{\vecgminio})$
has a decomposition such as (\ref{eq;13.1.22.20}).
In this case,
we have $|\vecgminio|=p_{\vecgminio}$.
We set 
$X^{\vecgminio}:=X^{\langle p_{\vecgminio}\rangle}$,
$\varphi_{\vecgminio}:=\varphi_{p_{\vecgminio}}$
and $z_{\vecgminio}:=z_{p_{\vecgminio}}$.
We have the decomposition
on $X^{\vecgminio}$:
\begin{equation}
\varphi_{\vecgminio}^{\ast}
 (\nbigp_{\ast}V_{\vecgminio},
 \theta_{\vecgminio})
=\bigoplus_{\gminia\in\vecgminio}
 (\nbigp_{\ast}V^{\vecgminio}_{\gminia},\theta^{\vecgminio}_{\gminia}).
\end{equation}
For any $\gminia\in\vecgminio$,
we have a natural isomorphism
$(\nbigp_{\ast}V_{\vecgminio},\theta_{\vecgminio})
\simeq
 \varphi_{\vecgminio\ast}
 (\nbigp_{\ast}V^{\vecgminio}_{\gminia},\theta^{\vecgminio}_{\gminia})$.
We set 
$m_{\vecgminio}:=(\ord_{z_{\vecgminio}^{-1}}\gminia)$
which is independent of $\gminia\in\vecgminio$.
In this paper,
if $(\nbigp_{\ast}V,\theta)
=(\nbigp_{\ast}V_{\vecgminio},\theta_{\vecgminio})$,
we say that
$(\nbigp_{\ast}V,\theta)$ has pure irregularity $\gminio$.

\vspace{.1in}
If $X$ is shrank appropriately,
we have the following decomposition,
which is a refinement of (\ref{eq;13.1.22.20}),
\[
 \varphi_p^{\ast}(\nbigp_{\ast}V,\theta)
=\bigoplus_{\gminia\in z_p^{-1}\cnum[z_p^{-1}]}
 \bigoplus_{\alpha\in\cnum}
 \bigl(
 \nbigp_{\ast}V^{\langle p\rangle}_{\gminia,\alpha},
 \theta^{\langle p\rangle}_{\gminia,\alpha}
 \bigr)
\]
such that
the eigenvalues of the residues 
$\Res\bigl(\theta^{\langle p\rangle}_{\gminia,\alpha}
-\bigl(d\gminia+p\alpha dz_p/z_p\bigr)
 \id_{V^{\langle p\rangle}_{\gminia,\alpha}}
 \bigr)$
are $0$.
Let 
$(\nbigp_{\ast}V_{\vecgminio,\alpha},\theta_{\vecgminio,\alpha})$
be the descent of 
$\bigoplus_{\gminia\in\vecgminio}
(\nbigp_{\ast}V^{\langle p\rangle}_{\gminia,\alpha},
 \theta^{\langle p\rangle}_{\gminia,\alpha})$
to $X$.
We obtain a decomposition:
\[
 (\nbigp_{\ast}V,\theta)
=\bigoplus_{\vecgminio\in\vecIrr(\theta)}
 \bigoplus_{\alpha\in\cnum}
 (\nbigp_{\ast}V_{\vecgminio,\alpha},
 \theta_{\gminio,\alpha})
\]
On $X^{\vecgminio}$,
we have a decomposition:
\[
 \varphi_{\vecgminio}^{\ast}
 (\nbigp_{\ast}V_{\vecgminio,\alpha},\theta_{\vecgminio,\alpha})
=\bigoplus_{\gminia\in\vecgminio}
  (\nbigp_{\ast}V^{\vecgminio}_{\gminia,\alpha},
 \theta^{\vecgminio}_{\gminia,\alpha})
\]

\begin{lem}
Let $(\nbigp_{\ast}V,\theta)$ be a good filtered Higgs bundle.
Let $\nbigp_{\ast}V'$ be a strict Higgs subbundle.
The restriction of $\theta$ to $V'$ is denoted by
$\theta'$.
Then, $(\nbigp_{\ast}V',\theta')$
is also good.
\end{lem}
\pf
Suppose $(\nbigp_{\ast}V,\theta)$ is unramified
with the decomposition
$(\nbigp_{\ast}V,\theta)
=\bigoplus (\nbigp_{\ast}V_{\gminia},\theta_{\gminia})$.
Because $\theta(V')\subset V'\otimes\Omega_X^1$,
we have
$V'=\bigoplus (V'\cap V_{\gminia})$.
By the strictness,
we obtain 
$\nbigp_{\ast}V'
=\bigoplus (V'\cap\nbigp_{\ast}V_{\gminia})$.
Hence, $(\nbigp_{\ast}V',\theta')$ is good.
The ramified case can be reduced to
the unramified case by the descent.
\hfill\qed

\vspace{.1in}
Take $p\in\seisuu_{>0}$ and $m\in\seisuu_{\geq 0}$
with $\gcd(p,m)=1$.
Let $\vecIrr(\theta,p,m):=
 \{\vecgminio\in\vecIrr(\theta)\,|\,p_{\vecgminio}/m_{\vecgminio}=p/m\}$.
We have
\[
 \nbigp_{\ast}V^{(p,m)}
=\bigoplus_{\vecgminio\in\vecIrr(\theta,p,m)}
 \nbigp_{\ast}V_{\vecgminio}.
\]
For any $\veco\in\nbigj(p,m)$,
we have
$\vecIrr(\theta,p,m,\veco)\subset\vecIrr(\theta,p,m)$
such that 
\[
  \nbigp_{\ast}V^{(p,m)}_{\veco}
=\bigoplus_{\vecgminio\in\vecIrr(\theta,p,m,\veco)}
 \nbigp_{\ast}V_{\vecgminio}.
\]
Take any $\alpha\in\veco$.
For each $\vecgminio\in\vecIrr(\theta,p,m,\veco)$,
we have $\gminia\in\vecgminio$
such that 
\[
 \nbigp_{\ast}V^{\langle p\rangle}_{\alpha}
=\bigoplus_{\vecgminio\in\vecIrr(\theta,p,m,\veco)}
 \varphi_{(p_{\vecgminio}/p)\,\ast}
 \bigl(
 \nbigp_{\ast}V^{\vecgminio}_{\gminia}
 \bigr).
\]
Here, 
$\varphi_{p_{\vecgminio}/p}$
is the ramified covering
$X^{\vecgminio}\lrarr X^{\langle p\rangle}$
given by
$\varphi_{p_{\vecgminio}/p}(z_{\vecgminio})=
 z_{\vecgminio}^{p_{\vecgminio}/p}$.
Let $c\in\real$.
We take a frame $\vecv_{\vecgminio}=(v_{\vecgminio,i})$
of $\nbigp_{p_{\vecgminio}c}
 V^{\vecgminio}_{\gminia}$
compatible with the parabolic structure.
Then, the tuple of the sections
\[
 \Bigl\{
 z_{\vecgminio}^jv_{\vecgminio,i}
\,\Big|\,
  \vecgminio\in\vecIrr(\theta,p,m,\veco),\,\,\,
 1\leq i\leq \rank V^{\vecgminio}_{\gminia},\,\,
 0\leq j< p_{\vecgminio}/p
\Bigr\}
\]
gives a frame of
$\nbigp_{pc}V^{\langle p\rangle}_{\alpha}$.

\subsubsection{Filtered bundles with an endomorphism}
\label{subsection;13.10.15.20}

Let $U_{\tau}$ be a small neighbourhood
of $0$ in $\cnum_{\tau}$.
Let $\nbigp_{\ast}V$ be a filtered bundle
on $(U_{\tau},0)$ with an endomorphism $g$.
We say that $(\nbigp_{\ast}V,g)$
has type $(p,m,\veco)$
(slope $(p,m)$),
if $(\nbigp_{\ast}V,-\tau^{-2}gd\tau)$
has type $(p,m,\veco)$
(resp. slope $(p,m)$).
The condition implies $p\geq m$.
We say that $(\nbigp_{\ast}V,g)$
is admissible,
if $(\nbigp_{\ast}V,-\tau^{-2}gd\tau)$
is admissible.
If $(\nbigp_{\ast}V,g)$ is admissible,
we have the type decomposition 
$(\nbigp_{\ast}V,g)
=\bigoplus(\nbigp_{\ast}V^{(p,m)}_{\veco},g^{(p,m)}_{\veco})$
and the slope decomposition
$(\nbigp_{\ast}V,g)
=\bigoplus(\nbigp_{\ast}V^{(p,m)},g^{(p,m)})$,
after $X$ is shrank appropriately.

Similarly, 
$(\nbigp_{\ast}V,g)$ is called good,
if $(\nbigp_{\ast}V,-\tau^{-2}gd\tau)$
is a good filtered Higgs bundle.

\begin{rem}
We regard $U_{\tau}$ as a neighbourhood of
$\infty$ in $\proj^1$.
Of course,
$-\tau^{-2}d\tau=dw$ for $w=\tau^{-1}$.
\hfill\qed
\end{rem}

\begin{rem}
We shall be interested in the case
that $(\nbigp_{\ast}V,g)$
is decomposed into
$\bigoplus_{\alpha\in\cnum}
 (\nbigp_{\ast}V_{\alpha},g_{\alpha})$
such that 
(i) $\Sp(g_{\alpha|0})=\{\alpha\}$,
(ii) $(\nbigp_{\ast}V_{\alpha},g_{\alpha}-\alpha)$
is admissible.
In that case,
in the slope decomposition
$(\nbigp_{\ast}V_{\alpha},g_{\alpha}-\alpha)
=\bigoplus\bigl(
 \nbigp_{\ast}V_{\alpha}^{(p,m)},g_{\alpha}^{(p,m)}
 \bigr)$,
we have $m/p<1$
for $V_{\alpha}^{(p,m)}\neq 0$.
\hfill\qed
\end{rem}

\subsection{Filtered bundles on $(T\times \proj^1,T\times\{\infty\})$}

\subsubsection{Local conditions}
\label{subsection;13.2.14.11}

Let $U\subset\proj^1$ be a small neighbourhood of $\infty$.
We introduce some conditions on filtered bundles
$\nbigp_{\ast}E$ on $(T\times U ,T\times\{\infty\})$.
\begin{description}
\item[(A1)]
 $\nbigp_cE_{|T\times \infty}$
 are semistable of degree $0$
 for any $c\in\real$.
\end{description}
The condition is equivalent to
that $\Gr^{\nbigp}_c(E)$ are semistable of degree $0$
for any $c\in\real$.
Let $\Sp_{\infty}(E)\subset T^{\lor}$
denote the spectrum of $\nbigp_cE_{|T\times\infty}$.
It is independent of $c$.
We fix its lift to 
$\Sptilde_{\infty}(E)\subset \cnum$.
Then, as observed in \S\ref{subsection;13.10.15.30},
for a small neighbourhood $U'$ of $\infty\in\proj^1$,
we obtain the corresponding filtered bundle
$\nbigp_{\ast}V$ with an endomorphism $g$
on $(U',\infty)$
such that
$\Sp(g_{|\infty})=\Sptilde_{\infty}(E)$.
We have the decomposition
\[
 (\nbigp_{\ast}V,g)
=\bigoplus_{P\in\Sp_{\infty}(E)}
 (\nbigp_{\ast}V_P,g_P)
\]
such that $\Sp(g_P)\cap(\cnum\times\{\infty\})=\{\Ptilde\}$
is the lift of $P$.
A filtered bundle satisfying {\bf(A1)}
is called admissible,
if it satisfies the following condition,
which is independent of the choice of
$\Sptilde_{\infty}(E)$.
\begin{description}
\item[(A2)]
$(\nbigp_{\ast}V_P,g_P-\Ptilde\,\id)$
is admissible in the sense of 
\S\ref{subsection;13.10.15.20}
for any $P\in\Sp_{\infty}(E)$.
\end{description}

We have the type decomposition
$(\nbigp_{\ast}V_P,g_P-\Ptilde\id)
=\bigoplus_{p,m,\veco}
 (\nbigp_{\ast}V^{(p,m)}_{P,\veco},g^{(p,m)}_{P,\veco})$.
We have the corresponding decomposition
$\nbigp_{\ast}E=\bigoplus_P
 \bigoplus_{p,m,\veco}
 \nbigp_{\ast}E_{P,\veco}^{(p,m)}$,
which is called the type decomposition
of $\nbigp_{\ast}E$.
The following lemma is clear.
\begin{lem}
If $\nbigp_{\ast}E$ satisfies the condition
{\bf (A1)}
(resp. the admissibility),
then the dual $\nbigp_{\ast}(E^{\lor})$
also satisfies the condition {\bf (A1)}
(resp. the admissibility).
\hfill\qed
\end{lem}

We also have the following condition.
\begin{description}
\item[(Good)]
Let $\nbigp_{\ast}E$ be a filtered bundle
on $(T\times U,T\times\{\infty\})$
satisfying {\bf (A1)}.
Take any lift $\Sptilde_{\infty}(E)\subset\cnum$
of $\Sp_{\infty}(E)$.
Then, the filtered bundle is called good, 
if the corresponding filtered bundle
$\nbigp_{\ast}V$
with an endomorphism $g$
 is good in the sense of \S\ref{subsection;13.10.15.20}.
\end{description}

\subsubsection{Some remarks on the cohomology}
\label{subsection;13.10.15.2}

Let $\nbigp_{\ast}E$ be a filtered bundle on
$(T\times\proj^1,T\times\{\infty\})$
satisfying {\bf (A1)}.
Let $\nbigu$ be any $\nbigo_{T\times\proj^1}$-submodule
of $\nbigp_cE$ for some $c\in\real$,
such that
(i) $\nbigu_{|T\times\cnum_w}
=\nbigp_cE_{|T\times\cnum_w}$,
(ii) $\nbigu_{|T\times\{\infty\}}$
 is semistable of degree $0$.
We give some remarks
on the cohomology groups of $\nbigu$.

\begin{lem}
\label{lem;13.1.21.30}
Suppose 
$0\not\in\Sp_{\infty}(\nbigp_{\ast}E)$.
Then, we have 
$H^j(T\times\proj^1,\nbigu)
=H^j(T\times\proj^1,\nbigp_cE)$.
\end{lem}
\pf
Let $\pi:T\times\proj^1\lrarr\proj^1$ be the projection.
By the assumption,
we obtain
$R\pi_{\ast}(\nbigu\otimes L)
 \simeq
 R\pi_{\ast}(\nbigp_cE\otimes L)$,
because both of them vanish around $\infty$.
Then, the claim of the lemma follows.
\hfill\qed

\vspace{.1in}
Suppose that $\nbigp_{\ast}E$ is admissible,
i.e., it satisfies {\bf (A1,2)},
and we give some refinement.
We have the decomposition 
$\nbigu=\bigoplus_P\bigoplus_{p,m,\veco}
 \nbigu^{(p,m)}_{P,\veco}$
around $T\times\{\infty\}$,
where 
$\nbigu^{(p,m)}_{P,\veco}:=
 \nbigu\cap\nbigp_c E^{(p,m)}_{P,\veco}$.
Let $\nbigu'\subset\nbigp_cE$ be a subsheaf
satisfying the above conditions (i,ii).
If $0\not\in\Sp_{\infty}(\nbigp_{\ast}E)$,
we have
$H^i(T\times\proj^1,\nbigu)
=H^i(T\times\proj^1,\nbigu')$
by Lemma \ref{lem;13.1.21.30}.

\begin{lem}
\label{lem;13.10.15.3}
Suppose $0\in \Sp_{\infty}(\nbigp_{\ast}E)$,
and that
$\nbigu^{(1,0)}_{0,0}
=\nbigu^{\prime(1,0)}_{0,0}$.
Then, 
we have natural isomorphisms
$H^i(T\times\proj^1,
 \nbigu)
\simeq
 H^i(T\times\proj^1,\nbigu')$
for $i=0,2$.
\end{lem}
\pf
We have only to consider the case that
$\nbigu\subset\nbigu'$,
and we shall prove that
the natural morphisms
$H^i(T\times\proj^1,\nbigu)
\lrarr
 H^i(T\times\proj^1,\nbigu')$
are isomorphisms.
Let $\varphi\in H^0(T\times\proj^1,\nbigu')$.
Around $T\times\{\infty\}$,
we have the decomposition
$\varphi=\sum_{P,p,m,\veco}\varphi_{P,\veco}^{(p,m)}$.
It is easy to observe that
$\varphi_{P,\veco}^{(p,m)}=0$
unless
$(P,p,m,\veco)=(0,1,0,0)$.
Hence, we obtain that
$H^0(T\times\proj^1,\nbigu)
\lrarr
 H^0(T\times\proj^1,\nbigu')$
is an isomorphism.
The duals $\nbigu$ and $(\nbigu^{\prime})^{\lor}$
are subsheaves of 
$\nbigp_{c'}(E^{\lor})$ for some $c'$,
and satisfy the conditions (i,ii).
Hence, by using the Serre duality,
we obtain that
$H^2(T\times\proj^1,\nbigu)
\lrarr
 H^2(T\times\proj^1,\nbigu')$
is an isomorphism.
\hfill\qed

\subsubsection{Vanishing condition}

Let $(\nbigp_{\ast}E,\theta)$
be a filtered bundle on 
$(T\times\proj^1,T\times\{0\})$.
We will be concerned with 
the following condition
on vanishing of the cohomology:
\begin{description}
\item[(A3)]
 $H^0(T\times \proj^1,\nbigp_0E\otimes p^{\ast}L)=0$
 and 
 $H^2(T\times\proj^1,\nbigp_{<-1}E\otimes p^{\ast}L)=0$
 for any line bundle $L$ of degree $0$ on $T$,
 where $p$ denotes the projection
 $T\times\proj^1\lrarr T$.
\end{description}
We shall often omit to denote $p^{\ast}$,
if there is no risk of confusion.

\begin{lem}
If $\nbigp_{\ast}E$ satisfies the condition
{\bf (A3)},
the dual $\nbigp_{\ast}(E^{\lor})$
also satisfies the condition {\bf (A3)}.
\end{lem}
\pf
Note 
$\nbigp_{a}(E^{\lor})^{\lor}
 \otimes\Omega^1_{\proj^1}
\simeq
 \nbigp_{<-a-1}(E)$.
Hence, by the Serre duality,
$H^0\bigl(T\times \proj^1,\nbigp_0(E^{\lor})\otimes L^{\lor}\bigr)$
is the dual space of
$H^2\bigl(T\times \proj^1,\nbigp_{<-1}E\otimes L\bigr)$,
and 
$H^2\bigl(T\times \proj^1,\nbigp_{<-1}(E^{\lor})\otimes L^{\lor}\bigr)$
is the dual space of
$H^0\bigl(T\times \proj^1,\nbigp_0E\otimes L\bigr)$.
The claim of the lemma follows.
\hfill\qed

\begin{lem}
Let $L$ be any holomorphic line bundle on $T$
of degree $0$.
If $\nbigp_{\ast}E$ satisfies {\bf (A3)},
we have 
$H^0(T\times\proj^1,\nbigp_cE\otimes L)=0$ for any $c\leq 0$
and 
$H^2(T\times\proj^1,\nbigp_{<c}E\otimes L)=0$ for any $c\geq -1$.
\end{lem}
\pf
We have only to consider the case $L=\nbigo_T$.
For $c\leq 0$,
we have
$H^0(T\times\proj^1,\nbigp_cE)
\subset
 H^0(T\times\proj^1,\nbigp_0E)=0$.
For $c\geq -1$,
the support of the quotient
$\nbigp_{<c}E/\nbigp_{<-1}E$
is one dimensional.
Hence, the morphism
$0=H^2(T\times\proj^1,\nbigp_{<-1}E)
\lrarr
 H^2(T\times\proj^1,\nbigp_{<c}E)$
is surjective.
\hfill\qed

\subsubsection{Stability condition}
\label{subsection;13.10.15.200}
We introduce a stability condition for filtered bundles
satisfying {\bf (A1)}
on $(T\times\proj^1,T\times\{\infty\})$,
by following \cite{Biquard-Jardim}.
Note that this is not the same as 
the standard slope stability condition
for filtered bundles
on projective varieties
\cite{Maruyama-Yokogawa}.

Let $\omega_{T}\in H^2(T\times\proj^1,\seisuu)$ 
denote the pull back of the fundamental class
of $T$ by the projection $T\times\proj^1\lrarr T$.
For any filtered torsion-free sheaf $\nbigp_{\ast}\nbige$
on $(T\times\proj^1,T\times\{\infty\})$,
we define the degree of $\nbigp_{\ast}\nbige$
by 
\[
 \deg(\nbigp_{\ast}\nbige):=
 \int_{T\times\proj^1}
 \parchern_1(\nbigp_{\ast}\nbige)
 \omega_T
=\int_{\{z\}\times\proj^1}
 \parchern_1(\nbigp_{\ast}\nbige).
\]
We set
$\mu(\nbigp_{\ast}\nbige):=
 \deg(\nbigp_{\ast}\nbige)/\rank\nbige$.
We say that a filtered bundle $\nbigp_{\ast}E$ is stable
(semistable),
if $\mu(\nbigp_{\ast}\nbige)<\mu(\nbigp_{\ast}E)$
(resp.
$\mu(\nbigp_{\ast}\nbige)\leq \mu(\nbigp_{\ast}E)$)
for any $\nbigp_{\ast}\nbige\subset\nbigp_{\ast}E$
such that 
(i) $0<\rank\nbige<\rank E$,
(ii) $\nbigp_{\ast}\nbige$ also satisfies {\bf (A1)}
around $T\times\{\infty\}$.
We say that a semistable filtered bundle 
$\nbigp_{\ast}E$ is polystable,
if it has a decomposition
$\nbigp_{\ast}E=\bigoplus\nbigp_{\ast}E_i$
such that each $\nbigp_{\ast}E_i$ is stable.
The following lemma is clear and standard.
\begin{lem}
Let $\nbigp_{\ast}E$ be a filtered bundle 
satisfying {\bf (A1)}
on $(T\times\proj^1,T\times\{\infty\})$.
If $\nbigp_{\ast}E$ is stable,
then $\nbigp_{\ast}E^{\lor}$ is also stable.
\hfill\qed
\end{lem}

It is standard
to obtain the vanishing of some cohomology groups
under the assumption of the stability and the degree $0$.

\begin{lem}
Let $\nbigp_{\ast}E$ be a filtered bundle 
satisfying {\bf (A1)} on
$(T\times\proj^1,T\times\{\infty\})$.
If $\nbigp_{\ast}E$ is stable with
$\deg(\nbigp_{\ast}E)=0$
and $\rank\nbigp_{\ast}E>1$,
it satisfies the condition {\bf (A3)}.
\end{lem}
\pf
Because $\nbigp_{\ast}E$ is stable of degree $0$,
we obtain that 
$H^0(T\times\proj^1,\nbigp_cE)=0$
for any $c\leq 0$.
Indeed, a non-zero section
induces a filtered strict subsheaf
$\nbigp_{\ast}\nbigo\subset\nbigp_{\ast}E$
with $\deg(\nbigp_{\ast}\nbigo)\geq 0$
and $0<\rank\nbigo<\rank E$.
Because $(\nbigp_{\ast}E)^{\lor}$
is also stable of degree $0$,
we obtain the vanishing
$H^2(T\times\proj^1,\nbigp_{<-1}E)$
by using the Serre duality.
\hfill\qed

\vspace{.1in}
\begin{rem}
\label{rem;13.10.20.1}
Filtered bundles
$\nbigp_{\ast}E$ satisfying {\bf (A1)}
with $\deg(\nbigp_{\ast}E)=0$
and $\rank\nbigp_{\ast}E=1$
naturally correspond to 
line bundles of degree $0$ on $T$.
Indeed,
there exists a line bundle of degree $0$
on $L$ on $T$ such that
$\nbigp_aE\simeq
 p^{\ast}L\otimes
 \nbigo\bigl([a](T\times\{\infty\})\bigr)$
for any $a\in\real$,
where
 $[a]:=\max\{n\in\seisuu\,|\,n\leq a\}$.
In this case,
{\bf (A3)} is not satisfied for $L^{-1}$.
\hfill\qed
\end{rem}

\section{Algebraic Nahm transforms}
\label{section;13.1.30.10}

\subsection{Local algebraic Nahm transforms}
\label{subsection;13.2.14.22}
\subsubsection{Complex}
\label{subsection;13.1.21.1}

Let $X:=\bigl\{
 z\in\cnum\,\big|\,|z|<\rho_0
 \bigr\}$
and $D:=\{0\}$.
In the following of this subsection,
we shall shrink $X$ without mentioning it.
We shall use the notation in \S\ref{subsection;13.1.21.200}.
We define a complex of sheaves associated to 
an admissible filtered Higgs bundle
$(\nbigp_{\ast}V,\theta)$ on $(X,D)$.
First, let us consider the case that
$(\nbigp_{\ast}V,\theta)$ has type $(p,m,\veco)$.
Suppose $(p,m,\veco)\neq (1,0,0)$.
For each $c\in\real$,
let $\nbigp_c(V\otimes\Omega_X^{\bullet},\theta)$
denote the complex 
\[
 \nbigp_{c-m/p}V
\lrarr
 \nbigp_{c+1}V\,dz
\]
where the first term sits in the degree $0$.
Take any $\alpha\in\veco$.
For each $c\in\real$,
let $\nbigp_c(V^{\langle p\rangle}_{\alpha}
 \otimes\Omega_{X^{\langle p\rangle}}^{\bullet},
 \theta^{\langle p\rangle}_{\alpha})$
denote the following complex 
on $X^{\langle p\rangle}$:
\[ 
\begin{CD}
 \nbigp_{c-m}V^{\langle p\rangle}_{\alpha}
@>{\theta^{\langle p\rangle}_{\alpha}}>>
 \nbigp_{c}V^{\langle p\rangle}_{\alpha}
 \otimes
 \frac{dz_p}{z_p}.
\end{CD}
\]
We have a natural isomorphism
$\nbigp_{c}(V\otimes\Omega^{\bullet},\theta)
\simeq
 \varphi_{p\ast}
 \nbigp_{cp}(V^{\langle p \rangle}_{\alpha}\otimes
 \Omega^1_{X^{\langle p\rangle}},\theta^{\langle p\rangle}_{\alpha})$.
It is also isomorphic to the descent of
$\bigoplus_{\alpha\in\veco}
\nbigp_{cp}(V^{\langle p\rangle}_{\alpha}
 \otimes\Omega^1_{X^{\langle p\rangle}},
 \theta^{\langle p\rangle}_{\alpha})$.
For $c\leq c'$,
the natural inclusion
$\nbigp_c(V\otimes\Omega^{\bullet},
 \theta)
\lrarr
 \nbigp_{c'}(V\otimes\Omega^{\bullet}, \theta)$
is a quasi-isomorphism.
We set
$\nbigc^{\bullet}(\nbigp_{\ast}V,\theta):=
 \nbigp_{-1/2}\bigl(V\otimes\Omega^{\bullet},
 \theta\bigr)$.
In the case $(p,m,\veco)=(1,0,0)$,
we set
$\nbigc^0(\nbigp_{\ast}V,\theta):=\nbigp_0V$
and 
\[
\nbigc^1(\nbigp_{\ast}V,\theta):=
 \nbigp_{<1}V
\otimes\Omega^1_X
+\theta\bigl(
 \nbigp_{0}V
 \bigr)
\subset
 \nbigp_1V\otimes\Omega^1_X.
\]
Thus, we obtain the complex
$\nbigc^{\bullet}(\nbigp_{\ast}V,\theta)$,
when $(\nbigp_{\ast}V,\theta)$ has type $(p,m,\veco)$.

For a general admissible filtered Higgs bundle
$(\nbigp_{\ast}V,\theta)$,
the complex 
$\nbigc^{\bullet}(\nbigp_{\ast}V,\theta)$
is defined as the extension of
the complex $(V\lrarr V\otimes\Omega_X^1)$
on $X\setminus D$
to a complex on $X$,
such that it is
$\bigoplus_{(m,p,\veco)}
 \nbigc^{\bullet}(\nbigp_{\ast}V^{(p,m)}_{\veco},\theta^{(p,m)}_{\veco})$
around $D$,
according to the type decomposition.

\begin{lem}
\label{lem;13.1.26.10}
If $(\nbigp_{\ast}V,\theta)$ comes from
a wild harmonic bundle
$(E,\delbar_E,\theta,h)$  on $(X,D)$,
then 
$\nbigc^{\bullet}
 (\nbigp_{\ast}V,\theta)$
is naturally quasi-isomorphic to
the complex of square-integrable sections of
the Higgs complex $E\otimes\Omega^{\bullet}$.
\end{lem}
\pf
By the descent, we have only to consider the unramified case.
We omit to denote $p$.
We have the naturally defined map
$\pi_{c}:\nbigp_cV\lrarr\Gr^{\nbigp}_c(V)$.
Let $W$ be the weight filtration of the nilpotent part of
the endomorphism
$\Res(\theta)$ on $\Gr^{\nbigp}_c(V)$.
We set 
$W_k\nbigp_cV:=
 \pi_c^{-1}(W_k\Gr^{\nbigp}_c(V))$.

Suppose that
$(\nbigp_{\ast}V,\theta)$
has type $(m,\veco)$.
If $(m,\veco)\neq (0,0)$,
let 
$\nbigc^{\bullet}_{L^2}(
 \nbigp_{\ast}V, \theta)$
be the following complex:
\[
 W_{-2}\nbigp_{-m}V
\lrarr
 W_{-2}\nbigp_0V\otimes\Omega^1_X(\log D)
\]
It is easy to check that the following natural morphisms
are quasi-isomorphisms:
\[
 \nbigc^{\bullet}_{L^2}(\nbigp_{\ast}V,\theta)
\lrarr
 \nbigp_0(V\otimes\Omega^{\bullet},\theta)
\llarr
 \nbigc^{\bullet}(\nbigp_{\ast}V,\theta)
\]
If $(m,\veco)=(0,0)$,
let 
$\nbigc^{\bullet}_{L^2}(\nbigp_{\ast}V,\theta)$
be the following complex:
\[
 W_{0}\nbigp_{0}V
\lrarr
 W_{-2}\nbigp_0V\otimes\Omega^1_X(\log D)
\]
It is easy to check that the natural inclusion
$\nbigc^{\bullet}_{L^2}(\nbigp_{\ast}V,\theta)
\lrarr
 \nbigc^{\bullet}(\nbigp_{\ast}V,\theta)$
is a quasi-isomorphism.

In general, we define
$\nbigc^{\bullet}_{L^2}(\nbigp_{\ast}V,\theta)
=\bigoplus
 \nbigc^{\bullet}_{L^2}(\nbigp_{\ast}V^{(m)}_{\veco},
 \theta^{(m)}_{\veco})$
by using the type decomposition.
According to the result in \S5.1 of \cite{Mochizuki-wild},
$\nbigc^{\bullet}_{L^2}(\nbigp_{\ast}V,\theta)$
is naturally quasi-isomorphic to
the complex of square-integrable sections
of the Higgs complex $E\otimes\Omega^{\bullet}$.
Thus, we obtain the claim of the lemma.
\hfill\qed

\subsubsection{Transform}
\label{subsection;13.1.20.20}

We shall construct some transformations for filtered Higgs bundles,
which are analogue to the local Fourier transform 
in \cite{Bloch-Esnault} and \cite{Garcia-Lopez}.

In the following,
for a variable $x$,
let $U_x$ denote a small neighbourhood of $0$
in $\cnum_x$.
For two variables $x$ and $y$,
let $U_{x,y}:=U_x\times U_y$,
and let $\pi_1:U_{x,y}\lrarr U_x$
and $\pi_2:U_{x,y}\lrarr U_y$ denote the projections.

Let $(\nbigp_{\ast}V,\theta)$ be an admissible filtered Higgs bundle
on $(U_{\zeta},0)$.
Let us define a filtered bundle
$\nbign_{\ast}^{0,\infty}(\nbigp_{\ast}V,\theta)$
with an endomorphism $g$
on $U_{\tau}$.
We consider the following complex
on $U_{\zeta,\tau}$:
\[
\begin{CD}
 \pi_1^{\ast}
 \nbigc^0(\nbigp_{\ast}V,\theta)
@>{\tau\theta+d\zeta}>>
 \pi_1^{\ast}
 \nbigc^1(\nbigp_{\ast}V,\theta)
\end{CD}
\]
Let $\nbigq$ be the quotient.
We define
\[
 \vecnbign^{0,\infty}(\nbigp_{\ast}V,\theta):=
 \pi_{2\ast}\nbigq,
\quad\quad\quad
 \nbign^{0,\infty}(\nbigp_{\ast}V,\theta):=
 \vecnbign^{0,\infty}(\nbigp_{\ast}V,\theta)(\ast\tau).
\]
Here, $(\ast \tau)$ means the localization
with respect to $\tau$.
If $U_{\tau}$ is sufficiently small,
then the support of $\nbigq$ is proper 
and relatively $0$-dimensional over $U_{\tau}$.
Indeed,
$\nbigq\cap\bigl(\{0\}\times U_{\zeta}\bigr)
=\{(0,0)\}$.
Hence,
$\vecnbign^{0,\infty}(\nbigp_{\ast}V,\theta)$ 
is coherent.
Let us check that
$\vecnbign^{0,\infty}(\nbigp_{\ast}V,\theta)$
is torsion-free.
Let $v$ be a section of
$\pi_1^{\ast}\nbigc^1(\nbigp_{\ast}V,\theta)$,
such that there exists a section $u$
of $\pi_1^{\ast}\nbigc^0(\nbigp_{\ast}V,\theta)$
satisfying 
$\tau v=(\tau\theta+d\zeta)u$.
We obtain that
$d\zeta \cdot u$ is contained in
$\tau\cdot \pi_1^{\ast}\nbigc^1(\nbigp_{\ast}V,\theta)$.
Then, we obtain that 
$u=\tau u'$
for some section $u'$ of
$\pi_1^{\ast}\nbigc^0(\nbigp_{\ast}V,\theta)$,
and we have
$v=(\tau\theta+d\zeta)u'$.
It implies that
$\vecnbign^{0,\infty}(\nbigp_{\ast}V,\theta)$
is torsion-free.
Hence,
$\vecnbign^{0,\infty}(\nbigp_{\ast}V,\theta)$
is a locally free $\nbigo_{U_{\tau}}$-module.
In particular,
$\nbign^{0,\infty}(\nbigp_{\ast}V,\theta)$
is a locally free $\nbigo_{U_{\tau}}(\ast\tau)$-module.
The multiplication of $\zeta$ induces the endomorphism $g$.
By setting $\psi:=-g\,\tau^{-2}d\tau$,
we obtain a Higgs field of
$\nbign^{0,\infty}(\nbigp_{\ast}V,\theta)$.
We shall introduce
a filtered bundle 
$\nbign_{\ast}^{0,\infty}(\nbigp_{\ast}V,\theta)
=\bigl(
 \nbign_a^{0,\infty}(\nbigp_{\ast}V,\theta)\,\big|\,
 a\in\real
 \bigr)$
over $\nbign^{0,\infty}(\nbigp_{\ast}V,\theta)$.

\vspace{.1in}

If $(\nbigp_{\ast}V,\theta)$ has type $(p,m,\veco)\neq (1,0,0)$,
we consider the following complexes on 
$U_{\zeta,\tau}$ for any $c\in\real$:
\begin{equation}
 \label{eq;13.1.29.20}
\begin{CD}
 \pi_1^{\ast}\nbigp_{c-m/p}(V)
 @>{\tau\theta+d\zeta}>>
 \pi_1^{\ast}\nbigp_c(V)(d\zeta/\zeta)
\end{CD}
\end{equation}
Let $\nbigq_c$ denote the quotient.
We define
\[
 \nbign^{0,\infty}_{\kappa_1(p,m,c)}(\nbigp_{\ast}V,\theta)
:=\pi_{2\ast}\nbigq_c,
\quad\quad
 \kappa_1(p,m,c):=\frac{2pc-m}{2(p+m)}.
\]
By construction,
we have 
$\nbign_{-1/2}^{0,\infty}(\nbigp_{\ast}V,\theta)
=\vecnbign(\nbigp_{\ast}V,\theta)$
in this case.
It is easy to check that 
$\nbign_a^{0,\infty}(\nbigp_{\ast}V,\theta)$
are locally free $\nbigo_{U_{\tau}}$-module
of finite rank.
We have a naturally induced map
$\nbign_{a'}^{0,\infty}(\nbigp_{\ast}V,\theta)
\lrarr
 \nbign_{a}^{0,\infty}(\nbigp_{\ast}V,\theta)$
for $a'\leq a$.
Its restriction to $\{\tau\neq 0\}$
is an isomorphism,
and hence it is injective.
We also obtain
$\nbign_a^{0,\infty}(\nbigp_{\ast}V,\theta)(\ast \tau)
=\nbign^{0,\infty}(\nbigp_{\ast}V,\theta)$.
For $c':=c-(1+m/p)$,
the images of
$\tau\cdot\pi_1^{\ast}\nbigp_{c}V(d\zeta/\zeta)$
and 
$\pi_1^{\ast}\nbigp_{c'}V(d\zeta/\zeta)$
are the same in the quotient of $\nbigq_c$.
It implies
$\tau\nbign_a^{0,\infty}(\nbigp_{\ast}V,\theta)
=\nbign_{a-1}^{0,\infty}(\nbigp_{\ast}V,\theta)$
for any $a\in\real$.
Hence,
$\nbign_{a}^{0,\infty}(\nbigp_{\ast}V,\theta)$ $(a\in\real)$
give a filtered bundle over
$\nbign^{0,\infty}(\nbigp_{\ast}V,\theta)$.

\vspace{.1in}
If $(\nbigp_{\ast}V,\theta)$ has type 
$(p,m,\veco)=(1,0,0)$,
we define
$\nbign_0^{0,\infty}(\nbigp_{\ast}V,\theta):=
 \vecnbign^{0,\infty}(\nbigp_{\ast}V,\theta)$.
We have the following natural morphisms:
\[
 \nbign_0^{0,\infty}(\nbigp_{\ast}V,\theta)_{|0}
\simeq
 \nbigc^1(\nbigp_{\ast}V,\theta)
 \big/
 \nbigc^0(\nbigp_{\ast}V,\theta)d\zeta
\lrarr
 (\nbigp_0V)_{|0}
\]
Here, the subscript ``$|0$'' means the fiber of
the vector bundle over $0$,
and the latter map is given by the residue,
which is injective.
Hence, the parabolic filtration of the right hand side induces
a parabolic filtration of 
$\nbign_0^{0,\infty}(\nbigp_{\ast}V,\theta)_{|0}$
indexed by $\openclosed{-1}{0}$.
It induces a filtered bundle
$\nbign_{\ast}^{0,\infty}(\nbigp_{\ast}V,\theta)$
over $\nbign^{0,\infty}(\nbigp_{\ast}V,\theta)$.

If $(\nbigp_{\ast}V,\theta)$ is admissible,
we replace $U_{\zeta}$ with smaller neighbourhoods
so that it has the type decomposition,
and we define
\[
 \nbign_{\ast}^{0,\infty}(\nbigp_{\ast}V,\theta)
:=\bigoplus_{p,m,\veco}
 \nbign_{\ast}^{0,\infty}(\nbigp_{\ast}V_{\veco}^{(p,m)},\theta).
\]
The construction $\nbign_{\ast}^{0,\infty}$
gives a functor from the category of 
the germs of admissible filtered Higgs bundles
to the category of the germs of filtered Higgs bundles.
We set
$\nbign^{0,\infty}_{<a}(\nbigp_{\ast}V,\theta):=
 \sum_{b<a}
 \nbign^{0,\infty}_b(\nbigp_{\ast}V,\theta)$.

\begin{lem}
Suppose $(\nbigp_{\ast}V,\theta)$ has type
$(p,m,\veco)$.
The rank of $\nbign^{0,\infty}(\nbigp_{\ast}V,\theta)$
is $(p+m)\rank V/p$ in the case $(p,m,\veco)\neq(1,0,0)$,
or 
$\rank V-\dim\Ker\Gr^{\nbigp}_0(\Res\theta)$
in the case $(p,m,\veco)=(1,0,0)$.
\end{lem}
\pf
The rank is equal to
the dimension of
$\nbigc^1(\nbigp_{\ast}V,\theta)
 \big/
 \nbigc^0(\nbigp_{\ast}V,\theta)d\zeta$
as a $\cnum$-vector space.
Then, the claim can be checked by a direct computation.
(See also the proof of Proposition \ref{prop;13.1.20.10} below
for the case 
$(p,m,\veco)\neq(1,0,0)$.)

\hfill\qed

\begin{prop}
\label{prop;13.1.20.10}
$\bigl(
 \nbign_{\ast}^{0,\infty}(\nbigp_{\ast}V,\theta),\psi
 \bigr)$
is admissible.
If $(\nbigp_{\ast}V,\theta)$ has type $(p,m,\veco)$,
then 
$\bigl(
 \nbign_{\ast}^{0,\infty}(\nbigp_{\ast}V,\theta),\psi
 \bigr)$
has type $(p+m,m,\veco')$
for some $\veco'$.
\end{prop}
\pf
We have only to consider the case that
$(\nbigp_{\ast}V,\theta)$ has type $(p,m,\veco)$.
Let us consider the case $(p,m,\veco)=(1,0,0)$.
For the expression
$\theta=f\,d\zeta/\zeta$,
$f$ gives an endomorphism of $\nbigp_cV$ for any $c$,
and $f_{|0}$ is nilpotent.
We have
$\psi=-\tau^{-1}g(d\tau/\tau)$,
and $-\tau^{-1}g$ is induced by $f$
by construction.
Hence, it preserves 
$\nbign_0^{0,\infty}(\nbigp_{\ast}V,\theta)$.
If we regard
$\nbign_0^{0,\infty}(\nbigp_{\ast}V,\theta)_{|0}$
as a subspace of
$\nbigp_0V_{|0}$ as above,
$(-\tau^{-1}g)_{|0}$
is the restriction of $f_{|0}$.
Hence, it is nilpotent,
and 
 it preserves the parabolic filtration,
i.e.,
$(\nbign^{0,\infty}(\nbigp_{\ast}V,\theta),\psi)$
is admissible of type $(1,0,0)$.

Let us consider the case
$(p,m,\veco)\neq (1,0,0)$.
Fix $\alpha\in\veco$.
We consider the following
on $U_{\zeta_p,\tau}$:
\begin{equation}
 \label{eq;13.1.20.1}
\begin{CD}
 \pi_1^{\ast}\nbigp_{pc-m}V^{\langle p\rangle}_{\alpha}
 @>{\tau \theta^{\langle p\rangle}_{\alpha}+d\zeta_p^p}>>
 \pi_1^{\ast}\nbigp_{pc}V^{\langle p\rangle}_{\alpha}(d\zeta_p/\zeta_p)
\end{CD}
\end{equation}
The quotient is denoted by $\nbigq'_c$.
We have a natural isomorphism
$\pi_{2\ast}\nbigq'_c
\simeq
 \nbign^{0,\infty}_{\kappa_1(p,m,c)}(\nbigp_{\ast}V,\theta)$.

The natural map
$\nbigq'_{c'}\lrarr\nbigq'_{c}$ $(c'\leq c)$
is injective.
We set
$\nbigq'_{<c}:=\bigcup_{b<c}\nbigq'_b$.
We have the following exact sequence:
\[
\begin{CD}
0@>>>
 \pi_1^{\ast}\Gr^{\nbigp}_{pc-m}(V^{\langle p\rangle}_{\alpha})
@>{\tau\theta^{\langle p\rangle}_{\alpha}}>>
 \pi_1^{\ast}\Gr^{\nbigp}_{pc}(V^{\langle p\rangle}_{\alpha})
 (d\zeta_p/\zeta_p)
@>>>
 \nbigq'_{c}/\nbigq'_{<c}
@>>>0
\end{CD}
\]
It induces the following isomorphism
of $\cnum$-vector spaces for any $c\in\real$:
\begin{equation}
 \label{eq;13.2.12.1}
 \Gr^{\nbigp}_{pc}(V^{\langle p\rangle}_{\alpha})
\simeq
\frac{
 \nbign_{\kappa_1(p,m,c)}^{0,\infty}(\nbigp_{\ast}V,\theta)}
 {\nbign_{<\kappa_1(p,m,c)}^{0,\infty}(\nbigp_{\ast}V,\theta)}
\end{equation}

Let $\vecv=(v_i)$ be a frame of
$\nbigp_{pc}V^{\langle p\rangle}_{\alpha}$.
We set 
$c_i:=\min\{a\in\real|\,
 v_i\in\nbigp_aV^{\langle p\rangle}_{\alpha}\}$.
We assume that $\vecv$ is compatible with the parabolic structure
in the sense that the induced tuple 
$\{[v_i]\,|\,c_i=d\}$ of elements in 
$\Gr^{\nbigp}_d(V^{\langle p\rangle}_{\alpha})$
is a base for any $d\in\openclosed{pc-1}{pc}$.
We set
$\nu_{ij}:=\zeta_p^{i}v_j\,(d\zeta_p/\zeta_p)$
for
$(0\leq i\leq p+m-1,\,\,\,
 1\leq j\leq \rank V^{\langle p\rangle}_{\alpha})$.
The induced sections of
$\nbign^{0,\infty}_{\kappa_1(p,m,c)}(\nbigp_{\ast}V,\theta)$
are also denoted by the same symbols.
Because they induce a base of
$\nbign^{0,\infty}_{\kappa_1(p,m,c)}(\nbigp_{\ast}V,\theta)
 \big/\tau\nbign^{0,\infty}_{\kappa_1(p,m,c)}(\nbigp_{\ast}V,\theta)$,
they give a frame
of 
$\nbign^{0,\infty}_{\kappa_1(p,m,c)}(\nbigp_{\ast}V,\theta)$
on a neighbourhood of $0$.
(In particular,
the rank of
$\nbign^{0,\infty}_{\kappa_1(p,m,c)}(\nbigp_{\ast}V,\theta)$
is $(p+m)\rank V^{\langle p\rangle}_{\alpha}=p^{-1}(p+m)\rank V$.)
Moreover, 
by the isomorphism (\ref{eq;13.2.12.1}),
the frame is compatible with the parabolic structure
of $\nbign^{0,\infty}_{\kappa_1(p,m,c)}(\nbigp_{\ast}V,\theta)$.

We take a ramified covering
$\varphi:U_{\eta}\lrarr U_{\tau}$
by $\varphi(\eta)=\eta^{p+m}$.
Let $\nbigp_{\ast}\nbigv$
be the filtered bundle on $(U_{\eta},0)$
obtained as the pull back of 
$\nbign_{\ast}^{0,\infty}(\nbigp_{\ast}V,\theta)$
by $\varphi$.
The tuple of the sections
$\nutilde_{ij}:=\eta^{-i}\varphi^{\ast}\nu_{ij}$
gives a frame $\vecnutilde$ of
$\nbigp_{pc-m/2}\nbigv$
which is compatible with the parabolic structure.
By the frames
$\vecv$ and $\vecnutilde$,
we obtain an isomorphism of 
$\nbigp_{pc-m/2}\nbigv_{|0}$
to 
$(\nbigp_{pc}V^{\langle p\rangle}_{\alpha})_{|0}
 \otimes\cnum^{p+m}$.

Let us show that
$\psi=-\tau^{-2}gd\tau$ has type $(p+m,m,\veco')$ for some $\veco'$.
Note that $g$ is induced by the multiplication of 
$\zeta=\zeta_p^p$.
Let $g_1$ be the endomorphism of
$\nbign^{0,\infty}(\nbigp_{\ast}V,\theta)$
induced by the multiplication of $\zeta_p$.
We have
$g_1\bigl(
 \nbign^{0,\infty}_a(\nbigp_{\ast}V,\theta)
 \bigr)
\subset
 \nbign^{0,\infty}_{a-1/(p+m)}(\nbigp_{\ast}V,\theta)$.
Hence, we obtain
$\eta^{-1}g_1$ gives an endomorphism
of $\nbigp_{pc-m/2}\nbigv$.
In particular,
$\eta^{-p}g$ gives an endomorphism of
$\nbigp_{pc-m/2}\nbigv$.
Let us show that the restriction $(\eta^{-p}g)_{|0}$ 
has a unique non-zero eigenvalue
modulo the action of $\Gal(\varphi)$.

We have the parabolic filtration $F$ of
$(\nbigp_{pc}V_{\alpha})_{|0}$,
indexed by
$\openclosed{pc-1}{pc}$.
Let $W$ denote the monodromy weight filtration
of the nilpotent part of
$\Res(\zeta_p^m\theta^{\langle p\rangle}_{\alpha})$
on $\Gr^{F}(\nbigp_{pc}V_{\alpha|0})$.
Let $\pi_a:F_a\bigl(\nbigp_{pc}V_{\alpha|0}\bigr)
\lrarr\Gr^F_a(\nbigp_{pc}V_{\alpha|0})$
denote the projection.
Let $M:=\min
 \bigl\{|a-b|\,\big|\,a,b\in\Par(\nbigp_{pc}V_{\alpha}),
 a\neq b\bigr\}$.
We take a small positive number $\delta$
such that 
$\delta\rank \nbigp_{pc}V_{\alpha}<M/100$.
We set
$\Ftilde_{a+\delta k}:=\pi_a^{-1}(W_k)$.
Then, we obtain a filtration
$\Ftilde$ of
$\nbigp_{pc}V_{\alpha|0}$
indexed by
$\openclosed{pc-1+\epsilon}{pc+\epsilon}$
for some small $\epsilon>0$.
Then, 
$\Ftilde$ is preserved by 
$\Res(\zeta_p^m\theta^{\langle p\rangle}_{\alpha})$,
and  the induced endomorphism on 
the associated graded space $\Gr^{\Ftilde}$
is semisimple.
We may assume that the frame $\vecv$
is compatible with $\Ftilde$.

Let $\Ftilde'$ be a filtration of
$\nbigp_{pc-m/2}\nbigv_{|0}
\simeq
 \nbigp_{pc}V^{\langle p\rangle}_{\alpha|0}
 \otimes\cnum^{p+m}$
indexed by 
$\openclosed{pc-1+\delta}{pc+\delta}$,
determined by the condition
$\deg^{\Ftilde'}(\nutilde_{ij})=
 \deg^{\Ftilde}(v_j)$.

The multiplication of $\eta^{-1}\zeta_p$
induces an endomorphism of
$\nbigp_{pc-m/2}\nbigv$.
We have
$(\eta^{-1}\zeta_p)
 \nutilde_{ij}
=\nutilde_{i+1,j}$
for $i<p+m-1$,
and 
$(\eta^{-1}\zeta_p)\nutilde_{p+m-1,j}$
is equal to the section $s$
induced by
\[
-p^{-1}\theta^{\langle p\rangle}_{\alpha}
 (\zeta_p^mv_j)
=\Bigl(
 -(\alpha/p)v_j
+\sum_{\deg^{\Ftilde}(v_k)<\deg^{\Ftilde}(v_j)}
 \gamma_k\cdot v_k
 +\zeta_p u
 \Bigr)(d\zeta_p/\zeta_p).
\]
Here, $\gamma_k$ are complex numbers,
and $u$ is a section of
$\nbigp_{pc}V^{\langle p\rangle}_{\alpha}$.
If $\deg^{\Ftilde}(v_{j|0})=a$,
then
$s_{|0}+(\alpha/p)\nutilde_{0,j|0}
 \in\Ftilde'_{<a}$.

The endomorphism $\eta^{-1}g$
of $\nbigp_{pc-m/2}\nbigv$
is induced by the multiplication of
the $p$-th power of $\eta^{-1}\zeta_p$.
Hence, $(\eta^{-p}g)_{|0}$
is compatible with $\Ftilde'$,
and the induced endomorphism
on $\Gr^{\Ftilde'}$ is represented by
the following matrix:
\[
 \sum_{i=1}^m I\otimes E_{p+i,i}
+\sum_{i=1}^p (-\alpha/p) I\otimes E_{i,m+i}
\]
Here, $I$ is the identity matrix 
and $E_{ij}$ denote the $(p+m)$-square matrix
whose $(k,\ell)$-entry is $1$ if $(k,\ell)=(i,j)$,
and $0$ otherwise.
Then, the set of the eigenvalues
are
$e^{2\pi\sqrt{-1}j/(p+m)}\alpha^p$
$(j=0,\ldots,p+m-1)$.
Thus, we are done.
\hfill\qed

\begin{cor}
\label{cor;13.10.13.2}
The construction
$\nbign_{\ast}^{0,\infty}$ gives a functor
from the category of the germs of
admissible filtered Higgs bundles
to the category of the germs of
admissible filtered Higgs bundles
whose slopes are strictly less than $1$.
\hfill\qed
\end{cor}

\subsubsection{Inverse transform}
\label{subsection;13.1.21.20}

Let $\nbigp_{\ast}V$ be a filtered bundle
on $(U_{\tau},0)$ with an endomorphism $g$,
which is admissible in the sense of
\S\ref{subsection;13.10.15.20}.
In this subsection,
we impose the following vanishing:
\begin{description}
\item[(C0)]
$V^{(1,0)}_{0}=0$,
and $V^{(p,m)}=0$ unless $p>m$.
\end{description}

Note that
the eigenvalues of $g(\tau)$ goes to $0$
when $\tau\to 0$
under the assumption {\bf (C0)}.

If $(\nbigp_{\ast}V,g)$ has slope $(p,m)$,
we consider the following complex
on $U_{\tau,\zeta}$:
\[
\begin{CD}
 \pi_1^{\ast}
 \nbigp_cV
 @>{g-\zeta}>>
 \pi_1^{\ast}
 \nbigp_cV
\end{CD}
\]
The quotient is denoted by $\nbigm_c$.
If $U_{\zeta}$ is sufficiently small,
the support of $\nbigm_c$ is proper over $U_{\zeta}$.
We define
\[
\nbign_{\kappa_2(p,m,c)+1}^{\infty,0}(\nbigp_{\ast}V,g)
:=
 \pi_{2\ast}\nbigm_c,
\quad\quad
 \kappa_2(p,m,c):=\frac{2pc+m}{2(p-m)}.
\]
They are locally free $\nbigo_{U_{\zeta}}$-modules.
For $a\leq a'$,
we have a naturally defined map
$\nbign_a^{\infty,0}(\nbigp_{\ast}V,g)
\lrarr
 \nbign_{a'}^{\infty,0}(\nbigp_{\ast}V,g)$
which induces
$\nbign_a^{\infty,0}(\nbigp_{\ast}V,g)(\ast \zeta)
\simeq
 \nbign_{a'}^{\infty,0}(\nbigp_{\ast}V,g)(\ast \zeta)$.
We have
$\nbign_{a-1}^{\infty,0}(\nbigp_{\ast}V,g)
=\zeta \nbign_{a}^{\infty,0}(\nbigp_{\ast}V,g)$
for any $a\in\real$.
Thus, we obtain a filtered bundle
$\nbign_{\ast}^{\infty,0}(\nbigp_{\ast}V,g)$
on $(U_{\zeta},0)$.
In the general case,
we define
$\nbign_{\ast}^{\infty,0}(\nbigp_{\ast}V,g):=
 \bigoplus
 \nbign_{\ast}^{\infty,0}(\nbigp_{\ast}V_{\veco}^{(p,m)},
 g^{(p,m)}_{\veco})$
by using the slope decomposition of $(\nbigp_{\ast}V,g)$.
The multiplication of $-\tau^{-1}$
gives a meromorphic endomorphism $f$.
We put $\theta=fd\zeta$.
The construction gives a functor
from the category of the germs of
admissible filtered Higgs bundles
satisfying {\bf(C0)}
to the category of germs of filtered Higgs bundles.

\begin{prop}
\label{prop;13.1.21.110}
$\bigl(
 \nbign_{\ast}^{\infty,0}(\nbigp_{\ast}V,g),\theta
 \bigr)$
is admissible.
If $(\nbigp_{\ast}V,g)$ has type $(p,m,\veco)$,
then $\nbign_{\ast}^{\infty,0}(\nbigp_{\ast}V,g)$
has type $(p-m,m,\veco')$
for some $\veco'$,
and the rank is $(p-m)\rank V/p$.
\end{prop}
\pf
We have only to consider the case that
$(\nbigp_{\ast}V,g)$ has type $(p,m,\veco)$.
Let $\varphi_p:U_{\eta}\lrarr U_{\tau}$
be given by $\varphi_p(\eta)=\eta^p$.
Let $\varphi:U_u\lrarr U_{\zeta}$ be given by
$\varphi(u)=u^{p-m}$.
Let $\nbigp_{\ast}\nbigv$
be the filtered bundle
on $U_{u}$
obtained as the pull back of 
$\nbign^{\infty,0}_{\ast}(\nbigp_{\ast}V,g)$
by $\varphi$.

We use the decomposition
$\varphi_p^{\ast}(\nbigp_{\ast}V,g)
=\bigoplus_{\alpha\in\veco}
 (\nbigp_{\ast}V^{\langle p\rangle}_{\alpha},g^{\langle p\rangle}_{\alpha})$.
We consider the following complex on $U_{\eta,\zeta}$:
\[
\begin{CD}
 \pi_1^{\ast}\nbigp_{pc}V^{\langle p\rangle}_{\alpha}
 @>{g^{\langle p\rangle}_{\alpha}-\zeta}>>
 \pi_1^{\ast}\nbigp_{pc}V^{\langle p\rangle}_{\alpha}
\end{CD}
\]
The quotient is denoted by $\nbigm'_c$.
We have 
$\pi_{2\ast}\nbigm'_c
\simeq
 \nbign^{\infty,0}_{\kappa_2(p,m,c)+1}(\nbigp_{\ast}V,g)$.
Because 
$g^{\langle p\rangle}_{\alpha}(\nbigp_aV^{\langle p\rangle}_{\alpha})
\subset
 \nbigp_aV^{\langle p\rangle}_{<\alpha}$,
we have the following exact sequence,
as in the case of $\nbign^{0,\infty}$
(see the proof of Proposition \ref{prop;13.1.20.10}):
\[
 \begin{CD}
 0@>>>
 \pi_1^{\ast}\Gr^{\nbigp}_{pc}(V^{\langle p\rangle}_{\alpha})
 @>{-\zeta}>>
 \pi_1^{\ast}\Gr^{\nbigp}_{pc}(V^{\langle p\rangle}_{\alpha})
 @>>>
 \nbigm'_c/\nbigm'_{<c}
 @>>>
 0
 \end{CD}
\]
It induces the following isomorphism
of $\cnum$-vector spaces:
\begin{equation}
\label{eq;13.2.12.2}
 \Gr^{\nbigp}_{pc}(V^{\langle p\rangle}_{\alpha})
\simeq
 \frac{\nbign^{\infty,0}_{\kappa_2(p,m,c)+1}(\nbigp_{\ast}V,g)}
 {\nbign^{\infty,0}_{<\kappa_2(p,m,c)+1}(\nbigp_{\ast}V,g)}
\end{equation}

We take a frame $\vecv$
of $\nbigp_{pc}V^{\langle p\rangle}_{\alpha}$
compatible with the parabolic structure.
We set $\nu_{ij}:=\eta^iv_j$.
By the isomorphism (\ref{eq;13.2.12.2}),
they induce a frame of 
$\nbign^{\infty,0}_{\kappa_2(p,m,c)+1}(\nbigp_{\ast}V,g)$
compatible with the parabolic structure.
We set
$\nutilde_{ij}:=u^{-i}\eta^iv_j$.
The tuple $\vecnutilde$ induces a frame 
of $\nbigp_{p(c+1)-m/2}\nbigv$
compatible with the parabolic structure.

We consider the endomorphism
$h:=\eta^{m-p}g^{\langle p\rangle}_{\alpha}$
on $\nbigp_{pc}V^{\langle p\rangle}_{\alpha}$,
which is invertible.
We have $\eta^{-p+m}u^{p-m}=h$ on $\nbigv$.
Let $k$ be the integer determined by the condition
$0\leq -p+k(p-m)<p-m$.
We set $a:=-p+k(p-m)$.
We have
$\eta^{-p}u^p=\eta^au^{-a}h^k=
 \eta^{a-(p-m)}u^{-a+(p-m)}h^{k-1}$.
We have
\[
 u^p\eta^{-p}
 \nu_{ij}
=\left\{
\begin{array}{ll}
\eta^{a+i}u^{-(a+i)}h^k(v_j)
 & (a+i<p-m)\\
 \mbox{{}}\\
\eta^{a+i-(p-m)}u^{-(a+i)+p-m}h^{k-1}(v_j)
 & (a+i\geq p-m)
\end{array}
\right.
\]
Hence,
$u^{p}\eta^{-p}$
preserves
$\nbigp_{p(c+1)-m/2}\nbigv$.

By the frames $\vecv$ and $\vecnutilde$,
we have an isomorphism
$\nbigp_{p(c+1)-m/2}\nbigv_{|0}$
and 
$\nbigp_{pc}V_{\alpha|0}\otimes \cnum^{p-m}$.
We take a refinement $\Ftilde$
of the parabolic filtration of $\nbigp_{pc}V_{\alpha|0}$
such that
(i) $\Ftilde$ is preserved by $h_{|0}$,
(ii) the induced endomorphism on
 $\Gr^{\Ftilde}$ is semisimple
 with a unique eigenvalue $\beta$.
It induces a filtration $\Ftilde'$ of
$\nbigp_{p(c+1)-m/2}\nbigv_{|0}$.
(See the proof of Proposition \ref{prop;13.1.20.10}
for a concrete construction.)
On $\Gr^{\Ftilde'}$,
$u^p\eta^{-p}$ is expressed by
the matrix
\[
 \sum E_{a+i,i}\otimes \beta^k I
+\sum E_{i,i+p-m-a}\otimes \beta^{k-1}I
\]
with respect to an appropriate base.
Then, we obtain that
$(\nbigp_{\ast}
 \nbign^{\infty,0}(\nbigp_{\ast}V,g),\theta)$
has type $(p-m,m,\veco')$
for some $\veco'$.
\hfill\qed

\begin{cor}
The construction $\nbign^{\infty,0}_{\ast}$
gives a functor from the category of
the germs of admissible filtered Higgs bundles 
satisfying {\bf(C0)}
to the category of the germs of 
admissible filtered Higgs bundles.
\hfill\qed
\end{cor}

\vspace{.1in}

We denote
$\bigl(\nbign^{0,\infty}_{\ast}(\nbigp_{\ast}V,\theta),g\bigr)$
in \S\ref{subsection;13.1.20.20}
by
$\nbign^{0,\infty}_{\ast}(\nbigp_{\ast}V,\theta)$
for simplicity.
We also denote 
$\bigl(\nbign^{\infty,0}_{\ast}(\nbigp_{\ast}V,g),\theta\bigr)$
by 
$\nbign^{\infty,0}_{\ast}(\nbigp_{\ast}V,g)$.
\begin{prop}
\label{prop;13.1.21.100}
\mbox{{}}
\begin{itemize}
\item
Suppose that $(\nbigp_{\ast}V,\theta)$ is admissible
such that $V^{(1,0)}_0=0$
in the type decomposition.
Then, we have a natural isomorphism
of germs of filtered Higgs bundles
$\nbign_{\ast}^{\infty,0}
 \nbign_{\ast}^{0,\infty}(\nbigp_{\ast}V,\theta)
\simeq
(\nbigp_{\ast}V,\theta)$.
\item
Suppose that 
 $(\nbigp_{\ast}V,g)$ is admissible 
and satisfies the condition {\bf(C0)}.
Then, we have a natural isomorphism
of germs of filtered bundles
with endomorphisms
$\nbign_{\ast}^{0,\infty}\nbign_{\ast}^{\infty,0}
 (\nbigp_{\ast}V,g)
\simeq
 (\nbigp_{\ast}V,g)$.
\end{itemize}
\end{prop}
\pf
Suppose that
$(\nbigp_{\ast}V,\theta)$ has type $(p,m)$.
Note that, if we set  $d:=\kappa_1(p,m,c)$,
then we have $\kappa_2(p+m,m,d)=c$.
Let $p_i$ be the projection of
$U_{\zeta}\times U_{\tau}\times U_{\zeta'}$
onto the $i$-th component.
We have the following diagram on
$U_{\zeta}\times U_{\tau}\times U_{\zeta'}$:
\[
 \begin{CD}
 p_1^{\ast}\nbigp_{c-m/p}(V)
 @>{\tau\theta+d\zeta}>>
 p_1^{\ast}\nbigp_c(V)d\zeta/\zeta \\
 @V{\zeta-\zeta'}VV @V{\zeta-\zeta'}VV \\
 p_1^{\ast}\nbigp_{c-m/p}(V)
 @>{\tau\theta+d\zeta}>>
 p_{1}^{\ast}\nbigp_c(V)\,d\zeta/\zeta
 \end{CD}
\]
We regard it as a double complex,
where the left upper $p_1^{\ast}\nbigp_{c-m/p}(V)$
sits in the degree $(0,0)$.
Let $C^{\bullet}$ denote the associated total complex.
By the construction,
$\nbign_{c+1}^{\infty,0}
 \nbign_{\ast}^{0,\infty}
 (\nbigp_{\ast}V,\theta)$
is obtained as 
$p_{3\ast}\nbigh^2(C^{\bullet})$.
We can observe that 
it is isomorphic to the push-forward of
$\nbigq_c$ in \S\ref{subsection;13.1.20.20}
by the projection
$U_{\zeta,\tau}\lrarr U_{\zeta}$,
which is naturally isomorphic to
$\nbigp_cV\,d\zeta/\zeta
\simeq
 \nbigp_{c+1}V$.
The action of $-\tau^{-1}$ is equal to $f$
for the expression $\theta=f\,d\zeta$.
Hence, we obtain the desired isomorphism
$\nbign_{\ast}^{\infty,0}
 \nbign_{\ast}^{0,\infty}(\nbigp_{\ast}V,\theta)
\simeq
(\nbigp_{\ast}V,\theta)$.

Suppose that $(\nbigp_{\ast}V,g)$ has type $(p,m)$
with $p>m$.
Let $p_i$ denote the projection of
$U_{\tau}\times U_{\zeta}\times U_{\tau'}$
onto the $i$-th component.
We have the following commutative diagram
of the sheaves on $U_{\tau}\times U_{\zeta}\times U_{\tau'}$:
\[
 \begin{CD}
 p_1^{\ast}\nbigp_{c-1}V
 @>{g-\zeta}>>
 p_1^{\ast}\nbigp_{c-1}V \\
 @V{(-\tau'(\tau^{-1})+1)d\zeta}VV 
 @VV{(-\tau'(\tau^{-1})+1)d\zeta}V \\
 p_1^{\ast}\nbigp_cV\,d\zeta
 @>{g-\zeta}>>
 p_1^{\ast}\nbigp_c V\,d\zeta\\
 \end{CD}
\]
We regard it as the double complex,
where the left upper $ p_1^{\ast}\nbigp_{c-1}V$
sits in the degree $(0,0)$.
Let $C^{\bullet}$ denote the associated total complex.
By the construction,
$\nbign_{c}^{0,\infty}\nbign^{\infty,0}_{\ast}(\nbigp_{\ast}V,g)$
is naturally isomorphic to $p_{3\ast}\nbigh^2(C^{\bullet})$.
We can observe that it is naturally isomorphic to
the push-forward of $\nbigm_c$ in \S\ref{subsection;13.1.21.20}
by the projection $U_{\tau}\times U_{\zeta}\lrarr U_{\tau}$,
which is naturally isomorphic to $\nbigp_cV$.
The action of $\zeta$ is given by $g$.
Hence, we obtain the desired isomorphism
$\nbign_{\ast}^{0,\infty}\nbign_{\ast}^{\infty,0}
 (\nbigp_{\ast}V,g)
\simeq
 (\nbigp_{\ast}V,g)$.
\hfill\qed

\subsubsection{Description of the functors}

Let $(\nbigp_{\ast}V,\theta)$ be a filtered Higgs bundle
with slope $(p,m)\neq (1,0)$ on $U_{\zeta}$.
Suppose that there exist a ramified covering
$\varphi_q:U_{\zeta_q}\lrarr U_{\zeta}$
and a filtered Higgs bundle 
$(\nbigp_{\ast}V',\theta')$ on $U_{\zeta_q}$
with an isomorphism
$\varphi_{q\ast}(\nbigp_{\ast}V',\theta')
\simeq (\nbigp_{\ast}V,\theta)$.
For $c\in\real$,
we consider the following morphism
on $U_{\zeta_q,\tau}$:
\[
\begin{CD}
 \nbigp_{q(c-m/p)}V'
@>{\tau\theta'+d\zeta_q^q}>>
 \nbigp_{qc}V'\,d\zeta_q/\zeta_q
\end{CD}
\]
The quotient is denoted by $\nbigq'_c$.
The following lemma is clear by construction.
\begin{lem}
$\pi_{2\ast}\nbigq'_c$
is naturally isomorphic to 
$\nbign^{0,\infty}_{\kappa_1(p,m,c)}(\nbigp_{\ast}V,\theta)$.
\hfill\qed
\end{lem}

Let $(\nbigp_{\ast}V,\psi)$ be a filtered Higgs bundle
with slope $(p,m)$ on $U_{\tau}$,
such that $(p,m)\neq (1,0)$ and $p>m$.
Suppose that there exist a ramified covering
$\varphi_q:U_{\tau_q}\lrarr U_{\tau}$
and a filtered Higgs bundle 
$(\nbigp_{\ast}V',\psi')$ on $U_{\tau_q}$
with an isomorphism
$\varphi_{q\ast}(\nbigp_{\ast}V',\psi')
\simeq(\nbigp_{\ast}V,\psi)$.
Let $\psi'=g'\,\varphi^{\ast}(-\tau^{-2}d\tau)$.
For $c\in\real$,
we consider the following morphism
on $U_{\tau_q,\zeta}$:
\[
\begin{CD}
 \nbigp_{qc}V'
@>{g'-\zeta}>>
 \nbigp_{qc}V'
\end{CD}
\]
Let $\nbigm'_c$ denote the quotient.
The following is clear by construction.
\begin{lem}
$\pi_{2\ast}\nbigm'_c$ is naturally isomorphic to
$\nbign^{\infty,0}_{\kappa_2(p,m,c)+1}(\nbigp_{\ast}V,g)$.
\hfill\qed
\end{lem}

\subsection{Algebraic Nahm transform for admissible filtered Higgs
  bundle}
\label{subsection;13.2.14.25}
\subsubsection{Construction of the transform}
\label{subsection;13.2.14.10}

Let $T^{\lor}:=\cnum/L^{\lor}$.
Let $D\subset T^{\lor}$ be an effective reduced divisor.
Let $(\nbigp_{\ast}\nbige,\theta)$
be a filtered Higgs bundle on $(T^{\lor},D)$.
Suppose that it is admissible around
each point of $D$
in the sense of \S\ref{subsection;13.1.21.200}.
We shall construct a filtered bundle
$\Nahm_{\ast}(\nbigp_{\ast}\nbige,\theta)$
on $(T\times\proj^1,T\times\{\infty\})$.
We begin with a construction of an object 
$\vecN(\nbigp_{\ast}\nbige,\theta)$ in 
$D^b(\nbigo_{T\times \proj^1})$.

For $I\subset\{1,2,3\}$,
let $p_I$ be the projections of
$T^{\lor}\times T\times\proj^1$
onto the product of the $i$-th components $(i\in I)$.
Let $\poincare$ be the Poincar\'e bundle on
$T^{\lor}\times T$.
Applying the construction
in \S\ref{subsection;13.1.21.1} around each point of $D$,
we extend 
$\nbige$ and $\nbige\otimes\Omega_X^1$
on $X\setminus D$
to $\nbigc^0(\nbigp_{\ast}\nbige,\theta)$
and $\nbigc^1(\nbigp_{\ast}\nbige,\theta)$,
respectively.
We set
\[
 \nbigctilde^0(\nbigp_{\ast}\nbige,\theta):=
 p_1^{\ast}\nbigc^0(\nbigp_{\ast}\nbige,\theta)
 \otimes
 p_{12}^{\ast}\poincare\otimes
 p_3^{\ast}\nbigo_{\proj^1}(-1),
\quad\quad
 \nbigctilde^1(\nbigp_{\ast}\nbige,\theta):=
 p_1^{\ast}\nbigc^1(\nbigp_{\ast}\nbige,\theta)
 \otimes
 p_{12}^{\ast}\poincare.
\]
Let $\zeta$ be the standard coordinate of $\cnum$,
which induces local coordinates of $T^{\lor}$.
We have the holomorphic $1$-form $d\zeta$ on $T^{\lor}$.
Let $w$ be the standard coordinate of $\cnum\subset\proj^1$,
which we can naturally regard as a section of
$\nbigo_{\proj^1}(1)$.
Then, we have the following morphism:
\begin{equation}
\label{eq;13.2.13.1}
 \theta+w\,d\zeta:
 \nbigctilde^0(\nbigp_{\ast}\nbige,\theta)
\lrarr
 \nbigctilde^1(\nbigp_{\ast}\nbige,\theta).
\end{equation}
Thus, we obtain a complex
$\nbigctilde^{\bullet}(\nbigp_{\ast}\nbige,\theta)$
on $T^{\lor}\times T\times\proj^1$.
We define
\[
 \vecN(\nbigp_{\ast}\nbige,\theta):=
 Rp_{23\ast}
 \bigl(
 \nbigctilde^{\bullet}(\nbigp_{\ast}\nbige,\theta)
 \bigr)[1].
\]

\begin{lem}
\label{lem;13.2.13.3}
There exists a neighbourhood $U$ of $\infty$ in $\proj^1$,
such that
$\nbigh^i\bigl(
 \vecN(\nbigp_{\ast}\nbige,\theta)
 \bigr)_{|T\times U}=0$
unless $i\neq 0$.
Moreover,
$\nbigh^0\bigl(
 \vecN(\nbigp_{\ast}\nbige,\theta)
 \bigr) _{|T\times\{P\}}$
are semistable bundles of degree $0$
for any $P\in U$.
\end{lem}
\pf
Let $\pi_i$ denote the projection of
$T^{\lor}\times\proj^1$ onto the $i$-th component.
We have the following complex
$\nbigctilde^{\bullet}_1(\nbigp_{\ast}\nbige,\theta)$
on $T^{\lor}\times \proj^1$:
\[
\begin{CD}
 \pi_1^{\ast}\nbigc^0(\nbigp_{\ast}\nbige,\theta)
 \otimes
 \pi_2^{\ast}\nbigo_{\proj^1}(-1)
@>{\theta+w\,d\zeta}>>
 \pi_1^{\ast}\nbigc^1(\nbigp_{\ast}\nbige,\theta)
\end{CD}
\]
By the construction,
$\vecN(\nbigp_{\ast}\nbige,\theta)$
is isomorphic to
$\widehat{\RFM}_+\bigl(
 \nbigctilde^{\bullet}_1(\nbigp_{\ast}\nbige,\theta)\bigr)[1]$.
If $U$ is sufficiently small,
$\theta+w\,d\zeta$ is injective
on $T^{\lor}\times U$,
and the support of the cokernel is
relatively $0$-dimensional over $U$.
Then, the claim of the lemma follows.
\hfill\qed

\vspace{.1in}

We consider the following vanishing condition.
\begin{description}
\item[(A0)]
$\hyperh^i\bigl(
 T^{\lor},
 \nbigc^{\bullet}(\nbigp_{\ast}\nbige\otimes L,
 \theta+wd\zeta)
 \bigr)=0$
unless $i=1$,
for any $w\in\cnum$
and any holomorphic line bundle $L$ 
of degree $0$ on $T^{\lor}$.

\end{description}
Under the assumption {\bf (A0)},
we naturally identify 
$\vecN(\nbigp_{\ast}\nbige,\theta)$
with the $0$-th cohomology sheaf
$\nbigh^0\bigl(\vecN(\nbigp_{\ast}\nbige,\theta)\bigr)$,
which is a locally free sheaf on $T\times\proj^1$.
Indeed,
$\nbigc^{\bullet}(\nbigp_{\ast}\nbige\otimes L,
 \theta+wd\zeta)$
is naturally identified with
the specialization of
$\nbigctilde^{\bullet}(\nbigp_{\ast}\nbige,\theta)$
to $T^{\lor}\times\{(L,w)\}$.
Note that we always have
$\hyperh^i\bigl(
 T^{\lor},
 \nbigc^{\bullet}(\nbigp_{\ast}\nbige\otimes L,
 d\zeta)
 \bigr)=0$
unless $i=1$,
for any $L$,
which corresponds to the specialization at $w=\infty$.
We define
\[
 \Nahm(\nbigp_{\ast}\nbige,\theta):=
 \vecN(\nbigp_{\ast}\nbige,\theta)
 \otimes
 \nbigo_{T\times\proj^1}\bigl(
 \ast (T\times\{\infty\})
 \bigr).
\]
We shall define a filtered bundle
$\Nahm_{\ast}(\nbigp_{\ast}\nbige,\theta)$
over 
$\Nahm(\nbigp_{\ast}\nbige,\theta)$.

By Lemma \ref{lem;13.2.13.3},
there exists a neighbourhood $U$
of $\infty\in\proj^1$
such that 
$\vecN(\nbigp_{\ast}\nbige,\theta)_{|T\times\{\tau_1\}}$
are semistable of degree $0$
for any $\tau_1\in U$.
Let $\gminis\subset T^{\lor}\times U$
denote the spectrum.
We have $\gminis\cap (T^{\lor}\times\{\infty\})\subset D$.
We fix a lift of $D$ to $\Dtilde\subset\cnum$.
Then, if $U$ is sufficiently small,
we may have a lift of $\gminis$
to $\gministilde\subset\cnum\times U$.
We obtain the corresponding holomorphic vector bundle
$\vecV$ with an endomorphism $g$
such that $\Sp(g)\subset\gministilde$.
(See \S\ref{subsection;13.1.21.10}.)
We have the decomposition
\[
 (\vecV,g)
=\bigoplus_{P\in D}
 (\vecV_{P},g_{P}),
\]
where $\Sp(g_P)\cap (\cnum\times\{\infty\})$
is the lift of $P$.
We have the induced decomposition
on $T\times U$:
\[
 \Nahm(\nbigp_{\ast}\nbige,\theta)
=\bigoplus_{P\in D}
 \Nahm(\nbigp_{\ast}\nbige,\theta)_P
\]

Let $U_P\subset T^{\lor}$ be a small neighbourhood of $P\in D$.
We use the coordinate $\zeta_P:=\zeta-\Ptilde$.
By the construction,
we have a natural isomorphism
$\vecV_P
\simeq
 \vecnbign^{0,\infty}\bigl(
 \nbigp_{\ast}(\nbige,\theta)_{|U_P}
 \bigr)$.
We have $g_P=g'_P+\Ptilde\,\id$,
where $g'_P$ is the endomorphism
induced by $\zeta_P$.
Thus, we obtain a filtered bundle 
$\Nahm_{\ast}(\nbigp_{\ast}\nbige,\theta)_P$
over $\Nahm(\nbigp_{\ast}\nbige,\theta)_P$,
by transferring 
$\nbign_{\ast}^{0,\infty}\bigl(
 \nbigp_{\ast}(\nbige,\theta)_{|U_P}\bigr)$.
By taking the direct sum,
we obtain a filtered bundle
$\Nahm_{\ast}(\nbigp_{\ast}\nbige,\theta)$
over $\Nahm(\nbigp_{\ast}\nbige,\theta)$.

\begin{rem}
We obtain a different transformation
by replacing
$\poincare$ and $wd\zeta$
with $\poincare^{\lor}$ and $-wd\zeta$,
respectively,
for which we can argue in a similar way.
\hfill\qed
\end{rem}

\begin{rem}
In {\rm\cite{Bonsdorff1}},
the Fourier transform for Higgs bundles on 
smooth projective curves are studied.
The algebraic Nahm transform in this paper
may be regarded as a filtered variant,
although we consider only the case
where the base space is an elliptic curve.
We also remark that this construction 
is an analogue of
the Fourier transform of the minimal extension
of algebraic meromorphic flat bundles on affine lines.
\hfill\qed
\end{rem}

\subsubsection{Some property}

Let $(\nbigp_{\ast}\nbige,\theta)$
be a filtered Higgs bundle
on $(T^{\lor},D)$
satisfying {\bf (A0)}.

\begin{prop}
\label{prop;13.10.18.2}
The filtered bundle
$\Nahm_{\ast}(\nbigp_{\ast}\nbige,\theta)$
is admissible
and satisfies the condition {\bf (A3)}.
\end{prop}
\pf
Let $(\nbigp_{\ast}\nbige,\theta)$ be an admissible Higgs bundle
on $(T^{\lor},D)$.
Clearly,
$\Nahm(\nbigp_{\ast}\nbige,\theta)$ satisfies {\bf(A1)}.
It satisfies {\bf(A2)} by Proposition \ref{prop;13.1.20.10}.

Let $L\in T^{\lor}$.
We set
$\vecN_L:=
 \vecN(\nbigp_{\ast}\nbige,\theta)\otimes L^{\lor}$.
We have the type decomposition
$\vecN_L=\bigoplus_P
 \bigoplus_{p,m,\veco}
 (\vecN_{L})^{(p,m)}_{P,\veco}$.
By the construction,
we have
\[
 (\vecN_L)^{(p,m)}_{P,\veco}
=\left\{
\begin{array}{ll}
 \nbigp_{-1/2}
 \Nahm(\nbigp_{\ast}\nbige,\theta)^{(p,m)}_{P\otimes L,\veco}
 \otimes L^{\lor}
 &
(\mbox{\rm if }(p,m,\veco)\neq (1,0,0))
 \\
  \nbigp_{0}
 \Nahm(\nbigp_{\ast}\nbige,\theta)^{(1,0)}_{P\otimes L,0}
 \otimes L^{\lor}
&
 (\mbox{\rm if } (p,m,\veco)=(1,0,0))
\end{array}
 \right.
\]
Here, $P\otimes L\in T^{\lor}$ denotes the multiplication 
of $P,L\in T^{\lor}$ in the group $T^{\lor}$.
We shall study 
the cohomology of $\vecN_L$ and its variant.
Let us consider the following complex on
$T^{\lor}\times T\times\proj^1$:
\[
 \begin{CD}
 \nbigctilde^0_L:=
 \nbigctilde^0\otimes p_2^{\ast}L^{\lor}
 @>{\theta+wd\zeta}>>
 \nbigctilde^1_L:=
 \nbigctilde^0\otimes p_2^{\ast}L^{\lor}
 \end{CD}
\]
By the construction,
we have
$\vecN_L\simeq
 R^1p_{23\ast}\nbigctilde^{\bullet}_L$.
We have 
$Rp_{12\ast}\nbigctilde_L^{\bullet}
\simeq
 R^1p_{12\ast}\nbigctilde_L^{\bullet}[-1]
\simeq
 \nbigc^1(\nbigp_{\ast}\nbige,\theta)
 \otimes\poincare \otimes L^{\lor}[-1]$
on $T^{\lor}\times T$.
For the projection
$\pi:T^{\lor}\times T\lrarr T^{\lor}$,
we have
$R\pi_{\ast}(\nbigc^1(\nbigp_{\ast}\nbige,\theta)
 \otimes \poincare\otimes L^{\lor})
\simeq
 \nbigc^1(\nbigp_{\ast}\nbige,\theta)
  \otimes 
 R^1\pi_{\ast}(
\poincare\otimes L^{\lor})[-1]$,
which is a skyscraper sheaf 
$\nbigc^1(\nbigp_{\ast}\nbige,\theta)_{|L}$
at $L$.
Hence, we have
\begin{equation}
 \label{eq;13.10.14.1}
 \hyperh^i\bigl(
 T^{\lor}\times T\times\proj^1,
 \nbigctilde_L^{\bullet}
 \bigr)
\simeq\left\{
 \begin{array}{ll}
 0 & (i\neq 2) \\
 \nbigc^1(\nbigp_{\ast}\nbige,\theta)_{|L}
 & (i=2)
 \end{array}
 \right.
\end{equation}
We obtain
$H^i\bigl(T\times\proj^1,
 \vecN_L
 \bigr)=0$ 
unless $i=1$,
and 
$H^1\bigl(T\times\proj^1,
 \vecN_L
 \bigr)
 \simeq
 \nbigc^1(\nbigp_{\ast}\nbige,\theta)_{|L}$.
We have
\[
 p_{12\ast}\bigl(
 \nbigctilde_L^{\bullet}\otimes\nbigo_{\proj^1}(-1)
 \bigr)
\simeq
 R^1p_{12\ast}\bigl(
 \nbigctilde_L^{\bullet}\otimes\nbigo_{\proj^1}(-1)
 \bigr)[-1]
\simeq
 \nbigc^0(\nbigp_{\ast}\nbige,\theta)
 \otimes\poincare\otimes L^{\lor}[-1]
\]
on $T\times T^{\lor}$.
Hence, we have
\[
 \hyperh^i\bigl(
 T^{\lor}\times T\times\proj^1,
 \nbigctilde_L^{\bullet}
 \otimes\nbigo_{\proj^1}(-1)
 \bigr)
\simeq\left\{
 \begin{array}{ll}
 0 & (i\neq 2) \\
 \nbigc^0(\nbigp_{\ast}\nbige,\theta)_{|L}
 & (i=2)
 \end{array}
 \right.
\]
We obtain
$H^i\bigl(T\times\proj^1,
 \vecN_L
 \otimes\nbigo_{\proj^1}(-1)\bigr)=0$
unless $i=1$,
and 
$H^1\bigl(T\times\proj^1,
 \vecN_L
 \otimes\nbigo_{\proj^1}(-1)\bigr)
\simeq
 \nbigc^0(\nbigp_{\ast}\nbige,\theta)_{|L}$.

\begin{lem}
The map
$H^1\bigl(T\times\proj^1,
 \vecN_L
 \otimes\nbigo_{\proj^1}(-1)\bigr)
\lrarr
 H^1\bigl(T\times\proj^1,
 \vecN_L\bigr)$
induced by the multiplication of $w$
is equal to the map 
$\nbigc^0(\nbigp_{\ast}\nbige,\theta)_{|L}
\lrarr\nbigc^1(\nbigp_{\ast}\nbige,\theta)_{|L}$
induced by $\theta$,
up to signatures.
\end{lem}
\pf
Let $V_i$ $(i=0,1)$ be vector spaces
with morphisms
$f_0,f_{\infty}\in\Hom(V_0,V_1)$.
Let $\alpha_{\kappa}:
 \nbigo_{\proj^1}(-1)\lrarr\nbigo_{\proj^1}$
be morphisms induced as
$\nbigo_{\proj^1}(-1)
\simeq
 \nbigo_{\proj^1}(-\{\kappa\})
\lrarr
 \nbigo_{\proj^1}$.
The induced morphisms
$\nbigo_{\proj^1}(-m-1)
\lrarr
 \nbigo_{\proj^1}(-m)$
are also denoted by $\alpha_{\kappa}$.

We consider a complex $C^{\bullet}$ on $\proj^1$
given as 
$C^0=V_0\otimes\nbigo_{\proj^1}(-1)$
and 
$C^1=V_1\otimes\nbigo_{\proj^1}$
with 
$f_0\alpha_0-f_{\infty}\alpha_{\infty}$.
The morphisms
$\alpha_0$ induce
$C^{\bullet}\otimes\nbigo_{\proj^1}(-1)
\lrarr
 C^{\bullet}$.
We have 
$\hyperh^1(\proj^1,C^{\bullet})
\simeq V_1$
and
$\hyperh^1(\proj^1,C^{\bullet}\otimes\nbigo(-1))\simeq V_0$.
Let us prove that the induced map
$a:\hyperh^1(\proj^1,C^{\bullet})
\lrarr
\hyperh^1(\proj^1,C^{\bullet}\otimes\nbigo(-1))$
is equal to $f_{\infty}$ up to signatures
under the identifications,
which implies the claim of the lemma.

Although we can check it by a direct computation,
we may also use the following argument.
We consider a double complex given as follows:
We set
$C^{00}=V_0\otimes\nbigo(-2)$,
$C^{01}=V_1\otimes\nbigo(-1)$,
$C^{10}=V_0\otimes\nbigo(-1)$
and 
$C^{11}=V_1\otimes\nbigo$.
The morphisms $C^{0i}\lrarr C^{1i}$
are given by $\alpha_0$,
and the morphisms
$C^{i0}\lrarr C^{i1}$ are given by
$f_0\alpha_0-f_{\infty}\alpha_{\infty}$.

For $i=0,1$,
we set 
$D_i^{ij}=C^{ij}$
and $D_i^{kj}=0$ for $k\neq i$.
Then, we have an exact sequence
of the double complexes
$0\lrarr D_1^{\bullet\bullet}
 \lrarr C^{\bullet\bullet}
 \lrarr D_0^{\bullet\bullet}\lrarr 0$.
Similarly, we set
$E_i^{ji}=C^{ji}$
and $E_i^{jk}=0$ for $k\neq i$.
Then, we have an exact sequence
$0\lrarr E_1^{\bullet\bullet}
\lrarr C^{\bullet\bullet}
\lrarr E_0^{\bullet\bullet}\lrarr 0$.
We set 
$F_0^{00}=C^{00}$
and $F_0^{ij}=0$ for $(i,j)\neq (0,0)$.
We set
$F_1^{ij}=C^{ij}$ for $(i,j)\neq (0,0)$
and $F_1^{00}=0$.
Then, we have an exact sequence
$0\lrarr F_1^{\bullet\bullet}
 \lrarr C^{\bullet\bullet}
 \lrarr F_0^{\bullet\bullet}\lrarr 0$.
We have the following commutative diagrams:
\[
 \begin{CD}
 D_1^{\bullet\bullet} @>>> 
 C^{\bullet\bullet} @>>>
 D_2^{\bullet\bullet}
 \\
 @VVV @VVV @VVV \\
 F_1^{\bullet\bullet} @>>> 
 C^{\bullet\bullet} @>>>
 F_2^{\bullet\bullet}
 \end{CD}
\quad\quad\quad\quad
 \begin{CD}
 E_1^{\bullet\bullet} @>>> 
 C^{\bullet\bullet} @>>>
 E_2^{\bullet\bullet}
 \\
 @VVV @VVV @VVV \\
 F_1^{\bullet\bullet} @>>> 
 C^{\bullet\bullet} @>>>
 F_2^{\bullet\bullet}
 \end{CD}
\]
The natural morphisms
$\hyperh^{\ast}(\proj^1,\Tot D_i^{\bullet\bullet})
\lrarr
 \hyperh^{\ast}(\proj^1,\Tot F_i^{\bullet\bullet})
\llarr
 \hyperh^{\ast}(\proj^1,\Tot E_i^{\bullet\bullet})$
are isomorphisms.
The map $a$ is regarded as the connecting homomorphism
of the long exact sequence associated to
$0\lrarr\Tot D_1^{\bullet\bullet}\lrarr
 \Tot C_1^{\bullet\bullet}\lrarr 
 \Tot D_0^{\bullet\bullet}\lrarr 0$.
The cokernel of
$C^{0i}\lrarr C^{1i}$ are 
the skyscraper sheaf at $\infty$,
whose fibers are $V_i$.
Hence, the connecting homomorphism
for $0\lrarr 
 \Tot E_1^{\bullet\bullet}
 \lrarr \Tot C^{\bullet\bullet}
 \lrarr \Tot E_0^{\bullet\bullet}
 \lrarr 0$
is $f_{\infty}$ up to signature.
Thus, we are done.
\hfill\qed

\vspace{.1in}
For any $Y$,
let $\iota_{\infty}:
 Y\times\{\infty\}\lrarr
 Y\times \proj^1$.
The morphism
$\vecN_L
 \lrarr
 \iota_{\infty\ast}\vecN_{L|T\times\{\infty\}}$
is obtained as the push-forward of
$\nbigctilde_L^{\bullet}
\lrarr
 i_{\infty\ast}\bigl(
 \nbigc^1/\nbigc^0
 \otimes\poincare
 \otimes L^{\lor}
 \bigr)$.
Hence, 
$H^1\bigl(T\times\proj^1,
 \vecN_L
 \bigr)
\lrarr
 H^1\bigl(T,
 \vecN_{L|T\times\{\infty\}}
 \bigr)$
is identified with
$\nbigc^1_{|L}\lrarr(\nbigc^1/\nbigc^0)_{|L}$.
By the construction,
the parabolic filtration of
$\bigl(
 (\vecN_L)_{0,0}^{(1,0)}
 \bigr)_{|T\times\{\infty\}}$
is induced by the isomorphism
$\bigl(
 (\vecN_L)_{0,0}^{(1,0)}
 \bigr)_{|T\times\{\infty\}}
\simeq
 (\nbigc^1/\nbigc^0)^{(1,0)}_{L,0}
 \otimes\nbigo_T$.

We have the following commutative diagram:
\[
 \begin{CD}
 H^1\bigl(T\times\proj^1,
 \vecN_L
 \otimes\nbigo(-1)
 \bigr)
 @>{b_1}>>
 H^1\Bigl(
 T\times\proj^1,
 \Gr^{\nbigp}_{-1}\bigl(
 (\vecN_L)^{1,0}_{P,0}
 \bigr)
 \Bigr)
 \\ 
 @V{b_2}VV @V{b_3}V{\simeq}V \\  
 H^1\bigl(T\times\proj^1,
 \vecN_L
 \bigr)
 @>{b_4}>>
H^1\Bigl(
 T\times\proj^1,
 \Gr^{\nbigp}_{0}\bigl(
 (\vecN_L)^{1,0}_{P,0}
 \bigr)
 \Bigr) 
 \end{CD}
\]
Here, $b_2$ and $b_3$ are induced by
the multiplication of $w$.
By the previous consideration,
the composite
$b_4\circ b_2$
is identified with
$\nbigc^1_{|L}
\lrarr
 \Gr^F_1\bigl((\nbigc^1/\nbigc^0)_L\bigr)$,
which is surjective.
Hence, $b_1$ is surjective.
Let 
$\nbign_L$ denote the kernel of
$\vecN_L\otimes\nbigo(-1)
\lrarr
 \Gr^{\nbigp}(\vecN_L)^{(1,0)}_{0,0}$.
We obtain
$H^2(T\times\proj^1,\nbign_L)=0$
from the surjectivity of $b_1$
and 
$H^2(T\times\proj^1,\vecN_L\otimes\nbigo(-1))=0$.
By the construction,
$\nbign_L\subset\nbigp_0\Nahm(\nbigp_{\ast}\nbige,\theta)
 \otimes L^{\lor}$
satisfies the conditions (i,ii) in 
\S\ref{subsection;13.10.15.2},
and we have
\[
 (\nbign_L)^{(1,0)}_{0,0}
=\nbigp_{<-1}
 \bigl(
 Nahm(\nbigp_{\ast}\nbige,\theta)
 \otimes L^{\lor}
 \bigr)^{(1,0)}_{0,0}
\]
Hence,
by using Lemma \ref{lem;13.10.15.3},
we obtain that
$\Nahm_{\ast}(\nbigp_{\ast}\nbige,\theta)$
satisfies the condition {\bf (A3)}.
\hfill\qed

\subsubsection{Characteristic number}
\label{subsection;13.1.24.10}

For compact complex manifolds $Z_i$ $(i=1,2)$,
let $\omega_{Z_i}\in H^{\ast}(Z_1\times Z_2)$ 
denote the pull back of the fundamental class of $Z_i$
by the projection.

Let $(\nbigp_{\ast}\nbige,\theta)$
be a filtered Higgs bundle on $(T^{\lor},D)$
satisfying the condition {\bf (A0)}.
We shall study the characteristic numbers
of $\Nahm(\nbigp_{\ast}\nbige,\theta)$.

\begin{lem}
We have
$\int_{T\times \proj^1}
 c_1\bigl(\Nahm_a(\nbigp_{\ast}\nbige,\theta)\bigr)
 \omega_{\proj^1}=0$
for any $a\in\real$.
\end{lem}
\pf
It follows from that
$\Nahm_a(\nbigp_{\ast}\nbige,\theta)
\big/
 \Nahm_{<a}(\nbigp_{\ast}\nbige,\theta)$
is of degree $0$
for any $a\in\real$.
\hfill\qed

\vspace{.1in}
The following lemma can be checked easily.

\begin{lem}
$c_2\bigl(
 \Nahm_a(\nbigp_{\ast}\nbige,\theta)\bigr)$
is independent of $a\in\real$.
We denote it by
$c_2\bigl(
 \Nahm_{\ast}(\nbigp_{\ast}\nbige,\theta)
 \bigr)$.
\hfill\qed
\end{lem}

We have the type decomposition
$(\nbigp_{\ast}\nbige,\theta)_{|U_P}
=\bigoplus_{(p,m,\veco)}
 (\nbigp_{\ast}\nbige^{(p,m)}_{P,\veco},\theta^{(p,m)}_{P,\veco})$
on a small neighbourhood $U_P$
of each $P\in D$.
We set
\[
 \ell_P:=
 \dim \Cok\Bigl(
 \Res(\theta):
 \Gr^{\nbigp}_{0}
 \bigl(
 \nbige^{(1,0)}_{P,0}
 \bigr)
\lrarr
  \Gr^{\nbigp}_{0}
 \bigl(
 \nbige^{(1,0)}_{P,0}
 \bigr)
 \Bigr).
\]
We put 
$r^{(p,m)}_{P,\veco}
=\rank\bigl(
 \nbige^{(p,m)}_{P,\veco}
 \bigr)\big/p$
and 
$r^{(p,m)}_{P}:=\sum_{\veco} r^{(p,m)}_{P,\veco}$.
We have
$\sum_{p,m} r^{(p,m)}_{P}p
=\rank\nbige$.

\begin{prop}
\label{prop;13.1.22.3}
We have the following equalities:
\begin{equation}
\label{eq;13.2.13.10}
\rank \Nahm(\nbigp_{\ast}\nbige,\theta)
= \sum_{P}\sum_{p,m}
 r^{(p,m)}_{P}(p+m)
-\sum_{P}\ell_P
\end{equation}
\begin{equation}
\label{eq;13.2.13.11}
 \int_{T\times\proj^1}
 c_1\bigl(\Nahm_{\ast}(\nbigp_{\ast}\nbige,\theta)\bigr)
 \cdot\omega_T
=\deg(\nbigp_{\ast}\nbige)
\end{equation}
\begin{equation}
 \label{eq;13.2.13.12}
 \int_{T\times\proj^1}
 c_2\bigl(\Nahm_{\ast}(\nbigp_{\ast}\nbige,\theta)
 \bigr)
=\rank \nbige
\end{equation}
\end{prop}
\pf
Let us prove (\ref{eq;13.2.13.10})
and (\ref{eq;13.2.13.11}).
We have only to consider
the rank and the degree of
$\Nahm_{\ast}(\nbigp_{\ast}\nbige,\theta)_{|\{0\}\times\proj^1}$.
Let $\nbigv\subset\nbigp_1\nbige$
be the subsheaf determined by the following conditions:
\begin{itemize}
\item
$\nbigv=\nbigp_1\nbige$ on the complement of $D$.
\item
It has a decomposition
$\nbigv=\bigoplus_P\nbigv^{(p,m)}_{P,\veco}$
around each $P\in D$.
\item
We have
 $\nbigv^{(p,m)}_{P,\veco}=
 \nbigp_{1/2}\nbige^{(p,m)}_{P,\veco}$
 for $(p,m,\veco)\neq (1,0,0)$,
and
$\nbigv^{(1,0)}_{P,0}=\nbigp_{1}\nbige^{(1,0)}_{P,0}$.
\end{itemize}
Let $\pi_i$ denote the projection of
$T^{\lor}\times\proj^1$ 
onto the $i$-th component.
We have the following $K$-theoretic description:
\begin{equation}
\label{eq;13.1.21.11}
 \Bigl(
 \nbigctilde^1(\nbigp_{\ast}\nbige,\theta)
-\nbigctilde^0(\nbigp_{\ast}\nbige,\theta)
 \Bigr)_{|T^{\lor}\times\{0\}\times\proj^1}
=\pi_1^{\ast}\Bigl(
 \nbigv
-\sum_{P\in D}\nbigo_P^{\oplus \ell_P}
 \Bigr)
-\pi_2^{\ast}\nbigo_{\proj^1}(-1)\cdot
 \pi_1^{\ast}\Bigl(
 \nbigv
-\sum_{P\in D}
 \sum_{p,m,\veco}
 \nbigo_P^{\oplus r^{(p,m)}_{P,\veco}(p+m)}
 \Bigr)
\end{equation}
The Chern character of (\ref{eq;13.1.21.11}) is equal to the following:
\begin{multline}
\pi_1^{\ast}\ch\bigl(\nbigv\bigr)
-\sum_{P\in D}\ell_P\omega_{T^{\lor}}
-\bigl(1-\omega_{\proj^1}\bigr)
 \Bigl(
 \pi_1^{\ast}\ch\bigl(\nbigv\bigr)
-\sum_P\sum_{p,m}
 r^{(p,m)}_{P}
 (p+m)\omega_{T^{\lor}}
 \Bigr)
 \\
=\Bigl(
 \sum_{P}\sum_{p,m}
 r^{(p,m)}_{P}(p+m)
-\sum_{P}\ell_P
 \Bigr)\omega_{T^{\lor}}
+\omega_{\proj^1}
 \pi_1^{\ast}\ch\bigl(\nbigv\bigr)
-\omega_{\proj^1}
 \sum_P\sum_{p,m}
 r^{(p,m)}_{P}(p+m)
 \omega_{T^{\lor}}
\end{multline}
Hence, the Chern character of 
$\vecN(\nbigp_{\ast}\nbige,\theta)_{|\{0\}\times\proj^1}$
is
\begin{multline}
 \sum_{P}\sum_{p,m}
 r^{(p,m)}_{P}(p+m)
-\sum_{P}\ell_P
+\omega_{\proj^1}
\Bigl(
 \deg\bigl(\nbigv\bigr)
-\sum_P\sum_{p,m}
 r^{(p,m)}_{P}(p+m)
\Bigr)
 \\
=\sum_{P}\sum_{p,m}
 r^{(p,m)}_{P}(p+m)
-\sum_{P}\ell_P
+\omega_{\proj^1}
\Bigl(
 \deg\bigl(\nbigv(-D)\bigr)
-\sum_P\sum_{p,m}
 r^{(p,m)}_{P}m
\Bigr).
\end{multline}
In particular,
we obtain (\ref{eq;13.2.13.10}).
We also obtain the following equality:
\[
\deg\bigl(\vecN(\nbigp_{\ast}\nbige,\theta)_{|\{0\}\times\proj^1}
 \bigr)
=\deg\bigl(\nbigv(-D)\bigr)
-\sum_P\sum_{p,m}r^{(p,m)}_Pm
\]

We set $a(p,m,\veco)=-1/2$ if
$(a,m,\veco)\neq (1,0,0)$,
and $a(1,0,0):=0$.
For the parabolic characteristic numbers,
we have the following expressions:
\[
 \deg(\nbigp_{\ast}\nbige)
=\deg\bigl(\nbigv(-D)\bigr)
-\sum_{P\in D}
 \sum_{p,m,\veco}
 \delta\bigl(\nbigp_{\ast}\nbige^{(p,m)}_{P,\veco},
 a(p,m,\veco)\bigr)
\]
\[
 \deg\bigl(\Nahm_{\ast}(\nbigp_{\ast}\nbige,\theta)
 _{|\{0\}\times\proj^1}
 \bigr)
=\deg\bigl(\vecN(\nbigp_{\ast}\nbige,\theta)_{|\{0\}\times\proj^1}
 \bigr)
-\sum_{P\in D}
 \sum_{p,m,\veco}
 \delta\bigl(
 \Nahm_{\ast}(\nbigp_{\ast}\nbige,\theta)^{(p+m,m)}_{P,\veco},
 a(p,m,\veco)
 \bigr)
\]
Here, $\delta(B,a(p,m,\veco))$ denote the contributions
of the locally given filtered bundles $B$ to the parabolic degree.
(See \S\ref{subsection;13.10.15.1}.)
In the following,
we omit $a(p,m,\veco)$.
In the case $(p,m,\veco)=(1,0,0)$,
we have
\[
 \delta\bigl(
 \Nahm_{\ast}(\nbigp_{\ast}\nbige,\theta)^{(1,0)}_{P,0}\bigr)
=\sum_{-1<c<0}
 c\dim \Gr^{\nbigp}_c
  \Nahm(\nbigp_{\ast}\nbige,\theta)^{(1,0)}_{P,0}
=\sum_{-1<c<0}
 c\dim \Gr^{\nbigp}_c(\nbige^{(1,0)}_{P,0})
= \delta(\nbigp_{\ast}\nbige^{(1,0)}_{P,0}).
\]
Let us consider the case $(p,m,\veco)\neq (1,0,0)$.
Let $\varphi_p:U_{u}\lrarr U_P$ be given by 
$\varphi_p(u)=u^p$.
We have the decomposition
$\varphi_p^{\ast}(\nbigp_{\ast}\nbige^{(p,m)}_{P,\veco},
 \theta^{(p,m)}_{P,\veco})
=\bigoplus_{\alpha\in\veco}
 (\nbigp_{\ast}V_{\alpha},\theta_{\alpha})$.
For any $c\in\real$, we put
\[
 r^{(p,m)}_{P,\veco,c}:=
 \dim\Gr^{\nbigp}_cV_{\alpha}.
\]
It is independent of the choice of $\alpha\in\veco$.
We have the following equality:
\[
 \delta\bigl(\nbigp_{\ast}\nbige^{(p,m)}_{O,\veco}\bigr)
=\sum_{\substack{-p/2-1<c\leq -p/2 \\
 0\leq j\leq p-1
 }}
 r^{(p,m)}_{P,\veco,c}\frac{c-j}{p}
=\sum_{-p/2-1<c\leq -p/2}
 r^{(p,m)}_{P,\veco,c}
 \left(
 c-\frac{1}{2}(p-1)
 \right)
\]
We also have the following equality
from the expression of the parabolic structure of 
$\nbign^{0,\infty}_{\ast}(\nbigp_{\ast}\nbige,\theta)$
in the proof of Proposition \ref{prop;13.1.20.10}:
\[
 \delta\Bigl(
 \Nahm(\nbigp_{\ast}\nbige,\theta)^{(p+m,m)}_{P,\veco}
 \Bigr)
=\sum_{\substack{-p/2-1<c\leq -p/2\\ 0\leq j\leq p+m-1}}
 r^{(p,m)}_{P,\veco,c}\frac{2c-2j-m}{2(p+m)}
=\sum_{-p/2-1<c\leq -p/2}
r^{(p,m)}_{P,\veco,c}
 \Bigl(
 c-m-\frac{1}{2}(p-1)
 \Bigr)
\]
Then, the equality (\ref{eq;13.2.13.11}) follows from
$\sum_{c,\veco}r^{(p,m)}_{P,c,\veco}
=r^{(p,m)}_{P}$.

Let us prove (\ref{eq;13.2.13.12}).
We have
$\int_{T\times\proj^1}
 c_2\bigl(
 \Nahm_{\ast}(\nbigp_{\ast}\nbige,\theta)
 \bigr)=
\int_{T\times\proj^1}
 c_2\bigl(\vecN(\nbigp_{\ast}\nbige,\theta)\bigr)$.
We also have the following:
\[
 \int_{T\times\proj^1}
 c_2\bigl(
 \vecN(\nbigp_{\ast}\nbige,\theta)
 \bigr)
=-\int_{T\times\proj^1}
 \ch_2\bigl(
 \vecN(\nbigp_{\ast}\nbige,\theta)
 \bigr)
=-\int_{T^{\lor}\times T\times\proj^1}
 \ch_3\bigl(
 \nbigctilde^1-\nbigctilde^0
 \bigr)
\]
We have
$\ch_3\bigl(
 \nbigctilde^1
 \bigr)=0$.
We have
$c_1(\poincare)^2
=-2\omega_{T}\omega_{T^{\lor}}$.
We also have
\[
 \int_{T^{\lor}\times T\times \proj^1}
 \ch_3(\nbigctilde^0)
=\int_{T^{\lor}\times T\times\proj^1}
\rank(\nbigv)
\omega_{T}\omega_{T^{\lor}}
\omega_{\proj^1}
=\rank(\nbigv).
\]
Hence, we obtain (\ref{eq;13.2.13.12}).
\hfill\qed

\subsubsection{Stable filtered Higgs bundles of degree $0$}

We consider the standard stability condition
for filtered Higgs bundles on $(T^{\lor},D)$.
For any filtered bundle 
$(\nbigp_{\ast}\nbige,\theta)$
on a projective curve $(X,D)$,
we define the slope
$\mu(\nbigp_{\ast}\nbige):=
 \int_X\parchern_1(\nbigp_{\ast}\nbige)
 \big/\rank\nbige$.
It is called stable (resp. semistable)
if $\mu(\nbigp_{\ast}\nbige')<\mu(\nbigp_{\ast}\nbige)$
(resp.
 $\mu(\nbigp_{\ast}\nbige')\leq \mu(\nbigp_{\ast}\nbige)$)
for any non-trivial filtered subbundle
$\nbigp_{\ast}\nbige'\subset\nbigp_{\ast}\nbige$
such that
$\theta(\nbige')\subset
 \nbige'\otimes\Omega^1$.
A semistable filtered Higgs bundle is called
poly-stable,
if it is a direct sum of stable ones.
The following lemma is easy to see.
\begin{lem}
If $(\nbigp_{\ast}\nbige,\theta)$
be a stable Higgs bundle 
on $(T^{\lor},D)$,
then its dual is also stable.
\hfill\qed
\end{lem}

The following proposition is standard.

\begin{prop}
\label{prop;13.10.20.3}
Let $(\nbigp_{\ast}\nbige,\theta)$ be 
a stable admissible filtered bundle on $(T^{\lor},D)$
with $\deg(\nbigp_{\ast}\nbige)=0$.
Then, if $\rank\nbige>1$,
it satisfies the condition {\bf (A0)}.
\end{prop}
\pf
Indeed, an element
of $\hyperh^0\bigl(T^{\lor},
 \nbigc^{\bullet}(\nbigp_{\ast}\nbige\otimes L,\theta+wd\zeta)
 \bigr)$ corresponds to
a morphism
$(\nbigo_{T^{\lor}}(\ast D),0)
\lrarr
 (\nbigp_{\ast}\nbige\otimes L,\theta)$.
By the stability with
$\deg(\nbigp_{\ast}\nbige)=0$
and $\rank\nbige>1$,
we obtain that such a morphism has to be $0$.
We obtain the vanishing of $\hyperh^2$
from the following lemma.

\begin{lem}
\label{lem;13.10.20.2}
$\hyperh^i\bigl(T^{\lor},
 \nbigc(\nbigp_{\ast}\nbige^{\lor},-\theta^{\lor})
 \bigr)$
is naturally isomorphic to
the dual space of
$\hyperh^{2-i}\bigl(T^{\lor},
 \nbigc(\nbigp_{\ast}\nbige,\theta)\bigr)$.
\end{lem}
\pf
We use the natural identification
$\Omega^1_{T^{\lor}}\simeq\nbigo_{T^{\lor}}$.
Let $P\in D$.
We have
$\bigl(
 \nbigp_0\nbige_{P,0}^{(1,0)}
 \bigr)^{\lor}
\simeq
 \nbigp_{<1}(\nbige^{\lor})_{P,0}^{(1,0)}
=:
 C^1(\nbigp_{\ast}\nbige^{\lor},-\theta^{\lor})^{(1,0)}_{P,0}$.
Let 
$\pi$ denote the projection
$\nbigp_0(\nbige^{\lor})^{(1,0)}_{P,0}
\lrarr
 \Gr^{\nbigp}_0((\nbige^{\lor})^{(1,0)}_{P,0})$.
We have a subspace
\[
 \Ker(\Gr^{\nbigp}_0(\theta^{(1,0)}_{P,0})^{\lor})
\subset
 \Gr^{\nbigp}_0((\nbige^{\lor})^{(1,0)}_{P,0}).
\]
We have a natural isomorphism
\[
 C^0(\nbigp_{\ast}\nbige^{\lor},-\theta^{\lor})^{(1,0)}_{P,0}:=
 \pi^{-1}\Bigl(
 \Ker(\Gr^{\nbigp}_0(\theta^{(1,0)}_{P,0})^{\lor})
 \Bigr)
\simeq
 \bigl(
 \nbigc^1(\nbigp_{\ast}\nbige,\theta)_{P,0}^{(1,0)}
 \bigr)^{\lor}
\]
The Higgs field $-\theta^{\lor}$
induces 
$C^0(\nbigp_{\ast}\nbige^{\lor},-\theta^{\lor})^{(1,0)}_{P,0}
 \lrarr 
 C^1(\nbigp_{\ast}\nbige^{\lor},-\theta^{\lor})^{(1,0)}_{P,0}$.
The complex 
$C^{\bullet}(\nbigp_{\ast}\nbige^{\lor},-\theta^{\lor})^{(1,0)}_{P,0}[1]$
is the dual of
$\nbigc^{\bullet}(\nbigp_{\ast}\nbige,-\theta^{\lor})^{(1,0)}_{P,0}$.
The natural inclusions induce a quasi isomorphism
$C^{\bullet}(\nbigp_{\ast}\nbige^{\lor},
 -\theta^{\lor})^{(1,0)}_{P,0}
\lrarr
 \nbigc^{\bullet}(\nbigp_{\ast}\nbige^{\lor},
 -\theta^{\lor})^{(1,0)}_{P,0}$.
For $(p,m,\veco)\neq (1,0,0)$,
the dual of the complex
$\nbigc^{\bullet}(\nbigp_{\ast}\nbige,\theta)^{(p,m)}_{P,\veco}$
is 
\[
\begin{CD}
\nbigp_{<1/2}(\nbige^{\lor})^{(p,m)}_{P,-\veco}
@>{-\theta^{\lor}}>>
 \nbigp_{<3/2+m/p}(\nbige^{\lor})^{(p,m)}_{P,-\veco} 
\end{CD}
\]
where 
the first term sits in the degree $-1$.
It is naturally quasi-isomorphic to
$\nbigc^{\bullet}(\nbigp_{\ast}\nbige^{\lor},-\theta^{\lor})
 ^{(p,m)}_{P,-\veco}[1]$.
Then, the claim of the lemma follows from Serre duality.
Thus, we complete 
the proof of Lemma \ref{lem;13.10.20.2}
and Proposition \ref{prop;13.10.20.3}.
\hfill\qed

\subsubsection{Filtered Higgs bundles of rank $1$
on $(T^{\lor},D)$}

Filtered Higgs bundles of rank $1$ are always
admissible and stable.
Let $(\nbigp_{\ast}\nbige,\theta)$
be a filtered Higgs field of rank $1$ on $(T^{\lor},D)$.
For each $P\in D$,
we have the complex number $\Res_P(\theta)$.
We also have $a(P)\in\real$
such that 
$\Par(\nbigp_{\ast}\nbige,P)=\{a(P)+n\,|\,n\in\seisuu\}$.
Such $a(P)$ is uniquely determined in $\real/\seisuu$.
We say that $P$ is a non-trivial singularity of
$(\nbigp_{\ast}E,\theta)$,
if $(\Res_P\theta,a(P))\neq (0,0)$ in $\cnum\times(\real/\seisuu)$.
If $P$ is a trivial singularity, i.e.,
$(\Res_P\theta,a(P))=(0,0)$,
we obtain a filtered Higgs bundle
on $(T^{\lor},D\setminus\{P\})$
by considering the lattice
$\nbigp_0(\nbige)$ around $P$.
The following proposition is clear.
\begin{prop}
\mbox{{}}
Let $(\nbigp_{\ast}\nbige,\theta)$
be a filtered Higgs bundle of rank $1$
on $(T^{\lor},D)$.
\begin{itemize}
\item
 If each $P\in D$ is a trivial singularity of
$(\nbigp_{\ast}\nbige,\theta)$,
 then 
 $(\nbigp_{\ast}\nbige,\theta)
\simeq
 (L(\ast D),\alpha\,d\zeta)$
for some $\alpha\in\cnum$
and some line bundle $L$ of degree $0$.
Here the parabolic structure of
$L(\ast D)$ is given in a typical way
as in {\rm\ref{subsection;13.10.27.1}}.
\item
If one of $P\in D$ is a non-trivial singularity of
$(\nbigp_{\ast}\nbige,\theta)$,
then $(\nbigp_{\ast}\nbige,\theta)$
satisfies {\bf (A0)}.
\hfill\qed
\end{itemize}
\end{prop}

\subsection{Algebraic Nahm transform for admissible filtered bundles}
\label{subsection;13.2.14.30}
\subsubsection{Construction of the transform}
\label{subsection;13.1.21.111}

For $I\subset\{1,2,3\}$,
let $p_I$ be the projection of
$T^{\lor}\times T\times \proj^1$
onto the product of the $i$-th components
$(i\in I)$.
Let $\poincare$ denote the Poincar\'e bundle
on $T^{\lor}\times T$.

Let $\nbigp_{\ast}E$ be an admissible filtered bundle
on $(T\times\proj^1,T\times\{\infty\})$
satisfying the conditions {\bf (A3)}.
We put $D:=\Sp_{\infty}(\nbigp_{\ast}E)$.
We define
\[
\Nahm(\nbigp_{\ast}E):=
 R^1p_{1\ast}
  \bigl(p_{12}^{\ast}\poincare^{\lor}\otimes 
 p_{23}^{\ast}\nbigp_{-1}E
 \bigr)
 \otimes\nbigo_{T^{\lor}}(\ast D).
\]
By {\bf(A3)},
$\Nahm(\nbigp_{\ast}E)$
is a locally free $\nbigo_{T^{\lor}}(\ast D)$-module.
By Lemma \ref{lem;13.1.21.30},
we have a natural isomorphism
\[
\Nahm(\nbigp_{\ast}E)\simeq
 R^1p_{1\ast}
  \bigl(p_{12}^{\ast}\poincare^{\lor}\otimes 
 p_{23}^{\ast}\nbigp_{0}E
 \bigr)
 \otimes\nbigo_{T^{\lor}}(\ast D).
\]
Let $w$ be the standard coordinate
of $\cnum\subset\proj^1$.
It naturally gives a section of
$\nbigo_{\proj^1}(1)$.
The multiplication of $-w$ induces
an endomorphism $f$ of 
$\Nahm(\nbigp_{\ast}E)$.
We obtain a Higgs field
$\theta:=fd\zeta$ of $\Nahm(\nbigp_{\ast}E)$.
We shall define a filtered bundle 
$\Nahm_{\ast}(\nbigp_{\ast}E)=\bigl(
 \Nahm_a(\nbigp_{\ast}E)\,\big|\,
 a\in\real\bigr)$
over $\Nahm(\nbigp_{\ast}E)$.

\vspace{.1in}
We have the type decomposition
$\nbigp_{\ast}E=
 \bigoplus_{P\in D}\bigoplus_{p,m,\veco}
 \nbigp_{\ast}E^{(p,m)}_{P,\veco}$
on a neighbourhood of $T\times\{\infty\}$.
Let $\nbigu\subset\nbigp_cE$ be 
an $\nbigo_{T\times\proj^1}$-submodule
for some large $c\in\real$,
satisfying the conditions (i,ii)
in \S\ref{subsection;13.10.15.2}.
We suppose 
$\nbigp_{<-1}E^{(1,0)}_{P,0}\subset
 \nbigu^{(1,0)}_{P,0}\subset 
 \nbigp_0E^{(1,0)}_{P,0}$
 for any $P\in D$.
We define 
$N(\nbigu):=
 R^1p_{1\ast}
 \bigl(p_{12}^{\ast}\poincare^{\lor}\otimes 
 p_{23}^{\ast}\nbigu
\bigr)$.
By Lemma \ref{lem;13.1.21.30}
and Lemma \ref{lem;13.10.15.3},
we have $R^ip_{1\ast}\bigl(
 \bigl(p_{12}^{\ast}\poincare^{\lor}\otimes 
 p_{23}^{\ast}\nbigu
 \bigr)=0$,
and $N(\nbigu)$ is a locally free sheaf
on $T^{\lor}$.

We have the following object in
$D^b(\nbigo_{T^{\lor}\times\proj^1})$:
\[
 \RFM_-(\nbigu):=
 Rp_{13\ast}\bigl(p_{12}^{\ast}\poincare^{\lor}\otimes 
 p_{23}^{\ast}\nbigu\bigr)[1]
\]
We can express $\RFM_-(\nbigu)$ as
a two term complex of
locally free $\nbigo_{T^{\lor}\times\proj^1}$
modules $\nbign_{-1}\stackrel{a}{\lrarr}\nbign_{0}$.
Because $a$ is generically isomorphism,
it is injective.
Hence, we have
$\RFM_-(\nbigu)
\simeq\nbigh^0\bigl(
 \RFM_-(\nbigu)\bigr)$.
We will not distinguish them.

Suppose $0\in D$.
Let $U_0\subset T^{\lor}$ 
denote a small neighbourhood of $0$.
Let $W_{\infty}\subset\proj^1$ be a small neighbourhood
of $\infty$.
We have the following decomposition:
\[
 \RFM_-(\nbigu)_{|U_0\times W_{\infty}}
=\bigoplus_{p,m,\veco}
 \RFM_-(\nbigu^{(p,m)}_{0,\veco}).
\]
If $(p,m,\veco)\neq (1,0,0)$,
the support of
$\RFM_-(\nbigu^{(p,m)}_{0,\veco})$
is proper over $U_0$.
Hence, we have the following decomposition:
\begin{equation}
 \label{eq;13.1.21.120}
 \RFM_-(\nbigu)_{|U_0\times\proj^1}
=\bigoplus_{(p,m,\veco)\neq (1,0,0)}
 \RFM_-(\nbigu^{(p,m)}_{0,\veco})
\oplus
 \nbigm(\nbigu).
\end{equation}
Here,
$\nbigm(\nbigu)_{|U_{0}\times W_{\infty}}
=\RFM_-(\nbigu^{(1,0)}_{0,0})$.
We have similar decompositions
for any $P\in D$.

\vspace{.1in}

We have
$N(\nbigu)(\ast D)
=\Nahm(\nbigp_{\ast}E)$.
We obtain the following decomposition
around any $P\in D$
induced by the decomposition (\ref{eq;13.1.21.120})
considered for $P$:
\[
 N(\nbigu)
=\bigoplus_{p,m,\veco}
 N(\nbigu)^{(p,m)}_{P,\veco}
\]
In particular,
we have the following decomposition
around any $P\in D$:
\begin{equation}
 \label{eq;13.1.21.40}
 \Nahm(\nbigp_{\ast}E)
=\bigoplus_{p,m,\veco}
 \Nahm(\nbigp_{\ast}E)^{(p,m)}_{P,\veco}
\end{equation}

We fix a lift $\Ptilde\in\cnum$ of any $P\in D$,
and we use a local coordinate 
$\zeta_P:=\zeta-\Ptilde$
around $P$.
Let $W_{\infty}$ be a small neighbourhood of $\infty$.
We have the filtered bundles with an endomorphism
$(\nbigp_{\ast}V_{P,\veco}^{(p,m)},g^{(p,m)}_{P,\veco})$
on $(W_{\infty},\infty)$,
as in \S\ref{subsection;13.2.14.11}.
If $(p,m,\veco)\neq (1,0,0)$,
we have a natural isomorphism
$\Nahm(\nbigp_{\ast}E)^{(p,m)}_{P,\veco}
\simeq
 \nbign^{(\infty,0)}\bigl(\nbigp_{\ast}V^{(p,m)}_{P,\veco},
 g^{(p,m)}_{P,\veco}\bigr)$.
Under the isomorphism,
we define
\[
 \Nahm_a(\nbigp_{\ast}E)^{(p,m)}_{P,\veco}
:=\nbign_a^{\infty,0}(\nbigp_{\ast}V^{(p,m)}_{P,\veco},
 g^{(p,m)}_{P,\veco}).
\]

Let us consider the case $(p,m,\veco)=(1,0,0)$.
First, we define
\[
 \Nahm_0(\nbigp_{\ast}E)^{(1,0)}_{P,0}:=
 N(\nbigp_{-1}E)^{(1,0)}_{P,0}.
\]
We set
$\gbigc_P:=\nbigp_{0}E^{(1,0)}_{P,0}
 \big/\nbigp_{-1}E^{(1,0)}_{P,0}$.
We have the following exact sequence
around $P$:
\[
 0\lrarr 
 N(\nbigp_{-1}E)^{(1,0)}_{P,0}
 \lrarr 
 N(\nbigp_{0}E)^{(1,0)}_{P,0}
 \lrarr
 R^1p_{1\ast}\bigl(
 p_{12}^{\ast}\poincare
 \otimes
 p_{23}^{\ast}\gbigc_P
 \bigr)
\lrarr 0
\]
We may regard $\gbigc_P$
as a locally free sheaf on $T$,
and then it is isomorphic to 
a direct sum of some copies of
the line bundle corresponding to $P$.
Hence,
the multiplication of $\zeta_P$ on 
$R^1p_{1\ast}\bigl(
 p_{12}^{\ast}\poincare
 \otimes
 p_{23}^{\ast}\gbigc_P
 \bigr)$ is $0$.
We obtain the induced surjection:
\[
  N(\nbigp_{0}E)^{(1,0)}_{P,0|0}:=
  N(\nbigp_{0}E)^{(1,0)}_{P,0}
 \otimes\nbigo_P
 \lrarr
 R^1p_{1\ast}\bigl(
 p_{12}^{\ast}\poincare
 \otimes
 p_{23}^{\ast}\gbigc_P
 \bigr)
\]
Let $K$ denote the kernel.
We have the following morphisms:
\[
 R^1p_{1\ast}\bigl(
 p_{12}^{\ast}\poincare
 \otimes
 p_{23}^{\ast}\gbigc_P
 \bigr)
\simeq
 N(\nbigp_{0}E)^{(1,0)}_{P,0|0}/K
\stackrel{h}{\lrarr}
 N(\nbigp_{-1}E)^{(1,0)}_{P,0|0}
\]
Here, $h$ is the injection induced by 
the multiplication of $\zeta_P$.
We have a natural isomorphism
of $\cnum$-vector spaces:
\[
 R^1p_{1\ast}\bigl(
 p_{12}^{\ast}\poincare
 \otimes
 p_{23}^{\ast}\gbigc_P
 \bigr)
\simeq
 \nbigp_0V^{(1,0)}_{P,0|\infty}
\]
Hence, for any $-1<c<0$,
we define
\[
 F_c\bigl(
 \Nahm_0(\nbigp_{\ast}E)^{(1,0)}_{P,0}
 \otimes\nbigo_P
 \bigr)
:=F_c(\nbigp_0V^{(1,0)}_{P,0|\infty}).
\]
We also set
$F_0\bigl(
 \Nahm_0(\nbigp_{\ast}E)^{(1,0)}_{P,0}
 \otimes\nbigo_P
 \bigr)
=\Nahm_0(\nbigp_{\ast}E)^{(1,0)}_{P,0}
 \otimes\nbigo_P$.
The filtration of
$\Nahm_0(\nbigp_{\ast}E)^{(1,0)}_{P,0}
 \otimes\nbigo_P$
indexed by $\openclosed{-1}{0}$
induces a filtered bundle
$\Nahm_{\ast}(\nbigp_{\ast}E)^{(1,0)}_{P,0}$
over $\Nahm(\nbigp_{\ast}E)^{(1,0)}_{P,0}$.
In all,
we obtain a filtered bundle
$\Nahm_{\ast}(\nbigp_{\ast}E)$
over $\Nahm(\nbigp_{\ast}E)$.

\begin{prop}
\label{prop;13.10.15.11}
$\Nahm_{\ast}(\nbigp_{\ast}E)$
with $\theta$ is admissible,
and satisfies the condition {\bf (A0)}.
Moreover,
the complex 
$N(\nbigp_{-1}E)
\stackrel{-w}{\lrarr}
 N(\nbigp_0E)$
is naturally identified with
$\nbigc^{\bullet}(\Nahm_{\ast}(\nbigp_{\ast}E))$.
\end{prop}
\pf
Suppose $(p,m,\veco)\neq (1,0,0)$.
By Proposition \ref{prop;13.1.21.110},
$\Nahm_{\ast}(\nbigp_{\ast}E)^{(p,m)}_{P,\veco}$
with $\theta$
is admissible.
Moreover, 
$N(\nbigp_{-1}E)^{(p,m)}_{P,\veco}
\stackrel{-w}{\lrarr}
 N(\nbigp_0E)^{(p,m)}_{P,\veco}$
is naturally identified with
$\nbigc^{\bullet}\bigl(
 \Nahm_{\ast}(\nbigp_{\ast}E)^{(p,m)}_{P,\veco}
 \bigr)$
by the construction.

\begin{lem}
\label{lem;13.10.15.10}
$\Nahm_{\ast}(\nbigp_{\ast}E)^{(1,0)}_{P,0}$
with $\theta$
is admissible of type $(1,0,0)$,
and $N(\nbigp_{-1}E)_{P,0}^{(1,0)}
\stackrel{-w}{\lrarr}
 N(\nbigp_0E)_{P,0}^{(1,0)}$
is naturally identified with
$\nbigc^{\bullet}
 \bigl(\Nahm(\nbigp_{\ast}E)_{P,0}^{(1,0)}\bigr)$.
\end{lem}
\pf
The morphism $f$ induces
$f_0:
 N(\nbigp_{-1}E)^{(1,0)}_{P,0|0}
\lrarr
 N(\nbigp_0E)^{(1,0)}_{P,0|0}/K$.
The endomorphism
$f_0\circ h$ on 
$N(\nbigp_0E)^{(1,0)}_{P,0|0}/K$
is identified with $-wg^{(1,0)}_{P,0}$ on 
$\nbigp_0V^{(1,0)}_{P,0|\infty}$.
It is nilpotent, and it preserves the parabolic filtration
on $\nbigp_0V^{(1,0)}_{P,0|\infty}$.
By the construction of the parabolic filtration,
$h\circ f_0$ preserves the parabolic filtration
on $N(\nbigp_{-1}E)^{(1,0)}_{P,0|0}$.
Thus, we obtain that
$\Nahm_{\ast}(\nbigp_{\ast}E)^{(1,0)}_{P,0}$
is admissible of type $(1,0,0)$.

Clearly, we have a natural isomorphism
$N(\nbigp_{-1}E)^{(1,0)}_{P,0}
\simeq
 \nbigc^0(\Nahm(\nbigp_{\ast}E)^{(1,0)}_{P,0})$
by the construction.
We have the natural morphism
$A:N(\nbigp_{0}E)^{(1,0)}_{P,0}
\lrarr
 N(\nbigp_{-1}E)^{(1,0)}_{P,0}\otimes\Omega^1(P)$.
Let $\rho$ denote the natural map
$N(\nbigp_0E)^{(1,0)}_{P,0}\lrarr
 N(\nbigp_0)^{(1,0)}_{P,0|0}/K
\simeq
 \nbigp_0V^{(1,0)}_{P,0|\infty}$.
By the construction
of $\Nahm_{\ast}(\nbigp_{\ast}E)^{(1,0)}_{P,0}$,
the image of $\rho^{-1}(F_{<0})$ by $A$
is equal to
$\nbigp_{<0}\Nahm_{\ast}(\nbigp_{\ast}E)^{(1,0)}_{P,0}
 \otimes\Omega^{1}(P)$.
By the construction of $\theta$,
we obtain that
$\Image(A)$ also contains
$\theta\bigl(
 N(\nbigp_{-1}E)^{(1,0)}_{P,0}\bigr)$.
Hence,
$\nbigc^1(\Nahm_{\ast}(\nbigp_{\ast}E)^{(1,0)}_{P,0})$
is contained in $\Image(A)$.
We remark that the following morphism is surjective,
because
$H^2\bigl(T\times\proj^1,\nbigp_{<-1}E\otimes L\bigr)=0$
for any holomorphic line bundle $L$ of degree $0$ on $T$:
\begin{equation}
 \label{eq;13.10.13.1}
\begin{CD}
 N(\nbigp_{-1}E)_{P,0}^{(1,0)}
@>>>
 R^1p_{1\ast}\bigl(
 p_{12}^{\ast}\poincare\otimes
 p_{23}^{\ast}\Gr^{\nbigp}_{-1}(E^{(1,0)}_{P,0})
 \bigr)
@>{w}>{\simeq}>
 R^1p_{1\ast}\bigl(
 p_{12}^{\ast}\poincare\otimes
 p_{23}^{\ast}\Gr^{\nbigp}_{0}(E^{(1,0)}_{P,0})
 \bigr)
\end{CD}
\end{equation}
It implies that
the morphism
$N(\nbigp_{-1}E)^{(1,0)}_{P,0}
\lrarr
 \nbigp_0V^{(1,0)}_{P,0|\infty}/F_{<0}$
induced by $\theta$
is surjective.
Then, we obtain that
$\nbigi=\nbigc^1(\Nahm_{\ast}(\nbigp_{\ast}E)^{(1,0)}_{P,0})$.
The proof of Lemma \ref{lem;13.10.15.10} is finished.
\hfill\qed

\vspace{.1in}
Let us prove that
$\Nahm_{\ast}(\nbigp_{\ast}E)$
with $\theta$ satisfies the condition {\bf (A0)}.
For $I\subset\{1,2,3\}$,
let $p_I$ denote the projection of
$T^{\lor}\times T\times \proj^1$
onto the product of the $i$-th components $(i\in I)$.
For any $a\in\cnum$
and a line bundle $L$ of degree $0$ on $T^{\lor}$,
we consider the complex
\[
\nbigctilde:=
 \Bigl(
\begin{CD}
  p_{23}^{\ast}\nbigp_{-1}E
 \otimes
 p_{12}^{\ast}\poincare
 \otimes p_1^{\ast}L
 @>{-w+a}>>
 p_{23}^{\ast}\nbigp_0E
 \otimes
 p_{12}^{\ast}\poincare
 \otimes p_1^{\ast}L
\end{CD}
\Bigr)
\]
where the first term sits in the degree $-1$.
Because $Rp_{1\ast}\nbigctilde$
is the complex
 $N(\nbigp_{-1}(E))\otimes L
 \stackrel{-w+a}\lrarr
 N(\nbigp_{0}E)\otimes L$
on $T^{\lor}$,
which is identified with
$\nbigc^{\bullet}(\Nahm_{\ast}(E)\otimes L,\theta+a d\zeta)$,
we have
\[
 \hyperh^i\bigl(
 T^{\lor}\times T\times \proj^1,
 \nbigctilde
 \bigr)
\simeq
 \hyperh^i\bigl(
 T^{\lor},
\nbigc^{\bullet}(\Nahm_{\ast}(E)\otimes L,\theta+a d\zeta)
 \bigr)
\]
Because
$Rp_{23\ast}\nbigctilde$
is quasi-isomorphic to
$\nbigp_{-1}E_{|\{L\}\times\proj^1}
 \lrarr
 \nbigp_{0}E_{|\{L\}\times\proj^1}$,
where the first term sits in the degree $0$,
we have
$\hyperh^i\bigl(T^{\lor}\times T\times\proj^1,
 \nbigctilde
 \bigr)=0$
unless $i=1$.
Thus, we obtain that
$\Nahm_{\ast}(\nbigp_{\ast}E)$
with $\theta$ satisfies {\bf (A0)},
and the proof of Proposition \ref{prop;13.10.15.11}
is finished.
\hfill\qed

\vspace{.1in}
We denote the filtered Higgs bundle
$\bigl(
 \Nahm_{\ast}(\nbigp_{\ast}E),\theta
\bigr)$
just by $\Nahm_{\ast}(\nbigp_{\ast}E)$.

\begin{rem}
We obtain a slightly different transformation
by replacing $\poincare$
with $\poincare^{\lor}$,
for which we can argue in a similar way.
\hfill\qed
\end{rem}

\subsubsection{Inversion}

\begin{prop}
\mbox{{}}\label{prop;13.1.22.2}
\begin{itemize}
\item
 Let $(\nbigp_{\ast}\nbige,\theta)$
 be an admissible filtered Higgs bundle
 on $(T^{\lor},D)$ satisfying {\bf(A0)}.
 Then,
we have a natural isomorphism
$\Nahm_{\ast}\bigl(
 \Nahm_{\ast}(\nbigp_{\ast}\nbige,\theta)
 \bigr)
\simeq
 (\nbigp_{\ast}\nbige,\theta)$.
\item
Let $\nbigp_{\ast}E$ be an admissible filtered bundle 
on $(T\times\proj^1,T\times\{\infty\})$
satisfying the conditions {\bf (A3)}.
Then,
we have a natural isomorphism
$\Nahm_{\ast}\bigl(
 \Nahm_{\ast}(\nbigp_{\ast}E)
 \bigr)
\simeq
 \nbigp_{\ast}E$.
\end{itemize}
\end{prop}
\pf
For any $I\subset\{1,2,3,4\}$,
let $p_I$ denote the projection of
$T^{\lor}\times T\times\proj^1\times T^{\lor}$
onto the product of the $i$-th components
$(i\in I)$.
We set
$\nbigc^i:=\nbigc^i(\nbigp_{\ast}\nbige,\theta)$.
We consider the following complex 
on $T^{\lor}\times T\times\proj^1$:
\[
\begin{CD}
 p_1^{\ast}\nbigc^0
\otimes
 p_{12}^{\ast}\poincare
\otimes
 \nbigo_{\proj^1}(-1)
\otimes
 p_{24}^{\ast}\poincare^{\lor}
@>{\theta+wd\zeta}>>
  p_1^{\ast}\nbigc^1
\otimes
 p_{12}^{\ast}\poincare
\otimes
 p_{24}^{\ast}\poincare^{\lor}
\end{CD}
\]
The complex is denoted by
$\overline{\nbigc}^{\bullet}$.
We can observe that
$Rp_{14\ast}\overline{\nbigc}^{\bullet}$
is quasi-isomorphic to
$p_1^{\ast}\nbigc^1\otimes\nbigo_{\Delta}[-2]$,
where $p_1:T^{\lor}\times T^{\lor}\lrarr T^{\lor}$
denotes the projection onto the first component,
and $\nbigo_{\Delta}$ denote the structure sheaf
of the diagonal.
Hence,
$Rp_{4\ast}\overline{\nbigc}^{\bullet}$
is naturally isomorphic to $\nbigc^1[-2]$.
We can also observe that
$Rp_{234\ast}\overline{\nbigc}^{\bullet}$
is quasi isomorphic to
$q_{12}^{\ast}\vecN(\nbigp_{\ast}\nbige,\theta)\otimes
 q_{13}^{\ast}\poincare^{\lor}[-1]$,
where $q_I$
denotes the projection of
$T\times\proj^1\times T^{\lor}$
onto the product of the $i$-th components
$(i\in I)$.
Hence, we have
$Rp_{4\ast}\overline{\nbigc}^{\bullet}(\ast D)$
is quasi isomorphic to
$\Nahm\bigl(
 \Nahm_{\ast}(\nbigp_{\ast}\nbige,\theta)
 \bigr)$.
We obtain
\[
 \Nahm(\Nahm(\nbigp_{\ast}\nbige,\theta))
 \otimes\nbigo(\ast D)
\simeq
 (\nbige,\theta).
\]
If $(p,m,\veco)\neq (1,0,0)$,
we obtain the comparison of the filtered bundles
over $\nbige^{(p,m)}_{P,\veco}$
from Proposition \ref{prop;13.1.21.100}.
We obtain the comparison of the filtered bundles
over $\nbige^{(1,0)}_{P,0}$
directly from the construction.
Thus, we obtain the first claim.

\vspace{.1in}

Let $\nbigp_{\ast}E$ be an admissible filtered bundle
on $(T\times\proj^1,T\times\{\infty\})$
satisfying {\bf (A3)}.
Let $\nbigv\subset\nbigp_0E$ be 
an $\nbigo_{T\times\proj^1}$-submodule
satisfying the conditions (i,ii)
in \S\ref{subsection;13.10.15.2}
and the following:
\[
 \nbigv^{(p,m)}_{P,\veco}=
 \left\{
 \begin{array}{ll}
 \nbigp_{0}E^{(1,0)}_{P,0},
 & ((p,m,\veco)=(1,0,0)),
 \\
 \nbigp_{-1/2}E^{(p,m)}_{P,\veco},
 & (\mbox{\rm otherwise}).
 \end{array}
 \right.
\]
By Proposition \ref{prop;13.10.15.11},
we have
$\nbigc^0(\Nahm_{\ast}(\nbigp_{\ast}E))
=N\bigl(\nbigv\otimes\nbigo_{\proj^1}(-1)\bigr)$
and 
$\nbigc^1(\Nahm_{\ast}(\nbigp_{\ast}E))
=N\bigl(\nbigv\bigr)\,d\zeta$.
The differential 
$\nbigc^0\lrarr\nbigc^1$
is induced by the multiplication of $-w$.
We shall rewrite the complex
\begin{equation}
 \label{eq;13.1.22.1}
\begin{CD}
 \nbigctilde^0(\Nahm_{\ast}(\nbigp_{\ast}E))
@>{\theta+wd\zeta}>>
 \nbigctilde^1(\Nahm_{\ast}(\nbigp_{\ast}E)).
\end{CD}
\end{equation}

For $I\subset \{1,2,3,4,5\}$,
let $p_I$ denote the projection of
$T\times\proj^1\times T^{\lor}\times T\times\proj^1$
onto the product of the $i$-th components $(i\in I)$.
We set
\[
 C_0:=p_{12}^{\ast}(
 \nbigv\otimes\nbigo_{\proj^1}(-1))
\otimes 
p_{13}^{\ast}\poincare^{\lor}
\otimes
p_{34}^{\ast}\poincare
\otimes
p_5^{\ast}\nbigo_{\proj^1}(-1)
\]
\[
C_1:=p_{12}^{\ast}\nbigv
\otimes 
p_{13}^{\ast}\poincare^{\lor}
\otimes
p_{34}^{\ast}\poincare
\]
We regard $\nbigo_{\proj^1}(1)=\nbigo_{\proj^1}(\{\infty\})$,
and let $\iota:\nbigo_{\proj^1}\lrarr \nbigo_{\proj^1}(\{\infty\})$
be the natural inclusion.
Let $G:C_0\lrarr C_1$ 
be induced by 
$-p_2^{\ast}w\otimes p_5^{\ast}\iota
+p_2^{\ast}\iota\otimes p_5^{\ast}w$.
Then, (\ref{eq;13.1.22.1})
is naturally isomorphic to
$R^1 p_{345\ast}\Bigl(
 C_0\stackrel{G}{\lrarr}C_1
 \Bigr)$.

For $I\subset \{1,2,3,4\}$
let $q_I$ denote the projection of 
$T\times T^{\lor}\times T\times\proj^1$
onto the product of the $i$-th components
$(i\in I)$.
The complex
$p_{1345\ast}\bigl(
 C_0\stackrel{G}{\lrarr} C_1
 \bigr)$
is quasi isomorphic to
\[
 q_{14}^{\ast}\nbigv\otimes
 q_{12}^{\ast}\poincare^{\lor}
 \otimes q_{23}^{\ast}\poincare[-1].
\]
For $I\subset\{1,2,3\}$,
let $s_I$ denote the projection
of $T\times T\times \proj^1$
onto the product of the $i$-th components
$(i\in I)$.
We have the following natural isomorphism
\[
q_{134\ast}\bigl(
 q_{14}^{\ast}\nbigv\otimes
 q_{12}^{\ast}\poincare^{\lor}
 \otimes q_{23}^{\ast}\poincare[-1]
\bigr)
\simeq
 s_{13}^{\ast}\nbigv
\otimes
 s_{12}^{\ast}\nbigo_{\Delta}[-2]
\]
Here,
$\nbigo_{\Delta}$ denote the structure sheaf
of the diagonal in $T\times T$.
Then, we obtain a natural isomorphism
\[
 \nbigv
 \simeq
 \vecN(\Nahm(\nbigp_{\ast}E,\theta))
\]
as $\nbigo_{T\times\proj^1} $-modules.
If $(p,m,\veco)\neq (1,0,0)$,
we obtain the comparison of the filtered bundles
over $\nbigv\bigl(\ast(T\times\{\infty\})\bigr)^{(p,m)}_{P,\veco}$
from Proposition \ref{prop;13.1.21.100}.
The comparison in the case 
$(p,m,\veco)=(1,0,0)$ follows directly
from the construction.
\hfill\qed

\begin{cor}
\label{cor;13.1.29.31}
Let $\nbigp_{\ast}E$ be an admissible filtered bundle on
$(T\times\proj^1,T\times\{\infty\})$
satisfying the conditions {\bf(A3)}.
We have 
\[
 \deg(\nbigp_{\ast}E)
=\deg(\Nahm(\nbigp_{\ast}E))
\]
\end{cor}
\pf
It follows from Proposition \ref{prop;13.1.22.3}
and Proposition \ref{prop;13.1.22.2}.
\hfill\qed

\subsection{Refinement for good filtered Higgs bundles}
\label{subsection;13.2.14.23}
\subsubsection{A stationary phase formula}
\label{subsection;13.11.4.1}

We have the following type of stationary phase formula
for the local Nahm transform,
which is analogue of the stationary phase formula
for the local Fourier transforms.
(See \cite{Fang}, \cite{Fu},
\cite{Graham-Squire},
\cite{Laumon},
\cite{Malgrange-book},
and \cite{Sabbah-local-Fourier}.)
We shall prove it in \S\ref{subsection;13.11.4.2}
after preliminaries in 
\S\ref{subsection;13.11.4.3}--\ref{subsection;13.11.4.4}.

\begin{thm}
\label{thm;13.2.3.2}
Let $U_{\zeta}$ be a small neighborhood of
$0$ in $\cnum_{\zeta}$.
Let $(\nbigp_{\ast}V,\theta)$
be an admissible filtered Higgs bundle on $U_{\zeta}$.
\begin{itemize}
\item
$(\nbigp_{\ast}V,\theta)$ is good,
if and only if 
$\nbign_{\ast}^{0,\infty}(\nbigp_{\ast}V,\theta)$
is good.
\item
Suppose $(\nbigp_{\ast}V,\theta)\simeq
\varphi_{p\ast}(\nbigp_{\ast}V',\theta')$,
where $\theta'-d\gminia\id$ is logarithmic
for some $\gminia\in\zeta_p^{-1}\cnum[\zeta_p^{-1}]$
with $\deg_{\zeta_p^{-1}}\gminia=m>0$.
Then, there exists
$(\nbigp_{\ast}W',\psi')$ on $U_{\tau_{p+m}}$
such that 
(i) $\psi'-d\gminib$ is logarithmic
for some $\gminib\in \tau_{p+m}^{-1}\cnum[\tau_{p+m}^{-1}]$
with $\deg_{\tau^{-1}_{p+m}}\gminib=m$,
(ii) we have an isomorphism
\[
  \varphi_{p+m\ast}(\nbigp_{\ast}W',\psi')
\simeq
 \nbign^{0,\infty}_{\ast}(\nbigp_{\ast}V,\theta).
\]
Moreover, we have an isomorphism
$\Gr^{\nbigp}_c(V')
\simeq
 \Gr^{\nbigp}_{c-m/2}(W')$
under which
$\Res(\varphi_p^{\ast}\theta)
=\Res(\varphi_{p+m}^{\ast}\psi')$.
(The choice of $\gminib$ will be explained in the proof.)
\item
If $(\nbigp_{\ast}V,\theta)$ is logarithmic,
$(\nbigp_{\ast}W,\psi):=
 \nbign^{0,\infty}_{\ast}(\nbigp_{\ast}V,\theta)$
is also logarithmic.
Moreover, we have an isomorphism
\[
 \Gr^{\nbigp}_c(W)\simeq
 \left\{
 \begin{array}{ll}
 \Gr^{\nbigp}_c(V) & (-1<c<0) \\
 \Image\bigl(
 \Res(\theta):
 \Gr^{\nbigp}_0(V)\lrarr
 \Gr^{\nbigp}_0(V)
 \bigr)
 & (c=0)
 \end{array}
 \right.
\]
Under the isomorphism,
we have
$\Res(\psi)=\Res(\theta)$.

\end{itemize}
\end{thm}

We obtain the following corollary from 
Theorem \ref{thm;13.2.3.2}.
(Recall the notion of good filtered bundle
in \S\ref{subsection;13.2.14.11}.)

\begin{cor}
\mbox{{}}\label{cor;13.1.24.21}
\begin{itemize}
\item
Let $(\nbigp_{\ast}\nbige,\theta)$ be a good filtered Higgs bundle
on $(T^{\lor},D)$ satisfying {\bf (A0)}.
Then, $\Nahm_{\ast}(\nbigp_{\ast}\nbige,\theta)$
is a good filtered bundle on $(T\times\proj^1,T\times\{\infty\})$.
\item
Let $\nbigp_{\ast}E$ be a good filtered bundle 
on $(T\times \proj^1,T\times\{\infty\})$
satisfying {\bf (A3)}
with $\Sp_{\infty}(E)=D$.
Then,
$\Nahm_{\ast}(\nbigp_{\ast}E)$
is a good filtered Higgs bundle
on $(T^{\lor},D)$.
\hfill\qed
\end{itemize}
\end{cor}

\subsubsection{Description of the parabolic structure of 
$\nbign^{0,\infty}_{\ast}(\nbigp_{\ast}V,\theta)$}
\label{subsection;13.11.4.3}

Let $(\nbigp_{\ast}V,\theta)$
be a good filtered Higgs bundle on
$(U_{\zeta},0)$.
For simplicity, we assume that
$(\nbigp_{\ast}V,\theta)$
has slope $(p,m)$.
We take $\gminia\in\vecgminio$
for each $\vecgminio\in\vecIrr(\theta)$.
Let $c\in\real$.
We take a frame 
$\vecv_{\vecgminio}=(v_{\vecgminio,i})$ of 
$\nbigp_{cp_{\vecgminio}}
 V^{\vecgminio}_{\gminia}$
compatible with the parabolic structure.
Each  $\zeta_{\vecgminio}^{j}v_{\vecgminio,i}dz_{\vecgminio}/z_{\vecgminio}$
induces a section of
$\nbign^{0,\infty}(\nbigp_{\ast}V,\theta)$,
denoted by 
$\bigl[
 \zeta_{\vecgminio}^{j}v_{\vecgminio,i}dz_{\vecgminio}/z_{\vecgminio}
 \bigr]$.
The following lemma is clear by the construction
of the filtered bundle $\nbign^{0,\infty}_{\ast}(\nbigp_{\ast}V,\theta)$.
(See the proof of Proposition \ref{prop;13.1.20.10}.)
\begin{lem}
\label{lem;13.1.23.30}
The tuple
\[
 \Bigl\{\bigl[
 \zeta_{\vecgminio}^jv_{\vecgminio,i}d\zeta_{\vecgminio}/\zeta_{\vecgminio}
 \bigr]
 \,\Big|\,
 \vecgminio\in\vecIrr(\theta),\,\,
 0\leq j<p_{\vecgminio}+m_{\vecgminio},\,\,
 1\leq i\leq \rank V^{\vecgminio}_{\gminia}
 \Bigr\}
\]
is a frame of 
$\nbign^{0,\infty}_{\kappa_1(p,m,c)}(\nbigp_{\ast}V,\theta)$,
compatible with the parabolic structure.
If the parabolic degree of $v_{\vecgminio,i}$ is $b$,
the parabolic degree of
$\bigl[
 \zeta_{\vecgminio}^jv_{\vecgminio,i}
 d\zeta_{\vecgminio}/\zeta_{\vecgminio}
 \bigr]$
is 
$\bigl(
b-j-m_{\vecgminio}/2
\bigr)
\bigl(p_{\vecgminio}+m_{\vecgminio}\bigr)^{-1}$.
\hfill\qed
\end{lem}

\subsubsection{Description of the parabolic structure of 
$\nbign^{\infty,0}_{\ast}(\nbigp_{\ast}V,g)$}
\label{subsection;13.11.4.4}

Let $(\nbigp_{\ast}V,g)$
be a filtered bundle with an endomorphism
on $(U_{\tau},0)$
such that
$\nbigp_{\ast}V$ with $\psi:=-\tau^{-2}g d\tau$
is a good filtered Higgs bundle.
For simplicity,
we assume that 
$(\nbigp_{\ast}V,g)$ has a slope $(p,m)$
with $p>m\neq 0$.

We take $\gminia\in\vecgminio$
for each $\vecgminio\in\vecIrr(\psi)$.
Let $c\in\real$.
We take a frame 
$\vecv_{\vecgminio}=(v_{\vecgminio,i})$
of $\nbigp_{cp_{\vecgminio}}V^{\vecgminio}_{\gminia}$
compatible with the parabolic structure.
Each $\tau_{\vecgminio}^jv_{\vecgminio,i}$ 
induces a section of
$\nbign^{\infty,0}(\nbigp_{\ast}V,g)$,
denoted by 
$\bigl[\tau_{\vecgminio}^jv_{\vecgminio,i}\bigr]$.
The following lemma is clear by the construction of
the filtered bundle $\nbign^{\infty,0}_{\ast}(\nbigp_{\ast}V,g)$.
(See the proof of Proposition \ref{prop;13.1.21.110}.)

\begin{lem}
\label{lem;13.1.23.102}
The tuple
\[
 \Bigl\{
 \bigl[
 \tau_{\vecgminio}^jv_{\vecgminio,i}
 \bigr]
 \,\Big|\,
 \vecgminio\in\vecIrr(\psi),\,\,
 0\leq j<p_{\vecgminio}-m_{\vecgminio},\,\,
 1\leq i\leq \rank V^{\vecgminio}_{\gminia}
 \Bigr\}
\]
is a frame of $\nbign^{\infty,0}_{\kappa_2(p,m,c)}(\nbigp_{\ast}V,g)$
compatible with the parabolic structure.
If the parabolic degree of $v_{\vecgminio,i}$ is $b$,
the parabolic degree of 
$\tau_{\vecgminio}^jv_{\vecgminio,i}$
is 
$\bigl(
 b-j+p_{\vecgminio}-m_{\vecgminio}/2
 \bigr)
 \bigl(
 p_{\vecgminio}-m_{\vecgminio}
 \bigr)^{-1}$.
\hfill\qed
\end{lem}

\subsubsection{Proof of Theorem \ref{thm;13.2.3.2}}
\label{subsection;13.11.4.2}

Let us return to the situation in \S\ref{subsection;13.11.4.1}.
Let us begin with the third claim.
We obtain the isomorphism of 
the associated graded vector spaces
by the construction of $\nbigp_{\ast}W$.
We have the expression
$\theta=f\,d\zeta/\zeta$,
where $f$ is an endomorphism of $\nbigp_{\ast}V$.
It naturally induces an endomorphism $f'$ of $\nbigp_{\ast}W$,
and we have $\psi=f'\,d\tau/\tau$ by the construction.
Thus, we obtain the third claim.

Let us consider the second claim.
Our argument is close to that in \cite{Fang}.
To simplify the notation,
we set $\eta:=\tau_{p+m}$ and $u:=\zeta_p$.
We set
$G(u):=u\del_u\gminia(u)
=\sum_{j=1}^{m}\alpha_j u^{-j}$.
Let $\omega:=e^{2\pi\sqrt{-1}/(p+m)}$.
We have  holomorphic functions $u^{(i)}(\eta)$
$(i=0,\ldots,p+m-1)$ on $U_{\eta}$,
satisfying  
$\del_{\eta}u^{(i)}(0)=\del_{\eta}u^{(0)}(0)\omega^i$ and
\[
 G\bigl(u^{(i)}(\eta)\bigr)
+pu^{(i)}(\eta)^p/\eta^{p+m}=0.
\]
For any   $c\in\real$,
we consider 
$\nbigp_{c-m/2}\nbigv:=
 \nbigp_{c-m/2}
 \varphi_{p+m}^{\ast}
 \nbign^{0,\infty}(\nbigp_{\ast}V,\theta)$.
We take a frame $\vecv$ of $\nbigp_{c}V'$
compatible with the parabolic structure.
We put $\nutilde_{ij}:=(\eta^{-1}u)^iv_j$
$(0\leq i\leq p+m-1,\,\,
 1\leq j\leq \rank V')$.
They induce a frame of $\nbigp_{c-m/2}\nbigv$,
which is compatible with the parabolic structure.
By the frames,
for $c-1<d\leq c$,
we obtain an isomorphism 
\begin{equation}
 \label{eq;13.2.2.1}
 \Gr^{\nbigp}_{d-m/2}(\nbigv)
\simeq
 \Gr^{\nbigp}_{d}(V')
\otimes
 \cnum^{p+m}.
\end{equation}

The following lemma can be checked 
by a direct computation.
\begin{lem}
$\eta^{-1}u$ gives an endomorphism $F$ of $\nbigp_{\ast}\nbigv$.
On $\Gr^{\nbigp}_{d-m/2}(\nbigv)$,
we have 
\[
 F(\nutilde_{i,j})=
 \left\{
 \begin{array}{ll}
 \nutilde_{i+1,j} & (i<p+m-1) \\
 \mbox{{}}\\
 -p^{-1}\alpha_m\nutilde_{0,j}
 & (i=p+m-1)
 \end{array}
 \right.
\]
The eigenvalues of
$F$ on $\Gr^{\nbigp}$
are $\del_{\eta}u^{(i)}(0)$
$(i=0,\ldots,p+m-1)$.
\hfill\qed
\end{lem}

By the lemma, we obtain the decomposition
$(\nbigp_{\ast}\nbigv,F)
=\bigoplus_{j=0}^{p+m-1}
 (\nbigp_{\ast}\nbigv^{(j)},F^{(j)})$
such that
$F^{(j)}_{|0}$ has a unique eigenvalue 
$\del_{\eta}u^{(j)}(0)$.
Note that 
$\nbign^{0,\infty}(\nbigp_{\ast}V,\theta)
\simeq
 \varphi_{p+m\ast}
 \bigl(\nbigp_{\ast}\nbigv^{(0)},
 \zeta(-\tau^{-2}d\tau) \bigr)$.
We also have an isomorphism
$\Gr^{\nbigp}_c(V')
\simeq
 \Gr^{\nbigp}_{c-m/2}(\nbigv^{(0)})$.

\vspace{.1in}

We have the expression
$\theta_{\gminia}
=(G(u)+f)\,du/u$,
where $f$ is an endomorphism of $\nbigp_{\ast}V'$.
On $\nbigp_{c-m/2}\nbigv$,
we have 
$\eta^m(G(u)+pu^p/\eta^{p+m})=-\eta^mf$.
We have the following decomposition:
\[
 \eta^{m}\bigl(
 G(u)
+pu^p/\eta^{p+m}
 \bigr)
=\bigl(\eta^{-1}u-\eta^{-1}u^{(0)}(\eta)\bigr)
\times
 p\prod_{i=1}^{p+m-1}
 \bigl(
 \eta^{-1}u-\eta^{-1}u^{(i)}(\eta)
 \bigr)
 (\eta^{-1}u)^{-m}
\]
Because 
$\eta^{-1}u-\eta^{-1}u^{(j)}(\eta)$
$(1\leq j<p-m)$ are invertible on
$\nbigp_{c-m/2}\nbigv^{(0)}$,
we obtain the following on
$\nbigp_{c-m/2}\nbigv^{(0)}$:
\[
 \eta^{-1}u-\eta^{-1}u^{(0)}(\eta)
=-p^{-1}\eta^m\cdot f \cdot
 \prod_{j=1}^{p+m-1}
 \bigl(\eta^{-1}u-\eta^{-1}u^{(j)}(\eta)\bigr)^{-1}
 (\eta^{-1}u)^m
\]
Let $Q_k(x,y)=\sum_{i+j=k} x^iy^j$.
We have
\begin{multline}
 \zeta/\tau
-\frac{u^{(0)}(\eta)^p}{\eta^{p+m}}
=\eta^{-m}
 \bigl(
 \eta^{-1}u-\eta^{-1}u^{(0)}(\eta)
 \bigr)
 \cdot
 Q_{p-1}\bigl(\eta^{-1}u,\eta^{-1}u^{(0)}(\eta)\bigr)
 \\
=-f\prod_{j=1}^{p+m-1}\bigl(\eta^{-1}u-\eta^{-1}u^{(j)}(\eta)\bigr)^{-1}
 (\eta^{-1}u)^{m}
  p^{-1}
 Q_{p-1}\bigl(\eta^{-1}u,\eta^{-1}u^{(0)}(\eta)\bigr)
\end{multline}
Hence, we obtain that
$\Bigl(
\zeta/\tau-u^{(0)}(\eta)^p\eta^{-p-m}
\Bigr)
\nbigp_{\ast}\nbigv^{(0)}\subset\nbigp_{\ast}\nbigv^{(0)}$.
On $\Gr_{a}^{\nbigp}(\nbigv^{(0)})$,
the endomorphisms
$u/\eta$ and $u^{(0)}(\eta)/\eta$
are the multiplication of 
$\del_{\eta}u^{(0)}(0)$.
Hence,
$\Bigl(
\zeta/\tau-u^{(0)}(\eta)^p\eta^{-p-m}
\Bigr)$
acts as $-(p+m)^{-1}f$
on $\Gr^{\nbigp}(\nbigv)$.
We set 
$\nbigp_{\ast}W':=\nbigp_{\ast}\nbigv^{(0)}$
and 
$\psi':=-\zeta \tau^{-2}d\tau
=-(\zeta/\tau)\,(p+m)d\eta/\eta$.
We have $\gminib\in \eta^{-1}\cnum[\eta^{-1}]$
uniquely determined by the condition that 
$\eta\del_{\eta}\gminib$
is equal to the polar part of
$-(p+m)u^{(0)}(\eta)^p\eta^{-p-m}$.
Then,
$\psi'-d\gminib$ is logarithmic.
The residue acts as $f$.
Hence, the second claim of Theorem \ref{thm;13.2.3.2} follows.
It also implies the ``only if'' part 
in the first claim.

\vspace{.1in}

Let us prove the ``if'' part of the first claim.
We use the inverse transform.
Let $(\nbigp_{\ast}W,\psi)$ be a good filtered Higgs bundle
on $(U_{\tau},0)$
which is isomorphic to 
$\varphi_{p\ast}
 (\nbigp_{\ast}W',\psi')$,
where
 $\psi'-d\gminib\id$ is logarithmic
 for some $\gminib\in \tau^{-1}_{p}\cnum[\tau_p^{-1}]$
with $\deg_{\tau_p^{-1}}\gminib=m<p$.
If $p=1$, we assume that any eigenvalue of $\Res(\psi')$
is not $0$.
The claim of Theorem \ref{thm;13.2.3.2}
follows from the next proposition.

\begin{prop}
\label{prop;13.2.3.1}
There exists 
$(\nbigp_{\ast}V',\theta')$ on $U_{\zeta_{p-m}}$
such that
(i) $\theta'-d\gminia\id$ is logarithmic
 for some $\gminia\in \zeta_{p-m}^{-1}\cnum[\zeta_{p-m}^{-1}]$,
(ii) we have an isomorphism 
$\varphi_{p-m\ast}
 (\nbigp_{\ast}V',\theta')
\simeq
\nbign_{\ast}^{\infty,0}
 (\nbigp_{\ast}W,\psi)$.
\end{prop}
\pf
To simplify the notation,
we set $\eta:=\tau_{p}$ and $u:=\zeta_{p-m}$.
We have the expression
\[
\psi'=
 \bigl(
 G(\eta)\id+\eta^{p}f
 \bigr)
 \varphi_p^{\ast}(-\tau^{-2}d\tau), 
\]
such that 
(i) $G(\eta)=\sum_{j=1}^{m} \beta_j\eta^{p-j}$
with $\beta_{m}\neq 0$,
(ii) $f$ is an endomorphism of $\nbigp_{\ast}W'$.
We fix a holomorphic function
$\eta^{(0)}(u)$ such that
$G(\eta^{(0)}(u))-u^{p-m}=0$
such that
$0<C_1\leq|\eta^{(0)}/u|\leq C_2$
for some constants $C_i$.

We set 
$\nbigp_{c+p-m/2}\nbigv:=
 \nbigp_{c+p-m/2}\varphi_{p-m}^{\ast}
 \nbign^{\infty,0}(\nbigp_{\ast}W,\psi)$.
Let $\vecv$ be a frame of $\nbigp_{c}W'$
compatible with the parabolic structure.
We set $\nutilde_{ij}=u^{-i}\eta^iv_j$
$(0\leq i\leq p-m-1,\,\,\,1\leq j\leq \rank W')$.
They induce a frame of 
$\nbigp_{c+p-m/2}\nbigv$
compatible with the parabolic structure.
By using the frame,
for any $c-1<d\leq c$,
we obtain an isomorphism
 $\Gr^{\nbigp}_{d+p-m/2}(\nbigv)
\simeq
 \Gr^{\nbigp}_{d}(W')\otimes \cnum^{p-m}$.
The following lemma can be checked directly.
\begin{lem}
$u^{-1}\eta$ gives an endomorphism $F$ of $\nbigp_{\ast}\nbigv$,
preserving the parabolic structure,
and the induced endomorphism
on $\Gr^{\nbigp}(\nbigv)$ is given by
$F(\nutilde_{ij})=\nutilde_{i+1,j}$ $(i=0,\ldots,p-m-2)$
and 
$F(\nutilde_{p-m-1,j})=-\beta_{m}^{-1}\nutilde_{0,j}$.
The eigenvalues are
$\omega^i\del_u\eta^{(0)}(0)$ $(i=0,\ldots,p-m-1)$,
where 
$\omega=e^{2\pi\sqrt{-1}/(p-m)}$.
\hfill\qed
\end{lem}

We obtain the decomposition
$(\nbigp_{\ast}\nbigv,F)
=\bigoplus_{i=0}^{p-m-1}
 (\nbigp_{\ast}\nbigv^{(i)},F^{(i)})$
such that
$F^{(i)}_{|0}$ has a unique eigenvalue
$\omega^i\del_u\eta^{(0)}(0)$.
We have an isomorphism
$\varphi_{p-m\ast}\bigl(
 \nbigp_{\ast}\nbigv^{(0)},
 -\tau^{-1}d\zeta
 \bigr)
\simeq
 \nbign^{\infty,0}_{\ast}(\nbigp_{\ast}W,\psi)$.
We also have an isomorphism
$\Gr^{\nbigp}_{c+p-m/2}(\nbigv^{(0)})
\simeq
 \Gr^{\nbigp}_c(W')$.

\vspace{.1in}

We have
$G(\eta)-u^{p-m}=-\eta^pf$
on $\nbigv$.
Note that
$u^{-(p-m-1)}\sum_{j=1}^m \beta_jQ_{p-j-1}(\eta^{(0)}(u),\eta)$
is invertible on $\nbigp_{c+p-m/2}\nbigv^{(0)}$.
Hence, we obtain the following
on $\nbigp_{c+p-m/2}\nbigv^{(0)}$:
\[
 u^{p-m-1}\bigl(
 \eta^{(0)}(u)-\eta
 \bigr)
=\eta^pf\cdot
 \Bigl(
 \sum_{j=1}^m \beta_jQ_{p-j-1}(\eta^{(0)}(u),\eta)
 \Bigr)^{-1} u^{p-m-1}
\]
We have the following:
\[
 u^{p-m}\bigl(
\eta^{-p}-\eta^{(0)}(u)^{-p}
 \bigr)
=f\eta^p
 Q_{p-1}\bigl(\eta^{(0)}(u)^{-1},\eta^{-1}\bigr)
 \,\eta^{(0)}(u)^{-1}\eta^{-1}
 \Bigl(
 \sum_{j=1}^m \beta_jQ_{p-j-1}(\eta^{(0)}(u),\eta)
 \Bigr)^{-1}u^{p-m}
\]
Hence, we obtain that 
$u^{p-m}\bigl(
 \eta^{-p}-\eta^{(0)}(u)^{-p}
 \bigr)$
is an endomorphism of
$\nbigp_{\ast}\nbigv^{(0)}$.
We set
$\nbigp_{\ast}V':=
 \nbigp_{\ast}\nbigv^{(0)}$
and 
$\theta':=-\tau^{-1}\varphi_{p-m}^{\ast}d\zeta
=-\eta^{-p}(p-m)u^{p-m}(du/u)$.
We have $\gminia\in u^{-1}\cnum[u^{-1}]$
uniquely determined by the condition
$u\del_u\gminia=
-\eta^{(0)}(u)^{-p}(p-m)u^{p-m}$.
Then, 
$\theta'-d\gminia$ is logarithmic.
Thus, the proof of Proposition \ref{prop;13.2.3.1}
and Theorem \ref{thm;13.2.3.2} are finished.
\hfill\qed

\section{Family of vector bundles on torus with small curvature}

\subsection{Small perturbation}
\label{subsection;12.11.5.1}

We use the notation in \S\ref{subsection;13.10.15.30}.
We use the metric $dz\,d\zbar$ of $T$.
For any finite dimensional vector space $V$,
let $L^p_k(V)$ be the space of $V$-valued
$L^p_k$-functions on $T$,
and let $L_k^p(V\otimes\Omega^{i,j})$ be 
the space of $V$-valued $L_k^p$-differential $(i,j)$-forms.
We have the linear map
$\int_T:L_k^p(V)\lrarr V$
given by $\int_Tf:=|T|^{-1}\int_Tf\,|dz\,d\zbar|$,
where $|T|$ denotes the volume of $T$.
The kernel is denoted by
$L_k^p(V)_0$.
We have a natural inclusion
$V\lrarr L_k^p(V)$ 
as constant functions.
We have the decomposition
$L_k^p(V)=L_k^p(V)_0\oplus V$
as topological vector spaces.

Suppose that $V$ is $r$-dimensional
and equipped with a hermitian metric $h_V$.
Let $p\geq 2$.
Let $\nbigg^p_{k}(V)$ be the space of $L_{k+2}^p$-maps
from $T$ to $\GL(V)$.
We set $\gbiga^p_{k}(V):=
\bigl\{\delbar_0+A\,\big|\,
 A\in L_{k+1}^p(\End(V)\otimes\Omega^{0,1})
 \bigr\}$,
i.e.,
the space of $(0,1)$-type differential operators
of the product bundle $\underline{V}$ of the form
$\delbar_0+A$ $(A\in L_{k+1}^p(\End(V)\otimes\Omega^{0,1}))$.
We have the natural right $\nbigg_k^p$-action
on $\gbiga_k^{p}(V)$ given by
$g\bullet \delbar:=
 g^{-1}\circ\delbar\circ g
=\delbar+g^{-1}\delbar g$.

Let $\Gamma$ be an endomorphism of $V$.
Let $U_1\subset L_{k+2}^p(\End(V))_0$ be 
a sufficiently small neighbourhood  of $0$
such that $1+U_1\subset\nbigg_k^p$.
Let $U_2$ be a neighbourhood of  $0$ in $\End(V)$.
We consider the map
$\Psi:U_1\times U_2\lrarr \gbiga_k^{p}(V)$
given by
\[
 \Psi(a,b):=
 (1+a)\bullet\Bigl(\delbar_0+(\Gamma+b)\,d\zbar\Bigr).
\]
We use the norm on $L_{k+2}^p(\End(V))$
such that $L_{k+2}^p(\End(V))
 \simeq L_{k+2}^p(\End(V))_0\oplus \End(V)$
is an isometry,
and the norm on $L_{k+1}^p(\End(V))$
such that
$L_{k+2}^p(\End(V))\lrarr L_{k+1}^p(\End(V))$,
$A\longmapsto
 \delbar_0A+\int_TA$
is an isometry.

\begin{prop}
\label{prop;11.12.18.10}
Fix $\delta>0$.
Suppose that
$\Gamma$ is decomposed as
$\Gamma=\Gamma_0+\Gamma_1$
satisfying the following conditions:
\begin{itemize} 
\item
 $\Gamma_0$ is commutative with its adjoint
 $\Gamma_0^{\dagger}$,
 i.e.,
 it is diagonalizable, 
 and the eigen spaces are orthogonal 
 with respect to $h_V$.
 Moreover, there exists $\zeta_0\in\cnum$
such that
$\Sp(\Gamma)$
is contained in 
\[
K_1(L,\zeta_0):=\bigl\{
 \zeta\in\cnum\,\big|\,
 0\leq\Image\bigl(\zeta-\zeta_0\bigr)\leq(1-\delta)\pi,\,\,
 0\leq\Image\bigl((\zeta-\zeta_0)\taubar\bigr)\leq(1-\delta)\pi
 \bigr\}.
\]
\item
$|\Gamma_1|_{h_V}\leq \delta/100$.
\end{itemize}

Then, there exist positive constants $C_i$ $(i=1,2)$,
independently from $\Gamma$ and $\zeta_0$,
such that the following holds:
\begin{itemize}
\item
For $\nbigb\in 
 L_{k+1}^p\bigl(\End(V)\otimes\Omega^{0,1}\bigr)$
with $|\nbigb|\leq C_1$,
there exists a unique
$(a,b)\in U_1\times U_2$ 
with $|a|+|b|\leq C_2\,|\nbigb|$
satisfying 
$\delbar_0+\Gamma\,d\zbar
+\nbigb\,d\zbar=\Psi(a,b)$.
\end{itemize}
\end{prop}
\pf
We set
$K(L):=\bigl\{
 \zeta\in\cnum\,\big|\,
 |\Image(\zeta)|\leq(1-\delta)\pi,\,\,
 |\Image(\zeta\taubar)|\leq(1-\delta)\pi
 \bigr\}$.
We have $\Sp(\ad(\Gamma_0))\subset K(L)$.
In the following,
$C_i$ will be positive constants
which are independent from $\Gamma$ and $\zeta_0$.

We have a morphism
$\Phi_{\Gamma}:L_{k+2}^p\bigl(\End(V)\bigr)=
 L_{k+2}^p\bigl(\End(V)\bigr)_0\oplus \End(V)
 \lrarr L_{k+1}^p\bigl(\End(V)\otimes\Omega^{0,1}\bigr)$
given by 
\[
 \Phi_{\Gamma}(A,B)
=\delbar A+[\Gamma,A]\,d\zbar
+B\,d\zbar,
\]
where
$A\in L_{k+2}^p\bigl(\End(V)\bigr)_0$
and $B\in \End(V)$.
We have
$\Phi_0(A,B)=\delbar A+B\,d\zbar$,
which is an isometry
by our choice of the norms.

\begin{lem}
$\Phi_{\Gamma}$ is a homeomorphism.
\end{lem}
\pf
Note that $\Phi_0$ is an isomorphism,
and that $\Phi_{\Gamma}-\Phi_0$ is compact.
Hence, the index of $\Phi_{\Gamma}$ is $0$.
Due to the condition for $\Gamma_0$,
we have
$\bigl\|
 \delbar_0A+[\Gamma_0,A]d\zbar
 \bigr\|_{L^2}
 \geq 
 \delta\pi |A|_{L^2}$
for any $A\in L_1^2(\End(V))$.
By the condition for $\Gamma_1$,
we obtain that
$\delbar_0A+[\Gamma,A]d\zbar\neq 0$
for any $A\in L_1^2(\End(V))$.
Then, we obtain that
$\Phi_{\Gamma}$ is injective.
\hfill\qed

\begin{lem}
We have 
$|\Phi_0^{-1}\circ\Phi_{\Gamma}|\leq C_3$
and 
$|\Phi^{-1}_{\Gamma}\circ\Phi_0|\leq C_3$,
independently from $\Gamma$,
where $|\cdot|$ denotes the operator norm.
\end{lem}
\pf
Let $\nbigs$ be the set of $\Gamma$
satisfying the conditions of the proposition.
It is compact.
For any fixed 
$(A,B)\in L_{k+2}^p(\End(V))_0\oplus \End(V)$,
the map
$\Gamma\longmapsto
\Phi_{\Gamma}^{-1}\circ\Phi_0(A,B)$
gives a continuous map
from $\nbigs$ to $L_{k+2}^p(\End(V))
 \oplus\End(V)$,
and hence bounded.
Then, we obtain the claim for
$\Phi_{\Gamma}^{-1}\circ\Phi_0$
by the uniform boundedness principle.
We obtain the claim for
$\Phi_{0}^{-1}\circ\Phi_{\Gamma}$
similarly.
\hfill\qed

\vspace{.1in}

We set
$\nbiga(a,b):=\Psi(a,b)-\Psi(0,0)
 \in L_{k+1}^p\bigl(\End(V)\otimes\Omega^{0,1}\bigr)$,
i.e.,
\[
 \nbiga(a,b)=
 (1+a)^{-1}\bigl(\delbar_0a+[\Gamma,a]\bigr)
+\Ad(1+a)\,b\,d\zbar.
\]
We have
$\bigl|\nbiga(a,b)\bigr|=
O\bigl(|a|+|b|\bigr)$,
independently from $\Gamma$.
The derivative $T_{(a,b)}\Psi$
of $\Psi$ 
at any $(a,b)\in U_1\times U_2$ is given by
\[
 T_{(a,b)}\Psi(X,Y)
=\Phi(X,Y)
+\bigl[\nbiga(a,b),(1+a)^{-1}X\bigr]
-\bigl[\Psi(0,0),(1+a)^{-1}a\,X\bigr]
+\bigl(\Ad(1+a)-1\bigr)Y.
\]
Hence, 
we obtain an estimate
$\bigl|
 \Phi_{\Gamma}^{-1}\circ T_{(a,b)}\Psi-\id
 \bigr|\leq
 C_4\bigl(|a|+|b|\bigr)$,
which is independent from $\Gamma$.
Then, the claim of Proposition \ref{prop;11.12.18.10}
follows from the classical inverse function theorem
(\cite{Lang}, for example).
\hfill\qed

\begin{cor}
$\Psi$ gives a diffeomorphism of
a neighbourhood of $(0,0)$ in $U_1\times U_2$
and a neighbourhood of $\delbar_0+\Gamma d\zbar$
in $\gbiga_k^p(V)$.
\hfill\qed
\end{cor}

\subsection{Frames}
\label{subsection;12.6.24.10}

\subsubsection{Preliminary}

We set
$U_1:=\bigl\{
 (x_1,x_2)\,\big|\,0\leq x_i\leq 1
 \bigr\}$
and $U_2:=\bigl\{
 (\xi_1,\ldots,\xi_{n-2})\,\big|\,
 |\xi_i|\leq 1
 \bigr\}$.
Let $T_0=\real^2/\seisuu^2$.
Let $U_1\times U_2\lrarr T_0\times U_2$
denote the natural projection.
We also use the variables
$t_i=x_i$ $(i=1,2)$ and  $t_i=\xi_{i-2}$ $(i=3,\ldots,n)$.
We also use $x=x_1$, $y=x_2$.

For any positive integer $k$,
we set $S_1(k):=\bigl\{(m_1,m_2)
 \,\big|\,
 m_1+m_2=k,m_i\geq 0
 \bigr\}$.
We also set
$S_2(k):=
 \bigl\{
 (m_1,\ldots,m_{n-2})
 \,\big|\,
 \sum m_i=k,m_i\geq 0
 \bigr\}$.
We set $S(k_1,k_2):=S_1(k_1)\times S_2(k_2)$.
We put 
$\del_{\vecx}^{\vecm}:=
 \prod\del_{x_i}^{m_i}$
and
$\del_{\vecxi}^{\vecm}:=
 \prod\del_{\xi_i}^{m_i}$.
We put $N_i(k):=\bigl|S_i(k)\bigr|$
and $N(k_1,k_2):=N_1(k_1)\times N_2(k_2)$.

Let $V$ be a vector space.
For $f\in C^{\infty}(U_1\times U_2,V)$,
we set
\[
 D_{\vecx}^{k_1}D_{\vecxi}^{k_2}(f):=\Bigl(
 \del_{\vecx}^{\vecm_1}\del_{\vecxi}^{\vecm_2}f\,\Big|\,
 (\vecm_1,\vecm_2)\in S(k_1,k_2)
 \Bigr)
 \in C^{\infty}\bigl(U_1\times U_2,V^{N(k_1,k_2)}\bigr).
\]

Formally,
we set $D^0f:=f\in C^{\infty}(U_1\times U_2,V)$.
We use similar notations
for the functions on $T_0\times U_2$
and $[0,1]\times U_2$.

\subsubsection{Orthonormal frame}

Let $E$ be a topologically trivial $C^{\infty}$-vector bundle on
$T_0\times U_2$
with a hermitian metric $h$ and a unitary connection
$\nabla$.
We set $r:=\rank E$.
Let $F$ denote the curvature of $\nabla$.
For any frame $\vecv$ of $E$,
let $A^{\vecv}=\sum_{i=1}^n A_i^{\vecv}\,dt_i$
denote the connection form of $\nabla$
with respect to $\vecv$.
We put $\lefttop{1}A^{\vecv}:=A_1^{\vecv}\,dt_1+A^{\vecv}_2\,dt_2$
and $\lefttop{2}A^{\vecv}:=
 \sum_{i=3}^nA^{\vecv}_idt_i$.
Similarly $F^{\vecv}=\sum F^{\vecv}_{ij}\,dt_i\,dt_j$
denote the curvature form with respect to $\vecv$.

Fix a positive number $M$.
Let $\epsilon$ be a small positive number.
Assume that $|D^{k_1}_{\vecx}D^{k_2}_{\vecxi}F|_{h}\leq \epsilon$
for any $k_1,k_2\leq M$.

\begin{lem}
\label{lem;11.12.28.10}
If $\epsilon$ is sufficiently small,
there exist an orthonormal frame $\vecv$ of $(E,h)$ 
on $T_0\times U_2$
and anti-hermitian matrices 
$\Lambda^{(x)},\Lambda^{(y)}$
such that the following holds:
\begin{description}
\item[(A1)]
 For $\kappa=x,y$,
 there exist $0\leq \theta_{\kappa}<2\pi$
 such that any eigenvalue $\sqrt{-1}\alpha$
 of $\Lambda^{(\kappa)}$ 
 satisfies
 $|\alpha-\theta_{\kappa}|\leq\pi(2r-1)/2r$. 
 They satisfy 
 $\bigl|
 \bigl[\Lambda^{(x)},\Lambda^{(y)}\bigr]
 \bigr|\leq C\epsilon$.
\item[(A2)]
 $\bigl|
 \lefttop{1}A^{\vecv}-\Lambda
 \bigr|\leq C\epsilon$,
 and
 $\bigl|D^{k_1}_{\vecx}D^{k_2}_{\vecxi}
 (\lefttop{1}A^{\vecv})\bigr|\leq C\epsilon$
 for any $0\leq k_1\leq M$ and $0\leq k_2\leq M$
 with $(k_1,k_2)\neq (0,0)$,
 where 
$\Lambda=\Lambda^{(x)}\,dx+\Lambda^{(y)}\,dy$.
\item[(A3)]
 $\bigl|
 D^{k_1}_{\vecx}D^{k_2}_{\vecxi}(\lefttop{2}A^{\vecv})
 \bigr|\leq C\epsilon$ for any $0\leq k_1,k_2\leq M$.
\end{description}
Here, the constant $C$ may depend only on $r$ and $M$.
\end{lem}
\pf
We shall indicate an outline of the construction,
although it is elementary.
We say that a quantity $\nbigp$ is $O(\epsilon)$,
if $\nbigp\leq C\epsilon$ for some
 constant $C$ which may depend only on 
$r$ and $M$.
Let $[a,b]_{\seisuu}$ denote the set of integers $k$
such that $a\leq k\leq b$.
For $j\geq 1$,
let $H_j$ be the subset of $U_1\times U_2$
determined by the condition $t_i=0$ $(i\in[1,j]_{\seisuu})$.
We set $H_0:=U_1\times U_2$.

Let $\vecu$ be an orthonormal frame of
$\pi^{\ast}(E,h)$ on $U_1\times U_2$
satisfying 
$\nabla_{t_i}\vecu=0$ on $H_{i-1}$
for any $i$.
We have $A^{\vecu}_p=0$ on $H_{p-1}$
by the construction.
For $j<p$,
we have $\del_{t_j}A^{\vecu}_p=F^{\vecu}_{jp}$
on $H_{j-1}$.
For a monomial $P$ of $\del_{t_{j+1}},\ldots,\del_{t_n}$,
we have 
$\del_{t_j}^{\alpha+1}P A^{\vecu}_p
=\del_{t_j}^{\alpha}PF^{\vecu}_{jp}$
on $H_{j-1}$.
Hence, 
for $j\leq p$
and for a monomial $\Ptilde=\prod_{i=j}^n\del_{t_i}^{m_i}$
satisfying
$m_1+m_2\leq M$ and $\sum_{i>2} m_i\leq M$,
we obtain 
$\Ptilde A^{\vecu}_{p}=O(\epsilon)$
on $H_{j-1}$ 
by a descending induction.
In particular, we obtain
$D_{\vecx}^{k_1}D_{\vecxi}^{k_2}A^{\vecu}_p=O(\epsilon)$
for any $(k_1,k_2)\in[0,M]^2_{\seisuu}$.

Let $G^{(x)}:H_{1}\lrarr U(r)$ be determined by 
$\vecu_{|(1,y,\vecxi)}
=\vecu_{|(0,y,\vecxi)}G^{(x)}(y,\vecxi)$,
where $U(r)$ denotes the $r$-th unitary group.
By the equation
\[
\del_{t_i}G^{(x)}(t_2,\ldots,t_n)
-G^{(x)}(t_2,\ldots,t_n)\,A^{\vecu}_{i|(1,t_2,\ldots,t_n)}
+A^{\vecu}_{i|(0,t_2,\ldots,t_n)}\,G^{(x)}(t_2,\ldots,t_n)=0, 
\]
we obtain 
$|D^1_{\vecx}G^{(x)}|
+|D^1_{\vecxi}G^{(x)}|
=O(\epsilon)$.
By an easy induction,
we obtain
$|D^{k_1}_{\vecx}D^{k_2}_{\vecxi}G^{(x)}|=O(\epsilon)$
for any $(k_1,k_2)\in[0,M]_{\seisuu}^2\setminus\{(0,0)\}$.
We also have
$\bigl|
 G^{(x)}(y,\vecxi)
-
 G^{(x)}(y',\vecxi')
 \bigr|
=O(\epsilon)$.

Let $G^{(y)}(x,\vecxi)$ be determined by
$\vecu_{|(x,1,\vecxi)}
=\vecu_{|(x,0,\vecxi)}\,G^{(y)}(x,\vecxi)$.
Similarly,
we have
$\bigl|D^{k_1}_{\vecx}D_{\vecxi}^{k_2}G^{(y)}\bigr|=O(\epsilon)$
for any $(k_1,k_2)\in[0,M]^2_{\seisuu}\setminus\{(0,0)\}$,
and $\bigl|G^{(y)}(x,\vecxi)-G^{(y)}(x',\vecxi')\bigr|=O(\epsilon)$.
Because 
$G^{(y)}(0,\vecxi)G^{(x)}(1,\vecxi)=G^{(x)}(0,\vecxi)G^{(y)}(1,\vecxi)$,
we have
$\bigl[G^{(y)}(0,0),G^{(x)}(0,0)\bigr]=O(\epsilon)$.
We set 
$\Gtilde^{(y)}:=G^{(y)}(0,0)$
and $\Gtilde^{(x)}:=G^{(x)}(0,0)$.

\vspace{.1in}
Let $\nbigi_{\kappa}$ denote the set of the eigenvalues of
$\Gtilde^{(\kappa)}$ for $\kappa=x,y$.
Let $d_{S^1}$ denote the standard distance on 
$S^1=\{e^{\sqrt{-1}\theta}\,|\,\theta\in\real\}$
induced by the metric $d\theta\,d\theta$.
There exist
$\gamma_{\kappa}\in S^1$
such that
$d_{S^1}(\gamma_{\kappa},\gamma)\geq\pi/(2r)$
for any $\gamma\in\nbigi_{\kappa}$.
Let $\theta_{\kappa}$ be determined by
$e^{\sqrt{-1}\theta_{\kappa}}=-\gamma_{\kappa}$
and $0\leq \theta_{\kappa}<2\pi$.
For any $\gamma_i\in\nbigi_{\kappa}$,
we can take $\alpha_i$ satisfying
(i) $e^{\sqrt{-1}\alpha_i}=\gamma_i$,
(ii) $|\theta_{\kappa}-\alpha_i|\leq\pi(2r-1)/2r$.
We remark that,
for any $\gamma_i,\gamma_j\in\nbigi_{\kappa}$,
we have
\begin{equation}
\label{eq;13.5.28.1}
 |\alpha_i-\alpha_j|
=O\bigl(
 |\gamma_i-\gamma_j|
 \bigr)
\end{equation}

We have the eigen decompositions
$\cnum^r=\bigoplus_{\gamma_i\in\nbigi_{\kappa}}
 V^{(\kappa)}_{\gamma_i}$
for $\Gtilde^{(\kappa)}$.
We set
$\Lambda^{(\kappa)}
=\bigoplus_{\gamma_i\in\nbigi_{\kappa}}
 \sqrt{-1}\alpha_i\,\id_{V^{(\kappa)}_{\alpha}}$.
By the construction,
we have
$\exp(\Lambda^{(\kappa)})=\Gtilde^{(\kappa)}$.

\begin{lem}
\label{lem;13.5.28.2}
We have
$\bigl[\Lambda^{(x)},
 \Lambda^{(y)}
 \bigr]=O(\epsilon)$.
\end{lem}
\pf
According to the decomposition
$\cnum^r=\bigoplus_{\gamma_i\in\nbigi_x} V^{(x)}_{\gamma_i}$,
we have the decomposition
$\Gtilde^{(y)}=\sum_{\gamma_i,\gamma_j\in\nbigi_{x}} 
\Gtilde^{(y)}_{\gamma_i,\gamma_j}$,
where
$\Gtilde^{(y)}_{\gamma_i,\gamma_j}
\in \Hom(V^{(x)}_{\gamma_j},V^{(x)}_{\gamma_i})$.
We have
$ (\gamma_i-\gamma_j)
 \Gtilde^{(y)}_{\gamma_i,\gamma_j}
=O(\epsilon)$.
By using (\ref{eq;13.5.28.1}),
we obtain
$\bigl[
 \Gtilde^{(y)},\Lambda^{(x)}
\bigr]=O(\epsilon)$.
By using a similar consideration,
we obtain
$\bigl[
 \Lambda^{(y)},\Lambda^{(x)}
 \bigr]=O(\epsilon)$.
\hfill\qed

\vspace{.1in}

Let us return to the proof of Lemma \ref{lem;11.12.28.10}.
We put 
$g^{(x)}(x):=\exp\bigl(-x\,\Lambda^{(x)}\bigr)$,
$g^{(y)}(y):=\exp\bigl(-y\,\Lambda^{(y)})$,
and
$g(x,y):=g^{(x)}(x)\,g^{(y)}(y)$.
We obtain an orthonormal frame
$\vecu':=\vecu\,g(x,y)$
of $\pi^{\ast}(E,h)$.
Let $A':=A^{\vecu'}$.
We have
$\bigl|A'-\Lambda\bigr|=O(\epsilon)$
and
$|D_{\vecx}^{k_1}D^{k_2}_{\vecxi}A'|=O(\epsilon)$ 
for any $(k_1,k_2)\in[0,M]_{\seisuu}^2\setminus\{(0,0)\}$.

Let $G^{\prime(x)}(y,\vecxi)$ 
and $G^{\prime\,(y)}(x,\vecxi)$
be determined by
\[
 \vecu'_{|(1,y,\vecxi)}
=\vecu'_{|(0,y,\vecxi)}
 \,G^{\prime\,(x)}(y,\vecxi),
\quad\quad
 \vecu'_{|(x,1,\vecxi)}
=\vecu'_{|(x,0,\vecxi)}
 \,G^{\prime\,(y)}(x,\vecxi).
\]
We have
$G^{\prime(x)}(y,\vecxi)=
 g^{(y)}(y)^{-1}\,G^{(x)}(y,\vecxi)\,
 (\Gtilde^{(x)})^{-1}\,g^{(y)}(y)$
and hence
$\bigl|
 G^{\prime\,(x)}-1\bigr|=O(\epsilon)$.
We have
\[
 dG^{\prime\,(x)}=
 g^{(y)}(y)^{-1}\,dG^{(x)}(y,\vecxi)\,(\Gtilde^{(x)})^{-1}\,g^{(y)}(y)
-\bigl[
 g^{(y)}(y)^{-1}dg^{(y)}(y),\,
 (G^{\prime\,(x)}-1)
 \bigr].
\]
Hence, we have
$\bigl|D^1_{y}G^{\prime(x)}\bigr|=O(\epsilon)$
and
$\bigl|D^1_{\vecxi}G^{\prime(x)}\bigr|=O(\epsilon)$.
By an easy induction,
we obtain
$\bigl|D_{y}^{k_1}D_{\vecxi}^{k_2}
 G^{\prime(x)}\bigr|=O(\epsilon)$
for $(k_1,k_2)\in[0,M]_{\seisuu}^2\setminus\{(0,0)\}$.
We have
$G^{\prime\,(y)}
=g^{(x)}(x)^{-1}G^{(y)}(x,\vecxi)g^{(x)}(x)
(\Gtilde^{(y)})^{-1}$.
We obtain
\[
 G^{\prime\,(y)}-1
=g^{(x)}(x)^{-1}
 \Bigl(
 G^{(y)}(x,\vecxi)
 (\Gtilde^{(y)})^{-1}-1
\Bigr)
g^{(x)}(x)
-g^{(x)}(x)^{-1}
 G^{(y)}(x,\vecxi)
 \bigl[
 (\Gtilde^{(y)})^{-1},\,
 g^{(x)}(x)
 \bigr]
=O(\epsilon).
\]
As in the case $\kappa=x$,
we also obtain
$\bigl|D_{y}^{k_1}D_{\vecxi}^{k_2}
 G^{\prime(x)}\bigr|=O(\epsilon)$
for $(k_1,k_2)\in[0,M]_{\seisuu}^2\setminus\{(0,0)\}$.

\vspace{.1in}

Let $\chi(x)$ be a non-negative valued
$C^{\infty}$-function on
$[0,1]$ such that 
$\chi(x)=0$ $(x\leq 1/3)$
and $\chi(x)=1$ $(x\geq 2/3)$.
We put 
$h_2(x,y,\vecxi):=\chi(x)\,
 \exp^{-1}\bigl(
 G^{\prime(x)}(y,\vecxi)
 \bigr)$.
By construction,
we have
$|D_{\vecx}^{k_1}D_{\vecxi}^{k_2}h_2|=O(\epsilon)$
for $(k_1,k_2)\in[0,M]_{\seisuu}^2$.

Let $g_2:=\exp(h_2)$,
and we set
$\vecu'':=\vecu'\,g_2$.
Let $A''=A^{\vecu''}$.
We have
$A''=g_2^{-1}A'g_2+g_2^{-1}dg_2$.
Hence, we have
$|\lefttop{1}A''-\Lambda|=O(\epsilon)$,
and
$|D_{\vecx}^{k_1}D_{\vecxi}^{k_2} (\lefttop{1}A'')|=O(\epsilon)$
for $(k_1,k_2)\in[0,M]_{\seisuu}^2\setminus\{(0,0)\}$.
We also have
$|D^{k_1}_{\vecx}D^{k_2}_{\vecxi} (\lefttop{2}A'')|=O(\epsilon)$
for $(k_1,k_2)\in[0,M]_{\seisuu}^2$.

We put
$G^{\prime\prime(y)}(x,\vecxi):=g_2(x,0,\vecxi)^{-1}\,
 G^{\prime(y)}(x,\vecxi)\,g_2(x,1,\vecxi)$.
We have 
$\vecu''_{|(x,1,\vecxi)}=
 \vecu''_{|(x,0,\vecxi)}\,G^{\prime\prime(y)}(x,\vecxi)$.
We have
$\bigl|
 G^{\prime\prime(y)}(x,\vecxi)-1\bigr|=O(\epsilon)$,
and
$|D_{\vecx}^{k_1}D_{\vecxi}^{k_2}G^{\prime\prime(y)}|=O(\epsilon)$
for $(k_1,k_2)\in[0,M]_{\seisuu}^2\setminus\{(0,0)\}$.

We put
$g_3:=\exp\Bigl(
 \chi(y)\,
 \exp^{-1}\bigl(
 G^{\prime\prime(y)}(x,\vecxi)
 \bigr)\Bigr)$,
and $\vecv:=\vecu''\,g_3$.
Then, it naturally gives an orthonormal frame of
$(E,h)$ on $T_0\times U_2$.
By construction,
we have the desired estimate
for the connection form $A^{\vecv}$.
Thus, the proof of Lemma \ref{lem;11.12.28.10}
is finished.
\hfill\qed

\subsubsection{Partially almost holomorphic frame}
\label{subsection;13.1.8.10}

We identify the $C^{\infty}$-manifolds
$T_0:=\real^2/\seisuu^2$ and $T$
by the diffeomorphism $T_0\simeq T$
given by $(x,y)\longmapsto x+\tau y=z$.
We have the description
$\Lambda=\Gamma d\zbar-\lefttop{t}\Gammabar\,dz$,
where $\Lambda$ is as in Lemma \ref{lem;11.12.28.10}.
Let $\nabla_{\zbar}:=\nabla(\del_{\zbar})$
and $\nabla_z:=\nabla(\del_z)$.
For any frame $\vecw$,
let $A^{\vecw}_z$ and $A^{\vecw}_{\zbar}$
be determined by
$\nabla_{z}\vecw=\vecw\,A^{\vecw}_{z}$
and
$\nabla_{\zbar}\vecw=\vecw\,A^{\vecw}_{\zbar}$,
respectively.
Let $H(h,\vecw)$ denote a function from $T\times U_2$
to the space of $r$-th positive definite hermitian matrices,
whose $(i,j)$-entries are $h(w_i,w_j)$.
When a function $f$ on $T\times U_2$ is regarded as a function 
$\ftilde:U_2\lrarr L_k^p(T)$,
we obtain an $\real_{\geq 0}$-valued function
 $\|f\|_{L^p_k}(\vecxi):=\|\ftilde(\vecxi)\|_{L_k^p(T)}$
on $U_2$.

\begin{prop}
\label{prop;12.6.24.3}
If $\epsilon>0$ is sufficiently small,
there exists a frame $\vecu$ of $E$ on $T\times U_2$
with the following property:
\begin{itemize}
\item
 $A^{\vecu}_{\zbar}$ is constant along the $T$-direction,
 and $|A_{\zbar}^{\vecu}-\Gamma|=O(\epsilon)$.
\item
 $\bigl\|A^{\vecu}_{z}+\lefttop{t}\Gammabar\bigr\|_{L_M^p}
 =O(\epsilon)$
and 
 $\bigl\|D^k_{\vecxi}A^{\vecu}_{z}\bigr\|_{L_M^p}
 =O(\epsilon)$ for $k\in[1,M]_{\seisuu}$.
\item
 $\bigl\|D_{\vecxi}^k (\lefttop{2}A^{\vecu})\bigr\|_{L_M^p}
 =O(\epsilon)$ for $k\in [0,M]_{\seisuu}$.
\end{itemize}
Moreover, 
$\|H(h,\vecu)-I\|_{L_{M+1}^p}=O(\epsilon)$
and
$\|D_{\vecxi}^kH(h,\vecu)\|_{L_{M+1}^p}=O(\epsilon)$
for $k\in[1,M]_{\seisuu}$,
where $I$ denotes the identity matrix.
\end{prop}
\pf
Let $\vecv$ be the orthonormal frame
as in Lemma \ref{lem;11.12.28.10}.
We have
$\nabla_{\zbar}\vecv
=\vecv
 (\Gamma+N)$,
where $\|D^k_{\vecxi}N\|_{L_M^p}=O(\epsilon)$.

\begin{lem}
\label{lem;13.10.12.1}
We have a decomposition
 $\Gamma=\Gamma_0+\Gamma_1$
 such that
 (i) $[\Gamma_0,\Gamma_0^{\dagger}]=0$
 and $\Sp(\Gamma_0)=\Sp(\Gamma)$,
 (ii) $|\Gamma_1|=O(\epsilon^{1/2})$.
Moreover, 
if $\delta>0$ is sufficiently smaller than $1/r$,
but independently from $\epsilon$,
there exists $\zeta_0\in\cnum$
such that
$\Sp(\Gamma)$ is contained in
$K_1(L,\zeta_0)$.
(See Proposition {\rm\ref{prop;11.12.18.10}}
 for $K_1(L,\zeta_0)$.)
\end{lem}
\pf
We give only an indication.
With an appropriate change of orthonormal base,
we may assume that $\Gamma$ is upper triangular.
By the base,
we identify matrices and endomorphisms.
Let $\Gamma_0$ be the diagonal part,
and we put $\Gamma_1:=\Gamma-\Gamma_0$.
By the construction, the condition (i)
is satisfied.
Let $\gamma_{ij}$ denote the $(i,j)$-entry of $\Gamma$.
Then, the $(k,k)$-entries of
$[\Gamma,\Gamma^{\dagger}]$
is 
$\sum_{i>k}|\gamma_{k,i}|^2
-\sum_{i<k}|\gamma_{k,i}|^2$.
Then, we obtain the desired estimate for $\Gamma_1$
from 
$[\Gamma,\Gamma^{\dagger}]=O(\epsilon)$,
which follows from the condition
{\bf (A1)} in Lemma \ref{lem;11.12.28.10}.
Thus, we obtain the first condition.
We also obtain the second condition
from {\bf (A1)} in Lemma \ref{lem;11.12.28.10}.
\hfill\qed

\vspace{.1in}

By Proposition \ref{prop;11.12.18.10},
if $\epsilon$ is sufficiently small,
there exist functions
$a:U_2\lrarr L_{M+1}^p\bigl(M_{r}(\cnum)\bigr)_0$
and 
$b:U_2\lrarr M_{r}(\cnum)$ satisfying the following:
\begin{itemize}
\item
 $\|D_{\vecxi}^ka\|_{L_{M+1}^p}=O(\epsilon)$ 
    for $k\in[0,M]_{\seisuu}$,
and
 $|D_{\vecxi}^k b|=O(\epsilon)$
 for $k\in [0,M]_{\seisuu}$.
\item
$(1+a)\bullet
 \bigl(
 \nabla_{\zbar,0}+(\Gamma+b)\,d\zbar
 \bigr)
=\nabla_{\zbar}$,
where $\nabla_{\zbar,0}$ is given by
$\nabla_{\zbar,0}\vecv=0$.
\end{itemize}
Let $\vecu:=\vecv(1+a)$.
By construction, we have
$\nabla_{\zbar}\vecu=
 \vecu\bigl(\Gamma+b\bigr)$.
The other estimates for
$A_{\zbar}^{\vecu}$
and $\lefttop{2}A^{\vecu}$
are also satisfied.
Because $H(h,\vecu)=\lefttop{t}(1+a)\overline{(1+a)}$,
we obtain the estimate for $H(h,\vecu)$.
\hfill\qed

\begin{rem}
If $A_{\zbar}^{\vecw}$ is constant along the $T$-direction,
such a frame $\vecw$ is called 
a partially almost holomorphic frame,
in this paper.
\hfill\qed
\end{rem}

\subsubsection{Spectra}

Let $E_{\vecxi}$ denote the holomorphic bundle on $T$
given by $E_{|T\times\vecxi}$ with 
$\nabla_{\zbar|T\times\vecxi}$.
According to Lemma \ref{lem;12.6.23.101},
if $\epsilon$ is sufficiently small,
$E_{\vecxi}$ are semistable of degree $0$
for any $\vecxi\in U_2$.
We have the spectrum
$\Sp(E_{\vecxi})\subset T^{\lor}$.
We regard it as a point in
$\Sym^rT^{\lor}$.
The point is denoted by $[\Sp(E_{\vecxi})]$.
Let $\Gamma$ be as in \S\ref{subsection;13.1.8.10}.
The eigenvalues of $\Gamma$
give a point in $\Sym^{r}\cnum$,
denoted by $[\Sp(\Gamma)]$.
The quotient map
$\Phi:\cnum\lrarr T^{\lor}$
induces
$\Sym^{r}\cnum\lrarr \Sym^rT^{\lor}$,
denoted by $\Phi$.
Recall that $\Sym^rT^{\lor}$
is naturally a smooth complex manifold.
Let $d_{\Sym^rT^{\lor}}$
be a distance induced by 
a $C^{\infty}$-Riemannian metric.

\begin{cor}
\label{cor;12.6.24.2}
There exist $\epsilon_0>0$ and $C>0$,
depending only on $r$,
such that the following holds if $\epsilon\leq \epsilon_0$:
\[
 d_{\Sym^r T^{\lor}}\Bigl(
 [\Sp(E_{\vecxi})],
 \Phi[\Sp(\Gamma)]\Bigr)
\leq
 C\epsilon
\]
In particular,
for $\vecxi,\vecxi'\in U_2$,
we have
$d_{\Sym^r T^{\lor}}\Bigl(
 [\Sp(E_{\vecxi})],
 [\Sp(E_{\vecxi'})]
 \Bigr)
\leq
 2C\epsilon$.
\end{cor}
\pf
Let $\vecu$ be a frame as in Proposition \ref{prop;12.6.24.3}.
Recall that $\Sym^r\cnum$ is naturally a complex manifold.
We take a distance $d_{\Sym^r\cnum}$
induced by a $C^{\infty}$-Riemannian metric.
We have 
$d_{\Sym^r\cnum}\bigl(
 [\Sp(\Gamma)],\,[\Sp(A_{\zbar}^{\vecu})]
 \bigr)
\leq C_1\epsilon$.
There exists $\zeta_0\in\cnum$
such that 
$\Sp(\Gamma)$
and $\Sp(A_{\zbar}^{\vecu})$
are contained in
$K_1(L,\zeta_0)$.
Note that
the restriction of $\Phi$
to $\Sym^r K_1(L,\zeta_0)$
is Lipschitz continuous,
and the Lipschitz constant is uniform for $\zeta_0$.
Then, the claim of the corollary follows.
\hfill\qed

\subsection{Estimates}
\label{subsection;12.11.5.10}

\subsubsection{Preliminary}

We continue to use the setting in \S\ref{subsection;12.6.24.10}.
We impose additional assumptions.

\begin{assumption}\mbox{{}}
\label{assumption;12.6.24.21}
\begin{itemize}
\item
 We are given a finite subset 
 $Z\subset \cnum$
 and a positive number $\rho>0$
 with the following property:
\begin{itemize}
\item
 $Z$ is contained in $K_1(L,\zeta_0)$,
 where $K_1(L,\zeta_0)$ is as in {\rm\S\ref{subsection;12.11.5.1}}
 for some appropriate $\delta>0$.
\item
 For any distinct points $\nu_1,\nu_2\in Z$,
 $d_{\cnum}(\nu_1,\nu_2)>100 r^2\rho$.
\item
 For any $\kappa\in \Sp(E_{\vecxi})$,
 there exists $\nu\in Z$ such that
 $d_{T^{\lor}}(\Phi(\nu),\kappa)<\rho$,
 where $\Phi:\cnum\lrarr T^{\lor}$
 denotes the projection.
\item
$\epsilon$ is sufficiently small so that
$E_{\vecxi}$ is semistable of degree $0$
for any $\vecxi\in U_2$.
We also assume that
$\epsilon$ is sufficiently smaller than $\rho^2$.
\hfill\qed
\end{itemize}
\end{itemize}
\end{assumption}

We have the spectral decomposition
$E_{\vecxi}=\bigoplus_{\nu'\in T^{\lor}} E_{\vecxi,\nu'}$.
Let $E_{\nu,\vecxi}$ be the direct sum of
 $E_{\vecxi,\nu'}$,
where $\nu'$ is contained in 
a $\rho$-ball of $\Phi(\nu)$.
We obtain a decomposition
$E_{\vecxi}=
 \bigoplus_{\nu\in Z}E_{\nu,\vecxi}$.
It induces a $C^{\infty}$-decomposition
$E=\bigoplus_{\nu\in Z}E_{\nu}$,
which is compatible with $\nabla_{\zbar}$.
We may assume that 
the partially almost holomorphic frame $\vecu$
in Proposition \ref{prop;12.6.24.3} is compatible
with the decomposition.

\vspace{.1in}

We have the decomposition
$\nabla_{\zbar}=\nabla_{\zbar,0}+f$
such that
(i) $(E,\nabla_{\zbar,0})_{|T\times\{\vecxi\}}$
are holomorphically trivial
 for any $\vecxi\in U_2$,
(ii) $\nabla_{\zbar}(f)=0$,
(iii) $\Sp(f)$ is contained in 
the union of the $\rho$-balls 
around $\nu\in Z$.
For each $\vecxi\in U_2$,
we obtain the vector space $\nbigv_{\vecxi}$
of the holomorphic global sections of
$(E,\nabla_{\zbar,0})_{|T\times\{\vecxi\}}$.
It is easy to see that 
$\nbigv_{\vecxi}$ $(\vecxi\in U_2)$ naturally gives 
a $C^{\infty}$-vector bundle $\nbigv$
on $U_2$,
and that we have a natural isomorphism
$p^{\ast}\nbigv\simeq E$ as $C^{\infty}$-bundles.
We identify them by the isomorphism.
A $C^{\infty}$-section $s$ of $p^{\ast}\nbigv$
is constant along the $T$-direction,
if and only if $\nabla_{\zbar,0}s=0$
under the identification.
It can be regarded as a section of $\nbigv$.
We have the decomposition
$\nbigv=\bigoplus _{\nu\in Z}\nbigv_z$,
corresponding to $E=\bigoplus_{\nu\in Z}E_{\nu}$.

\subsubsection{Spaces of functions}

Let $C^M_{\vecxi}L_{M,\vecx}^p$ denote the space of
$C^M$-functions $U_2\lrarr L_M^p(T)$.
Let $C_{\vecxi}^ML_{M,\vecx}^p(E)$ denote the 
sections $f=\sum f_iu_i$ of $E$
such that $f_i\in C_{\vecxi}^ML_{M,\vecx}^p$,
where $\vecu=(u_i)$ is a frame as
in Proposition \ref{prop;12.6.24.3}.
It is independent of the choice of $\vecu$.
We have the naturally defined integration
$\int_T:
 C_{\vecxi}^ML_{M,\vecx}^p(E)
 \lrarr C^M(U_2,\nbigv)$.
The kernel is denoted by 
$C_{\vecxi}^ML_{M,\vecx}^p(E)_0$.
Similar spaces are defined for
$\End(E)$ and $\Hom(E_i,E_j)$.
We set
\[
 C^M_{\vecxi}L_{M,\vecx}^p\bigl(\End(E)\bigr)^{\circ}:=
 \bigoplus_{\nu}
 C^M\bigl(U_2,\End(\nbigv_{\nu})\bigr),
\]
\[
 C^M_{\vecxi}L_{M,\vecx}^p\bigl(\End(E)\bigr)^{\bot}:=
\bigoplus_{\nu}
 C^M_{\vecxi}L_{M,\vecx}^p\bigl(\End(E_{\nu})\bigr)_0
\oplus
 \bigoplus_{\nu\neq \mu}
 C^M_{\vecxi}L_{M,\vecx}^p\bigl(\Hom(E_{\nu},E_{\mu})\bigr).
\]
We have a decomposition
$C^M_{\vecxi}L_{M,\vecx}^p\bigl(\End(E)\bigr)=
C^M_{\vecxi}L_{M,\vecx}^p\bigl(\End(E)\bigr)^{\circ}
\oplus
C^M_{\vecxi}L_{M,\vecx}^p\bigl(\End(E)\bigr)^{\bot}$.
For any $s\in C^M_{\vecxi}L_{M,\vecx}^p\bigl(\End(E)\bigr)$,
the corresponding decomposition is denoted by
$s=s^{\circ}+s^{\bot}$.
We use similar notations for sections of
$\End(E)\otimes \Omega^{i,j}_T$.

\subsubsection{Some estimates}
\label{subsection;12.6.24.22}

Let $\vecu$ be a frame as in Proposition \ref{prop;12.6.24.3}.
We set $H(h,\vecu)_{i,j}:=h(u_i,u_j)$,
and we obtain a function $H(h,\vecu)$
from $T\times U_2$ to the space $\nbigh$ of 
positive definite hermitian $r$-th matrices.
Each entry is $C^M_{\vecxi}L_{M,\vecx}^p$-class.
Let $H_1$ be a function of $U_2$ to $\nbigh$
determined by
$(H_1)^2=\int_{T}H(h,\vecu)$.
Then, we have 
$|H_1-I|=O(\epsilon)$
and
$|D^k_{\vecxi}H_1|=O(\epsilon)$
for $k\in[1,M]_{\seisuu}$.
Note that $\vecu':=\vecu\,H_1$
also has the property in Proposition \ref{prop;12.6.24.3}.
So, we may assume that $\int_TH(h,\vecu)=I$
from the beginning.

We set $\gtilde:=\overline{H(h,\vecu)}$.
We have 
$\|\gtilde-I\|_{L_{M+1}^p}=O(\epsilon)$,
$\|D_{\vecxi}^k\gtilde\|_{L_{M+1}^p}=O(\epsilon)$
$(k\in [1,M]_{\seisuu})$,
and $\int_T\gtilde=I$.

\begin{lem}
\label{lem;11.12.29.2}
There exist $C>0$ and $\epsilon_0>0$,
such that 
$\|\gtilde-I\|_{L_{M+2}^p}\leq
 C\,\|F_{z\zbar}^{\bot}\|_{L_M^p}$
holds if $\epsilon<\epsilon_0$.
In particular,
$\sup_{T\times\{\vecxi\}}\bigl|\gtilde-I\bigr|
 \leq
 C'\|F_{z\zbar}^{\bot}\|_{L^2}$
for some $C'>0$.
\end{lem}
\pf
We put $B:=A^{\vecu}_{\zbar}$.
Let $\Gamma_2$ be the diagonal matrix
whose $(i,i)$-entry $\nu_i$
is determined by $u_i\in E_{\nu_i}$.
Let $\Gamma$ be as in \S\ref{subsection;13.1.8.10},
which is decomposed 
$\Gamma=\Gamma_0+\Gamma_1$
as in Lemma \ref{lem;13.10.12.1}.
We have
$|\Gamma_0-\Gamma_2|\leq  r\rho$
and $|\Gamma_1|=O(\epsilon^{1/2})$.
We have $|A^{\vecu}_{\zbar}-\Gamma|=O(\epsilon)$.
Hence, 
if $\epsilon$ is sufficiently small,
we may have
$|A^{\vecu}_{\zbar}-\Gamma_2|\leq 
2r\rho$.

We have
$A^{\vecu}_z=-\gtilde^{-1}
 \bigl(\lefttop{t}\overline{B}\bigr)
 \gtilde
+\gtilde^{-1}\del_z\gtilde$.
Let $\nbigb_{z\zbar}$ be the matrix-valued function
determined by
$F_{z\zbar}\vecu=\vecu\,\nbigb_{z\zbar}$.
We have
$\nbigb_{z\zbar}
=\del_zA^{\vecu}_{\zbar}-\del_{\zbar}A^{\vecu}_z
+\bigl[A^{\vecu}_z,A^{\vecu}_{\zbar}\bigr]$.
Hence, we have the following equation:
\begin{equation}
\label{eq;11.12.29.3}
 \nbigb_{z\zbar}=
 \bigl[
 \gtilde^{-1}\lefttop{t}\overline{B}\gtilde,
 \gtilde^{-1}\del_{\zbar}\gtilde
 \bigr]
-\gtilde^{-1}\del_{\zbar}\del_z(\gtilde)
+(\gtilde^{-1}\del_{\zbar}\gtilde)\,
 (\gtilde^{-1}\del_{z}\gtilde)
-\bigl[\gtilde^{-1}\lefttop{t}\overline{B}\gtilde,B\bigr]
-\bigl[B,\gtilde^{-1}\del_z\gtilde\bigr]
\end{equation}
Let $b:=\gtilde-I$.
We have a polynomial 
$Q(t_1,t_2,t_3,t_4,t_5,t_6) 
=\sum
 Q_{j_1,\ldots,j_m}
 t_{j_1}t_{j_2}\cdots t_{j_m}$
in non-commutative variables $t_i$
such that (i)
if $Q_{j_1,\ldots,j_m}\neq 0$
then $m_1+m_2+m_3\geq 2$,
where $m_i=\{k\,|\,j_k=i\}$,
(ii) the following equality holds:
\begin{equation}
 \bigl(\del_{\zbar}+\ad(B)\bigr)\circ
 \bigl(\del_z-\ad(\lefttop{t}\Bbar)\bigr)b
=-\gtilde\nbigb_{z\zbar}
-\bigl[\lefttop{t}\Bbar,B\bigr]
+Q\bigl(b,\del_zb,\del_{\zbar}b,
 (1+b)^{-1},B,\lefttop{t}\overline{B}\bigr).
\end{equation}
By taking the $\bot$-part,
we obtain the following
\begin{equation}
 \label{eq;12.6.24.20}
 \bigl(\del_{\zbar}+\ad(B)\bigr)\circ
 \bigl(\del_z-\ad(\lefttop{t}\Bbar)\bigr)b
=-(\gtilde\nbigb_{z\zbar})^{\bot}
+Q\bigl(b,\del_zb,\del_{\zbar}b,
 (1+b)^{-1},B,\lefttop{t}\overline{B}
 \bigr)^{\bot}.
\end{equation}
We obtain
\[
\|b\|_{L^p_{m+2}}
\leq
 C_2\|F^{\bot}_{z\zbar}\|_{L^p_m}
+C_2\epsilon
\|b\|_{L^p_{m+2}}
\]
Hence, we obtain
$\|b\|_{L_{m+2}^p}
\leq
 C_3\,\|F^{\bot}_{z\zbar}\|_{L_{m}^p}$.
\hfill\qed

\begin{lem}
\label{lem;11.12.29.3}
Let $a_1$ and $a_2$ be sections of $\End(E)_{|T\times\{\vecxi\}}$.
Assume that $a_1=a_1^{\bot}$
and $a_2=a_2^{\circ}$.
Then, we have 
\[
 \Bigl|
\int_{T}h(a_1,a_2)
 \Bigr|
\leq
 \|a_1\|_{L^2}\,\|a_2\|_{L^2}\,
 \bigl\|(F_{z\zbar}^{\bot})_{|T\times\{\vecxi\}}\bigr\|_{L^2}.
\]
\end{lem}
\pf
It follows from Lemma \ref{lem;11.12.29.2}
and $\overline{H(h,\vecu)}=\gtilde$.
\hfill\qed

\begin{lem}
\label{lem;11.12.30.2}
Let $\nbigp$ be an endomorphism of $E$,
and let $\nbigp^{\dagger}$ denote the adjoint
with respect to $h$.
Let $\nbigr$ (resp. $\nbigr^{\dagger}$)
be the matrix representing $\nbigp$ 
(resp. $\nbigp^{\dagger}$)
with respect to $\vecu$.
Then, we have
\[
 \bigl(\nbigr^{\dagger}\bigr)^{\circ}
=\bigl(
 \lefttop{t}\overline{\nbigr}
 \bigr)^{\circ}
+O\bigl(
 |\nbigr^{\bot}|\,\|F_{z\zbar}^{\bot}\|_{L^2}
 \bigr)
+O\bigl(
 \|F_{z\zbar}^{\bot}\|^2_{L^2}
 \bigl|
 \lefttop{t}\overline{\nbigr}^{\circ}
 \bigr|
 \bigr)
\]
\[
 (\nbigr^{\dagger})^{\bot}
=
 \bigl(
 \lefttop{t}\overline{\nbigr}
 \bigr)^{\bot}
+O\bigl(|\nbigr^{\bot}|\,
 \|F_{z\zbar}^{\bot}\|_{L^2}\bigr)
+O\bigl(
 \|F_{z\zbar}^{\bot}\|_{L^2}\,
 |\nbigr^{\circ}|\bigr)
\]
In particular,
we have
$|(\nbigr^{\dagger})^{\bot}|=
 |\nbigr^{\bot}|+O\bigl(
 |\nbigr|\,\|F^{\bot}_{z\zbar}\|_{L^2}
 \bigr)$.
\end{lem}
\pf
Let $H=H(h,\vecu)$.
We have
$\nbigr^{\dagger}=\overline{H^{-1}}
 (\lefttop{t}\overline{\nbigr})\,\overline{H}$.
Then, the claim follows from the estimate for $H$.
\hfill\qed

\begin{lem}
\label{lem;11.12.30.10}
For $k\in[1,M]_{\seisuu}$,
we have
$\|D_{\vecxi}^k\gtilde\|_{L_{m+2}^p}
=O\Bigl(
 \sum_{j=0}^k
 \|D_{\vecxi}^jF_{z\zbar}^{\bot}\|_{L_m^p}
 \Bigr)$.
\end{lem}
\pf
We obtain the estimate from (\ref{eq;12.6.24.20})
by a standard inductive argument.
\hfill\qed

\section{Estimates for $L^2$-Instantons}

\subsection{Preliminary}

Let $\tau$ be a complex number such that $\Image\tau>0$.
Let $T$ be a complex torus
obtained as the quotient of $\cnum$
by a lattice $\seisuu+\seisuu\,\tau$.
Let $z$ be the standard coordinate of $\cnum$.
It also gives a local coordinate of 
a small open subset in $T$,
once we fix a lift of the open subset in $\cnum$.
We shall use the metric $dz\,d\zbar$ 
for $\cnum$ and $T$
unless otherwise specified.

For any open subset $W\subset \cnum_w$,
we use the metric $dw\,d\wbar$ on $W$,
and the metric $dz\,d\zbar+dw\,d\wbar$
on $T\times W$ unless otherwise specified.
Let $\omega$ denote the associated K\"ahler form.
For $w\in W$,
we put $T_w:=T\times\{w\}\subset T\times W$.

Let $E$ be a complex $C^{\infty}$-vector bundle on $T\times W$
with a hermitian metric $h$ and a unitary connection $\nabla$.
Let $F(\nabla)$ denote the curvature of $\nabla$.
We shall often denote it simply by $F$.
The $(0,1)$-part and the $(1,0)$-part of $\nabla$ are 
denoted by $\delbar_E$ and $\del_E$,
respectively.
The restrictions of $(E,h)$ to $T_w$
are denoted by  $(E_w,h_w)$.

Recall that $(E,\nabla,h)$ is called an instanton,
if $\Lambda_{\omega} F(\nabla)=0$.
For the expression
$F(\nabla)=F_{z\zbar}dz\,d\zbar+F_{z\wbar}dz\,d\wbar
+F_{w\zbar}dw\,d\zbar+F_{w\wbar}dw\,d\wbar$,
the equation is $F_{z\zbar}+F_{w\wbar}=0$.
We have the following equalities:
\begin{equation}
 \label{eq;12.6.22.1}
 (\nabla_z\nabla_{\zbar}+\nabla_w\nabla_{\wbar})F_{w\wbar}
=-(\nabla_z\nabla_{\zbar}+\nabla_w\nabla_{\wbar})F_{z\zbar}
=[F_{z\wbar},F_{w\zbar}]
\end{equation}
\begin{equation}
 \label{eq;12.6.22.2}
 (\nabla_z\nabla_{\zbar}+\nabla_w\nabla_{\wbar})F_{z\wbar}
=2[F_{w\wbar},F_{z\wbar}],
\quad
 \bigl(
 \nabla_z\nabla_{\zbar}+\nabla_w\nabla_{\wbar}
 \bigr)F_{w\zbar}
=2[F_{w\zbar},F_{w\wbar}]
\end{equation}

\subsubsection{Hitchin's equivalence}

Let us recall the relation between harmonic bundles 
on an open subset $W\subset \cnum_w$
and instantons on $T\times W$ due to Hitchin.
Let $(E,\delbar_E,h,\theta)$ be a harmonic bundle on $W$.
Let $\nabla^0:=\delbar_E+\del_E$ be the Chern connection.
Let $\theta^{\dagger}$ be the adjoint of $\theta$.
Let $p:T\times W\lrarr W$ be the projection.
The pull back $p^{\ast}(E,\nabla^0,h)$ 
is denoted by $(E_1,\nabla^1,h_1)$.
We set $\nabla:=\nabla^1+f\,d\zbar-f^{\dagger}dz$.
Then, $(E_1,\nabla,h_1)$ is an instanton on $T\times W$.

Conversely,
let $(E_2,\nabla^2,h_2)$ be a $T$-equivariant
instanton on $T\times W$.
By considering $T$-equivariant sections,
we obtain a vector bundle $E$ on $W$
such that $p^{\ast}E\simeq E_2$.
It is naturally equipped with a connection $\nabla^0$
such that
$p^{\ast}\nabla^0_v
=\nabla^2_v$,
where $v$ denotes the natural horizontal lift
of vector fields on $W$.
By using the $T$-equivariance of $\nabla^2$,
we have the expression
$\nabla^2-p^{\ast}\nabla^0=
 p^{\ast}f\,d\zbar-p^{\ast}f^{\dagger}dz$,
where $f$ is a section of $\End(E)$.
Then, 
$(E,\delbar_E,h,f\,dz)$ is a harmonic bundle.
In summary, we have the following.

\begin{prop}[Hitchin]
Harmonic bundles on $W$ naturally correspond
to $T$-equivariant instantons on $T\times W$.
\hfill\qed
\end{prop}

\subsection{Local estimate}
\label{subsection;12.6.22.30}

Let $U$ be a closed disc
$\bigl\{w\,\big|\,|w-w_0|\leq 1\bigr\}$
of $\cnum$.
Let $(E,\nabla,h)$ be an instanton on $T\times U$.

\begin{assumption}
\label{assumption;13.2.26.2}
We assume that 
$|F(\nabla)|\leq \epsilon$
for a given positive small number $\epsilon$.
We also impose Assumption {\rm\ref{assumption;12.6.24.21}}.
\hfill\qed
\end{assumption}

We use the notation in \S\ref{subsection;12.6.24.10}.
Note that we have
$\bigl|D_{\vecx}^{k_1}D_{w}^{k_2}F\bigr|\leq 
C_{k_1,k_2}\epsilon$,
where $C_{k_1,k_2}$ are constants depending only on 
$(k_1,k_2)$.

\subsubsection{Estimates of the $\bot$-part of the connection form}

Let $\vecu$ be a partially almost holomorphic frame
as in Proposition \ref{prop;12.6.24.3}.
We also assume that $\int_TH(h,\vecu)=I$,
as in \S\ref{subsection;12.6.24.22}.
Let $A$ be the connection form of $\nabla$
with respect to $\vecu$.
Let $\nbigb_{z\zbar}$ represent
$F_{z\zbar}$ with respect to $\vecu$.
We use $\nbigb_{z\wbar}$
and $\nbigb_{w\zbar}$ in similar meanings.

We prepare a notation
in a general situation.
Let $V$ be any vector bundle with a hermitian metric $h_V$
on $U$.
Let $\pi:T\times U\lrarr U$ be the projection.
Let $p\geq 2$.
For any section $f$ of $\pi^{\ast}V$ on $T\times U$,
let $\|f\|_p$ denote the function on $U$
given by 
$\|f\|_p(w)=
 \Bigl(
 \int_{T\times\{w\}}|f|_{h_V}^p
 \Bigr)^{1/p}$.

\begin{lem}
\label{lem;11.12.29.200}
We have
$\|A_{w}^{\bot}\|_p=
 O\bigl(\|F_{w\zbar}^{\bot}\|_p\bigr)$
and
$\bigl\|\del_{\wbar}A^{\bot}_{w}\bigr\|_p
=O\bigl(
 \|\nabla_{\wbar}F^{\bot}_{w\zbar}\|_p
 \bigr)
+O\bigl(
 \epsilon\,\|F_{w\zbar}^{\bot}\|_p
 \bigr)$.
\end{lem}
\pf
Because
$\del_wA_{\zbar}-\del_{\zbar}A_w+[A_{w},A_{\zbar}]
=\nbigb_{w\zbar}$,
we have the following equalities:
\begin{equation}
 \label{eq;11.12.29.210}
 \del_{\zbar}A^{\bot}_w
+\bigl[A_{\zbar},\,A_w^{\bot}\bigr]
=-\nbigb^{\bot}_{w\zbar}
\end{equation}
Then, we obtain the first estimate.
We also obtain the following equation:
\[
 \del_{\zbar}\del_{\wbar}A_w^{\bot}
+\bigl[A_{\zbar},\,\del_{\wbar}A_w^{\bot}\bigr]
=-\del_{\wbar}\nbigb^{\bot}_{w\zbar}
-\bigl[\del_{\wbar}A_{\zbar},\,A_w^{\bot}\bigr]
\]
Because
$\del_{\wbar}A_{\zbar}=O(\epsilon)$,
we obtain the second estimate.
\hfill\qed

\begin{lem}
\label{lem;11.12.30.100}
We have
$\|A_{\wbar}^{\bot}\|_p
=O\Bigl(
 \|F_{w\zbar}^{\bot}\|_p
+\|F_{z\zbar}^{\bot}\|_p
+\|\nabla_{\wbar}F_{z\zbar}^{\bot}\|_p
 \Bigr)$.
We also have
\[
 \|\del_wA_{\wbar}^{\bot}\|_p
=O\Bigl(
 \|\nabla_{\wbar}F_{w\zbar}^{\bot}\|_p
+\|F_{w\zbar}^{\bot}\|_p
+\|F_{z\zbar}^{\bot}\|_p
+\|\nabla_{\wbar}F_{z\zbar}^{\bot}\|_p
 \Bigr).
\]
\end{lem}
\pf
We set $\gtilde:=\overline{H(h,\vecu)}$.
We have 
$A_{\wbar}=
-\gtilde^{-1}\,\overline{(\lefttop{t}{A_w})}\,\gtilde
+\gtilde^{-1}\del_{\wbar}\gtilde$.
Hence, the first claim follows from Lemma \ref{lem;11.12.29.200},
Lemma \ref{lem;11.12.29.2}
and Lemma \ref{lem;11.12.30.10}.
We have
$\del_{w}A_{\wbar}-\del_{\wbar}A_{w}
+[A_{w},A_{\wbar}]=\nbigb_{w\wbar}$.
Hence, we have
\[
 \|\del_wA_{\wbar}^{\bot}\|_p
=O\bigl(
 \|\del_{\wbar}A_w^{\bot}\|_p
 \bigr)
+O\bigl(
 \|A_{\wbar}^{\bot}\|_p+\|A_{w}^{\bot}\|_p
 \bigr)
+\|F_{w\wbar}^{\bot}\|_p.
\]
Then, the second claim follows.
\hfill\qed

\subsubsection{Estimate of the $\bot$-part of the curvature}
\label{subsection;13.1.10.4}

We prepare a notation in a general situation.
Let $V$ be any vector bundle with a hermitian metric $h_V$
on $T\times U$.
Let $\pi:T\times U\lrarr U$ be the projection.
For any section $f$ of $V$ on $T\times U$,
let $\|f\|$ denote the function on $U$
given by $\Bigl( \int_T|f|_{h_V}^2\Bigr)^{1/2}$.
For any sections $f$ and $g$ of $V$,
let $(\!(f,g)\!)$ denote the function on $U$
given by $\int_T h_V(f,g)$.

\begin{prop}
\label{prop;11.12.29.220}
We have the following:
\begin{multline}
 \Delta_w\|F_{z\zbar}^{\bot}\|^2
\leq
-\|\nabla_{\zbar}F^{\bot}_{z\zbar}\|^2
-\|\nabla_zF^{\bot}_{z\zbar}\|^2
-\|\nabla_{\wbar}F^{\bot}_{z\zbar}\|^2
-\|\nabla_{w}F^{\bot}_{z\zbar}\|^2
 \\
+O\Bigl(
 \epsilon\|F_{z\zbar}^{\bot}\|^2
+\epsilon\|F_{w\zbar}^{\bot}\|\,
 \|F_{z\zbar}^{\bot}\|
+\epsilon\|\nabla_{\wbar}F_{w\zbar}^{\bot}\|\,
 \|F_{z\zbar}^{\bot}\|
+\epsilon\|\nabla_{\zbar}F_{z\zbar}^{\bot}\|\,
 \|F_{z\zbar}^{\bot}\|
 \Bigr)\\
+O\Bigl(
 \epsilon\|\nabla_{\wbar}F_{z\zbar}^{\bot}\|\,
 \|F_{z\zbar}^{\bot}\|
+\epsilon\|F_{w\zbar}^{\bot}\|^2
+\epsilon\|F_{w\zbar}^{\bot}\|\,
 \|\nabla_{\wbar}F_{z\zbar}^{\bot}\|
 \Bigr)
\end{multline}
\end{prop}
\pf
We have the following equation:
\[
 \Delta_w|F_{z\zbar}^{\bot}|^2
=-(\nabla_{w}\nabla_{\wbar}F_{z\zbar}^{\bot},F_{z\zbar}^{\bot})
  -(F_{z\zbar}^{\bot},\nabla_{\wbar}\nabla_wF_{z\zbar}^{\bot})
-(\nabla_wF_{z\zbar}^{\bot},\nabla_{w}F_{z\zbar}^{\bot})
-(\nabla_{\wbar}F_{z\zbar}^{\bot},\nabla_{\wbar}F_{z\zbar}^{\bot})
\]
We have
\[
-(\nabla_w\nabla_{\wbar}F_{z\zbar}^{\bot},F_{z\zbar}^{\bot})=
-\bigl(\nabla_w\nabla_{\wbar}F_{z\zbar},\,F_{z\zbar}^{\bot}\bigr)
+\bigl(\nabla_w\nabla_{\wbar}F_{z\zbar}^{\circ},\,
 F_{z\zbar}^{\bot}\bigr)
\]
Let us consider the estimate of
$(\nabla_{w}\nabla_{\wbar}F_{z\zbar}^{\circ},F_{z\zbar}^{\bot})$.
The endomorphism $\nabla_{w}\nabla_{\wbar}F_{z\zbar}^{\circ}$
is represented by the following with respect to $\vecu$:
\[
 \del_{w}\del_{\wbar}\nbigb^{\circ}_{z\zbar}
+[A_{w},\del_{\wbar}\nbigb^{\circ}_{z\zbar}]
+\del_{w}[A_{\wbar},\nbigb_{z\zbar}^{\circ}]
+[A_w,[A_{\wbar},\nbigb^{\circ}_{z\zbar}]]
\]
Recall Lemma \ref{lem;11.12.29.3}.
We have the following estimates:
\begin{equation}
\label{eq;11.12.29.300}
\bigl(\!\bigl(
 \del_{w}\del_{\wbar}\nbigb_{z\zbar}^{\circ},
 \nbigb_{z\zbar}^{\bot}\bigr)\!\bigr)_h
=O\Bigl(
 \bigl\|\del_{w}\del_{\wbar}\nbigb^{\circ}_{z\zbar}\bigr\|_h\,
 \bigl\|\nbigb^{\bot}_{z\zbar}\bigr\|_h\,
 \bigl\|F_{z\zbar}^{\bot}\bigr\|_h
 \Bigr)
\end{equation}
\begin{equation}
\bigl(\!\bigl(
[A_{w},\del_{\wbar}\nbigb_{z\zbar}^{\circ}],
\nbigb_{z\zbar}^{\bot}
\bigr)\!\bigr)_h
=O\Bigl(
 \|A_{w}^{\bot}\|_h\,\|\del_{\wbar}\nbigb^{\circ}_{z\zbar}\|_h\,
 \|\nbigb_{z\zbar}^{\bot}\|_h
 \Bigr)
+O\Bigl(
 \bigl\|[A^{\circ}_{w},\del_{\wbar}\nbigb^{\circ}_{z\zbar}]\bigr\|_h\,
 \bigl\|\nbigb^{\bot}_{z\zbar}\bigr\|_h\,
 \bigl\|F_{z\zbar}^{\bot}\bigr\|_h
 \Bigr)
\end{equation}
\begin{multline}
 \bigl(\!\bigl(
 [\del_{w}A_{\wbar},\nbigb_{z\zbar}^{\circ}],\,
 \nbigb_{z\zbar}^{\bot}
 \bigr)\!\bigr)_h 
=\bigl(\!\bigl(
 [\del_{w}A_{\wbar}^{\bot},\nbigb_{z\zbar}^{\circ}],\,
 \nbigb_{z\zbar}^{\bot}
\bigr)\!\bigr)_h
+\bigl(\!\bigl(
 [\del_{w}A^{\circ}_{\wbar},\nbigb_{z\zbar}^{\circ}],
 \nbigb^{\bot}_{z\zbar}
\bigr)\!\bigr)_h
 \\
=O\bigl(
 \|\nbigb_{z\zbar}^{\circ}\|_h\,
 \|\del_{w}A^{\bot}_{\wbar}\|_h\,
 \|\nbigb_{z\zbar}^{\bot}\|_h
 \bigr)
+O\bigl(
 \bigl\|[\del_{w}A^{\circ}_{\wbar},\nbigb_{z\zbar}^{\circ}]\bigr\|_h\,
 \|\nbigb_{z\zbar}^{\bot}\|_h\,
 \|F_{z\zbar}^{\bot}\|
 \bigr)
\end{multline}
\begin{equation}
 \bigl(\!\bigl(
 [A_{\wbar},\del_{w}\nbigb_{z\zbar}^{\circ}],\,\nbigb^{\bot}_{z\zbar}
 \bigr)\!\bigr)_h
=O\bigl(
 \|\del_{w}\nbigb_{z\zbar}^{\circ}\|_h\,
 \|A_{\wbar}^{\bot}\|_h\, 
 \|\nbigb_{z\zbar}^{\bot}\|_h
 \bigr)
+O\bigl(
 \bigl\|[A_{\wbar}^{\circ},\del_{w}\nbigb^{\circ}_{z\zbar}]\bigr\|_h\,
 \|\nbigb_{z\zbar}^{\bot}\|_h\, 
 \|F_{z\zbar}^{\bot}\|_h
 \bigr) 
\end{equation}
\begin{multline}
\label{eq;11.12.29.301}
 \bigl(\!\bigl(
 \bigl[A_{w},[A_{\wbar},\nbigb_{z\zbar}^{\circ}]\bigr],\,
 \nbigb_{z\zbar}^{\bot}\bigr)\!\bigr)_h
=O\bigl(
 \|A_{w}^{\bot}\|_h\,\|A_{\wbar}^{\bot}\|_h\,
 \|\nbigb_{z\zbar}^{\circ}\|_h\,\|\nbigb_{z\zbar}^{\bot}\|_h
 \bigr)
+O\bigl(
 \|A_{w}^{\bot}\|_h\,
 \|A_{\wbar}^{\circ}\|_h\,\|\nbigb_{z\zbar}^{\circ}\|_h\,
 \|\nbigb_{z\zbar}^{\bot}\|_h
 \bigr)
 \\ 
+O\bigl(
 \|A_{\wbar}^{\bot}\|_h\,
 \|A_{w}^{\circ}\|_h\,
 \|\nbigb_{z\zbar}^{\circ}\|_h\,
 \|\nbigb_{z\zbar}^{\bot}\|_h
 \bigr) 
+O\Bigl(
 \bigl\|A_{w}^{\circ}\bigr\|_h\,
 \bigl\| A_{\wbar}^{\circ}\bigr\|_h\,
\bigl\|
\nbigb_{z\zbar}^{\circ}
\bigr\|_h\,
 \bigl\|\nbigb_{z\zbar}^{\bot}
 \bigr\|_h\,
 \bigl\|F_{z\zbar}^{\bot}\bigr\|
\Bigr)
\end{multline}

We obtain the following estimate for
$(\nabla_w\nabla_{\wbar}F_{z\zbar}^{\circ},F_{z\zbar}^{\bot})$
from (\ref{eq;11.12.29.300})--(\ref{eq;11.12.29.301})
with Lemma \ref{lem;11.12.29.200}:
\begin{equation}
(\!(\nabla_w\nabla_{\wbar}F_{z\zbar}^{\circ},F_{z\zbar}^{\bot})\!)
=O\Bigl(
 \epsilon\,
 \|F_{z\zbar}^{\bot}\|^2
+\epsilon\,
 \|F_{w\zbar}^{\bot}\|\,
 \|F_{z\zbar}^{\bot}\|
+\epsilon\,
 \|\nabla_{\wbar}F_{w\zbar}^{\bot}\|\,
 \|F_{z\zbar}^{\bot}\|
+\epsilon\,
 \|\nabla_{\wbar}F_{z\zbar}^{\bot}\|\,
 \|F_{z\zbar}^{\bot}\|
 \Bigr)
\end{equation}

We have
\begin{equation}
 -(\!(\nabla_{w}\nabla_{\wbar}F_{z\zbar},F_{z\zbar}^{\bot})\!)
=\bigl(\!
 \bigl(
 \nabla_z\nabla_{\zbar}F_{z\zbar},\,F_{z\zbar}^{\bot}
 \bigr)\!\bigr)
+\bigl(\!\bigl([F_{z\wbar},F_{w\zbar}],\,
 F_{z\zbar}^{\bot}\bigr)\!\bigr)
 \\
=-(\!(\nabla_{\zbar}F_{z\zbar},\,\nabla_{\zbar}F_{z\zbar}^{\bot})\!)
+(\!([F_{z\wbar},F_{w\zbar}],F_{z\zbar}^{\bot})\!)
\end{equation}
We have
\begin{multline}
 -(\!(\nabla_{\zbar}F_{z\zbar},\nabla_{\zbar}F_{z\zbar}^{\bot})\!) 
=-(\!(\nabla_{\zbar}F_{z\zbar}^{\bot},\nabla_{\zbar}F_{z\zbar}^{\bot})\!)
-(\!(\nabla_{\zbar}F_{z\zbar}^{\circ},\nabla_{\zbar}F_{z\zbar}^{\bot})\!)
 \\
=-(\!(\nabla_{\zbar}F_{z\zbar}^{\bot},\nabla_{\zbar}F_{z\zbar}^{\bot})\!)
+O\bigl(
 \bigl\|
 \nabla_{\zbar}F_{z\zbar}^{\circ}
\bigr\|\,
 \bigl\|\nabla_{\zbar}F_{z\zbar}^{\bot}\bigr\|\,
 \bigl\|F_{z\zbar}^{\bot}\bigr\|
 \bigr)
\end{multline}
We also have the following:
\begin{multline}
 \bigl(\!\bigl(
 [F_{z\wbar},F_{w\zbar}],\,
 F_{z\zbar}^{\bot} 
 \bigr)\!\bigr)
=O\bigl(
 \bigl\|[F_{z\wbar}^{\circ},\,F_{w\zbar}^{\circ}]\bigr\|\,
 \bigl\|F_{z\zbar}^{\bot}\bigr\|\,
 \bigl\|F_{z\zbar}^{\bot}\bigr\|
 \bigr)
+O\bigl(
 \|F_{z\wbar}^{\circ}\|\,\|F_{w\zbar}^{\bot}\|\,
 \|F_{z\zbar}^{\bot}\|
 \bigr)
 \\
+O\bigl(
 \|F_{w\zbar}^{\circ}\|\,\|F_{z\wbar}^{\bot}\|\,
 \|F_{z\zbar}^{\bot}\|
 \bigr)
+O\bigl(
 \|F_{z\wbar}^{\bot}\|\,\|F_{w\zbar}^{\bot}\|\,
 \|F_{z\zbar}^{\bot}\|
 \bigr)
\end{multline}
We have a similar estimate for the contribution of
$-\bigl(\!\bigl(
 F_{z\zbar}^{\bot},\nabla_{\wbar}\nabla_wF_{z\zbar}^{\bot}
 \bigr)\!\bigr)$.
In all, 
we obtain the claim of Proposition 
\ref{prop;11.12.29.220}.
\hfill\qed

\begin{prop}
\label{prop;12.6.22.20}
We have the following inequality:
\begin{multline}
\label{eq;11.12.30.20}
 \Delta_w\|F_{z\wbar}^{\bot}\|^2
\leq
-\|\nabla_{\zbar}F_{z\wbar}^{\bot}\|^2
-\|\nabla_zF_{z\wbar}^{\bot}\|^2
-\|\nabla_wF_{z\wbar}^{\bot}\|^2
-\|\nabla_{\wbar}F_{z\wbar}^{\bot}\|^2
 \\
+O\Bigl(
 \epsilon\,\|F_{z\wbar}^{\bot}\|\,\|F_{z\zbar}^{\bot}\|
+\epsilon\|\nabla_{\wbar}F_{w\zbar}^{\bot}\|\,
 \|F_{z\wbar}^{\bot}\|
+\epsilon\|F_{w\zbar}^{\bot}\|\,\|F_{z\wbar}^{\bot}\|
+\epsilon
 \|\nabla_{\wbar}F_{z\zbar}^{\bot}\|\,
 \|F_{z\wbar}^{\bot}\|
\Bigr)
 \\
+O\Bigl(
\epsilon
 \|\nabla_{\zbar}F_{z\wbar}^{\bot}\|\,
 \|F_{z\zbar}^{\bot}\|
+\epsilon\|F_{z\wbar}^{\bot}\|\,\|F^{\bot}_{w\zbar}\|
 \Bigr)
\end{multline}
\end{prop}
\pf
We have the following:
\begin{equation}
-\del_w\del_{\wbar}\bigl|F_{z\wbar}^{\bot}\bigr|^2
=-\bigl|\nabla_{\wbar}F_{z\wbar}^{\bot}\bigr|^2
  -\bigl|\nabla_wF_{z\wbar}^{\bot}\bigr|^2
-\bigl(\nabla_w\nabla_{\wbar}F_{z\wbar}^{\bot},\,
 F_{z\wbar}^{\bot}\bigr)
-\bigl(F_{z\wbar}^{\bot},\,
 \nabla_{\wbar}\nabla_wF_{z\wbar}^{\bot}
 \bigr)
\end{equation}

We have
\begin{equation}
-(\nabla_w\nabla_{\wbar}F_{z\wbar}^{\bot},\,F_{z\wbar}^{\bot})
=-\bigl(\nabla_w\nabla_{\wbar}F_{z\wbar},\,F_{z\wbar}^{\bot}\bigr)
+\bigl(\nabla_w\nabla_{\wbar}F_{z\wbar}^{\circ},F_{z\wbar}^{\bot}\bigr)
\end{equation}
Let us look at the contribution of
$\bigl(\nabla_w\nabla_{\wbar}F^{\circ}_{z\wbar},\,
 F_{z\wbar}^{\bot}\bigr)$.
Let $\nbigb_{z\wbar}$ express
$F_{z\wbar}$ with respect to $\vecu$
as in the proof of Proposition \ref{prop;11.12.29.220}.
Then, $\nabla_w\nabla_{\wbar}F_{z\wbar}^{\circ}$
is represented by the following:
\[
 \del_w\del_{\wbar}\nbigb_{z\wbar}^{\circ}
+\bigl[\del_wA_{\wbar},\,\nbigb_{z\wbar}^{\circ}\bigr]
+\bigl[A_{\wbar},\del_{w}\nbigb_{z\wbar}^{\circ}\bigr]
+\bigl[A_{w},\del_{\wbar}\nbigb_{z,\wbar}^{\circ}\bigr]
+\bigl[A_w,\,[A_{\wbar},\nbigb_{z\wbar}^{\circ}]\bigr]
\]
We have the following estimates:
\begin{equation}
 -\bigl(\!\bigl(
 \del_w\del_{\wbar}\nbigb_{z\wbar}^{\circ},\,
 \nbigb_{z\wbar}^{\bot}\bigr)\!\bigr)_h
=O\Bigl(
 \bigl\|\del_w\del_{\wbar}\nbigb_{z\wbar}^{\circ}\bigr\|_h\,
 \bigl\|\nbigb_{z\wbar}^{\bot}\bigr\|_h\,
 \bigl\|F_{z\zbar}^{\bot}\bigr\|
 \Bigr)
\end{equation}
\begin{equation}
 \bigl(\!\bigl(
 [\del_wA_{\wbar},\,\nbigb_{z\wbar}^{\circ}],\,
 \nbigb_{z\wbar}^{\bot}\bigr)\!\bigr)_h
=O\Bigl(
 \bigl\|[\del_wA^{\circ}_{\wbar},\nbigb_{z\wbar}^{\circ}]\bigr\|_h\,
 \bigl\|\nbigb_{z\wbar}^{\bot}\bigr\|_h\,
 \bigl\|F_{z\zbar}^{\bot} \bigr\|
 \Bigr)
+O\Bigl(
 \bigl\|\del_wA_{\wbar}^{\bot}\bigr\|_h\,
 \bigl\|\nbigb_{z\wbar}^{\circ}\bigr\|_h\,
 \bigl\|\nbigb_{z\wbar}^{\bot}\bigr\|_h
 \Bigr)
\end{equation}
\begin{equation}
 \bigl(\!\bigl(
 [A_{\wbar},\del_w\nbigb_{z\wbar}^{\circ}],\,
 \nbigb_{z\wbar}^{\bot}
 \bigr)\!\bigr)_h
=O\Bigl(
 \bigl\|A_{\wbar}^{\circ}\bigr\|_h\,
 \bigl\|\del_w\nbigb_{z\wbar}^{\circ}\bigr\|_h\,
 \bigl\|\nbigb_{z\wbar}^{\bot}\bigr\|_h\,
 \bigl\|F_{z\zbar}^{\bot}\bigr\|
 \Bigr)
+O\Bigl(
 \bigl\|A_{\wbar}^{\bot}\bigr\|_h\,
 \bigl\|\del_w\nbigb_{z\wbar}^{\circ}\bigr\|_h\,
 \bigl\|\nbigb_{z\wbar}^{\bot}\bigr\|_h
 \Bigr)
\end{equation}
\begin{equation}
\bigl(\!\bigl(
 [A_w,\del_{\wbar}\nbigb^{\circ}_{z\wbar}],\,
 \nbigb_{z\wbar}^{\bot}
 \bigr)\!\bigr)_h
=O\Bigl(
\|A_w^{\circ}\|_h\|\del_{\wbar}\nbigb^{\circ}_{z\wbar}\|_h
 \|\nbigb_{z\wbar}^{\bot}\|_h\,\|F_{z\zbar}^{\bot}\|
 \Bigr)
+O\Bigl(
\|A_w^{\bot}\|_h\|\del_{\wbar}\nbigb^{\circ}_{z\wbar}\|_h
 \|\nbigb_{z\wbar}^{\bot}\|_h
 \Bigr)
\end{equation}
\begin{multline}
 \bigl(\!\bigl(
 \bigl[A_w,[A_{\wbar},\nbigb_{z\wbar}^{\circ}]\bigr],\,
 \nbigb_{z\wbar}^{\bot}
 \bigr)\!\bigr)_h
=O\Bigl(
 \bigl\|A_w^{\bot}\bigr\|_h\,
 \bigl\|A_{\wbar}^{\bot}\bigr\|_h\,
 \bigl\|\nbigb_{z\wbar}^{\circ}\bigr\|_h\,
 \bigl\|\nbigb_{z\wbar}^{\bot}\bigr\|_h
 \Bigr)
+O\Bigl(
 \bigl\|A_w^{\circ}\bigr\|_h\,
 \bigl\|A_{\wbar}^{\bot}\bigr\|_h\,
 \bigl\|\nbigb_{z\wbar}^{\circ}\bigr\|_h\,
 \bigl\|\nbigb_{z\wbar}^{\bot}\bigr\|_h
 \Bigr)
 \\
+O\Bigl(
 \bigl\|A_w^{\bot}\bigr\|_h\,
 \bigl\|A_{\wbar}^{\circ}\bigr\|_h\,
 \bigl\|\nbigb_{z\wbar}^{\circ}\bigr\|_h\,
 \bigl\|\nbigb_{z\wbar}^{\bot}\bigr\|_h
 \Bigr)
+O\Bigl(
 \bigl\|A_w^{\circ}\bigr\|_h\,
 \bigl\|A_{\wbar}^{\circ}\bigr\|_h\,
 \bigl\|\nbigb_{z\wbar}^{\circ}\bigr\|_h\,
 \bigl\|\nbigb_{z\wbar}^{\bot}\bigr\|_h\,
 \bigl\|F_{z\zbar}^{\bot}\bigr\|_h
 \Bigr)
\end{multline}
Hence, we obtain the following:
\begin{equation}
\bigl(\!\bigl(
 \nabla_w\nabla_{\wbar}F_{z\wbar}^{\circ},F_{z\wbar}^{\bot}
\bigr)\!\bigr) 
=O\Bigl(
 \epsilon\,\|F_{z\wbar}^{\bot}\|\,\|F_{z\zbar}^{\bot}\|_h
+\epsilon\|\nabla_{\wbar}F_{w\zbar}^{\bot}\|\,\|F_{z\wbar}^{\bot}\|_h
+\epsilon\|F_{w\zbar}^{\bot}|\,|F_{z\wbar}^{\bot}\|_h
+\epsilon
 \|\nabla_{\wbar}F_{z\zbar}^{\bot}\|_h\,
 \|F_{z\wbar}^{\bot}\|_h
 \Bigr)
\end{equation}
We have
\begin{multline}
 -(\!(\nabla_w\nabla_{\wbar}F_{z\wbar},F_{z\wbar}^{\bot})\!)
=(\!(\nabla_z\nabla_{\zbar}F_{z\wbar},F_{z\wbar}^{\bot})\!)
-2\bigl(\!\bigl(
 [F_{w\wbar},F_{z\wbar}],\,F_{z\wbar}^{\bot}
 \bigr)\!\bigr) 
 \\
=-(\!(\nabla_{\zbar}F_{z\wbar},\nabla_{\zbar}F_{z,\wbar}^{\bot})\!)
-2\bigl(\!\bigl(
 [F_{w,\wbar},F_{z\wbar}],\,F_{z\wbar}^{\bot}\bigr)\!\bigr)
 \end{multline}
We have
\begin{multline}
 -(\!(\nabla_{\zbar}F_{z\wbar},\nabla_{\zbar}F_{z\wbar}^{\bot})\!)
=-(\!(\nabla_{\zbar}F_{z\wbar}^{\bot},\nabla_{\zbar}F_{z\wbar}^{\bot})\!)
-(\!(\nabla_{\zbar}F_{z\wbar}^{\circ},\nabla_{\zbar}F_{z\wbar}^{\bot})\!)
 \\
=-(\!(\nabla_{\zbar}F_{z\wbar}^{\bot},\nabla_{\zbar}F_{z\wbar}^{\bot})\!)
+O\Bigl(
 \bigl\|[\nabla_{\zbar}F_{z\wbar}^{\circ}]\bigr\|\,
 \bigl\|\nabla_{\zbar}F_{z\wbar}^{\bot}\bigr\|\,
 \bigl\|F_{z\zbar}^{\bot}\bigr\|
 \Bigr)
\end{multline}
We also have
\begin{multline}
 \bigl(\!\bigl(
 [F_{w\wbar},F_{z\wbar}],\,F^{\bot}_{z\wbar}
 \bigr)\!\bigr) 
=O\Bigl(
 \bigl\|F_{w\wbar}^{\bot}\bigr\|\,
 \bigl\|F_{z\wbar}^{\bot}\bigr\|\,
 \bigl\|F_{z\wbar}^{\bot}\bigr\|
 \Bigr) 
+O\Bigl(
 \bigl\|F_{w\wbar}^{\circ}\bigr\|\,
 \bigl\|F_{z\wbar}^{\bot}\bigr\|\,
 \bigl\|F_{z\wbar}^{\bot}\bigr\|
 \Bigr)
 \\
+O\Bigl(
 \bigl\|F_{w\wbar}^{\bot}\bigr\|\,
 \bigl\|F_{z\wbar}^{\circ}\bigr\|\,
 \bigl\|F_{z\wbar}^{\bot}\bigr\| 
 \Bigr)
+O\Bigl(
 \bigl\|[F_{w\wbar}^{\circ},F_{z\wbar}^{\circ}]\bigr\|\,
 \bigl\|F_{z\wbar}^{\bot}\bigr\|\,
 \bigl\|F_{z\zbar}^{\bot}\bigr\|
 \Bigr)
\end{multline}
We have a similar estimate for the contribution of
$-(F^{\bot}_{z\wbar},\,\nabla_{\wbar}\nabla_wF_{z\wbar}^{\bot})$.
In all,
we obtain the desired estimate
(\ref{eq;11.12.30.20}).
\hfill\qed

\begin{prop}
\label{prop;12.6.22.102}
There exist $C>0$ and $\epsilon_0>0$
such that
the following inequality holds
if $\epsilon<\epsilon_0$:
\begin{multline}
 \Delta_w\bigl(
 \|F_{z\zbar}^{\bot}\|^2
+\|F_{z\wbar}^{\bot}\|^2
 \bigr)
\leq
 -C\Bigl(
 \|F_{z\zbar}^{\bot}\|^2
+ \|F_{z\wbar}^{\bot}\|^2
 \Bigr)
-C\Bigl(
 \|\nabla_zF_{z\zbar}^{\bot}\|^2
+\|\nabla_{\zbar}F_{z\zbar}^{\bot}\|^2
+\|\nabla_wF_{z\zbar}^{\bot}\|^2
+\|\nabla_{\wbar}F_{z\zbar}^{\bot}\|^2
 \Bigr)
 \\
-C\Bigl(
 \|\nabla_zF_{z\wbar}^{\bot}\|^2
+\|\nabla_{\zbar}F_{z\wbar}^{\bot}\|^2
+\|\nabla_wF_{z\wbar}^{\bot}\|^2
+\|\nabla_{\wbar}F_{z\wbar}^{\bot}\|^2
\Bigr)
\end{multline}
\end{prop}
\pf
There exist $C_1>0$
such that
$\|\nabla_z s\|\geq C_1\|s\|$
and
$\|\nabla_{\zbar}s\|\geq C_1\|s\|$
for any section of $\End(E)$ such that 
$s=s^{\bot}$.
Then, the claim follows from 
Proposition \ref{prop;11.12.29.220}
and Proposition \ref{prop;12.6.22.20}.
\hfill\qed

\subsubsection{Higher derivative}

Assume that 
$\|F_{z\zbar}^{\bot}\|^2
+\|F_{z\wbar}^{\bot}\|^2
\leq \delta^2$
for some $\delta<\!<\epsilon$.
For $\rho<1$,
we set $U(\rho)=\bigl\{
 w\,\big|\,|w-w_0|\leq \rho
 \bigr\}\subset U$.

\begin{prop}
\label{prop;12.6.22.200}
For any $k,p$, there exists $C>0$ such that
\[
 \bigl\|
 F_{z\zbar}^{\bot}
 \bigr\|_{L_k^p(T\times U(\rho))}
\leq C\delta,
\quad
\bigl\|
 F_{z\wbar}^{\bot}
 \bigr\|_{L_k^p(T\times U(\rho))}
\leq C\delta.
\]
\end{prop}
\pf
It can be shown by a standard bootstrapping argument.
We give only an indication.
We take $\rho<\rho'<1$.
In the following, we shall replace $\rho'$ with smaller one.
Let $\kappa$ denote $z$, $\zbar$, $w$ and $\wbar$.
By Proposition \ref{prop;12.6.22.102},
we obtain 
$\|\nabla_{\kappa}F_{z\zbar}^{\bot}
 \|_{L^2(T\times U(\rho'))}
=O(\delta)$
and
$\|\nabla_{\kappa}F_{z\wbar}^{\bot}
  \|_{L^2(T\times U(\rho'))}
=O(\delta)$.

With respect to the frame $\vecu$,
the endomorphism
$-\nabla_w\nabla_{\wbar}F_{z\zbar}$
is represented by
\begin{equation}
\label{eq;11.12.27.1}
 -\del_{w}\del_{\wbar}\nbigb_{z\zbar}
-\bigl[\del_{w}A_{\wbar},\,\nbigb_{z\zbar}\bigr]
-\bigl[A_{\wbar},\del_{w}\nbigb_{z\zbar}\bigr]
+\bigl[A_{w},\del_{\wbar}\nbigb_{z\zbar}\bigr]
+\bigl[A_{w},\,[A_{\wbar},\nbigb_{z\zbar}]\bigr],
\end{equation}
and the endomorphism $-\nabla_z\nabla_{\zbar}F_{z\zbar}$
is represented by
\begin{equation}
\label{eq;11.12.27.2}
 -\del_{z}\del_{\zbar}\nbigb_{z\zbar}
-\bigl[\del_{z}A_{\zbar},\nbigb_{z\zbar}\bigr]
-\bigl[A_{\zbar},\,\del_{z}\nbigb_{z\zbar}\bigr]
+\bigl[A_{z},\,\del_{\zbar}\nbigb_{z\zbar}\bigr]
+\bigl[A_{z},\,[A_{\zbar},\nbigb_{z\zbar}]\bigr].
\end{equation}
The sum of (\ref{eq;11.12.27.1}) and (\ref{eq;11.12.27.2})
is equal to $[\nbigb_{z\wbar},\,\nbigb_{w\zbar}]$.
By looking at the $\bot$-part of the equation,
we obtain the following equation:
\begin{equation}
\label{eq;12.6.22.110}
 \mbox{The $\bot$-part of (\ref{eq;11.12.27.1})}
+\mbox{The $\bot$-part of (\ref{eq;11.12.27.2})}
=[\nbigb_{z\wbar},\,\nbigb_{w\zbar}]^{\bot}
\end{equation}
By using Lemma \ref{lem;11.12.29.200}
and Lemma \ref{lem;11.12.30.100},
we obtain 
$\bigl\|
 (\del_w\del_{\wbar}+\del_z\del_{\zbar})\nbigb^{\bot}_{z\zbar}
\bigr\|_{L^2(T\times U(\rho'))}
=O(\delta)$.
Similarly,
we obtain
$\bigl\|
 (\del_w\del_{\wbar}+\del_z\del_{\zbar})
 \nbigb^{\bot}_{z\wbar}
\bigr\|_{L^2(T\times U(\rho'))}
=O(\delta)$.
It follows that
\[
 \|F^{\bot}_{z\zbar}\|_{L^4(T\times U(\rho'))}
+\|F^{\bot}_{z\wbar}\|_{L^4(T\times U(\rho'))}
=O(\delta)
\]
\[
  \|\nabla_{\kappa}F^{\bot}_{z\zbar}\|_{L^4(T\times U(\rho'))}
+\|\nabla_{\kappa}F^{\bot}_{z\wbar}\|_{L^4(T\times U(\rho'))}
=O(\delta)
\]
By using Lemma \ref{lem;11.12.29.200},
Lemma \ref{lem;11.12.30.100}
and (\ref{eq;12.6.22.110}),
we obtain 
$\bigl\|
 (\del_w\del_{\wbar}+\del_z\del_{\zbar})
 \nbigb^{\bot}_{z\zbar}
\bigr\|_{L^4(T\times U(\rho'))}
=O(\delta)$.
Similarly,
we obtain
$\bigl\|
 (\del_w\del_{\wbar}+\del_z\del_{\zbar})
 \nbigb^{\bot}_{z\wbar}
\bigr\|_{L^4(T\times U(\rho'))}
=O(\delta)$.
By the same argument,
we obtain the following for any $p$:
\[
 \|F^{\bot}_{z\zbar}\|_{L^p(T\times U(\rho'))}
+\|F^{\bot}_{z\wbar}\|_{L^p(T\times U(\rho'))}
+\|\nabla_{\kappa}F^{\bot}_{z\zbar} \|_{L^p(T\times U(\rho'))}
+\|\nabla_{\kappa}F^{\bot}_{z\wbar}\|_{L^p(T\times U(\rho'))}
=O(\delta)
\]
\[
  \bigl\|
  (\del_w\del_{\wbar}+\del_z\del_{\zbar})\nbigb^{\bot}_{z\zbar}
\bigr\|_{L^p(T\times U(\rho'))}
+\bigl\|
 (\del_w\del_{\wbar}+\del_z\del_{\zbar})\nbigb^{\bot}_{z\wbar}
\bigr\|_{L^p(T\times U(\rho'))}
=O(\delta)
\]
Namely, 
we obtain 
$\|F^{\bot}_{z\zbar}\|_{L^p_2(T\times U(\rho'))}
+\|F^{\bot}_{z\wbar}\|_{L^p_2(T\times U(\rho'))}
=O(\delta)$.

By the argument in Lemma \ref{lem;11.12.29.200},
we obtain 
$\|A_w^{\bot}\|_{L^p_2}=O(\delta)$.
By the argument in Lemma \ref{lem;11.12.30.100},
we obtain 
$\|A_{\wbar}^{\bot}\|_{L^p_1}=O(\delta)$.
By the relation
$\del_wA_{\wbar}-\del_{\wbar}A_w+[A_{w},A_{\wbar}]
=\nbigb_{w\wbar}$,
we obtain
$\|\del_wA^{\bot}_{\wbar}\|_{L_1^p}=O(\delta)$.
We also have
$\|A^{\bot}_z\|_{L_2^p}=O(\delta)$,
which follows from Lemma \ref{lem;11.12.30.10}.
Then, we obtain
\[
  \bigl\|
 (\del_w\del_{\wbar}+\del_z\del_{\zbar})\nbigb^{\bot}_{z\zbar}
\bigr\|_{L^p_1(T\times U(\rho'))}
=O(\delta)
\]
Hence, we obtain
$\bigl\|\nbigb^{\bot}_{z\zbar}
 \bigr\|_{L^p_3(T\times U(\rho'))}=O(\delta)$.
Similarly, we obtain
$\bigl\|
 \nbigb^{\bot}_{z\wbar} 
 \bigr\|_{L^p_3(T\times U(\rho'))}=O(\delta)$.
By the inductive argument,
we obtain
$\bigl\|
 \nbigb^{\bot}_{z\zbar}
 \bigr\|_{L^p_k(T\times U(\rho'))}
+\bigl\|\nbigb^{\bot}_{z\wbar}
 \bigr\|_{L^p_k(T\times U(\rho'))}=O(\delta)$
for any $k$.
\hfill\qed

\begin{cor}
\label{cor;12.6.22.201}
For any $k$ and $p$,
there exists $C>0$ such that
$\|H(h,\vecu)^{\bot}\|_{L_k^p(T\times U(\rho))}
\leq C\delta$.
\end{cor}
\pf
It follows from Proposition \ref{prop;12.6.22.200}
and Lemma \ref{lem;11.12.30.10}.
\hfill\qed

\subsection{Global estimate}
\label{subsection;13.2.26.1}

\subsubsection{Preliminary}

For $R>0$,
we set
$Y_R:=\bigl\{w\in\cnum\,\big|\,|w|\geq R\bigr\}$
and $X_R:=T\times Y_R$.
An instanton $(E,\nabla,h)$ is called $L^2$,
if the curvature $F:=F(\nabla)$ is $L^2$.
We study the behaviour of $L^2$-instantons around infinity.
We suppose that $(E,\nabla,h)$ is an $L^2$-instanton
in this subsection.
For $w_0\in Y_{R}$ and $a>0$,
let $B_{w_0}(a):=\bigl\{
 w\in\cnum\,\big|\,|w-w_0|\leq a
 \bigr\}$.

Let $\epsilon>0$ be sufficiently small.
There exists $R_1$
such that
$\bigl\|F_{|X_{R_1}}\bigr\|_{L^2}<\epsilon$.
Let $w_0\in Y_{2R_1}$.
By the theorem of Uhlenbeck \cite{Uhlenbeck1},
for any $(x,w)\in T\times B_{w_0}(1)$,
we have 
$\bigl|F(x,w)\bigr|
 =O\bigl(
 \bigl\|F_{|T\times B_{w_0}(2)}\bigr\|_{L^2}
 \bigr)
=O(\epsilon)$.
In particular, we may assume that
$E_{w}$ are semistable if $w\in Y_{2R_1}$.
Because we are interested in
the behaviour around infinity,
we may assume that
$(E_w,\delbar_{E_w})$ are semistable of degree $0$
for any $w\in Y_R$,
from the beginning.

\subsubsection{Prolongation of the spectral curve}

We consider the relative Fourier-Mukai transform
$\RFM_-(E,\delbar_E)$,
which is a coherent sheaf on $T^{\lor}\times Y_R$.
The support is relatively $0$-dimensional over $Y_R$,
denoted by $\Sp(E)$.
It is called the spectral curve of $(E,\delbar_E)$.
Let $\Ybar_R$ be the closure of $Y_R$ in $\proj^1$,
i.e.,
$\Ybar_R=Y_R\cup\{\infty\}$.

\begin{thm}
\label{thm;12.6.22.101}
$\Sp(E)$ is extended to a closed subvariety
$\Spbar(E)$ in $T^{\lor}\times \Ybar_R$.
\end{thm}
\pf
Let $\rho$ denote the rank of $E$.
We have the holomorphic map
$\varphi:Y_R\lrarr \Sym^{\rho} T^{\lor}$
induced by $\Sp(E)$.
We have only to show that 
it is extended to a holomorphic map
$\Ybar_R\lrarr\Sym^{\rho}T^{\lor}$.
We fix a closed immersion
$\Sym^{\rho}T^{\lor}\subset\proj^N$
for a sufficiently large $N$,
and we regard $\varphi$ as a  holomorphic map
$Y_R\lrarr \proj^N$.
Let $d_{\proj^N}$ denote the distance of $\proj^N$,
induced by the Fubini-Study metric.

Take any $w_0\in Y_{2R}$.
By Corollary \ref{cor;12.6.24.2},
for any $w_1,w_2\in B_{w_0}(1/2)$,
we have
\begin{equation}
\label{eq;13.5.28.10}
 d_{\proj^N}(\varphi(w_1),\varphi(w_2))
=O\bigl(\|F_{|T\times B_{w_0}(2)}\|_{L^2}\bigr).
\end{equation}
Note that $\varphi$ is holomorphic.
We can also regard it as a harmonic map
between K\"ahler manifolds.
Let $T_w\varphi$ be the derivative of $\varphi$,
and 
$\bigl|
 T_w\varphi
 \bigr|$
denote the norm of $T_w\varphi$
with respect to the Euclidean metric $dw\,d\wbar$
and the Fubini-Study metric of $\proj^N$.
For any $w\in B_{w_0}(1/4)$,
we obtain the following estimate
from (\ref{eq;13.5.28.10})
by using Cauchy's formula for differentiation
in complex analysis:
\begin{equation}
 \bigl|
 T_w\varphi
 \bigr|
=O\bigl(
 \|F_{|T\times B_{w_0}(2)}\|
 \bigr)
\end{equation}
Hence, we obtain the following finiteness
of the energy of the harmonic map $\varphi$:
\[
 \int_{Y_{2R_1}}
 \bigl|T_w\varphi\bigr|^2\,|dw\,d\wbar|
<
 C\bigl\|
 F_{|Y_R}
 \bigr\|_{L^2}^2
<\infty
\]
Then, $\varphi$ is extended on $\Ybar_R$,
according to Theorem 3.6 in \cite{Sacks-Uhlenbeck}.
\hfill\qed

\vspace{.1in}
The intersection
$\Spbar(E)\cap (T^{\lor}\times\{\infty\})$
is denoted by
$\Sp_{\infty}(E)$.

\subsubsection{Asymptotic decay}

We fix a lift of $\Spbar(E)$ to 
a closed subvariety 
$\Spbar(E)_1\subset \Ybar_R\times \cnum_{\zeta}$,
which induces an action of $\zeta$
on $\RFM_-(E,\delbar_E)$.
(See \S\ref{subsection;13.10.15.30}.)
Let $f_{\zeta}$ be the corresponding holomorphic endomorphism of $E$.
We set $\delbar_0:=\delbar_E-f_{\zeta}d\zbar$,
which gives a holomorphic structure of $E$.
For each $w$,
the restriction of 
$\nbige'=(E,\delbar_0)$ 
to $T_z\times\{w\}$
is holomorphically trivial.
It is naturally isomorphic to
$p^{\ast}p_{\ast}(\nbige')$,
where $p:X_R\lrarr Y_R$ denotes the natural projection.
We obtain the decomposition $h=h^{\circ}+h^{\bot}$
as in \S\ref{subsection;12.11.5.10}.

\begin{thm}
\label{thm;12.8.7.1}
For any polynomial $P(t_1,t_2,t_3,t_4)$
of non-commutative variables,
there exists $C>0$
such that 
\[
 P(\nabla_z,\nabla_{\zbar},\nabla_w,\nabla_{\wbar})
 h^{\bot}=O\bigl(\exp(-C|w|)\bigr). 
\]
\end{thm}
\pf
Let $\epsilon>0$ be any sufficiently small number.
We may assume that
$\|F_{|X_{R_1}}\|<\epsilon$
for some $R_1>0$.
By Theorem \ref{thm;12.6.22.101},
we may assume that Assumption \ref{assumption;13.2.26.2}
is satisfied for the restriction of $(E,\nabla,h)$
to any disc contained in $X_{R_1}$.
In particular, we can apply 
Proposition \ref{prop;12.6.22.102}
to $(E,\nabla,h)_{|X_{R_1}}$.
We obtain the following inequality
for some $C_1>0$:
\[
 \Delta_w\bigl(
 \|F^{\bot}_{z\zbar}\|^2
+\|F^{\bot}_{z\wbar}\|^2
\bigr)
\leq
-C_1\bigl(
 \|F^{\bot}_{z\zbar}\|^2
+\|F^{\bot}_{z\wbar}\|^2
 \bigr)
\]
We obtain the following lemma
by a standard argument.
\begin{lem}
\label{lem;13.1.10.3}
We have
$\|F^{\bot}_{z\zbar}\|^2
+\|F^{\bot}_{z\wbar}\|^2
=O\Bigl(
 \exp\bigl(-C_2|w|\bigr)
 \Bigr)$
for some $C_2>0$.
\end{lem}
\pf
This is a variant of a lemma of Ahlfors
(\cite{a}, \cite{Simpson90}).
We give only an indication.
We put 
$G:=\|F^{\bot}_{z\zbar}\|^2+\|F^{\bot}_{z\wbar}\|^2$.
We put
$f_{\epsilon}:=C_3\exp(-2C_1^{1/2}|w|)+\epsilon$,
where $\epsilon>0$ and $C_3>0$.
We have the inequality
$\Delta_wf_{\epsilon}\geq -C_1f_{\epsilon}$.
If $C_3$ is sufficiently large,
we have $f_{\epsilon}>G$
on $\{|w|=R_1\}$.
For each $\epsilon>0$,
we have $f_{\epsilon}>G$
outside a compact subset.
We put
$U:=\bigl\{
 w\,\big|\,
 f_{\epsilon}(w)<G(w)
 \bigr\}$.
Then, $U$ is relatively compact,
and we have $f_{\epsilon}=G$
on the boundary of $U$.
On $U$,
we have
$\Delta_w(G-f_{\epsilon})
\leq
 -C(G-f_{\epsilon})
\leq 0$.
By the maximum principle,
we have
$\sup_U(G-f_{\epsilon})
\leq
 \max_{\del U}(G-f_{\epsilon})=0$.
Hence, we obtain that $U$ is empty.
It means $G\leq f_{\epsilon}$ on $Y_R$
for any $\epsilon$.
We obtain the desired inequality
by taking the limit $\epsilon\to 0$.
\hfill\qed

\vspace{.1in}
Then, the claim of Theorem \ref{thm;12.8.7.1}
follows from Corollary \ref{cor;12.6.22.201}.
\hfill\qed

\subsubsection{Reduction to asymptotic harmonic bundles}

Let $p:X_R\lrarr Y_R$ denote the projection.
By using the push-forward of $\nbigo$-modules,
we obtain a holomorphic vector bundle
$V:=p_{\ast}\nbige'$ on $Y_R$.
It is equipped with a Higgs field 
$\theta_V:=f_{\zeta}dw$.
For any $s_i\in V_{|w}$ $(i=1,2)$,
the corresponding holomorphic  section
of $\nbige'_{|T_w}$
is denoted by $\stilde_i$.
We set
$h_V(s_1,s_2):=\int_{T}h(\stilde_1,\stilde_2)$.
We have the Chern connection $\delbar_V+\del_V$
with respect to $h_V$.
Let $\theta^{\dagger}_V$ denote the adjoint of
$\theta_V$.

\begin{prop}
\label{prop;13.1.11.21}
There exists $C>0$
such that the following holds:
\begin{equation}
 \label{eq;12.6.22.202}
 F(h_V)+\bigl[\theta_V,\theta_V^{\dagger}\bigr]=
O\bigl(\exp(-C|w|)\bigr). 
\end{equation}
\end{prop}
\pf
We identify $p^{\ast}V=\nbige'$.
According to Theorem \ref{thm;12.8.7.1},
the difference $h-p^{\ast}h_V$ and its derivatives are
$O\bigl(\exp(-C_1|w|)\bigr)$.
(The constant $C_1$ may depend on the order 
of derivatives.)
We also have
$\delbar_E=p^{\ast}\delbar_V+f_{\zeta}d\zbar$.
Hence,
$(p^{\ast}V,p^{\ast}\delbar_V+f_{\zeta}d\zbar,p^{\ast}h_V)$
satisfies 
$\Lambda_{\omega}F(p^{\ast}h_V)=O\bigl(\exp(-C_2|w|)\bigr)$,
which is equivalent to (\ref{eq;12.6.22.202}).
\hfill\qed

\subsubsection{Estimate of the curvature}

\begin{thm}
\label{thm;13.1.30.2}
There exists $\rho>0$
such that the following holds:
\[
 F(h)=O\Bigl(
 \frac{dz\,d\zbar}{|w|^{2}(-\log|w|)^2}
 \Bigr)
+O\Bigl(
 \frac{dw\,d\wbar}{|w|^{2}(-\log|w|)^2}
 \Bigr)
+O\Bigl(
 \frac{dw\,d\zbar}{|w|^{1+\rho}}
 \Bigr)
+O\Bigl(
 \frac{dz\,d\wbar}{|w|^{1+\rho}}
 \Bigr)
\]
\end{thm}
\pf
We shall use an estimate for asymptotic harmonic bundles
explained in \S\ref{subsection;13.1.11.10}.
Let $\varphi:\Delta_{u}=\{|u|<R^{-1/e}\}\lrarr \Ybar_R$
be given by
$\varphi(u)=u^e$.
For the expression 
$\theta=f_{\zeta}dw
=f_{\zeta}(-eu^{-e-1}du)$,
according to Theorem \ref{thm;12.6.22.101},
the spectral curve $\Sp(f_{\zeta})\subset \cnum\times Y_R$
is contained in 
$\{|\zeta|\leq R'\}\times Y_R$
for some $R'$,
and the closure in $\cnum\times\Ybar_R$
is a complex variety.
Hence, we may assume that
$\varphi^{\ast}(E,\delbar_E,\theta)$
has the following decomposition as in (\ref{eq;13.1.11.3}):
\[
 \varphi^{\ast}E
=\bigoplus_{\gminia\in\Irr(\varphi^{\ast}\theta_V)}
 E_{\gminia}
\]
Moreover, we have $\deg_{u^{-1}}\gminia\leq e$
for any $\gminia\in\Irr(\varphi^{\ast}\theta_V)$.

We set
$(V',\delbar_{V'},\theta_{V'},h'):=
 \varphi^{-1}
 (V,\delbar_V,\theta_V,h_V)$.
According to Proposition \ref{prop;13.1.11.21},
it satisfies (\ref{eq;13.1.11.1}).
By Corollary \ref{cor;13.1.11.22},
we have
\[
 |F(h_V)|_h=O\bigl(|u|^{-2}(\log|u|)^{-2}du d\ubar\bigr)
\]
Hence, 
we have
$|F(h)_{w\wbar}|_h=|F(h)_{z\zbar}|_h
=O\bigl(
 |w|^{-2}(\log|w|)^{-2}\bigr)$.

We take a frame $\vecv$ of $\nbigp_aV'$
as in \S\ref{subsection;13.1.11.22}.
Let $\Theta$ be determined by
$\varphi^{\ast}f_{\zeta}\vecv=
 \vecv\Theta$.
Let $C_w$ be determined by
$\varphi^{\ast}(\del_w)\vecv=\vecv C_w$.
We have
$\varphi^{\ast}(\del_wf_{\zeta})\vecv
=\vecv\,\bigl(\varphi^{\ast}(\del_w)\Theta+[C_w,\Theta]\bigr)$.
We have the expression
\[
 \Theta=\bigoplus\bigl(
 (\varphi^{\ast}\del_w\gminia
-e^{-1}\alpha u^e )\,I_{\gminia,\alpha}
-e^{-1}u^{e}\Theta_{\gminia,\alpha}
 \bigr),
\]
where the entries of $\Theta_{\gminia,\alpha}$
are holomorphic at $u=0$.
The norm of the endomorphism determined by $\vecv$
and $\Theta_{\gminia,\alpha}$ is 
$O\bigl((\log|w|)^{-1}\bigr)$
by Proposition \ref{prop;12.8.7.3}.
Note that
$\varphi^{\ast}(\del_w)=-e^{-1}u^{e+1}\del_u$
and $\varphi^{\ast}(\del^2_w)\gminia
=O(|\varphi^{\ast}(w)|^{-1-\rho})$ for some $\rho>0$.
Hence, the contribution of 
$\varphi^{\ast}(\del_w)\Theta$
to $\varphi^{\ast}(\del_wf)$ is dominated as 
$O(\varphi^{\ast}|w|^{-1-\rho})$ for some $\rho>0$.
Let $G_w$ be the endomorphism
determined by $\vecv$ and $C_w$.
By using Lemma \ref{lem;13.1.11.4},
we obtain
$[G_w,\varphi^{\ast}f_{\zeta}]=
 O\bigl(\varphi^{\ast}|w|^{-2}\bigr)$.
Hence, we obtain
$|\del_wf_{\zeta}|_{h_V}=O(|w|^{-1-\rho})$
for some $\rho>0$.
Then, we obtain
$|F_{z\wbar}|_h=|F_{w\zbar}|_h
=O(|w|^{-1-\rho})$ for some $\rho>0$.
\hfill\qed

\begin{cor}
$(E,\delbar_E,h)$ is acceptable,
i.e.,
the curvature $F(h)$ is bounded
with respect to
$h$ and 
the Poincar\'e metric
$|w|^{-2}(\log|w|)^{-2}dw\,d\wbar
+dz\,d\zbar$
on $X_R$
around $T\times\{\infty\}$.
\hfill\qed
\end{cor}

\subsubsection{Prolongation to a filtered bundle}

We set $\Xbar_R:=T\times\Ybar_R$.

\begin{cor}
\label{cor;13.10.15.40}
The holomorphic vector bundle $(E,\delbar_E)$
is naturally extended to
a filtered bundle $\nbigp_{\ast}E$
on $(\Xbar_R,T\times\{\infty\})$.
(See {\rm\S\ref{subsection;13.2.26.3}}
for filtered bundles.)
Moreover, the filtered bundle is good
in the sense of {\rm\S\ref{subsection;13.2.14.11}}.
\end{cor}
\pf
Because $(E,\delbar_E,h)$ is acceptable,
we obtain the first claim from
Theorem 21.31 of \cite{Mochizuki-wild}.
As explained in \ref{subsection;13.1.11.22},
we obtain a filtered bundle
$\nbigp_{\ast}V$ on $(\Ybar_R,\infty)$
from the Higgs bundle with hermitian metric
$(V,\delbar_V,h_V,\theta_V)$.
By Proposition \ref{prop;12.8.7.3},
we obtain that
the filtered Higgs bundle $(\nbigp_{\ast}V,\theta_V)$
is good.
It implies the claim of the corollary.
\hfill\qed

\vspace{.1in}

We obtain the spectral curve
$\Sp(\nbigp_a E)\subset 
 T^{\lor}\times \Ybar_R$ of $\nbigp_aE$.
It is equal to $\Spbar(E)$ in 
Theorem \ref{thm;12.6.22.101},
and independent of the choice of $a\in\real$.

\subsection{An estimate in a variant case}
\label{subsection;13.2.26.10}

We continue to use the notation in
\S\ref{subsection;13.2.26.1}.
Let $(E,\nabla,h)$ be an instanton on $X_R$.
Let $F=F(\nabla)$ be its curvature.
We suppose the following:
\begin{itemize}
\item
 $|F(z,w)|\to 0$ when $|w|\to\infty$,
 i.e.,
 for any $\delta>0$,
 there exists $R_{\delta}>0$
 such that 
 $|F(z,w)|_h\leq \delta$
 for any $|w|\geq R_{\delta}$.
 In particular, we obtain 
 $\Sp(E,\delbar_E)
 \subset T^{\lor}\times Y_{R_{\delta}}$,
 if $\delta$ is sufficiently small.
\item
 The closure of $\Sp(E)$ in $T^{\lor}\times \Ybar_{R_{\delta}}$
 is a complex subvariety.
\end{itemize}
We denote the closure by $\Spbar(E)$,
and we set $\Sp_{\infty}(E):=\Spbar(E)\cap(T^{\lor}\times\{\infty\})$.
We obtain the following theorem.

\begin{thm}
\label{thm;13.2.26.3}
Under the assumption,
$(E,\nabla,h)$ is an $L^2$-instanton.
\end{thm}
\pf
By the assumption,
there exists $R_1>0$,
such that Assumption \ref{assumption;13.2.26.2}
is satisfied for $(E,\nabla,h)_{|X_{R_1}}$.
In particular, we can apply 
Proposition \ref{prop;12.6.22.102}
to $(E,\nabla,h)_{|X_{R_1}}$.
We obtain the estimate as in Theorem \ref{thm;12.8.7.1}
by the same argument.
Then, we obtain estimates
as in Proposition \ref{prop;13.1.11.21}
and Theorem \ref{thm;13.1.30.2}
by the same arguments.
In particular,
$(E,\nabla,h)$ is an $L^2$-instanton.
\hfill\qed

\vspace{.1in}
Theorem \ref{thm;13.2.26.3} implies that 
we can replace the $L^2$-condition with a weaker one,
under the assumption that the spectral curve
is extended in a complex analytic way.

\subsection{Asymptotic harmonic bundles}
\label{subsection;13.1.11.10}
In this subsection,
we explain that some of the results for the asymptotic behaviour of
wild harmonic bundles 
are naturally extended for
Higgs bundles with a hermitian metric
satisfying the Hitchin equation
up to an exponentially small term.
It is used in the proof of Theorem \ref{thm;13.1.30.2}.

We put 
$X:=\Delta_z=\bigl\{z\in\cnum\,\big|\,|z|<1\bigr\}$,
$\Xbar:=\{|z|\leq 1\}$,
and $D:=\{0\}$.
Let $g_{\Poin}$ be the Poincar\'e metric of $X\setminus D$.
Let $(E,\delbar_E,\theta)$ be a Higgs bundle
on $\Xbar\setminus D$.
We suppose that there exists a decomposition
\begin{equation}
 \label{eq;13.1.11.3}
 (E,\theta)=\bigoplus_{
 \substack{\gminia\in z^{-1}\cnum[z^{-1}]\\
 \alpha\in\cnum}}
 (E_{\gminia,\alpha},\theta_{\gminia,\alpha})
\end{equation}
such that,
for the expression
$\theta_{\gminia,\alpha}=
 d\gminia+\alpha\,dz/z+f_{\gminia,\alpha}dz/z$,
the eigenvalues of $f_{\gminia,\alpha}(z)$
goes to $0$ when $z\to 0$.
We put
$\Irr(\theta):=\bigl\{\gminia\,|\,\exists\alpha\,\mbox{\rm such that }
 E_{\gminia,\alpha}\neq 0\bigr\}$.

For any $\gminia(z)=\sum_{j\geq -N}\gminia_jz^j$
with $\gminia_{-N}\neq 0$,
we set $\ord(\gminia):=-N$.
We also set $\ord(0):=0$.
We take a negative number $p$ satisfying
$p<\min\bigl\{\ord(\gminia-\gminib)\,\big|\,
 \gminia,\gminib\in\Irr(\theta),\,
 \gminia\neq\gminib\bigr\}$.

Let $h$ be a hermitian metric of $E$.
Let $\theta^{\dagger}$ denote the adjoint of $\theta$
with respect to $h$.
Let $F(h)$ denote the curvature of $(E,\delbar_E,h)$.
We impose the following condition
for some $C_0>0$ and $\epsilon_0>0$:
\begin{equation}
 \label{eq;13.1.11.1}
 \Bigl|
 F(h)+[\theta,\theta^{\dagger}]
 \Bigr|_{h,g_{\Poin}}
\leq
C_0\exp\bigl(-\epsilon_0|z|^{p}\bigr)
\end{equation}

\subsubsection{Asymptotic orthogonality and acceptability}

We have the following version of Simpson's main estimate.

\begin{prop}
\label{prop;12.8.7.3}
Suppose that $(E,\delbar_E,\theta,h)$
satisfies {\rm(\ref{eq;13.1.11.1})}.
\begin{itemize}
\item
If $\gminia\neq\gminib$,
there exists $\epsilon>0$ such that
$E_{\gminia,\alpha}$ and $E_{\gminib,\beta}$
are $O\bigl(
 \exp(-\epsilon|z|^{\ord(\gminia-\gminib)})\bigr)$-asymptotically
orthogonal,
i.e., there exists $C>0$ such that,
for any $u,v\in E_{|Q}$,
we have
$\bigl|h(u,v)\bigr|
\leq C_1\exp\bigl(-\epsilon|z(Q)|^{\ord(\gminia-\gminib)}\bigr)$.
\item
If $\alpha\neq\beta$,
there exists $\epsilon>0$ such that 
$E_{\gminia,\alpha}$
and $E_{\gminia,\beta}$ 
are $O(|z|^{\epsilon})$-asymptotically orthogonal. 
\item
 $\theta_{\gminia,\alpha}
 -(d\gminia+\alpha\,dz/z)\,\id_{E_{\gminia,\alpha}}$
 is bounded with respect to
 $h$ and the Poincar\'e metric $g_{\Poin}$.
\end{itemize}
\end{prop}
\pf
By considering the tensor product
with a harmonic bundle of a rank one,
we may assume 
$p<\min\bigl\{
 \ord(\gminia)\,\big|\,
 \gminia\in\Irr(\theta)
 \bigr\}$.
We have a map
$\eta_{\ell}:
 z^{-1}\cnum[z^{-1}]\lrarr 
 \nbigi_{\ell}:=
 z^{-\ell}\cnum[z^{-1}]$
by forgetting the terms
$\sum_{j\geq -\ell+1}\gminia_jz^{j}$.
For each $\gminib\in \nbigi_{\ell}$,
we set
$E^{(\ell)}_{\gminib}:=
 \bigoplus_{\eta_{\ell}(\gminia)=\gminib}
 \bigoplus_{\alpha\in\cnum}E_{\gminia,\alpha}$.
Let $\pi^{(\ell)}_{\gminia}$ denote the projection of $E$
onto $E^{(\ell)}_{\gminib}$
with respect to the decomposition
$E=\bigoplus E^{(\ell)}_{\gminib}$.
In the case $\ell=1$,
we omit the superscript $(1)$

Let $\Irr(\theta,\ell)$ be the image of
$\Irr(\theta)$ by $\eta_{\ell}$.
We take a total order $\leq'$ on  $\Irr(\theta,\ell)$
for each $\ell$
such that the induced map
$\Irr(\theta,1)\lrarr\Irr(\theta,\ell)$
is order-preserving.
Let $E^{\prime(\ell)}_{\gminib}$
be the orthogonal complement of
$\bigoplus_{\gminic<'\gminib}E_{\gminic}$
in 
$\bigoplus_{\gminic\leq'\gminib}E_{\gminic}$.
Let $\pi^{\prime(\ell)}_{\gminib}$ be the orthogonal projection
onto $E^{\prime(\ell)}_{\gminib}$.
In the case $\ell=1$,
we omit the superscript $(1)$.
We have 
$\pi^{\prime(\ell)}_{\gminib}
=\sum_{\eta_{\ell}(\gminia)=\gminib}
 \pi^{\prime}_{\gminia}$.

We put $\zeta_{\ell}:=\eta_{\ell}-\eta_{\ell+1}$.
We have the expression $\theta=f\,dz$.
We put
$f^{(\ell)}:=
 f-\sum_{\gminia}\del_z\eta_{\ell+1}(\gminia)\pi_{\gminia}$,
$\mu^{(\ell)}:=
 f^{(\ell)}-\sum_{\gminia}
 \del_z\zeta_{\ell}(\gminia)\pi^{\prime}_{\gminia}$
and 
$\nbigr_{\gminib}^{(\ell)}:=
 \pi_{\gminib}^{(\ell)}-\pi_{\gminib}^{\prime(\ell)}$.
We consider the following claims.
\begin{description}
\item[($P_\ell$)]
 $|f^{(\ell')}|_h=O(|z|^{-\ell'-1})$
 for $\ell'\geq\ell$.
\item[($Q_\ell$)]
 $|\mu^{(\ell')}|_h=O(|z|^{-\ell'})$
 for $\ell'\geq\ell$.
\item[($R_\ell$)]
 $|\nbigr_{\gminib}^{(\ell')}|_h=O\bigl(
 \exp(-C|z|^{-\ell'})
 \bigr)$
for $\ell'\geq \ell$
and for $\gminib\in\Irr(\theta,\ell')$.
\end{description}
The asymptotic orthogonality of
$E_{\gminia,\alpha}$ and $E_{\gminib,\beta}$ 
$(\gminia\neq\gminib)$
follows from $(R_{1})$.

In the proof of Theorem 7.2.1 of \cite{Mochizuki-wild},
we proved the claims for a wild harmonic bundle
by using a descending induction on $\ell$.
The essentially same argument can work.
We give an indication for a modification in this situation.

We have the expression
$\theta^{\dagger}=f^{\dagger}d\zbar$.
Let $\Delta:=-\del_z\del_{\zbar}$.
If a holomorphic section $s$ of $\End(E)$
satisfies $[s,f]=0$,
we obtain the following inequality from
(\ref{eq;13.1.11.1}):
\begin{equation}
 \label{eq;13.1.11.2}
\Delta\log|s|^2_h
\leq
-\frac{\bigl|[f^{\dagger},s]\bigr|_h^2}{|s|_h^2}
+C_0\exp(-\epsilon_0|z|^p)
\end{equation}
Let $f^{(\ell)\dagger}$ denote the adjoint of $f^{(\ell)}$
with respect to $h$.
Suppose $P_{\ell+1}$, $Q_{\ell+1}$ and $R_{\ell+1}$.
By applying (\ref{eq;13.1.11.2}) to $f^{(\ell)}$,
we obtain the following, as in (99) of \cite{Mochizuki-wild}:
\[
 \Delta\log|f^{(\ell)}|^2_h
\leq
-\frac{\bigl|[f^{(\ell)\dagger},f^{(\ell)}]\bigr|_h^2}{|f^{(\ell)}|_h^2}
+C_{1}
\]
Then, by the same argument
as that in \S7.3.2--\S7.3.3 of \cite{Mochizuki-wild},
we obtain $P_{\ell}$ and $Q_{\ell}$.
We put
\[
 k^{(\ell)}_{\gminib}:=
\log\bigl(|\pi^{(\ell)}_{\gminib}|_h^2\big/
 |\pi^{\prime(\ell)}_{\gminib}|_h^2\bigr)
=\log\bigl(1+|\nbigr^{(\ell)}_{\gminib}|_h^2\big/
 |\pi^{\prime(\ell)}_{\gminib}|_h^2\bigr).
\]
By applying (\ref{eq;13.1.11.2})
to $\pi^{(\ell)}_{\gminib}$,
we obtain
\[
 \Delta\log k^{(\ell)}_{\gminib}
\leq
 -\frac{\bigl|[f^{\dagger},\pi_{\gminib}^{(\ell)}]\bigr|_h^2}
 {|\pi_{\gminib}^{(\ell)}|_h^2}
+C_0\exp\bigl(
 -\epsilon_0|z|^{p}
 \bigr).
\]
There exists $C_1>0$ and $R_1>0$
such that
the following holds for any $|z|<R_1$:
\begin{multline}
 \Delta \exp(-A|z|^{-\ell})
\geq
 -\exp(-A|z|^{-\ell})\,
\Bigl(
 \frac{\ell^2}{4}A^2|z|^{-2(\ell+1)}
 \Bigr)
 \\
\geq
 -\exp(-A|z|^{-\ell})\,
 \frac{\ell^2}{4}A^2C_1|z|^{-2(\ell+1)}
+C_0\exp(-\epsilon_0|z|^p)
\end{multline}
Hence, 
we obtain $R_{\ell}$
by using the argument in \S7.3.4 of \cite{Mochizuki-wild}.
Similarly,
we obtain the asymptotic orthogonality
of $E_{\gminia,\alpha}$ and $E_{\gminia,\beta}$ $(\alpha\neq\beta)$,
and the boundedness of
$\theta_{\gminia,\alpha}-(d\gminia+\alpha dz/z)\id_{E_{\gminia,\alpha}}$
by using the argument in \S7.3.5--\S7.3.7 of \cite{Mochizuki-wild}
with (\ref{eq;13.1.11.2}).
\hfill\qed

\vspace{.1in}
We obtain the following corollary.
(See \S7.2.5 of \cite{Mochizuki-wild} for the argument.)
\begin{cor}
\label{cor;13.1.11.22}
$(E,\delbar_E,h)$ is acceptable,
i.e., the curvature
$F(h)$ is bounded with respect to $h$ and $g_{\Poin}$.
\hfill\qed
\end{cor}

\subsubsection{Prolongation and the norm estimate}
\label{subsection;13.1.11.22}

For any $U\subset X$ and for any $a\in\real$,
let $\nbigp_aE(U)$ denote the space of
holomorphic sections $s$ of $E_{|U\setminus D}$
such that $|s|_h=O(|z|^{-a-\epsilon})$ $(\forall \epsilon)$
locally around any point of $U$.
(See \S\ref{subsection;13.10.28.3}.)
According to a general theory of acceptable bundles,
we obtain a locally free $\nbigo_X$-module
$\nbigp_aE$,
and a filtered bundle 
$\nbigp_{\ast}E=(\nbigp_{a}E\,|\,a\in\real)$.
(See \S\ref{subsection;13.2.26.3}
for a review of filtered bundles.)
The decomposition (\ref{eq;13.1.11.3})
is extended to a decomposition of
$\nbigp_aE$:
\[
 \nbigp_aE
=\bigoplus \nbigp_aE_{\gminia,\alpha}
\]
We set 
$\nbigp E:=
 \bigcup_{a\in\real}\nbigp_aE$
and 
$\nbigp E_{\gminia,\alpha}:=
 \bigcup_{a\in\real}\nbigp_aE_{\gminia,\alpha}$.
We set
$\Gr^{\nbigp}_a(E):=\nbigp_aE/\nbigp_{<a}E$,
which we naturally regard as $\cnum$-vector spaces.

By Proposition \ref{prop;12.8.7.3},
$\theta$ gives a section of $\End(\nbigp E)\otimes\Omega^1_X$,
which preserves the decomposition
$\nbigp E=\bigoplus \nbigp E_{\gminia,\alpha}$.
By the estimate in Proposition \ref{prop;12.8.7.3},
$\theta_{\gminia,\alpha}
-(d\gminia+\alpha dz/z)\id_{E_{\gminia,\alpha}}$
is logarithmic with respect to the lattice
$\nbigp_aE_{\gminia,\alpha}$.
Hence, we have the induced endomorphism
$\Res(\theta_{\gminia,\alpha})$
of $\Gr^{\nbigp}_aE_{\gminia,\alpha}$,
which has a unique eigenvalue $\alpha$.
We set
$\Res(\theta)=\bigoplus\Res(\theta_{\gminia,\alpha})$.
Let $W\Gr^{\nbigp}_a(E)$ be the monodromy weight filtration
of the nilpotent part of $\Res(\theta)$.

For each section $s$ of $\nbigp E$,
let $\deg^{\nbigp}(s):=\min\{a\,|\,s\in\nbigp_a E\}$.
For any $g\in \Gr^{\nbigp}_aE$,
let $\deg^W(g):=\min\{m\,|\,g\in W_m\}$.
Let $\vecv=(v_i)$ be a frame of $\nbigp_aE$
which is compatible with
the decomposition $\nbigp_aE=\bigoplus \nbigp_a E_{\gminia,\alpha}$,
the parabolic filtration and the weight filtration,
i.e., 
each $v_i$ is a section of a direct summand
$E_{\gminia,\alpha}$,
the tuple $\vecv^{(b)}:=(v_i\,|\,\deg^{\nbigp}v_i=b)$
induces a base $[\vecv^{(b)}]:=([v^{(b)}_i])$
of $\Gr^{\nbigp}_bE$
for any $a-1<b\leq a$,
and the tuple
$[\vecv^{(b),m}]:=([v^{(b)}_i]\,|\,\deg^Wv^{(b)}_i=m)$
induces a base of
$\Gr^W_m\Gr^{\nbigp}_bE$.
We set
$a_i:=\deg^{\nbigp}(v_i)$
and $k_i:=\deg^W(v_i)$.
Let $h_0$ be the metric of $E$ determined by
$h_0(v_i,v_i)=|z|^{-2a_i}(-\log|z|)^{k_i}$
and $h_0(v_i,v_j)=0$ $(i\neq j)$.
The following proposition can be shown by the argument
in \S8.1.2 of \cite{Mochizuki-wild}.
\begin{prop}
\label{prop;13.1.26.1}
$h$ and $h_0$ are mutually bounded.
\hfill\qed
\end{prop}

\subsubsection{Connection form}

Let $\vecv$ be a frame of $\nbigp_aE$,
which is compatible with the decomposition 
$\nbigp_aE=\bigoplus \nbigp_aE_{\gminia,\alpha}$,
the parabolic filtration
and the weight filtration.
Let $G$ be the endomorphism of $E$
determined by $G(v_i)\,dz=\del v_i$ for $i=1,\ldots,\rank E$.
We can show the following by the arguments
of Lemma 7.5.5, Lemma 10.1.3 and Proposition 10.3.3
of \cite{Mochizuki-wild}.
\begin{lem}
\label{lem;13.1.11.4}
We have $|G|_h=O(|z|^{-1})$.
For the decomposition
$G=\sum G_{(\gminia,\alpha),(\gminib,\beta)}$
according to
$E=\bigoplus E_{\gminia,\alpha}$,
we have the following estimate for some $\epsilon>0$:
\[
 \bigl|
 G_{(\gminia,\alpha),(\gminib,\beta)}
 \bigr|_h
=\left\{
 \begin{array}{ll}
 O\bigl(\exp(-\epsilon|z|^{\ord(\gminia-\gminib)})\bigr)
 & (\gminia\neq\gminib)\\
 \mbox{{}}\\
O(|z|^{-1+\epsilon}) &(\gminia=\gminib,\alpha\neq\beta)
 \end{array}
 \right.
\]
\hfill\qed
\end{lem}

We have the expression $\theta=f\,dz$.
Let us consider $\del_hf$.
Let $\Theta$ be determined by
$f\vecv=\vecv\Theta$.
Let $C$ be determined by
$\del_h\vecv=\vecv C$.
We have
$(\del_hf)\vecv=
 \vecv\bigl(
 \del_z\Theta\,dz+[C,\Theta]
 \bigr)$
and
$[G,f]\vecv=\vecv[C,\Theta]$.
We have the decompositions
$\del_hf=\sum (\del_hf)_{(\gminia,\alpha),(\gminib,\beta)}$
and
$\delbar f^{\dagger}=
 \sum (\delbar f^{\dagger})_{(\gminia,\alpha),(\gminib,\beta)}$
according to
$E=\bigoplus E_{\gminia,\alpha}$.

\begin{cor}
\label{cor;13.2.15.11}
Let $m:=\min\{\ord(\gminia)\,|\,\gminia\in\Irr(\theta)\}$.
If $m<0$,
we have
$\del_{h} f=O\bigl(
 |z|^{-2+m}dz
 \bigr)$ 
with respect to $h$ and $dz\,d\zbar$.
We have
\[
 \bigl|
 (\del_{h} f)_{(\gminia,\alpha),(\gminib,\beta)}
 \bigr|_h
=\left\{
 \begin{array}{ll}
 O\bigl(
 \exp(-\epsilon|z|^{\ord(\gminia-\gminib)})
 \bigr)
 & (\gminia\neq\gminib)\\
 \mbox{{}}\\
O(|z|^{\epsilon-2}) & (\gminia=\gminib,\,\alpha\neq\beta)
 \end{array}
 \right.
\]
We also have the following:
\[
 \bigl|
 (\delbar_E f^{\dagger})_{(\gminia,\alpha),(\gminib,\beta)}
 \bigr|_h
=\left\{
 \begin{array}{ll}
 O\bigl(
 \exp(-\epsilon|z|^{\ord(\gminia-\gminib)})
 \bigr)
 & (\gminia\neq\gminib)\\
 \mbox{{}}\\
O(|z|^{\epsilon-2}) & (\gminia=\gminib,\,\alpha\neq\beta)
 \end{array}
 \right.
\]
\end{cor}
\pf
It follows from Lemma \ref{lem;13.1.11.4}.
\hfill\qed

\subsubsection{Some estimate}

Let $t$ be a $C^{\infty}$-endomorphism of $E$.
According to the decomposition
$E=\bigoplus E_{\gminia,\alpha}$,
we have the decomposition
$t=\sum t_{(\gminia,\alpha)(\gminib,\beta)}$,
where
$t_{(\gminia,\alpha),(\gminib,\beta)}
 \in \Hom(E_{\gminib,\beta},E_{\gminia,\alpha})$.
Let $\nbigc$ be the set of $C^{\infty}$-endomorphisms $t$
such that the following holds for some $\epsilon>0$
which may depend on $t$:
\[
 \bigl|
 t_{(\gminia,\alpha),(\gminib,\beta)}
 \bigr|_h
=\left\{
 \begin{array}{ll}
 O\bigl(
 |z|^{\epsilon}\exp(-\epsilon|z|^{\ord(\gminia-\gminib)})
 \bigr)
 & ((\gminia,\alpha)\neq(\gminib,\beta))
 \\
 \mbox{{}}\\
O(1)
 & (\mbox{\rm otherwise})
 \end{array}
 \right.
\]
Note that $\nbigc$ is closed under the addition and 
the composition.

\begin{prop}
\label{prop;13.2.15.1}
Suppose $t$ and $|z|^2\del_z\del_{\zbar}t$
are contained in $\nbigc$.
Then,
$z\del_zt$ and $\zbar\del_{\zbar}t$
are also contained in $\nbigc$.
\end{prop}
\pf
Let $\Psi:\hyperh:=
 \{u\in\cnum\,|\,\Image u>0\}\lrarr 
 \bigl\{z\in\cnum\,\big|\,0<|z|<1\bigr\}$
be given by $\Psi(u)=\exp(\sqrt{-1}u)$.
Because $\Psi^{\ast}t$ and $\del_u\del_{\ubar}\Psi^{\ast}(t)$
are bounded,
we obtain that $\del_u\Psi^{\ast}t$ and $\del_{\ubar}\Psi^{\ast}t$
are also bounded.

In the following argument,
positive constants $\epsilon$ can change.
We use the notation in the proof of Proposition \ref{prop;12.8.7.3}.
We clearly have $\del_{\zbar}\pi^{(\ell)}_{\gminib}=0$.
We have
$\del_z\pi^{(\ell)}_{\gminib}=O\bigl(\exp(-\epsilon|z|^{-\ell})\bigr)$
by Lemma \ref{lem;13.1.11.4}.
We also have
$\del_{\zbar}\del_z\pi^{(\ell)}_{\gminib}
=\bigl[F(h),\pi^{(\ell)}_{\gminib}\bigr]
=O\bigl(\exp(-\epsilon|z|^{-\ell})\bigr)$.

We have the decomposition
$t=\sum t^{(\ell)}_{\gminia,\gminib}$
according to the decomposition
$E=\bigoplus E^{(\ell)}_{\gminia}$.
We have 
$t^{(\ell)}_{\gminia,\gminib}
=O\bigl(
 \exp(-\epsilon|z|^{-\ell})
 \bigr)$
if $\gminia\neq\gminib$.
Hence,
we have
\[
 [t,\pi^{(\ell)}_{\gminib}]=
\sum_{\gminia\neq\gminib} t^{(\ell)}_{\gminia,\gminib}
-\sum_{\gminia\neq\gminib} t^{(\ell)}_{\gminib,\gminia}
=O\bigl(\exp(-\epsilon|z|^{-\ell})\bigr)
\]
We also have
$|z|^2\del_{\zbar}\del_z[t,\pi^{(\ell)}_{\gminib}]
=\bigl[
 |z|^2\del_{\zbar}\del_zt,\pi^{(\ell)}_{\gminib}
 \bigr]
+\bigl[\zbar\del_{\zbar}t,z\del_z\pi^{(\ell)}_{\gminib}\bigr]
+\bigl[t,|z|^2\del_{\zbar}\del_z\pi^{(\ell)}_{\gminib}\bigr]
=O\bigl(
 \exp\bigl(
 -\epsilon|z|^{-\ell}
 \bigr)
 \bigr)$.
Hence, we obtain
$z\del_z[t,\pi^{(\ell)}_{\gminib}]=
 O\bigl(
 \exp(-\epsilon|z|^{-\ell})
 \bigr)$
and 
$\zbar\del_{\zbar}[t,\pi^{(\ell)}_{\gminib}]=
 O\bigl(
 \exp(-\epsilon|z|^{-\ell})
 \bigr)$.
Therefore, we obtain
$z\del_z t^{(\ell)}_{\gminia,\gminib}
=O\bigl(
 \exp(-\epsilon|z|^{-\ell})
 \bigr)$
and 
$\zbar\del_{\zbar} t^{(\ell)}_{\gminia,\gminib}
=O\bigl(
 \exp(-\epsilon|z|^{-\ell})
 \bigr)$
for $\gminia\neq\gminib$.

We have 
$z\del_z\pi_{\gminia,\alpha}=O(|z|^{\epsilon})$
and 
$|z|^2\del_{\zbar}\del_z\pi_{\gminia,\alpha}
=O(|z|^{\epsilon})$
by Lemma \ref{lem;13.1.11.4}.
Then, we obtain
$z\del_z t_{(\gminia,\alpha),(\gminib,\beta)}
=|z|^{\epsilon}$
and
$\zbar\del_{\zbar} t_{(\gminia,\alpha),(\gminia,\beta)}
=|z|^{\epsilon}$
for $\alpha\neq\beta$.
If $\gminia\neq\gminib$ with $\ell=\ord(\gminia-\gminib)$,
we obtain the desired estimate
by using
$t_{(\gminia,\alpha),(\gminib,\beta)}
=\pi_{\gminia,\alpha}\circ 
 t^{(\ell)}_{\eta_{\ell}(\gminia),\eta_{\ell}(\gminib)}
 \circ
\pi_{\gminib,\beta}$.
\hfill\qed

\subsubsection{Refined asymptotic orthogonality}
\label{subsection;12.8.7.4}

We obtain an asymptotic orthogonality of the derivative
by assuming the following 
with respect to $h$ and $dz\,d\zbar$,
in addition to (\ref{eq;13.1.11.1}):
\begin{equation}
\del_{\zbar}\del_z
 \bigl(F(h)+[\theta,\theta^{\dagger}]\bigr)
=O\bigl(
 \exp(-\epsilon_0|z|^p)
 \bigr).
\end{equation}

Let $\vecv$ be a holomorphic frame of $\nbigp_0E$,
compatible with the decomposition
$\nbigp_0E=
 \bigoplus \nbigp_0E_{\gminia,\alpha}$,
the parabolic filtration and the weight filtration.
Let $(\gminia_i,\alpha_i)$
be determined by
$v_i\in \nbigp_0E_{\gminia_i,\alpha_i}$.
We say that a matrix valued function $B=(B_{ij})$ 
satisfies the condition $\nbigc_1$,
if the following holds for some $\epsilon>0$
which may depend on $B$:
\[
 B_{ij}=
\left\{
 \begin{array}{ll}
 O\bigl(
 |z|^{\epsilon}
 \exp\bigl(-\epsilon|z|^{\ord(\gminia_i-\gminia_j)}\bigr)
 \bigr)
 &
 ((\gminia_i,\alpha_i)\neq (\gminia_j,\alpha_j))
 \\
 \mbox{{}}
 \\
 O\bigl(
 |v_i|_h |v_j|_h
 \bigr)
 & (\mbox{\rm otherwise})
 \end{array}
\right.
\]
Let $H$ be the matrix valued function
determined by $H_{ij}=h(v_i,v_j)$.
Lemma \ref{lem;13.1.11.4} implies that
$z\del_z H$
and $\zbar\del_{\zbar}H$ satisfy the condition $\nbigc_1$.

\begin{prop}
\label{prop;13.2.15.3}
$(|z|^2\del_{\zbar}\del_z)^2H$
satisfies the condition $\nbigc_1$.
\end{prop}
\pf
Let $G(A)$ denote the endomorphism determined
by $\vecv$ and a matrix-valued function $A$.
By Lemma \ref{lem;13.1.11.4},
we have 
\[
G(H^{-1}z\del_zH),
G(H^{-1}\zbar\del_{\zbar}H),
G(\Hbar^{-1}z\del_z\Hbar),
G(\Hbar^{-1}\zbar\del_{\zbar}\Hbar)
\in\nbigc.
\]
Because
$G(\zbar\del_{\zbar}(\Hbar^{-1}z\del_z\Hbar))=|z|^2F(h)
\in\nbigc$,
we have
$G(\Hbar^{-1}|z|^2\del_{\zbar}\del_z\Hbar)
\in\nbigc$.

We have the expression $\theta=f\,dz$.
We have
$\del_{\zbar}\del_z[f,f^{\dagger}]
=\bigl[
 [F(h)_{\zbar,z},f],f^{\dagger}
 \bigr]
+\bigl[
 \del_zf,\delbar_zf^{\dagger}
 \bigr]$.
It gives an estimate for 
$\del_{\zbar}\del_z[f,f^{\dagger}]$
by Corollary \ref{cor;13.2.15.11},
from which we can deduce that
$|z|^2\del_z\del_{\zbar}(|z|^2F(h))\in\nbigc$.
By Proposition \ref{prop;13.2.15.1},
we obtain
$z\del_z(|z|^2F(h))\in\nbigc$
and 
$\zbar\del_{\zbar}(|z|^2F(h))\in\nbigc$.
We obtain
$G(\zbar\del_{\zbar}(\zbar\del_{\zbar}(\Hbar^{-1}z\del_z\Hbar))),
G(z\del_{z}(\zbar\del_{\zbar}(\Hbar^{-1}z\del_z\Hbar)))
\in\nbigc$.
We obtain
$G\bigl(\Hbar^{-1}(\zbar\del_{\zbar})^2z\del_z\Hbar\bigr),
G\bigl(\Hbar^{-1}\zbar\del_{\zbar}(z\del_z)^2\Hbar\bigr)
\in\nbigc$.
Then, we obtain
$G\bigl(
 \Hbar^{-1}(z\del_z)^2(\zbar\del_{\zbar}^2\Hbar)
 \bigr)
\in \nbigc$
from $|z|^2\del_z\del_{\zbar}(|z|^2F(h))\in\nbigc$.
It implies the claim of the lemma.
\hfill\qed

\begin{cor}
\label{cor;13.2.15.10}
$(z\del_z)^2H$ satisfies the condition $\nbigc_1$.
\hfill\qed
\end{cor}

\begin{rem}
The estimate as in Corollary {\rm\ref{cor;13.2.15.10}}
will be used in the study for the extension of
the associated twistor family,
which will be discussed elsewhere.
\hfill\qed
\end{rem}

\section{$L^2$-instantons on $T\times\cnum$}

\subsection{Some standard property}

\subsubsection{Instantons of rank one}

Let $(E,\nabla,h)$ be an $L^2$-instanton on $T\times\cnum$
with $\rank E=1$.
\begin{lem}
$(E,\nabla,h)$ is a unitary flat bundle.
\end{lem}
\pf
Because $\rank E=1$,
we have
$\bigl(
 \nabla_{z}\nabla_{\zbar}+\nabla_w\nabla_{\wbar}
\bigr)F_{z\zbar}=0$
and 
$\bigl(
 \nabla_{z}\nabla_{\zbar}+\nabla_w\nabla_{\wbar}
\bigr)F_{z\wbar}=0$.
We obtain the following inequalities:
\[
 -(\del_w\del_{\wbar}+\del_{z}\del_{\zbar})
 \bigl|F_{z\zbar}\bigr|^2
\leq
 0,
\quad
 -(\del_w\del_{\wbar}+\del_{z}\del_{\zbar})
 \bigl|F_{z\wbar}\bigr|^2
\leq
 0.
\]
We use the notation 
in \S\ref{subsection;13.1.10.4}.
By applying the fiber integral for 
$T\times\cnum\lrarr\cnum$,
we obtain
$-\del_w\del_{\wbar}
 \|F_{z\zbar}\|^2\leq 0$
and
$ -\del_w\del_{\wbar}
 \|F_{z\wbar}\|^2\leq 0$.
Because the functions
$\|F_{z\zbar}\|^2$ 
and $\|F_{z\wbar}\|^2$
are $L^1$ on $\cnum_w$,
they are $0$.
\hfill\qed

\begin{cor}
\label{cor;13.1.29.32}
Let $(E,\nabla,h)$ be an $L^2$-instanton on
$T\times\cnum$ of an arbitrary rank.
Then,
$\det(E,\nabla,h)$ is a flat unitary bundle,
i.e.,
we have $\Tr F(\nabla)=0$.
\hfill\qed
\end{cor}

If we do not impose the $L^2$-property,
there exist much more instantons of rank one
on $T\times\cnum$.

{\bf (i)} Let $\gminia$ be any holomorphic function on $\cnum$.
Then, the trivial holomorphic line bundle $\nbigo_{\cnum}$
with the trivial metric and the Higgs field $d\gminia$
gives a harmonic bundle $\nbigl(\gminia)$ on $\cnum$.
By the equivalence of Hitchin,
we have the associated instanton on $T\times\cnum$.

{\bf (ii)} 
Let $\rho$ be an $\real$-valued harmonic function
on $T\times\cnum$.
Then, the trivial holomorphic line bundle
$\nbigo_{T\times\cnum}\,e$ with the metric
$h_{\rho}$ given by $\log h_{\rho}(e,e)=\rho$
gives an instanton $\nbigl(\rho)$ on $T\times\cnum$.
Note that there exist many harmonic functions
which is not the real part of a holomorphic function on
$T\times\cnum$.
We can construct such a function by using a Bessel function
$I_0(r)=\int_{-1}^1\cosh(rt)(t^2-1)^{-1/2}dt$
which satisfies $I_0''+r^{-1}I_0'-I_0=0$.
It is a $C^{\infty}$-function on $\real$,
satisfying $I_0(r)=I_0(-r)$.
In particular,
$\kappa(w):=I_0(|w|)$ gives a $C^{\infty}$-function on $\cnum$
satisfying
$(-\del_w\del_{\wbar}+4)\kappa=0$.
By using the Fourier expansion on $T\times\cnum$ in a standard way,
we can construct a harmonic function $\rho$ on $T\times\cnum$
from $\kappa$ such that $\rho$ is not constant along $T$.
(See \cite{Jardim2}.)
It is not the real part of any holomorphic function.

In general,
any instanton of rank one $(E,\delbar_E,h)$
can be expressed
as the tensor product of instantons of types (i) and (ii).
Indeed,
by considering the support $\RFM_-(E,\delbar_E)$,
we obtain a holomorphic function $\cnum\lrarr T^{\lor}$.
Because $\cnum$ is simply connected,
it is lifted to a holomorphic function
$\gminib:\cnum\lrarr \cnum$.
We have a holomorphic function $\gminia$
such that $\del_w\gminia=\gminib$.
Then, we can observe that
$(E,\delbar_E,h)$ is isomorphic to
$\nbigl(\gminia)\otimes \nbigl(\rho)$
for a harmonic function $\rho$ on $T\times\cnum$.

\subsubsection{Polystability of the associated filtered bundle}

Let $(E,\nabla,h)$ be an $L^2$-instanton on $T\times\cnum$.
Let $(E,\delbar_E)$ be the underlying holomorphic 
vector bundle on $T\times\cnum$.
For a saturated $\nbigo_{T\times\cnum}$-subsheaf $\nbigf\subset E$,
let $h_{\nbigf}$ denote the induced hermitian metric
of the smooth part of $\nbigf$.
Let $F(h_{\nbigf})$ denote the curvature.
As in \cite{Biquard-Jardim} and \cite{Simpson88},
we set
\[
 \deg(\nbigf,h):=
 \sqrt{-1}
 \int_{T\times\cnum}
 \Tr\bigl(\Lambda F(h_{\nbigf})\bigr)
 \dvol_{T\times\cnum}.
\]
Let $\pi_{\nbigf}$ denote the orthogonal projection
of $E$ to $\nbigf$,
where it is considered only on the smooth part of $\nbigf$.
By the Chern-Weil formula \cite{Simpson88},
we have
\[
 \deg(\nbigf,h)
=-\int_{T\times\cnum}|\delbar \pi|_h^2\dvol_{T\times\cnum}.
\]
\begin{lem}
\label{lem;13.1.29.1}
$\deg(\nbigf,h)$ is finite,
if and only if
(i) the degree of $\nbigf_{|T\times\{w\}}$ are $0$ 
 for any $w\in\cnum$,
(ii) $\nbigf$ is extended to a subsheaf 
$\nbigp_0\nbigf$ of $\nbigp_0E$.
In that case,
we have 
$\deg(\nbigf,h)=
 \int_{z\times\proj^1}
 \parchern_1(\nbigp_{\ast}\nbigf)$,
where $\nbigp_{\ast}\nbigf$
denotes $\nbigp_0\nbigf$
with the induced parabolic structure.
\end{lem}
\pf
This type of claim is standard 
in the study of Kobayashi-Hitchin correspondence
for parabolic objects,
and well established in \cite{Li-Narasimhan},
based on the fundamental results
in \cite{Simpson88}, \cite{Simpson90}
and \cite{Siu-extension}.
We give only an indication in our situation.

By Lemma 10.5 and Lemma 10.6 of \cite{Simpson88},
$\nbigf_{|\{z\}\times \cnum}$
is extended to a parabolic subsheaf,
if and only if
$\int_{\cnum} |\delbar\pi_{|\{z\}\times \cnum}|^2<\infty$.
In that case,
$\frac{\sqrt{-1}}{2\pi}\int_{\cnum}
 Tr(F(h_{\nbigf}))_{|z\times\cnum}$
is equal to the parabolic degree
of the parabolic subsheaf.

If the conditions (i) and (ii) are satisfied, 
then we have
\[
 \deg(\nbigf,h)=
 \int_{T}\dvol_T
 \Bigl(
 \int_{z\times\cnum}
 \sqrt{-1}\Tr\bigl(F(h_{\nbigf})\bigr)
 \Bigr)
=2\pi|T|\int_{z\times\proj^1}
 \parchern_1(\nbigp_{\ast}\nbigf)>-\infty.
\]
Conversely, suppose $\deg(\nbigf,h)$ is finite.
Because
\[
 \deg(\nbigf,h)
\leq
 \int_{\cnum}
 \dvol_{\cnum}
 \Bigl(
 -\int_{T}
 \bigl|
 \nabla_{\zbar}\pi
 \bigr|^2
 \dvol_{T}
 \Bigr)
=2\pi\int_{\cnum}
 \deg(\nbigf_{|T\times \{w\}})
 \dvol_{\cnum},
\]
we have $\deg(\nbigf_{|T\times\{w\}})=0$
for any $w$.
We obtain
$-\deg(\nbigf,h)
=\int_T
 \dvol_{T}
 \Bigl(
 \int_{\cnum}
\bigl|
 \delbar\pi_{|\{z\}\times\cnum}
 \bigr|^2
 \Bigr)
<\infty$.
Hence, 
there exists a thick subset $A\subset T^{\lor}$
such that 
$\nbigf_{|z\times\cnum}$
is extendable
for any $z\in A$.
(A subset is thick,
if it is not contained
in a countable union of a complex
analytically closed subsets.)
Then, we obtain that
$\nbigf$ is extendable
according to Theorem 4.5 of \cite{Siu-extension}.
\hfill\qed

\begin{prop}
\label{prop;13.1.24.3}
$\nbigp_{\ast}E$ is polystable.
We have
$\deg(\nbigp_{\ast}E)=0$.
(See {\rm\S\ref{subsection;13.2.14.11}}
for the stability condition in this case.)
\end{prop}
\pf
The second claim directly follows from
Lemma \ref{lem;13.1.29.1}
and Corollary \ref{cor;13.1.29.32}.
Let $\nbigp_{\ast}\nbigf$ be a filtered subsheaf
$\nbigp_{\ast}E$ satisfying {\bf(A1,2)}
in \S\ref{subsection;13.2.14.11}.
Let $\nbigf$ be its restriction to $X\times\cnum$.
By Lemma \ref{lem;13.1.29.1},
we have 
$\mu(\nbigp_{\ast}\nbigf)
=\mu(\nbigf,h)\leq 0$.
Moreover, if it is $0$,
the orthogonal projection onto $\nbigf$ is holomorphic.
Hence, the orthogonal decomposition
$E=\nbigf\oplus \nbigf^{\bot}$ is holomorphic.
It is extended to a decomposition
$\nbigp_{\ast}E=
 \nbigp_{\ast}\nbigf\oplus
 \nbigp_{\ast}\nbigf^{\bot}$.
Both $\nbigf$ and $\nbigf^{\bot}$
with the induced metrics are $L^2$-instantons.
Hence, we obtain the first claim of the corollary
by an easy induction on the rank.
\hfill\qed

\subsubsection{Uniqueness of the $L^2$-instanton 
adapted to a filtered bundle}

Let $(E,\nabla,h)$ be an $L^2$-instanton on $T\times\cnum$.
We have the associated filtered bundle $\nbigp_{\ast}E$
on $(T\times\proj^1,T\times\{\infty\})$.
Let $h'$ be a hermitian metric of $E$,
and let $\nabla_{h'}$ be a unitary connection
of $(E,h')$
such that 
(i) $(E,\nabla_{h'},h')$ is an $L^2$-instanton,
(ii) the $(0,1)$-parts of $\nabla_{h'}$ and $\nabla_{h}$ are equal,
(iii) $h'$ is adapted to $\nbigp_{\ast}E$.
(See \S\ref{subsection;13.10.28.3}
for adaptedness.)

\begin{prop}
\label{prop;13.1.29.20}
We have a holomorphic decomposition
$(E,\delbar_E)=
 \bigoplus_i (E_{i},\delbar_{E_i})$
such that
(i) it is orthogonal with respect to both $h$ and $h'$,
(ii) for each $i$, there exists $\alpha_i>0$
 such that $h_{|E_i}=\alpha_i\,h'_{|E_i}$.
In particular,
we have $\nabla_h=\nabla_{h'}$.
\end{prop}
\pf
Let $s$ be the self-adjoint endomorphism
determined by $h'=h\,s$.
According to \cite{Simpson88},
we have the following inequality
(see p.876 of \cite{Simpson88}):
\[
-(\delbar_z\del_z+\delbar_w\del_w) \Tr(s)
+\bigl|
 \delbar (s)\,s^{-1/2}
 \bigr|_h^2
\leq 0
\]
By taking the fiber integral for
$T\times\cnum\lrarr\cnum$,
we obtain
\[
-\del_{\wbar}\del_w
 \int_{T}\Tr(s)
+\int_{T}\bigl|
 \delbar (s)\,s^{-1/2}
 \bigr|_h^2
 \leq 0
\]
It implies that 
$\int_{T}\Tr(s)$ is a subharmonic function on $\cnum_w$.
By using the norm estimate for asymptotically harmonic bundle
(Proposition \ref{prop;13.1.26.1}),
we obtain that $h$ and $h'$ are mutually bounded,
i.e.,
$s$ and $s^{-1}$ are bounded
with respect to both of $h$ and $h'$.
Hence, we obtain that $\int_T\Tr(s)$ is constant.
We obtain
$\int_{T}\bigl|
 \delbar (s)\,s^{-1/2}
 \bigr|_h^2=0$,
which implies
$\delbar(s)=0$.
Then, the claim of the proposition follows.
\hfill\qed

\subsubsection{Instanton number}

Let $(E,\nabla,h)$ be an $L^2$-instanton
on $T\times\cnum$.
We have the associated filtered bundle
$\nbigp_{\ast}E$ on $(T\times\proj^1,T\times\{0\})$.
Note that the second Chern class of
$\nbigp_aE$ is independent of $a\in\real$.

\begin{prop}
\label{prop;13.10.15.110}
We have the following equality
for any $a\in\real$.
\[
 \frac{1}{8\pi^2}
 \int_{T\times\cnum} \Tr\bigl(F(h)^2\bigr)
=\int_{T\times\proj^1}c_2(\nbigp_aE)
\]
\end{prop}
\pf
Let $U\subset\proj^1$ be a small neighbourhood of $\infty$
such that 
$\nbigp_aE_{|T\times w}$ are semistable of degree $0$
for any $w\in U$.
In the following argument,
we will shrink $U$.
We have the filtered Higgs bundle 
$(\nbigp_{\ast}V,g\,dw)$ on $(U,\infty)$
corresponding to $\nbigp_{\ast}E$.

Let $p:T\times U\lrarr U$ be the projection.
We have a natural $C^{\infty}$-isomorphism
$\nbigp_aE\simeq
 p^{\ast}(\nbigp_aV)$,
and the holomorphic structure of $\nbigp_aE$
is described as 
$p^{\ast}(\delbar_{\nbigp_aV})+g\,d\zbar$.

We take a holomorphic structure 
$\delbar_{\nbigp_aE}'$ of $\nbigp_aE$
such that 
its restriction to
$T\times (\proj^1\setminus U)$ is $\delbar_E$,
and 
its restriction to
$T\times U$ is 
$p^{\ast}\bigl(
 \delbar_{\nbigp_aV}\bigr)$.

We take a holomorphic frame $\vecv$ of $\nbigp_aV$
which is compatible with the parabolic structure.
It induces a $C^{\infty}$-frame $\vecu$ of $\nbigp_aE$
on $T\times U$.
We take a $C^{\infty}$-metric $h_0$ of $\nbigp_aE$
such that $\vecu$ is orthonormal 
with respect to $h_{0|T\times U}$.
Let $\nabla_0$ be the Chern connection 
determined by $h_0$ and 
$\delbar_{\nbigp_aE}'$.
We have $\nabla_0\vecu=0$
on $T\times U$.
The curvature $F(\nabla_0)$ vanishes
on $T\times U$.
We set $A:=\nabla-\nabla_0$.

Let $J$ be the endomorphism
of $E_{|T\times(U\setminus\{\infty\})}$
determined by
$\nabla_wu_i=J(u_i)$
$(i=1,\ldots,\rank E)$.
According to Lemma \ref{lem;13.1.11.4}
and Theorem \ref{thm;12.8.7.1},
we have $J=O(|w|^{-1})$
with respect to $h$.
On $T\times U$, we have 
\[
A=J\,dw+g\,d\zbar-g^{\dagger}_h\,dz.
\]
Here, $g^{\dagger}_h$ denotes
the adjoint of $g$ with respect to $h$.
We have $|g|_h=|g^{\dagger}|_h=O(1)$.
According to Theorem \ref{thm;12.8.7.1}
and Proposition \ref{prop;12.8.7.3},
we have
$[g,g^{\dagger}_h]=
O\bigl(|w|^{-2}(\log|w|)^{-2}\bigr)$
with respect to $h$.
According to Lemma \ref{lem;13.1.11.4}
and Theorem \ref{thm;12.8.7.1},
we have
$[g,J]=O(|w|^{-2})$
and 
$[g^{\dagger}_h,J]=O(|w|^{-2})$
with respect to $h$.
Hence, we have
\[
 A^2=O(|w|^{-2})\,dw\,d\zbar
+O(|w|^{-2})\,dw\,dz
+O(|w|^{-2})\,dz\,d\zbar.
\]

We set
$\nabla_t:=t\nabla+(1-t)\nabla_0$
for $0\leq t\leq 1$.
On $T\times (U\setminus\{\infty\})$,
we have the following estimate
for some $\rho>0$,
which is uniform for $t$:
\[
 F(\nabla_t)
=tF(\nabla)+(t^2-t)A^2
=O(|w|^{-2})\,dw\, d\wbar
+O(|w|^{-1-\rho})\,dw\,d\zbar
+O(|w|^{-1-\rho})\,d\wbar\,dz
+O(|w|^{-2})dz\,d\zbar
\]
We obtain the following estimate,
which is uniform for $t$:
\[
 \Tr\bigl(
 F(\nabla_t)\,
 A
 \bigr)
=O(|w|^{-2})\,dw\,d\wbar\,dz
+O(|w|^{-2})\,dw\,d\wbar\,d\zbar
+O(|w|^{-1-\rho})\,dw\,dz\,d\zbar
+O(|w|^{-1-\rho})\,d\wbar\,dz\,d\zbar
\]
We obtain the following:
\[
 -\frac{1}{8\pi^2}
 \int_{T\times\cnum}
 \Tr(F(h)^2)
=-\frac{1}{8\pi^2}
 \int_{T\times\proj^1}
 \Tr(F(\nabla_0)^2)
=\int_{X\times\proj^1}\ch_2(\nbigp_aE)
=-\int_{X\times\proj^1}c_2(\nbigp_aE)
\]
Thus, we are done.
\hfill\qed

\subsection{Cohomology}
\label{subsection;13.1.24.1}

Let $(E,\nabla,h)$ be an $L^2$-instanton on
$X:=T\times \cnum$.
The $(0,1)$-part of $\nabla$ is denoted by $\delbar_E$.
Let $\Xbar:=T\times\proj^1$.
We put $D:=T\times\{\infty\}$.
Let $A_{c}^{0,i}(E)$ denote the space of
$C^{\infty}$-sections of $E\otimes \Omega^{0,i}$ on $X$
with compact supports.
Its cohomology group is denoted by
$H^{0,i}_c(X,E)$.
Let $A^{0,i}(\nbigp_aE)$ denote the space of
$C^{\infty}$-sections of
$\nbigp_aE\otimes\Omega^{0,i}$ on $\Xbar$.
Its cohomology group is
$H^i(\Xbar,\nbigp_aE)$.
In this subsection, we suppose that
\[
0\not\in\Sp_{\infty}(E).
\]

\begin{prop}
\label{prop;13.2.13.100}
The natural map
$H_c^{0,i}(X,E)\lrarr
 H^{0,i}(\Xbar,\nbigp_aE)$
is an isomorphism
for any $a\in\real$.
\end{prop}
\pf
There exists $R>0$ such that,
if $|w|>R$,
$E_{|T\times \{w\}}$ is semistable of degree $0$,
and 
$0\not\in \Sp(E_{|T\times\{w\}})$.
We have two consequences
for a $C^{\infty}$-section $s$ of $\nbigp_aE$ on $X_R$.
\begin{itemize}
\item
There exists a $C^{\infty}$-section $t$ of
$\nbigp_aE$ on $X_R$
such that $\nabla_{\zbar}t=s$.
\item
If $\nabla_{\zbar}s=0$, then $s=0$.
\end{itemize}
Then, the claim can be shown easily.
\hfill\qed

\vspace{.1in}

Let $A_{L^2}^{0,i}(E)$ be the space of
$L^2$-sections $s$ of
$E\otimes\Omega^{0,i}$ on $X$
such that $\delbar_E s$ is also $L^2$.
The cohomology group of the complex
$\bigl(A^{0,\bullet}_{L^2}(E),\delbar_E\bigr)$
is denoted by $H^{0,i}_{L^2}(X,E)$.

\begin{prop}
\label{prop;13.1.11.30}
The natural map
$H^{0,i}_{c}(X,E)\lrarr
 H^{0,i}_{L^2}(X,E)$
is an isomorphism.
\end{prop}
\pf
Let $A_{L^2,c}^{0,i}(E)\subset A_{L^2}^{0,i}(E)$ 
be the subspace of the sections with compact supports.
It gives a subcomplex,
and its cohomology is denoted by
$H^{0,i}_{L^2,c}(X,E)$.
\begin{lem}
The natural map
$H^{0,i}_{L^2,c}(X,E)
\lrarr
 H^{0,i}_{L^2}(X,E)$
is an isomorphism.
\end{lem}
\pf
For any $L^2$-section $s$ of $E$ on $X_R$,
there exists an $L^2$-section $t$ of $E$ on $X_{R'}$
$(R'>R)$ 
such that
$\nabla_{\zbar}t=s$ on $X_{R'}$.
If an $L^2$-section $s$ of $E$ on $X_R$
satisfies $\nabla_{\zbar}s=0$,
then we have $s=0$.
Then, the claim of the lemma can be shown.
\hfill\qed

\vspace{.1in}
We take a smooth K\"ahler metric of $\Xbar$.
Let $B^{0,i}_{L^2}(\nbigp_aE)$ be the space of
$L^2$-sections $\omega$ of $\nbigp_aE$ on $\Xbar$
such that $\delbar \omega$ is $L^2$.
Let $B^{0,i}_{L^2,c}(\nbigp_aE)\subset
 B^{0,i}_{L^2}(\nbigp_aE)$
denote the subspace of the sections
whose support is contained in $X$.
By the same argument,
the natural map
$B^{0,\bullet}_{L^2,c}(\nbigp_aE)
\lrarr
 B^{0,\bullet}_{L^2}(\nbigp_aE)$
is a quasi isomorphism.
We have a natural identification
$B^{0,\bullet}_{L^2,c}(\nbigp_aE)
=A^{0,\bullet}_{L^2,c}(E)$
as $\cnum$-linear spaces.
By the $L^2$-Dolbeault theorem,
the cohomology group of
$B^{0,\bullet}_{L^2}(\nbigp_aE)$
is naturally isomorphic to 
$H^i(\Xbar,\nbigp_aE)$.
Then, the claim of Proposition \ref{prop;13.1.11.30}
follows.
\hfill\qed

\begin{cor}
$H^{0,i}_{L^2}(\Xbar,E)$ is finite dimensional.
\hfill\qed
\end{cor}

\begin{prop}
We have
$H^0(\Xbar,\nbigp_aE)=H^2(\Xbar,\nbigp_aE)=0$.
\end{prop}
\pf
Clearly $H^0(\Xbar,\nbigp_a(E))=0$.
Let $p:\Xbar\lrarr \proj^1$ be the projection.
We have $p_{\ast}E=0$,
and the support of $R^1p_{\ast}E$
is $0$-dimensional.
Then, we obtain $H^2(\Xbar,\nbigp_aE)=0$.
\hfill\qed

\subsection{Exponential decay of harmonic sections}

\subsubsection{Statement}
\label{subsection;13.1.10.11}

Let $(E,\nabla,h)$ be an $L^2$-instanton
on $T\times\cnum$.
Let $\delbar_E$ denote the $(0,1)$-part of $\nabla$,
and let $\delbar_E^{\ast}$ denote the formal adjoint
with respect to
$h$ and $dz\,d\zbar+dw\,d\wbar$.
We set
$\Delta_E:=\delbar_E^{\ast}\delbar_E
+\delbar_E\delbar_E^{\ast}$.

\begin{prop}
\label{prop;13.1.10.12}
Assume that $0\not\in\Sp_{\infty}(E)$.
Let $\omega$ be an $L^2$-section of
$E\otimes\Omega^{0,1}$ on $T\times\cnum$
such that
$\Delta_E\omega=0$.
Then, we have 
$|\omega|=O\bigl(
 \exp(-C|w|)
 \bigr)$
for some $C>0$.
\end{prop}

\subsubsection{An estimate}
\label{subsection;13.2.16.1}

Take $R>0$,
and put  $Y_R:=\{|w|\geq R\}$
and $X_R:=T\times Y_R$.
Let $(E,\nabla,h)$ be an $L^2$-instanton on $X_R$.

\begin{lem}
\label{lem;13.1.10.20}
Assume that $0\not\in \Sp_{\infty}(E)$.
Suppose that $\omega$ is an $L^2$-section of
$E\otimes\Omega^{0,1}_{Y_R}$ such that
$\delbar_E\omega=\delbar_E^{\ast}\omega=0$.
Then, there exists $C>0$ such that
$\bigl|\omega\bigr|_h
=O\Bigl( \exp(-C|w|) \Bigr)$.
\end{lem}
\pf
Let $\omega=f\,d\zbar+g\,d\wbar$ be a harmonic section.
We have 
$-\nabla_{\wbar}f+\nabla_{\zbar}g=0$
and
$\nabla_zf+\nabla_wg=0$.
We have the following equalities:
\[
 \nabla_w\nabla_{\wbar}f
=\nabla_w(\nabla_{\zbar}g)
=\nabla_{\zbar}\nabla_wg+F_{w\zbar}g
=-\nabla_{\zbar}\nabla_zf+F_{w\zbar}g
=-\nabla_z\nabla_{\zbar}f+F_{z\zbar}f+F_{w\zbar}g
\]
\[
 \nabla_w\nabla_{\wbar}g
=F_{w\wbar}g+\nabla_{\wbar}\nabla_wg
=F_{w\wbar}g+\nabla_{\wbar}(-\nabla_{z}f)
=F_{w\wbar}g-\nabla_z\nabla_{\wbar}f+F_{z\wbar}f
=F_{w\wbar}g-\nabla_z\nabla_{\zbar}g+F_{z\wbar}f
\]
We obtain the following
\begin{multline}
-(\del_w\del_{\wbar}+\del_z\del_{\zbar})(f,f)
\leq
 -(\nabla_{\zbar}f,\nabla_{\zbar}f)
-2\Re\bigl((\nabla_w\nabla_{\wbar}+\nabla_z\nabla_{\zbar})f,f\bigr)
 \\
=-(\nabla_{\zbar}f,\nabla_{\zbar}f)
-2\Re\bigl(
 F_{z\zbar}f+F_{z\wbar}g,\,f
 \bigr)
\end{multline}
Using the notation in \S\ref{subsection;13.1.10.4},
we obtain
\[
 -\del_w\del_{\wbar}\|f\|^2
\leq
 -\|\nabla_{\zbar}f\|^2
+O\bigl(
 \|F\|\,(\|f\|^2+\|g\|^2)
 \bigr)
\]
Similarly, we obtain
\[
 -\del_w\del_{\wbar}\|g\|^2
\leq
 -\|\nabla_{\zbar}g\|^2
+O\bigl(
 \|F\|\,(\|f\|^2+\|g\|^2)
 \bigr)
\]
By the assumption
$0\not\in\Sp_{\infty}(E)$,
there exist $R_1>R$ and $C_1>0$
such that,
if $|w|\geq R_1$,
we have
$\|\del_{\zbar}g\|
\geq
 C_1\|g\|$
and
$\|\del_{\zbar}f\|
\geq
 C_1\|f\|$.
Hence, there exist
$\epsilon>0$ and $R_2>R$ such that
the following holds if $|w|>R_2$:
\begin{equation}
 \label{eq;13.1.10.2}
 -\del_w\del_{\wbar}\bigl(
 \|f\|^2+\|g\|^2
 \bigr)
\leq
 -\epsilon\bigl(
 \|f\|^2+\|g\|^2
 \bigr)
\end{equation}

In general,
if $\varphi$ is a positive $L^1$-subharmonic function on $Y_{R_2}$,
Then $\varphi(w)=O(|w|^{-2})$.
Indeed, by the mean value property,
we have
\[
\varphi(w)
\leq
 \frac{4}{\pi (|w|-R_2)^2}
 \int_{|w-w'|\leq (|w|-R_2)/2}
 \varphi(w')
\leq
 \frac{C_2}{(|w|-R_2)^2}.
\]
Hence,
we have
$\|f\|^2+\|g\|^2=O(|w|^{-2})$.
Then, by a standard argument with (\ref{eq;13.1.10.2}),
we obtain
$\|f\|^2+\|g\|^2=O\bigl(
 \exp(-C_3|w|)
 \bigr)$.
(See the proof of Lemma \ref{lem;13.1.10.3}.)
By a bootstrapping argument,
we obtain $|f(z,w)|=O\bigl(\exp(-C_4|w|)\bigr)$
and $|g(z,w)|=O\bigl(\exp(-C_4|w|)\bigr)$.
\hfill\qed

\subsubsection{Finiteness}

We continue to use the notation in \S\ref{subsection;13.2.16.1}.
Let $\omega$ be a $C^{\infty}$-section of
$E\otimes\Omega^{0,1}$ on 
$X_R$.
Suppose that the support of $\omega$
is contained in $T\times \{|w|\geq R+1\}$.
We set $\nbigd:=\delbar_E+\delbar_E^{\ast}$.
Let $\dvol$ denote the volume form 
induced by the Euclidean metric.

\begin{lem}
\label{lem;13.1.10.10}
Assume that $\omega$ and 
$\Delta_E\omega$ are $L^2$.
Then,
$\delbar_E^{\ast}\omega$
and 
$\delbar_E\omega$
are $L^2$,
and we have
\[
 \int h(\omega,\Delta_E\omega)\dvol
=\int\bigl|\nbigd\omega\bigr|_h^2
 \dvol
\]
\end{lem}
\pf
Let $g:=dz\,d\zbar+dw\,d\wbar$.
Let $|\cdot|_{h,g}$ denote the norm of
sections of $E\otimes\Omega^{\bullet}$
induced by $h$ and $g$.
Let $\chi(t)$ be a non-negative valued $C^{\infty}$-function
such that
$\chi(t)=1$ $(t\leq 0)$
and $\chi(t)=0$ $(t\geq 1)$,
and that
$\del_t(\chi)\big/\chi^{1/2}$ is also $C^{\infty}$.
For a large $N$, we put
$\chi_N(w):=\chi(\log|w|-N)$.
There exists $C_1>0$
such that
$|\del_w\chi_N|\leq C_1|w|^{-1}$,
$|\del_{\wbar}\chi_N|\leq C_1|w|^{-1}$,
and
$|\del_{w}\del_{\wbar}\chi_N|\leq C_1|w|^{-2}$.
We have
\begin{multline}
\left|
 \int\chi_Nh(\omega,\Delta_E\omega)\dvol
-\int\chi_N|\nbigd\omega|_{h}^2\dvol
\right|
\leq 
 \left(
 \int|\delbar\chi_N|_g^2\chi_N^{-1}
 |\omega|_{h,g}^2
 \dvol
 \right)^{1/2}
 \left(
 \int\chi_N|\delbar\omega|_{h,g}^2\dvol
 \right)^{1/2}
 \\
+ \left(
 \int|\del\chi_N|_g^2\,\chi_N^{-1}
 |\omega|_{h,g}^2\dvol
 \right)^{1/2}
 \left(
 \int \chi_N|\delbar\omega|_{h,g}^2\dvol
 \right)^{1/2}
\end{multline}
There exist $C_i>0$ $(i=2,3)$
such that the following holds:
\[
 \int\chi_N|\nbigd\omega|_{h,g}^2\dvol
\leq
C_2\left(
 \int\chi_N|\nbigd\omega|_{h,g}^2\dvol
 \right)^{1/2}
+C_3
\]
Then, the first claim of Lemma \ref{lem;13.1.10.10}
follows.
We have
\[
\left|
\int h(\chi_N\omega,\Delta_E\omega)\dvol
-\int\chi_N|\nbigd\omega|^2_{h,g}\dvol
\right|
\leq
 C_4\int 
 |\delbar\chi_N|_g\,|\omega|_{h,g}
 \,|\delbar\omega|_{h,g}
 \dvol
+C_4\int
 |\del\chi_N|_g\,
 |\omega|_{h,g}\,
 |\delbar^{\ast}\omega|_{h,g}
\]
for some $C_4>0$.
By the first claim,
the integrands of the right hand side
are dominated by some integrable functions,
independently from $N$.
By taking the limit,
we obtain the second claim.
\hfill\qed

\subsubsection{Proof of Proposition \ref{prop;13.1.10.12}}

Let us return to the setting in \S\ref{subsection;13.1.10.11}.
According to Lemma \ref{lem;13.1.10.20},
we have only to show the following lemma
to establish Proposition \ref{prop;13.1.10.12}.

\begin{lem}
$\delbar_E\omega=\delbar_E^{\ast}\omega=0$.
\end{lem}
\pf
By the first claim of Lemma \ref{lem;13.1.10.10},
we obtain that
$\nbigd\omega$ is $L^2$.
By the argument in the proof of the second claim
of the same lemma,
we obtain
$\int|\nbigd\omega|_{h,g}^2\dvol=0$,
i.e.,
$\nbigd\omega=0$.
\hfill\qed

\subsection{Nahm transform for $L^2$-instantons}
\label{subsection;13.2.3.10}
Let $(E,\nabla,h)$ be an $L^2$-instanton on $T\times\cnum$
with $\rank E>1$.
Let $D:=\Sp_{\infty}(E)$.
We shall construct a harmonic bundle on $T^{\lor}\setminus D$
with the method in 
\cite{Donaldson-Kronheimer}
and \cite{Jardim2}.
For any $\zeta\in T^{\lor}\setminus D$,
let $\nbigl_{-\zeta}=(\underline{\cnum},\delbar_T-\zeta d\zbar)$ 
denote the corresponding line bundle on $T$
with the natural hermitian metric.
Let $\Nahm(E,\nabla)_{\zeta}$ denote the space of
$L^2$-harmonic sections of 
$E\otimes \nbigl_{-\zeta}\otimes\Omega^{0,1}$.
It is finite dimensional,
and naturally isomorphic to 
$H^1(T\times\proj^1,\nbigp_{-1}E\otimes \nbigl_{-\zeta})
\simeq
 H^1(T\times\proj^1,\nbigp_0E\otimes \nbigl_{-\zeta})$.
(See \S\ref{subsection;13.1.24.1}.)
The Euclidean metric $dz\,d\zbar+dw\,d\wbar$
of $T\times\cnum$
and the hermitian metric $h$ of $E$
induce a metric $h_1(\zeta)$
of $\Nahm(E,\nabla)_{\zeta}$.
The multiplication of  $-w\in\nbigo_{\proj^1}(1)$
induces an endomorphism $F_w(\zeta)$ of
$\Nahm(E,\nabla)_{\zeta}$.
It is also described as $-P_{\zeta}\circ w$,
where $P_{\zeta}$ denotes the orthogonal projection
of the space of $L^2$-sections of
$E\otimes \nbigl_{-\zeta}\otimes\Omega^{0,1}_{T\times\cnum}$
onto $\Nahm(E,\nabla)_{\zeta}$.
(Note Proposition \ref{prop;13.1.10.12}.)

\vspace{.1in}

Let $A^{p,q}(E\otimes \nbigl_{-\zeta})$ denote 
the space of $L^2$-sections of
$E\otimes \nbigl_{-\zeta}\otimes\Omega_{T\times\cnum}^{p,q}$.
Let $\delbar_{E,\zeta}$ denote the $\delbar$-operator 
of $E\otimes \nbigl_{-\zeta}$.
Let $\delbar_{E,\zeta}^{\ast}$ denote its adjoint.
Let 
$\nbigd_{\zeta}:=\delbar_{E,\zeta}+\delbar_{E,\zeta}^{\ast}$
be the closed operator
$A^{0,0}(E\otimes \nbigl_{-\zeta})\oplus A^{0,2}(E\otimes \nbigl_{-\zeta})
\lrarr A^{0,1}(E\otimes \nbigl_{-\zeta})$,
and let  $\nbigd^{\ast}_{\zeta}:=\delbar_{E,\zeta}+\delbar_{E,\zeta}^{\ast}$
denote its adjoint
$A^{0,1}(E\otimes \nbigl_{-\zeta})\lrarr 
 A^{0,0}(E\otimes \nbigl_{-\zeta})\oplus A^{0,2}(E\otimes \nbigl_{-\zeta})$.
By the results in \S\ref{subsection;13.1.24.1},
we obtain that 
$\nbigd_{\zeta}^{\ast}$ is surjective.
We have
$\Ker(\nbigd_{\zeta})=\Nahm(E,\nabla)_{\zeta}$.
We obtain that the family
$\bigcup_{\zeta}\Nahm(E,\nabla)_{\zeta}$
gives a $C^{\infty}$-bundle 
$\Nahm(E,\nabla)$ on $T^{\lor}\setminus D$.
It is equipped with a $C^{\infty}$-metric $h_1$
and a $C^{\infty}$-endomorphism $F_w$.
It is also equipped with the induced unitary connection
$\nabla_1$.
The $C^{\infty}$-bundle $\Nahm(E,\nabla)$
is also constructed as a family of the cohomology of
the complexes of the closed operators
$\bigl(
 A^{0,\bullet}(E\otimes \nbigl_{-\zeta}),\delbar_{E,\zeta}
 \bigr)$.
It induces a holomorphic structure of $\Nahm(E,\nabla)$
as a bundle on $T^{\lor}\setminus D$,
and $F_w$ is holomorphic.
We set $\theta_1:=F_w\,d\zetabar$.
The $(0,1)$-part of $\nabla_1$ is equal 
to the $\delbar$-operator of $\Nahm(E,\nabla)$.

\begin{prop}
\label{prop;13.1.29.10}
$(E_1,\delbar_{E_1},\theta_1,h_1)$
is a wild harmonic bundle.
\end{prop}
\pf
Because the argument is rather standard,
we give only an indication
for the convenience of the readers.
For $I\subset \{1,2,3\}$,
let $p_{I}$ denote the projection of
$T^{\lor}\times T\times \proj^1$
onto the product of the $i$-th components.
By the construction,
we have a natural isomorphism
$(E_1,\delbar_{E_1})
\simeq
 Rp_{1\ast}
 \bigl(
 p_{23}^{\ast}\nbigp_0E\otimes
 p_{12}^{\ast}\poincare^{-1}
 \bigr)_{|T^{\lor}\setminus D}$.
The endomorphism $F_w$
is equal 
the multiplication of 
$-w:
 Rp_{1\ast}
 \bigl(
 p_{23}^{\ast}\nbigp_{-1}E\otimes
 p_{12}^{\ast}\poincare^{-1}
 \bigr)_{|T^{\lor}\setminus D}
\lrarr
 Rp_{1\ast}
 \bigl(
 p_{23}^{\ast}\nbigp_{0}E\otimes
 p_{12}^{\ast}\poincare^{-1}
 \bigr)_{|T^{\lor}\setminus D}$.
Hence, we obtain that
$\theta$ is a wild Higgs field
in the sense that,
for the local expression $\theta=f\,d\zeta$
around $P\in D$,
the coefficients of 
the characteristic polynomial $\det(t\id-f)$
are meromorphic at $P$.

Let us prove that $(E_1,\delbar_{E_1},\theta_1,h_1)$
is a harmonic bundle.
Although we follow a standard argument,
we give rather details for the convenience of readers.
Let $\Delta_E$ denote the Laplacian on $A^{0,0}(E)$,
i.e., $\Delta_E=\delbar_E^{\ast}\delbar_E=
-\sqrt{-1}\Lambda\del_E\delbar_E$.
We have
\[
\Delta_E\psi
=-2\bigl(\nabla_z\nabla_{\zbar}+\nabla_w\nabla_{\wbar}\bigr)\psi.
\]
On $A^{0,2}(E)$,
the Laplacian is given by
$\delbar_E\delbar_E^{\ast}
=(-\sqrt{-1})\delbar_E\Lambda\del_E$.
We have
\[
 \delbar_E\delbar_E^{\ast}\bigl(\psi\,d\zbar\,d\wbar\bigr)
=-2(\nabla_{\zbar}\nabla_z+\nabla_{\wbar}\nabla_w)\psi\,d\zbar\,d\wbar.
\]
Because $F_{z\zbar}+F_{w\wbar}=0$,
it is equal to $\Delta_E(\psi)\,d\zbar\,d\wbar$.
Hence, under the natural identification
$A^{0,0}(E)\oplus A^{0,2}(E)
\simeq
 A^{0,0}(E)\otimes\langle\!\langle
 1,\,d\zbar\,d\wbar
 \rangle\!\rangle$,
the Laplacian $\nbigd^{\ast}\nbigd$
acts as $\Delta_E\otimes\id$,
where $\langle\!\langle a,b\rangle\!\rangle$
denotes the $2$-dimensional vector space
generated by $a,b$.
The Green operator of $\nbigd^{\ast}\nbigd$
acts as $\nbigg_E\otimes\id$,
where $\nbigg_E$ denotes the Green operator
for $\Delta_E$ on $A^{0,0}(E)$.

We naturally identify $A^{p,q}(E\otimes \nbigl_{-\zeta})$
with $A^{p,q}(E)$.
For a differential form $\tau$,
let $\mu(\tau)$ be an endomorphism of
$\bigoplus A^{p,q}(E)$ given by
$\mu(\tau)(\varphi)=\tau\wedge \varphi$.
We have
$\delbar_{E,\zeta}=\delbar_E-\zeta \mu(d\zbar)$
and
$\delbar_{E,\zeta}^{\ast}
=\delbar_E^{\ast}+(\sqrt{-1})\zetabar\Lambda\circ\mu(dz)$.
Let $d_{T^{\lor}}$ denote the trivial connection
of the product vector bundle
$A^{0,1}(E)\times (T^{\lor}\setminus D)$
over $T^{\lor}\setminus D$.
We have the following relation
for the operators
on the space of the sections
$T^{\lor}\setminus D\lrarr 
 A^{0,1}(E)\times (T^{\lor}\setminus D)$:
\[
 \bigl[
 d_{T^{\lor}},\,\delbar+\zeta\,d\zbar
 \bigr]
=d\zeta\,\mu(d\zbar),
\quad\quad
 \bigl[
 d_{T^{\lor}},\,(\delbar+\zeta\,d\zbar)^{\ast}
 \bigr]
=\sqrt{-1}d\zetabar\Lambda\circ\mu(dz).
\]
We set 
$\Omega:=d\zeta\,\mu(d\zbar)
 +d\zetabar\sqrt{-1}\Lambda\circ\mu(dz)$.
Let $P_{\zeta}$ denote the orthogonal projection
of $A^{0,1}(E)$ onto 
the kernel of $\nbigd^{\ast}_{\zeta}$.
Let $\Delta_{\zeta}=\delbar_{E,\zeta}^{\ast}\delbar_{E,\zeta}$
denote the Laplacian on $A^{0,0}(E)$
for $E\otimes \nbigl_{-\zeta}$.
Let $\nbigg_{\zeta}$ denote the Laplacian for $\Delta_{\zeta}$
on $A^{0,0}(E)$,
i.e.,
$\nbigg_{\zeta}\Delta_{\zeta}=\id_{A^{0,0}(E)}$.
The Green operator $G_{\zeta}$
for $\nbigd_{\zeta}^{\ast}\nbigd_{\zeta}$
on $A^{0,0}\oplus A^{0,2}$
is given by 
$\nbigg_{\zeta}\otimes\id$.
We have 
$P_{\zeta}=1-\nbigd_{\zeta}\circ
 G_{\zeta}\circ \nbigd_{\zeta}^{\ast}$.
Let $\nbigg_{\zeta}\otimes \id$
also denote the naturally induced operator 
on $A^{0,1}
\simeq 
A^{0,0}\otimes\langle\!\langle
 d\zbar,d\wbar
 \rangle\!\rangle$.

Let $\langle\cdot,\cdot\rangle$
denote the inner product of $A^{0,\bullet}(E)$
induced by $h$ and $dz\,d\zbar+dw\,d\wbar$.
By a standard computation,
the curvature $F$ of the connection $\nabla_1$
is described as follows,
for any sections $\psi_i$ $(i=1,2)$ of 
$\Nahm(E,\nabla)$:
\begin{multline}
 \bigl\langle
 \psi_1,F\psi_2
 \bigr\rangle
=\bigl\langle
 \psi_1,
 d_{T^{\lor}}\circ
 P_{\zeta}(d_{T^{\lor}}\psi_2)
 \bigr\rangle
=-\bigl\langle
 \psi_1,
 d_{T^{\lor}}\circ
 \nbigd_{\zeta}\circ G_{\zeta}\circ 
 \nbigd_{\zeta}^{\ast}
 (d_{T^{\lor}}\psi_2)
 \bigr\rangle
\\
=\bigl\langle
d_{T^{\lor}}\psi_1,
 \nbigd_{\zeta}\circ G_{\zeta}\circ 
 \nbigd_{\zeta}^{\ast}
 (d_{T^{\lor}}\psi_2)
 \bigr\rangle
=\bigl\langle
\nbigd_{\zeta}^{\ast}
d_{T^{\lor}}\psi_1,
G_{\zeta}\circ 
 \nbigd_{\zeta}^{\ast}
 (d_{T^{\lor}}\psi_2)
 \bigr\rangle
=\bigl\langle
 \Omega\psi_1,\,G_{\zeta}\Omega\psi_2
 \bigr\rangle
\\
=d\zeta\,d\zetabar
 \Bigl(
 \bigl\langle
 d\zbar\,\psi_1,\,d\zbar\,(\nbigg_{\zeta}\otimes\id)\psi_2
 \bigr\rangle
-\bigl\langle
 \Lambda (dz\,\psi_1),\,
 \Lambda\bigl(dz\,(\nbigg_{\zeta}\otimes\id)\psi_2\bigr)
 \bigr\rangle
 \Bigr)
\end{multline}
We have $\theta(\psi)=P_{\zeta}(w\psi)d\zeta$
and $\theta^{\dagger}(\psi)=P_{\zeta}(\wbar\psi)d\zetabar$.
We have
\begin{multline}
 \bigl\langle
 \psi_1,\,(P_{\zeta}w\circ P_{\zeta}\wbar
 -P_{\zeta}\wbar\circ P_{\zeta}w)\psi_2
 \bigr\rangle\,d\zetabar\,d\zeta
=-\bigl\langle
 \psi_1,\,
 \bigl(
 w(P_{\zeta}-1)\wbar-\wbar(P_{\zeta}-1)w
 \bigr)\psi_2
 \bigr\rangle d\zeta\,d\zetabar
 \\
=\Bigl(
 \bigl\langle
 \wbar\psi_1,
 \nbigd_{\zeta}
 G_{\zeta}
 \nbigd_{\zeta}^{\ast}\wbar\psi_2
 \bigr\rangle
- \bigl\langle
 w\psi_1,\nbigd_{\zeta} G_{\zeta}\nbigd_{\zeta}^{\ast}w\psi_2
 \bigr\rangle
 \Bigr)\,d\zeta\,d\zetabar
 \\
=\Bigl(
 \bigl\langle
 \nbigd_{\zeta}^{\ast}(\wbar\psi_1),
 G_{\zeta}\nbigd_{\zeta}^{\ast}(\wbar\psi_2)
 \bigr\rangle
-\bigl\langle
 \nbigd_{\zeta}^{\ast}(w\psi_1),
 G_{\zeta}\nbigd_{\zeta}^{\ast}(w\psi_2)
 \bigr\rangle
 \Bigr)
 d\zeta\,d\zetabar
 \\
=\Bigl(
 \bigl\langle
 [\nbigd_{\zeta}^{\ast}, \wbar]\psi_1,
 G_{\zeta}[\nbigd_{\zeta}^{\ast},\wbar]\psi_2
 \bigr\rangle
-\bigl\langle
 [\nbigd_{\zeta}^{\ast},w]\psi_1,
 G_{\zeta}[\nbigd_{\zeta}^{\ast},w]\psi_2
 \bigr\rangle
 \Bigr)
 d\zeta\,d\zetabar
\end{multline}
We have 
$[\nbigd_{\zeta}^{\ast},\wbar]=\mu(d\wbar)$
and
$[\nbigd_{\zeta}^{\ast},w]=-\sqrt{-1}\Lambda\circ\mu(dw)$.
Hence, 
we obtain the following expression:
\begin{equation}
\bigl\langle
 \psi_1,\,(P_{\zeta}w\circ P_{\zeta}\wbar
 -P_{\zeta}\wbar\circ P_{\zeta}w)\psi_2
 d\zeta\,d\zetabar
 \bigr\rangle
=\Bigl(
 \bigl\langle
 d\wbar\psi_1,\,d\wbar(\nbigg_{\zeta}\otimes\id)\psi_2
 \bigr\rangle
-\bigl\langle
 \Lambda (dw\psi_1),\,
 \Lambda (dw(\nbigg_{\zeta}\otimes\id)\psi_2)
 \bigr\rangle
 \Bigr)
 d\zeta\,d\zetabar
\end{equation}
By using
$\langle
 d\zbar d\wbar,d\zbar d\wbar\rangle
=\langle
 \Lambda dw d\wbar,\,
 \Lambda dw d\wbar\rangle
=\langle
 \Lambda dz d\zbar,\,
 \Lambda dz d\zbar\rangle$
for the metric on
$\Omega^{\bullet}_{T\times\cnum}$,
we obtain the following:
\[
 \Bigl\langle
 \psi_1,\,
 \bigl(
 F+(P\circ w\,d\zeta)\circ (P\circ \wbar d\zetabar)
 \bigr)\psi_2
 \Bigr\rangle=0
\]
Namely, the Hitchin equation is satisfied.
Thus, the proof of Proposition \ref{prop;13.1.29.10}
is finished.
\hfill\qed

\begin{rem}
We obtain a different transformation
by replacing $\nbigl_{-\zeta}$
with $\nbigl_{\zeta}$,
for which we do not need any essential change.
\end{rem}

\begin{rem}
We use the operators,
which looks natural in the complex geometry,
instead of the Dirac operator itself.
\end{rem}

\section{$L^2$-instantons and wild harmonic bundles}

\subsection{Nahm transform for wild harmonic bundles on $T^{\lor}$}
\label{subsection;13.1.27.10}

\subsubsection{Construction}
\label{subsection;13.2.14.1}

Let $D$ be a non-empty finite subset of $T^{\lor}$.
We fix a K\"ahler metric $g_{T^{\lor}\setminus D}$ of 
$T^{\lor}\setminus D$,
which is Poincar\'e like around $D$.
Let $\harmonicbundle$ be a wild harmonic bundle
on $(T^{\lor},D)$.
We assume that $\harmonicbundle$
has a singularity at each point of $D$,
i.e., $\theta$ has a pole, or the parabolic structure is non-trivial.
We shall construct an $L^2$-instanton
from $\harmonicbundle$
with the method in \cite{Donaldson-Kronheimer}
and \cite{Jardim1}.
Let $H^i_{L^2}\harmonicbundle$ 
denote the $i$-th $L^2$-cohomology group of 
$\harmonicbundle$.
As recalled in Lemma \ref{lem;13.1.26.10},
they are isomorphic to the hypercohomology groups
of the complex
$\nbigc^{\bullet}\bigl(
 \nbigp_{\ast}E\otimes\Omega^{\bullet},\theta\bigr)$,
where $(\nbigp_{\ast}E,\theta)$
is the associated good filtered Higgs bundle.
In particular,
they are finite dimensional,
and isomorphic to
the space of $L^2$-harmonic $i$-forms
of $\harmonicbundle$.
We have $H^0_{L^2}\harmonicbundle=
H^2_{L^2}\harmonicbundle=0$
unless $\harmonicbundle$
is $\nbigo_{T^{\lor}}$ with 
the trivial metric and the trivial Higgs field.

\begin{rem}
If $D$ is empty,
$(E,\delbar_E,\theta,h)$
is isomorphic to a direct sum
$\bigoplus (L_i,\delbar_{L_i},\theta_i,h_i)$
such that $\rank L_i=1$.
So we exclude the case $D=\emptyset$.
\hfill\qed
\end{rem}

For any $(z,w)$,
let $\nbigl_{z,w}$ denote the harmonic bundle of rank one
given by $(\underline{\cnum},\delbar+zd\zetabar)$
with trivial metric and the Higgs field $w\,d\zeta$.
Let 
$(E,\delbar_{E,z},\theta_w,h)$
denote
$\harmonicbundle\otimes \nbigl_{z,w}$.
Let $\Nahm\harmonicbundle_{(z,w)}$
be the space of $L^2$-harmonic $1$-forms
of $(E,\delbar_{E,z},\theta_w,h)$.
It is independent of the choice of 
the Poincar\'e like metric $g_{T^{\lor}\setminus D}$.
It is finite dimensional, and naturally isomorphic to
$\Nahm(\nbigp_{\ast}E,\theta)_{(z,w)}$.
It is naturally equipped with the metric $h_1$
induced by $h$.

Let $A^{p,q}(E)$
denote the space of $L^2$-sections of 
$E\otimes\Omega^{p,q}_{T^{\lor}\setminus D}$.
Let $\delbar_{E,z}^{\ast}:A^{p,q}\lrarr A^{p,q-1}$ denote 
the adjoint of the closed operator
$\delbar_{E,z}:A^{p,q}\lrarr A^{p,q+1}$.
Let $\theta_w^{\dagger}:A^{p,q}\lrarr A^{p-1,q}$ denote
the adjoint of 
$\theta_w:A^{p,q}\lrarr A^{p+1,q}$.
We have
$\delbar_{E,z}^{\ast}:=
 -\sqrt{-1}\bigl[\Lambda,\del_{E}-\zbar d\zeta\bigr]$
and 
$\theta_w^{\ast}=
 -\sqrt{-1}\bigl[\Lambda,\theta_w^{\dagger}\bigr]
=
-\sqrt{-1}\bigl[
 \Lambda,\theta^{\dagger}+\wbar\,d\zetabar
 \bigr]$.

We set $S^+:=A^{0,0}(E)
 \oplus A^{1,1}(E)$
and $S^-:=A^1(E)=A^{0,1}(E)\oplus A^{1,0}(E)$.
Let $\nbigd_{z,w}:=
 \delbar_{E,z}+\theta_w
+\delbar_{E,z}^{\ast}+\theta_w^{\ast}$
be the closed operator
$S^+\lrarr S^{-}$,
and let 
$\nbigd_{z,w}^{\ast}:=
 \delbar_{E,z}+\theta_w
+\delbar_{E,z}^{\ast}+\theta_w^{\ast}$
denote its adjoint
$S^-\lrarr S^+$.
We have
$\Ker\nbigd_{z,w}^{\ast}
=\Nahm\harmonicbundle_{(z,w)}$.
(See \cite{Mochizuki-wild}.)
By the vanishing 
$H^i_{L^2}(E,\delbar_{E,z},\theta_w,h)$ $(i=0,2)$,
we obtain that $\nbigd^{\ast}$ is surjective.
Hence, the family 
$\bigcup_{(z,w)}
 \Nahm\harmonicbundle_{(z,w)}$
gives a $C^{\infty}$-vector bundle
on $T\times\cnum$.
(See \cite{Donaldson-Kronheimer}.)
It is equipped with an induced $C^{\infty}$-metric $h_1$
and an induced unitary connection $\nabla_1$.
Because the $C^{\infty}$-bundle is also constructed
as a family of the cohomology of the complexes
$\bigl(
 A^{\bullet}(E),
 \delbar_{E,z}+\theta_w\,d\zeta
 \bigr)$,
it is equipped with a naturally induced holomorphic structure, 
which is equal to the $(0,1)$-part of $\nabla_1$. 
By the construction, 
the holomorphic bundle 
is naturally isomorphic to 
$\Nahm(\nbigp_{\ast}E,\theta)_{|T\times\cnum}$.
(See \S\ref{subsection;13.3.5.10} for more details
on this isomorphism.)
We shall give the proof of the following theorem
in \S\ref{subsection;13.2.26.20}
after preliminaries.

\begin{thm}
\label{thm;13.1.25.10}
$\bigl(
\Nahm\harmonicbundle,h_1,\nabla_1\bigr)$
is an $L^2$-instanton.
\end{thm}

We give a remark on the proof.
It is rather easy and standard to prove that
$\bigl(
 \Nahm(E,\delbar_E,\theta,h),h_1,\nabla_1
\bigr)$
is an instanton
by using the twistor property of instantons and harmonic bundles.
But, we do not give such an argument in the following.
Instead, we follow another standard argument to use a description of 
the curvature $F(\nabla_1)$ in terms of the Green operator.
Because we need an estimate for the decay of $F(\nabla_1)$,
we need the description, anyway.

\subsubsection{Preliminary}

Let $X$ be a torus $\cnum_{\zeta}/L$.
Let $D\subset X$ be a finite set.
Let $g_{A}=Ad\zeta\,d\zetabar$ be a K\"ahler metric
of $X\setminus D$ 
for some positive valued function $A$,
which is Poincar\'e like around $D$.
Let $\harmonicbundle$ be a wild harmonic bundle
on $X\setminus D$.
We set $\nbigd:=\delbar_E+\theta$.
Let $\nbigd_A^{\ast}$ (resp $\nbigd_1^{\ast}$)
denote the formal adjoint of
$\nbigd$ with respect to $h$ and $g_{A}$
(resp. $d\zeta\,d\zetabar$).
We set
$\Delta_A=\nbigd_A^{\ast}\nbigd$
and 
$\Delta_1=\nbigd_1^{\ast}\nbigd$.
We have $\Delta_A=A^{-1}\Delta_1$.

\begin{lem}
\label{lem;13.1.25.120}
Let $\varphi$ be a section of $E$ on $X\setminus D$
such that 
\[
\int |\varphi|_h^2A|d\zeta\,d\zetabar|
+\int|\Delta_1\varphi|_h^2|d\zeta\,d\zetabar|<\infty.
\]
Then, we have the following finiteness:
\begin{equation}
\label{eq;13.1.25.100}
 \int|\varphi|_h^2 |d\zeta\,d\zetabar|
+\int|\Delta_A\varphi|_h^2
 \,A|d\zeta\,d\zetabar|<\infty
\end{equation}
\begin{equation}
\label{eq;13.1.25.101}
\int h(\varphi,\Delta_1\varphi)
 |d\zeta\,d\zetabar|
=\int h(\varphi,\Delta_A\varphi)
 A|d\zeta\,d\zetabar|
=\int |\nbigd \varphi|_h^2
<\infty
\end{equation}
\end{lem}
\pf
The finiteness (\ref{eq;13.1.25.100}) is clear.
In (\ref{eq;13.1.25.101}),
the first equality is trivial.
The second equality and finiteness
can be shown by an argument
in the proof of Lemma \ref{lem;13.1.10.10}.
\hfill\qed

\vspace{.1in}
We set
$\nbigd^{\dagger}:=\del_E+\theta^{\dagger}$.
Let $(\nbigd^{\dagger})^{\ast}_A$ 
(resp. $(\nbigd^{\dagger})^{\ast}_1$)
denote the formal adjoint of $\nbigd$
with respect to $g_A$
(resp. $d\zeta d\zetabar$).
We have
$\Delta_1=(\nbigd^{\dagger})_1^{\ast}\nbigd^{\dagger}$
and 
$\Delta_A=(\nbigd^{\dagger})_A^{\ast}\nbigd^{\dagger}$.
\begin{lem}
\label{lem;13.1.26.20}
Let $\varphi$ be as in Lemma {\rm\ref{lem;13.1.25.120}}.
Then, we have
\begin{equation}
\int h(\varphi,\Delta_1\varphi)
 |d\zeta\,d\zetabar|
=\int h(\varphi,\Delta_A\varphi)
 A|d\zeta\,d\zetabar|
=\int |\nbigd^{\dagger} \varphi|_h^2
<\infty.
\end{equation}
\end{lem}
\pf
The first equality is trivial.
For the second,
we have only to apply Lemma \ref{lem;13.1.25.120}
to a harmonic bundle
$(E,\del_E,\theta^{\dagger},h)$ on $X\setminus D$.
\hfill\qed

\subsubsection{Estimate}
\label{subsection;13.2.28.1}

Let $X$ be a torus $\cnum_{\zeta}/L$
with a non-empty finite subset $D$.
We use the Euclidean metric
$d\zeta\,d\zetabar$ of $X$.
Let $\dvol_X=|d\zeta\,d\zetabar|$ 
denote the associated volume form.
Let $\harmonicbundle$ be a wild harmonic bundle
on $(X,D)$.
Assume that the harmonic bundle
has a singularity at each point of $D$.

Let $\nabla_{h}=\delbar_E+\del_E$
be the Chern connection.
Let $\nbigh_{w}$ be the space of the sections of $E$
on $X\setminus D$
such that
\[
\int_{X}
 \bigl|\varphi\bigr|_h^2
 \dvol_{X}
+\int_{X}\Bigl(
\bigl|\nabla_{h}\varphi\bigr|_h^2
 +\bigl|
 (\theta+w\,d\zeta)\varphi
 \bigr|_h^2
 \Bigr)
<\infty.
\]

\begin{prop}
\label{prop;13.2.27.20}
There exist positive constants
$R>0$, $C>0$ and $\rho>0$ such that,
if $|w|>R$,
the following holds for any 
$\varphi\in\nbigh_{w}$:
\[
  \int_{X}\Bigl(
 \bigl|\nabla_h\varphi\bigr|_h^2
 +\bigl|
 (\theta+w\,d\zeta)\varphi
 \bigr|^2_h
 \Bigr)
\geq
 C\,|w|^{\rho}\int_{X}
 |\varphi|_h^2\dvol_{X}
\]
(See also a refined estimate in Proposition
{\rm\ref{prop;12.7.19.4}} below.)
\end{prop}
\pf
We use an argument in \S2.4 of \cite{Szabo}
with an adjustment to our situation. 
We use the standard distance on $X$.
We take small neighbourhoods $B_P$ of $P\in D$.
There exists $R_1>0$ and $C_1>0$
such that, if $|w|\geq R_1$, then
we have
$\bigl|
 (\theta+w\,d\zeta)\varphi
 \bigr|_h^2
\geq
 C_1\,|w|^2\,|\varphi|_h^2\,\dvol_X$
on $X\setminus
 \bigcup_{P\in D}B_P$.
We have only to show the estimate on
each $B_P$.
We may assume $P=0$,
and $B_P$ is an $\epsilon$-ball
$B_{\epsilon}=\{|\zeta|\leq \epsilon\}$.

We have
a ramified covering
$\psi:(B'_{\epsilon},0)\lrarr (B_{\epsilon},0)$
given by $\psi(u)=u^p$
such that
$\psi^{\ast}(E,\delbar_E,\theta,h)$
is unramified,
i.e.,
we have the decomposition
\[
 \psi^{\ast}(E,\delbar_E,\theta)
=\bigoplus_{\gminia\in u^{-1}\cnum[u^{-1}]}
 (E_{\gminia},\delbar_{E_{\gminia}},\theta_{\gminia}),
\]
where 
the Higgs field
$\theta_{\gminia}-d\gminia\,\id_{E_{\gminia}}$
are tame.
Let $\ell:=\max\bigl\{
 \deg_{u^{-1}}\gminia\,\big|\,
 E_{\gminia}\neq 0 \bigr\}$.

\begin{lem}
\label{lem;13.2.27.10}
There exists $R'>0$, $C'_i>0$ $(i=1,2)$
such that
$|\theta\varphi|_h\geq 
C_1'|w|\,|d\zeta|\,|\varphi|_h$
on 
$B_{\epsilon}\setminus
 \bigl\{|\zeta|<C_2'|w|^{-p/(\ell+p)}\bigr\}$,
if $|w|\geq R'$.
\end{lem}
\pf
We have only to estimate each $\theta_{\gminia}$
on $B_{\epsilon}'$.
Let us consider the case
$\gminia\neq 0$.
We set $n:=\deg_{u^{-1}}\gminia$.
For each $w$,
we have the solutions $b_i(w)$ $(i=0,\ldots,n+p-1)$
of the following equation:
\[
 \del_u\gminia(u)+pw u^{p-1}=0
\]
We have the equality
$u^{-p+1}\del_u\gminia(u)+pw
=\alpha\prod_{i=0}^{n+p-1}\bigl(u^{-1}-b_i(w)^{-1}\bigr)$
for some $\alpha\in\cnum\setminus\{0\}$.
We have
\[
 \theta_{\gminia}
=\del_u\gminia\,\id_{E_{\gminia}}\,du
+g_{\gminia}\,du,
\]
where $|g_{\gminia}|_{h}\leq C_1|u|^{-1}$.
We have $R_{2}>0$ and $C_{2}>0$ such that 
the following holds if $|w|>R_2$:
\[
 C_{2}^{-1}\leq 
 |b_i(w)|\,|w|^{1/(n+p)}\leq C_{2}
\]
We take $C_{3}>\!>C_{2}$.
We set
$\nbigw_1:=\bigl\{
 |u|\leq C_{3}|w|^{-1/(n+p)}
 \bigr\}$.

On $B_{\epsilon}'\setminus\nbigw_1$,
we have
$|g_{\gminia}|_{h}\leq (C_1/C_3)|w|^{1/(n+p)}$.
We also have
\[
\bigl| u^{-1}-b_i(w)^{-1} \bigr|
\geq
 |b_i(w)^{-1}|-|u^{-1}|
\geq (C_2^{-1}-C_3^{-1})
|w|^{1/(n+p)}
\]
for any $i$,
and hence
$\Bigl|
 u^{-p+1}\del_{u}\gminia+pw
 \Bigr|
\geq 
|\alpha|(C_2^{-1}-C_3^{-1})^{n+p}|w|$.
Hence, if $C_3$ is sufficiently larger than $C_2$,
there exist $R_4>0$ and $C_4>0$
such that the following holds 
if $|w|>R_4$:
\[
 \bigl|
 \bigl(
 \del_u\gminia+pw u^{p-1}
 \bigr)\id_{E_{\gminia}}+g_{\gminia}
 \bigr|_{h}
\geq
C_{4}|w|
 |u|^{p-1}
\]
Hence, we obtain the desired inequality for the integral over
$B_{\epsilon}'\setminus \nbigw_1$
in the case $\gminia\neq 0$.

Let us consider the case $\gminia=0$.
We have the expression $\theta_0=g_0\,du$,
and $|g_0|_h\leq C_{10}|u|^{-1}$ for some $C_{10}>0$.
We take $C_{11}>C_{10}$,
and we consider
$\nbigw:=\{|u|\leq C_{11}|w|^{-1/p}\}$.
On $B_{\epsilon}'\setminus \nbigw$,
we have
$|wu^{p-1}|\geq C_{11}^{p-1}|w|^{1/p}$.
We also have
$|g_0|_h\leq (C_{10}/C_{11})||w|^{1/p}$.
Hence, if $C_{11}$ is sufficiently larger than $C_{10}$,
we have the following for some $C_{12}>0$:
\[
 \bigl|pwu^{p-1}\id_{E_0}+g_0\bigr|_h
\geq
 C_{12}\bigl|wu^{p-1}\bigr|
\]
Hence, we obtain the desired inequality in the case
$\gminia=0$.
Thus, the proof of Lemma \ref{lem;13.2.27.10} is finished.
\hfill\qed

\vspace{.1in}

Let $\varphi$ be an $L^2$-section of $E$
on $B_{\epsilon}$ with respect to $\dvol_X$,
such that
\[
 \int_{B_{\epsilon}}
 \Bigl(
 \bigl|\nabla_{h}\varphi\bigr|_{h}^2
+\bigl|(\theta+w d\zeta)
 \,\varphi\bigr|_{h}^2
 \Bigr)
 \dvol_X
<\infty.
\]
We set
$\nbigw_1:=\{|\zeta|<2C_2'|w|^{-p/(\ell+p)}\}$
and 
$\nbigw_2:=\{|\zeta|<C_2'|w|^{-p/(\ell+p)}\}$.
We have the following type of Poincar\'e inequality,
i.e., there exist $C''>0$ and $R''>0$ 
such that the following holds if $|w|>R''$
(see \cite{biquard-boalch}, and (2.12) of \cite{Szabo}):
\[
 |w|^{2p/(\ell+p)}
 \int_{\nbigw_1}|\varphi|_{h}^2
 |d\zeta d\zetabar|
\leq
 C''\left(
 \int_{\nbigw_1}\bigl|
 d|\varphi|_{h}
 \bigr|^2
+|w|^{2p/(n+p)}
 \int_{\nbigw_1\setminus\nbigw_2}
 |\varphi|_{h}^2\,|d\zeta\,d\zetabar|
 \right)
\]
There exists $C'''$ such that 
the right hand side is dominated by
\[
 C'''\left(
 \int_{\nbigw_1}
 \bigl|\nabla_h\varphi\bigr|_h^2
+\int_{\nbigw_1\setminus\nbigw_2}
\bigl|(\theta+wd\zeta)\varphi\bigr|_h^2
 \right)
\]
Thus, the proof of Proposition \ref{prop;13.2.27.20}
is finished.
\hfill\qed

\vspace{.1in}

Let $\nbigd:=\delbar_E+\theta$.
Let $\nbigd^{\ast}_1$ 
denote the adjoint with respect to 
the Euclidean metric $d\zeta\,d\zetabar$.
Let $\Delta_1:=\nbigd_1^{\ast}\circ\nbigd$.
Let $g_{X\setminus D}$
be a K\"ahler metric of $X\setminus D$
which is Poincar\'e like around $D$.
Let $\dvol_{X\setminus D}$
be the volume form associated to $g_{X\setminus D}$.

\begin{cor}
\label{cor;13.2.27.22}
There exist $\rho>0$ and $C>0$ such that the following holds:
\begin{itemize}
\item
Let $\varphi$ be a section of $E$
such that 
\begin{equation}
 \label{eq;13.2.27.30}
 \int |\varphi|_h^2\dvol_{X\setminus D}
+\int |\Delta_{1}\varphi|_h^2
 \dvol_{X}
<\infty.
\end{equation}
Then, we have the following inequality:
\begin{equation}
 \label{eq;13.2.27.31}
  C|w|^{\rho}
 \Bigl(
 \int |\varphi|_h^2\dvol_{X}
 \Bigr)^{1/2}
\leq
 \Bigl(
\int|\Delta_{1}\varphi|_{h}^2\dvol_{X}
 \Bigr)^{1/2}
\end{equation}
\end{itemize}
(See Corollary {\rm\ref{cor;12.7.20.3}}
for a refinement.)
\end{cor}
\pf
Let $\nbigd^{\dagger}=\del_E+\theta^{\dagger}$.
From (\ref{eq;13.2.27.30}),
Lemma \ref{lem;13.1.25.120}
and Lemma \ref{lem;13.1.26.20},
we obtain 
$\int |\nbigd \varphi|_h^2<\infty$
and $\int|\nbigd^{\dagger}\varphi|_h^2<\infty$.
By using the same lemmas and
Proposition \ref{prop;13.2.27.20},
we obtain 
\[
 C|w|^{\rho}
 \int |\varphi|_h^2\dvol_{X}
\leq
 \int|\nbigd\varphi|_h^2
+\int|\nbigd^{\dagger}\varphi|_h^2
=2\int
 h\bigl(
 \varphi,\Delta_{1}\varphi
 \bigr)\dvol_X.
\]
Then, the claim of the corollary follows.
\hfill\qed

\subsubsection{Proof of Theorem \ref{thm;13.1.25.10}}
\label{subsection;13.2.26.20}

Let $\omega_{T^{\lor}\setminus D}$ be the K\"ahler form
associated to the metric $g_{T^{\lor}\setminus D}$.
The multiplication of $\omega_{T^{\lor}\setminus D}$
induces an isomorphism $A^{0,0}(E)\simeq A^{1,1}(E)$.
It gives an identification
$S^+\simeq
 A^{0,0}(E)\otimes\langle\!\langle
1,\omega_{T^{\lor}\setminus D}\rangle\!\rangle$,
where 
$\langle\!\langle 1,
 \omega_{T^{\lor}\setminus D}\rangle\!\rangle$
denotes the $2$-dimensional vector space
generated by $1$ and $\omega_{T^{\lor}\setminus D}$.
By a general theory of harmonic bundles,
the Laplacian $\nbigd_{zw}^{\ast}\nbigd_{zw}$ on $S^+$
is identified with 
$\Delta_{zw}\otimes\id$
on $A^{0,0}(E)\otimes
 \langle\!\langle 1,\omega_{T^{\lor}\setminus D}
 \rangle\!\rangle$,
where $\Delta_{zw}:=(\delbar_{E,z}^{\ast}+\theta_w^{\ast})\circ
 (\delbar_{E,z}+\theta_w)$ on $A^{0,0}(E)$.
(See \cite{s5}. 
In this case, it can be easily checked directly.)
The Green operator $G_{zw}$
for $\nbigd_{zw}^{\ast}\nbigd_{zw}$
is identified with $\nbigg_{zw}\otimes\id$,
where $\nbigg_{zw}$ is the Green operator
of $\Delta_{zw}$ on $A^{0,0}(E)$.

\vspace{.1in}

For a differential form $\tau$ on $T^{\lor}$,
let $\mu(\tau)$ be an endomorphism of
$\bigoplus A^{p,q}(E)$ given by
$\mu(\tau)(\varphi)=\tau\wedge \varphi$.
Let $d_{T\times\cnum}$ denote the trivial connection
of the product vector bundle
$S^{-}\times (T\times\cnum)$ over $T\times\cnum$.
We have the following relation
 for the operators on the space of the sections
$T\times\cnum\lrarr
 S^-\times(T\times\cnum)$:
\[
 \bigl[d_{T\times\cnum},\,\delbar+z\,d\zetabar\bigr]
=dz\,\mu(d\zetabar),
\quad
 \bigl[
 d_{T\times\cnum},\,\theta+w\,d\zeta
 \bigr]
=dw\,\mu(d\zeta),
\]
\[
 \bigl[d_{T\times\cnum},(\delbar+z\,d\zetabar)^{\ast}\bigr]
=d\zbar\,\bigl(\sqrt{-1}\Lambda\circ\mu(d\zeta)\bigr),
\quad
 \bigl[d_{T\times\cnum},\,(\theta+w\,d\zeta)^{\ast}\bigr]
=d\wbar\bigl(-\sqrt{-1}\Lambda\,\mu(d\zetabar)\bigr)
\]
We set 
$\Omega:=
 dz\,\mu(d\zetabar)+dw\,\mu(d\zeta)
+d\zbar\bigl(\sqrt{-1}\Lambda\,\mu(d\zeta)\bigr)
+d\wbar\bigl(-\sqrt{-1}\Lambda\, \mu(d\zetabar)\bigr)$.

Let $F(\nabla_1)$ be the curvature of the transformed bundle
$\Nahm\harmonicbundle$
with the metric and the unitary connection.
Let $P_{zw}$ denote the orthogonal projection of
$S^-$ onto $\Nahm\harmonicbundle_{(z,w)}$.
Let $\psi_i$ be sections of $\Nahm\harmonicbundle$.
Let $\langle\cdot,\cdot\rangle$
denote the hermitian pairing on
$A^{p,q}(E)$ induced by $h$ and $\omega_{T^{\lor}\setminus D}$.
We have the following standard computation:
\begin{multline}
\label{eq;12.7.20.1}
 \bigl\langle
 \psi_1,F(\nabla_1)\psi_2
 \bigr\rangle
=\bigl\langle
 \psi_1,P_{zw}\circ d\circ P_{zw}d\psi_2
 \bigr\rangle
=\bigl\langle
 \psi_1,d\circ P_{zw}\circ d\psi_2
 \bigr\rangle
=\bigl\langle
 \psi_1,d\circ(P_{zw}-1)\circ d\psi_2
 \bigr\rangle
 \\
=-\bigl\langle
 d\psi_1,(P_{zw}-1)\circ d\psi_2
 \bigr\rangle
=\bigl\langle
 d\psi_1,\nbigd_{zw}\circ G_{zw}\circ \nbigd_{zw}^{\ast}d\psi_2
 \bigr\rangle
=\bigl\langle
 \nbigd_{zw}^{\ast}d\psi_1,\,
 G_{z,w}\nbigd_{zw}^{\ast}d\psi_2
 \bigr\rangle
 \\
=\Bigl\langle
 [d,\nbigd_{zw}^{\ast}]\psi_1,\,
 G_{zw}[d,\nbigd_{zw}^{\ast}]\psi_2
 \Bigr\rangle
=\bigl\langle
 \Omega\psi_1,\, G_{zw}\Omega\psi_2
 \bigr\rangle
\end{multline}

We have the expression
$\psi_i=\psi_{i1}d\zeta+\psi_{i2}d\zetabar$.
We have
\[
 \Omega\psi_1
=dz\,\psi_{11}d\zetabar d\zeta
+dw\, \psi_{12}d\zeta d\zetabar
-\sqrt{-1}d\wbar\, \psi_{11}\Lambda(d\zetabar d\zeta)
+\sqrt{-1}d\zbar\, \psi_{12}\Lambda(d\zetabar d\zeta).
\]
Let $A$ be determined by
$g_{T^{\lor}\setminus D}=A\,d\zeta d\zetabar$.
We have the following:
\begin{multline}
 G_{zw}\Omega\psi_2
=dz\,\nbigg_{zw}(A^{-1}\psi_{21})Ad\zetabar d\zeta
+dw\, \nbigg_{zw}(A^{-1}\psi_{22})Ad\zeta d\zetabar
 \\
-\sqrt{-1}
 d\wbar\, \nbigg_{zw}\bigl(\psi_{21}\Lambda(d\zetabar d\zeta)\bigr)
+\sqrt{-1}
 d\zbar\, \nbigg_{zw}\bigl(\psi_{22}\Lambda(d\zetabar d\zeta)\bigr).
\end{multline}
We have the following:
\begin{equation}
 \bigl\langle
 \psi_{11}d\zetabar d\zeta,\,
 A\nbigg_{zw}(A^{-1}\psi_{21})d\zetabar d\zeta
 \bigr\rangle
=\bigl\langle
 \psi_{11}\Lambda(d\zetabar d\zeta),\,
 \nbigg_{zw}\bigl(\psi_{21}\,\Lambda(d\zetabar d\zeta)\bigr)
 \bigr\rangle
\\
=4\int
 \bigl(\psi_{11},\,\nbigg_{zw}(A^{-1}\psi_{21})\bigr)
 \dvol_{T^{\lor}}
\end{equation}
\begin{equation}
 \bigl\langle
 \psi_{12}d\zetabar d\zeta,\,
 A\nbigg_{zw}(A^{-1}\psi_{22})d\zetabar d\zeta
 \bigr\rangle
=\bigl\langle
 \psi_{12}\Lambda(d\zetabar d\zeta),\,
 \nbigg_{zw}\bigl(\psi_{22}\,\Lambda(d\zetabar d\zeta)\bigr)
 \bigr\rangle
\\
=4\int
 \bigl(\psi_{12},\,\nbigg_{zw}(A^{-1}\psi_{22})\bigr)
 \dvol_{T^{\lor}}
\end{equation}
From these equalities,
we obtain
$ (dz\,d\zbar+dw\,d\wbar)\wedge
 \langle \psi_1,F(\nabla_1)\psi_2\rangle=0$,
which means that
$\Nahm\harmonicbundle$
with the induced metric $h_1$ and connection $\nabla_1$
is an instanton.

\vspace{.1in}
Let us show that it is an $L^2$-instanton.
Let $(\delbar_{E,z}+\theta_w)_1^{\ast}$
denote the formal adjoint of
$\delbar_{E,z}+\theta_w$
with respect to $h$ and $d\zeta\,d\zetabar$.
We set
$\Delta_{zw,1}:=
 (\delbar_{E,z}+\theta_w)_1^{\ast}
 (\delbar_{E,z}+\theta_w)$.
Because $\Delta_{zw,1}=A\Delta_{zw}$,
we have
$\Delta_{zw,1}\bigl(
 \nbigg_{zw}(A^{-1}\psi_{21})
 \bigr)=\psi_{21}$.
We have
\[
 \int \bigl|\nbigg_{zw}(A^{-1}\psi_{21})\bigr|_h^2
 \dvol_{T^{\lor}\setminus D}
+\int \bigl|
 \psi_{21}
 \bigr|_h^2\dvol_{T^{\lor}}<\infty
\]
By Corollary \ref{cor;13.2.27.22},
we have the following for some $\rho>0$ and $C>0$:
\[
 C|w|^{2\rho}
 \int_{T^{\lor}}
 \bigl|
 \nbigg_{zw}(A^{-1}\psi_{21})
 \bigr|_h^2\dvol_{T^{\lor}}
<
 \int_{T^{\lor}}
 \bigl|
 \psi_{21}
 \bigr|^2\dvol_{T^{\lor}}
\]
Hence, we obtain
\begin{multline}
 \label{eq;13.1.26.30}
\left|
 \bigl\langle
 \psi_{11}d\zetabar d\zeta,\,
 A\nbigg_{zw}(A^{-1}\psi_{21})d\zetabar d\zeta
 \bigr\rangle
\right|
=\left|
 \bigl\langle
 \psi_{11}\Lambda(d\zetabar d\zeta),\,
 \nbigg_{zw}\bigl(\psi_{21}\,\Lambda(d\zetabar d\zeta)\bigr)
 \bigr\rangle
 \right|
 \\
<C|w|^{-\rho}
\left(
\int \bigl|\psi_{11}d\zeta\bigr|^2_h
\right)^{1/2}
\left(
\int \bigl|\psi_{21}d\zeta\bigr|^2_h
\right)^{1/2}
\end{multline}
We have a similar estimate for
$\bigl|
 \langle
 \psi_{12}d\zetabar d\zeta,\,
 \nbigg_{zw}(\psi_{22})d\zetabar d\zeta
 \rangle
 \bigr|
=|\bigl\langle
 \psi_{12}\Lambda(d\zetabar d\zeta),\,
 \nbigg_{zw}\bigl(\psi_{22}\,\Lambda(d\zetabar d\zeta)\bigr)
 \bigr\rangle|$.
From those estimates,
we obtain
$|F(\nabla_1)|=O\bigl(|w|^{-\rho}\bigr)$
for some $\rho>0$.
Because
$\Nahm\harmonicbundle
\simeq
 \Nahm(\nbigp_{\ast}E,\theta)_{|T\times\cnum}$,
we can apply Theorem \ref{thm;13.2.26.3},
and hence we obtain that $F(\nabla_1)$ is $L^2$.
Thus, the proof of Theorem \ref{thm;13.1.25.10}
is finished.
\hfill\qed

\begin{rem}
We can directly prove that the curvature is $L^2$ 
by using Corollary {\rm\ref{cor;12.7.20.3}} below.
\hfill\qed
\end{rem}

\subsubsection{Refined estimates (Appendix)}

We refine the estimates in 
\S\ref{subsection;13.2.28.1},
i.e.,
we show that $\rho$ can be replaced with $1+\rho$.
Although we do not use it in this paper,
this type of argument
seems useful in the study of a different type of
Nahm transform,
and so we would like to keep it.

\begin{prop}
\label{prop;12.7.19.4}
There exist positive constants
$R>0$, $C>0$ and $\rho>0$ such that,
if $|w|>R$,
the following holds for any 
$\varphi\in\nbigh_{w}$:
\[
  \int_{X}\Bigl(
 \bigl|\nabla_h\varphi\bigr|_h^2
 +\bigl|
 (\theta+w\,d\zeta)\varphi
 \bigr|^2_h
 \Bigr)
\geq
 C\,|w|^{1+\rho}\int_{X}
 |\varphi|_h^2\dvol_{X}
\]
\end{prop}
\pf
We again use the argument in \S2.4 of \cite{Szabo}
with an adjustment to our situation. 
We use the standard distance on $X$.
We take small neighbourhoods $B_P$ of $P\in D$.
There exists $R_1>0$ and $C_1>0$
such that, if $|w|\geq R_1$, then
we have
$\bigl|
 (\theta+w\,d\zeta)\varphi
 \bigr|_h^2
\geq
 C_1\,|w|^2\,|\varphi|_h^2\,\dvol_X$
on $X\setminus
 \bigcup_{P\in D}B_P$.
We have only to show the estimate on
each $B_P$.
We may assume $P=0$,
and $B_P$ is an $\epsilon$-ball
$B_{\epsilon}=\{|\zeta|\leq \epsilon\}$.

We have a ramified covering
$\psi:(B'_{\epsilon},0)\lrarr (B_{\epsilon},0)$
given by $\psi(u)=u^p$
such that
$\psi^{\ast}(E,\delbar_E,\theta,h)$
is unramified,
i.e.,
we have the decomposition
\begin{equation}
 \label{eq;13.1.9.1}
 \psi^{\ast}(E,\delbar_E,\theta)
=\bigoplus_{\gminia\in u^{-1}\cnum[u^{-1}]}
 (E_{\gminia},\delbar_{E_{\gminia}},\theta_{\gminia}),
\end{equation}
where 
the Higgs field
$\theta_{\gminia}-d\gminia\,\id_{E_{\gminia}}$
are tame.
Let $h'=\bigoplus h_{|E_{\gminia}}$,
and let $\nabla_{h'}$ denote the unitary connection
associated to $\psi^{\ast}(E,\delbar_E)$ with $h'$.
By the asymptotic orthogonality of 
the decomposition (\ref{eq;13.1.9.1})
with respect to $h$
(see \cite{Mochizuki-wild}),
we have the following inequality:
\[
 \int_{B_{\epsilon}'}
 \Bigl(
 \bigl|\nabla_{h}\varphi
 \bigr|_{h}^2
+\bigl|
 (\theta+w\,d\zeta)\varphi
 \bigr|_h^2
 \Bigr)
\geq
 C_2
 \int_{B_{\epsilon}'}
 \Bigl(
 \bigl|
 \nabla_{h'}\varphi
 \bigr|_{h'}^2
+\bigl|
 (\theta+w\,d\zeta)\varphi
 \bigr|_{h'}^2
 \Bigr)
\]
\[
 \int_{B_{\epsilon}'(P)}
 \bigl|\varphi\bigr|_h^2
 \psi^{\ast}\dvol_X
\leq C_3\int_{B_{\epsilon}'(P)}
 \bigl|
 \varphi
 \bigr|_{h'}^2
 \psi^{\ast}\dvol_X
\]
Hence, we need only the estimate
with respect to the metric $h'$.

\vspace{.1in}

Let us begin with the estimate for sections of
$E_{\gminia}$ with $\gminia\neq 0$.
We set $n:=\deg_{u^{-1}}\gminia$.
\begin{lem}
\label{lem;12.7.19.1}
There exist constants 
$R'>0$ and $C'>0$
such the following holds if $|w|\geq R'$:
\begin{itemize}
\item
Let $\varphi$ be an $L^2$-section of $E_{\gminia}$
on $B_{\epsilon}'$ with respect to
$\psi^{\ast}\dvol_X$,
such that
\[
 \int_{B_{\epsilon}'}
 \Bigl(
 \bigl|\nabla_{h'}\varphi\bigr|_{h'}^2
+\bigl|(\theta_{\gminia}+w d\zeta)
 \,\varphi\bigr|_{h'}^2
 \Bigr)
 \psi^{\ast}\dvol_X
<\infty.
\]
Then, we have
\[
 |w|^{e}\int_{B_{\epsilon}'}
 \bigl|\varphi\bigr|_{h'}^2\,
 \psi^{\ast}\dvol_X
<C'\int_{B_{\epsilon}'}
 \Bigl(
 \bigl|\nabla_{h'}\varphi\bigr|_{h'}^2
+\bigl|(\theta_{\gminia}+w d\zeta)
 \,\varphi\bigr|_{h'}^2
 \Bigr)
 \psi^{\ast}\dvol_X.
\]
Here, $e=1+p/(n+p)>1$.
\end{itemize}
\end{lem}
\pf
For each $w$,
we have the solutions $b_i(w)$ $(i=0,\ldots,n+p-1)$
of the following equation:
\[
 \del_u\gminia(u)+pw u^{p-1}=0
\]
We have the equality
$u^{-p+1}\del_u\gminia(u)+pw
=\alpha\prod_{i=0}^{n+p-1}\bigl(u^{-1}-b_i(w)^{-1}\bigr)$
for some $\alpha\in\cnum\setminus\{0\}$.
We have
\[
 \theta_{\gminia}
=\del_u\gminia\,\id_{E_{\gminia}}\,du
+g_{\gminia}\,du,
\]
where $|g_{\gminia}|_{h'}\leq C_1|u|^{-1}$.
We have $R_{2}>0$ and $C_{2}>0$ such that 
the following holds if $|w|>R_2$:
\[
 C_{2}^{-1}\leq 
 |b_i(w)|\,|w|^{1/(n+p)}\leq C_{2}
\]
We take $C_{3}>\!>C_{2}$.
We set
$\nbigu_1:=\bigl\{
 |u|\leq C_{3}^{-1}|w|^{-1/(n+p)}
 \bigr\}$
and
$\nbigu_2:=\bigl\{
 |u|\leq C_{3}|w|^{-1/(n+p)}
 \bigr\}$.

Let us consider the estimate on $B_{\epsilon}'\setminus\nbigu_2$.
We have
$|g_{\gminia}|_{h'}\leq (C_1/C_3)|w|^{1/(n+p)}$.
We also have
\[
\bigl| u^{-1}-b_i(w)^{-1} \bigr|
\geq
 |b_i(w)^{-1}|-|u^{-1}|
\geq (C_2^{-1}-C_3^{-1})
|w|^{1/(n+p)}
\]
for any $i$,
and hence
$\Bigl|
 u^{-p+1}\del_{u}\gminia+pw
 \Bigr|
\geq 
|\alpha|(C_2^{-1}-C_3^{-1})^{n+p}|w|$.
Hence, if $C_3$ is sufficiently larger than $C_2$,
we have the following for some $C_4>0$:
\[
 \bigl|
 \bigl(
 \del_u\gminia+pw u^{p-1}
 \bigr)\varphi+g_{\gminia}(\varphi)
 \bigr|_{h'}
\geq
C_{4}|w|\,|\varphi|_{h'}
 |u|^{p-1}
\]
Hence, we obtain the inequality for the integral over
$B_{\epsilon}'\setminus \nbigu_2$.

\vspace{.1in}
Let us consider the estimate on $\nbigu_1$.
There exist $C_5>0$ and $R_5>0$
such that 
\[
 \bigl|
 (\del_u\gminia+pw u^{p-1})
 \varphi\,du
 \bigr|_{h'}
\geq
 C_5|u|^{-n-1}|\varphi|_{h'}\,|du|. 
\]
We also have
$\bigl|
 g_{\gminia}\varphi\,du
 \bigr|_{h'}
\leq
 C_{1}\,|u^{-1}|\,|\varphi|_{h'}\,|du|$.
Hence, there exists $C_6>0$ such that
\[
 \Bigl|
 \bigl( \theta_\gminia+w p u^{p-1}\,du\bigr) 
 \varphi
 \Bigr|^2_{h'}
\geq
 C_{6}|\varphi|_{h'}^2|u|^{-2(n+p)}\,|u|^{2(p-1)}|du\,d\ubar|
\geq
 C_6C_3|\varphi|_{h'}^2|w|^2\,|u|^{2(p-1)}|du\,d\ubar|.
\]
Therefore, we have the desired inequality
for the integral over $\nbigu_1$.

\vspace{.1in}
Let us consider the estimate on $\nbigu_2\setminus \nbigu_1$.
For each $i=0,\ldots,n+p-1$,
we set
$\nbigvtilde_i:=\bigl\{
 |u-b_i(w)|\leq \epsilon_1|w|^{-1/(n+p)}
 \bigr\}$ for some $\epsilon_1>0$.
Let $u\in \nbigu_2\setminus(\nbigu_1\cup\bigcup_i\nbigvtilde_i)$.
We have
\[
 \bigl|
 u^{-p+1}\del_u\gminia+pw
 \bigr|
=|pw|\,|u|^{-p-n}
\prod_{i=0}^{n+p-1}
 \bigl|u-b_i(w)\bigr|
\geq
 pC_3^{-1}|w|^2\prod_{i=0}^{n+p-1}
 \bigl|u-b_i(w)\bigr|
\geq 
 pC_3^{-1}\epsilon_1^{p+n} |w|
\]
We also have the following:
\[
 |g_{\gminia}\varphi|_{h'}
\leq
 C_1|u|^{p-1}|\varphi|_{h'}\cdot |u|^{-p}
\leq
 C_1|u|^{p-1}|\varphi|_{h'}
 \cdot C_3|w|^{p/(n+p)}
=C_1|u|^{p-1}|\varphi|_{h'}\,|w|\,
 \cdot C_3|w|^{-n/(n+p)}
\]
Hence, there exists $C_7>0$ and $R_7>0$
such that the following holds on
$\nbigu_2\setminus(\nbigu_1\cup\bigcup_i\nbigvtilde_i)$,
if $|w|\geq R_7$:
\[
 \bigl|
 (\del_u\gminia+pw u^{p-1})\varphi\,du
+g_{\gminia}\varphi\,du
 \bigr|_{h'}
\geq
 C_{7}|w|\,|\varphi|_{h'}\,|u|^{p-1}|du|
\]
We set 
$a:=(n+2)/2(n+p)$.
We put
$\nbigv_i:=\bigl\{
 |u-b_i(w)|\leq \epsilon_1|w|^{-a}
 \bigr\}$
and 
$\nbigv_i':=\bigl\{
 |u-b_i(w)|\leq \epsilon_1|w|^{-a}/2
 \bigr\}$.
On $\nbigvtilde_i\setminus\nbigv_i'$,
we have
\begin{multline}
 \bigl|
 u^{-p+1}\del_u\gminia+pw
 \bigr|
=pC_3^{-1}|w|^2
 \prod_{i=0}^{n+p-1}
 |u-b_i(w)|
\geq
 pC_3^{-1}|w|^2
 (\epsilon_1|w|^{-1/(n+p)})^{n+p-1}
\times
 (\epsilon_1|w|^{-a}/2)
\\
\geq
 pC_3^{-1}\epsilon_1^{p+n}
|w|^{1+1/(n+p)-a}
\end{multline}
We also have
$|g_{\gminia}|\leq C_1C_3|w|^{1/(n+p)}$.
Because 
$-(p-1)/(n+p)+1+1/(n+p)-a>1/(n+p)$,
there exist $C_8>0$ and $R_8>0$
such that the following holds on
$\nbigvtilde_i\setminus\nbigv_i'$,
if $|w|\geq R_8$:
\[
 \Bigl|
 \bigl(
 \theta_{\gminia}
+pw u^{p-1}du
\bigr)
 \varphi
 \Bigr|_{h'}
\geq
 C_8|w|^{1+1/(n+p)-a}|u|^{p-1}\,|du|\,|\varphi|_{h'}
=C_8|w|^{(n+2p)/2(n+p)}|u|^{p-1}\,|du|\,|\varphi|_{h'}
\]
We have the following type of Poincar\'e inequality,
i.e., there exist $C_9>0$ and $R_9>0$ 
such that the following holds on $\nbigv_i$,
if $|w|>R_9$
(see \cite{biquard-boalch} and (2.12) of \cite{Szabo}):
\[
 |w|^{(n+2)/(n+p)}
 \int_{\nbigv_i}|\varphi|_{h'}^2|du\,d\ubar|
\leq
 C_{9}\left(
 \int_{\nbigv_i}\bigl|
 d|\varphi|_{h'}
 \bigr|^2
+|w|^{(n+2)/(n+p)}
 \int_{\nbigv_i\setminus\nbigv_i'}
 |\varphi|_{h'}^2\,|du\,d\ubar|
 \right)
\]
We also have the following inequalities:
\begin{equation}
 |w|^{(n+2p)/(n+p)}\int_{\nbigv_i}
 |\varphi|_{h'}^2|u|^{2(p-1)}\,|du\,d\ubar|
\leq
 C_3^{2(p-1)}|w|^{(n+2)/(n+p)}
 \int_{\nbigv_i}|\varphi|_{h'}^2\,|du\,d\ubar|
\end{equation}
\begin{multline}
 |w|^{(n+2)/(n+p)}
 \int_{\nbigv_i\setminus\nbigv_i'}|\varphi|_{h'}^2\,|du\,d\ubar|
\leq
 C_3^{2(p-1)}
 |w|^{(n+2p)/(n+p)}
 \int_{\nbigv_i\setminus\nbigv_i'}|\varphi|_{h'}^2\,
 |u|^{2(p-1)}\, |du\,d\ubar|
 \\
 \leq
 C_8 C_3^{2(p-1)}
 \int_{\nbigv_i\setminus\nbigv_i'}
 \bigl|(\theta_{\gminia}+pw u^{p-1}du)
 \varphi\bigr|_{h'}^2
\end{multline}
Then, we obtain the desired inequality
for the integral over $\nbigv_i$.
Thus, the proof of Lemma \ref{lem;12.7.19.1}
is finished.
\hfill\qed

\vspace{.1in}

Let us consider the case $\gminia=0$.
Because this part is essentially contained in \cite{Szabo},
we give just an indication.
We take a positive number $C_{10}$ 
which is sufficiently larger than
$|\alpha|$ for any eigenvalues $\alpha$
of the residue of $\theta_0$.
We may assume
$|g_0|\leq (C_{10}/10)|\zeta|^{-1}$
on $B_{\epsilon}$.
Take $R_{10}>0$ sufficiently larger than $C_{10}$.
For $|w|\geq R_{10}$,
let $\nbigu:=\bigl\{
 |\zeta|\leq C_{10}|w|^{-1}
 \bigr\}$
and 
$\nbigu':=\bigl\{
 |\zeta|\leq C_{10}|w|^{-1}/2
 \bigr\}$.
On $B_{\epsilon}\setminus \nbigu'$,
we have 
\begin{equation}
 \label{eq;12.7.19.2}
 \bigl|
 (\theta_{0}+w\,d\zeta)\varphi
 \bigr|_{h'}
\geq
 |w|\,|\varphi|_{h'}\,|d\zeta|
-|g_0|_{h'}\,|\varphi|_{h'}\,|d\zeta|
\geq
\frac{4}{5}|w|\,|\varphi|_{h'}|d\zeta|
\end{equation}
There exist $C_{11}>0$ and $R_{11}>0$
such that the following holds on $\nbigu$,
if $|w|\geq R_{11}$:
\begin{equation}
 \label{eq;12.7.19.3}
 |w|^2\int_{\nbigu}|\varphi|_{h'}^2\,|d\zeta\,d\zetabar|
\leq
 C_{11}\int_{\nbigu}\bigl|d|\varphi|_{h'}\bigr|^2
+|w|^2\int_{\nbigu\setminus\nbigu'}
 |\varphi|_{h'}^2\,|d\zeta\,d\zetabar|
\leq
 \int_{\nbigu}
 \Bigl(
 C_{11} |\nabla_{h'}\varphi|_{h'}^2
+4\bigl|(\theta_0+w d\zeta)\varphi\bigr|_{h'}^2
\Bigr)
\end{equation}
We obtain the desired inequality
for sections of $E_0$
from (\ref{eq;12.7.19.2}) 
and (\ref{eq;12.7.19.3}).
Thus, the proof of Proposition \ref{prop;12.7.19.4}
is finished.
\hfill\qed

\vspace{.1in}
The following is a refinement of Corollary
\ref{cor;13.2.27.22}.

\begin{cor}
\label{cor;12.7.20.3}
There exist $\rho>0$ and $C>0$ such that the following holds:
\begin{itemize}
\item
Let $\varphi$ be a section of $E$
such that 
\begin{equation}
 \label{eq;13.1.25.110}
 \int |\varphi|_h^2\dvol_{X\setminus D}
+\int |\Delta_{1}\varphi|_h^2
 \dvol_{X}
<\infty.
\end{equation}
Then, we have the following inequality:
\begin{equation}
 \label{eq;13.1.25.111}
  C|w|^{1+\rho}
 \Bigl(
 \int |\varphi|_h^2\dvol_{X}
 \Bigr)^{1/2}
\leq
 \Bigl(
\int|\Delta_{1}\varphi|_{h}^2\dvol_{X}
 \Bigr)^{1/2}
\end{equation}
\end{itemize}
\end{cor}
\pf
It is shown by the argument for
Corollary \ref{cor;13.2.27.22},
by using Proposition \ref{prop;12.7.19.4},
instead of Proposition \ref{prop;13.2.27.20}.
\hfill\qed

\subsection{Comparison with the algebraic Nahm transform}

\subsubsection{Statements}

Let $\harmonicbundle$ be a wild harmonic bundle
on $(T^{\lor},D)$.
Let $\nbigp_{\ast}E$ be the associated filtered bundle 
on $(T^{\lor},D)$.
Let $(E_1,h_1,\nabla_1)$ be the $L^2$-instanton on 
$T\times\cnum$
obtained as the Nahm transform of $\harmonicbundle$
(see \S\ref{subsection;13.1.27.10}).
Let $\nbigp_{\ast}E_1$ be the associated filtered bundle
on $(T\times\proj^1,T\times\{\infty\})$.

\begin{thm}
\label{thm;13.1.23.20}
There is a natural isomorphism
of the filtered bundles
$\nbigp_{\ast}E_1\simeq
 \Nahm_{\ast}\bigl(\nbigp_{\ast}E,\theta\bigr)$.
\end{thm}

Conversely, let $(E_1,\nabla_1,h_1)$ be an $L^2$-instanton 
on $T\times\cnum$.
Let $\nbigp_{\ast}E_1$ be the associated filtered bundle
on $(T\times\proj^1,T\times\{\infty\})$.
Let $\harmonicbundle$ be the wild harmonic bundle
on $(T^{\lor},D)$
obtained as the Nahm transform of
$(E_1,\nabla_1,h_1)$
(see \S\ref{subsection;13.2.3.10}).
Let $(\nbigp_{\ast}E,\theta)$ be the associated filtered
Higgs bundle.

\begin{thm}
\label{thm;13.1.23.100}
There is a natural isomorphism of
the filtered Higgs bundles
$(\nbigp_{\ast}E,\theta)\simeq
\Nahm_{\ast}(\nbigp_{\ast}E_1)$.
\end{thm}

We obtain the involutivity of the Nahm transform
in the following sense.

\begin{cor}
\label{cor;13.1.30.31}
For an $L^2$-instanton
$(E_1,\nabla_1,h_1)$ on $T\times\cnum$,
we have an isomorphism
\[
 \Nahm\bigl(
 \Nahm(E_1,\nabla_1,h_1)
 \bigr)
\simeq
 (E_1,\nabla_1,h_1).
\]
For a wild harmonic bundle 
$\harmonicbundle$ on $T^{\lor}$,
we have an isomorphism
\[
 \Nahm\bigl(
 \Nahm\harmonicbundle
 \bigr)
\simeq \harmonicbundle.
\]
\end{cor}
\pf
It follows from
Proposition \ref{prop;13.1.22.2},
Theorem \ref{thm;13.1.23.20},
Theorem \ref{thm;13.1.23.100},
and the uniqueness
of the harmonic metric or Hermitian-Einstein metric
adapted to the filtered bundle.
(See Proposition \ref{prop;13.1.29.20}
for the uniqueness of Hermitian-Einstein metric.
See \cite{biquard-boalch} for the uniqueness of
the harmonic metric.
See also \cite{Mochizuki-wild}.
See \S\ref{subsection;13.10.28.3}
for adaptedness of metrics
and filtered bundles.)
\hfill\qed

\subsubsection{Proof of Theorem \ref{thm;13.1.23.20}}
\label{subsection;13.3.5.10}

Let us construct an isomorphism
$(E_1,\delbar_{E_1})\simeq
 \Nahm(\nbigp_{\ast}E,\theta)_{|T\times\cnum}$.
We recall the monad construction of
$E_1=\Nahm\harmonicbundle$ \cite{Donaldson-Kronheimer}.
We use the notation in \S\ref{subsection;13.2.14.1}.
Let $g_{T^{\lor}\setminus D}$ be a Poincar\'e like K\"ahler metric
of $T^{\lor}\setminus D$.
Let $\nbiga^{i}\harmonicbundle$ denote the space of
sections $\varphi$
of $E\otimes\Omega^i$ on $T^{\lor}\setminus D$
such that
$\varphi$ and $(\delbar_E+\theta)\varphi$
are $L^2$ with respect to $h$
and $g_{T^{\lor}\setminus D}$.
Note that the conditions also imply
$(\delbar_{E,z}+\theta_w)\varphi$
are $L^2$ for any $(z,w)\in T\times \cnum$.
Let $\underline{\nbiga}^i$
denote the sheaf of holomorphic sections of
the product bundle
$\nbiga^{i}\harmonicbundle\times(T\times\cnum)$
over $T\times\cnum$.
We have the morphisms
$\delta^i:
 \underline{\nbiga}^i\lrarr
 \underline{\nbiga}^{i+1}$
induced by
$\delbar_{E,z}+\theta_w$,
and the sheaf of holomorphic sections of $E_1$
is isomorphic to
$\Ker\delta^1/\Image\delta^0$.

Applying the construction in the proof of 
Lemma \ref{lem;13.1.26.10} around each point of $D$,
we extend $E$ and $E\otimes\Omega^1$
to $\nbigc_{L^2}^0(\nbigp_{\ast}E,\theta)$
and $\nbigc_{L^2}^1(\nbigp_{\ast}E,\theta)$.
Let $\nbigc_{L^2}^{i,\bullet}(\nbigp_{\ast}E,\theta)$
denote the Dolbeault resolution of 
$\nbigc_{L^2}^i(\nbigp_{\ast}E,\theta)$.

For $I\subset \{1,2,3\}$,
let $p_I$ denote the projection of
$T^{\lor}\times T\times\cnum$
onto the product of the $i$-th components
$(i\in I)$.
On $T^{\lor}\times T\times\cnum$,
we set
\[
 \nbigctilde^{i}_{L^2}:=
 \bigoplus_{k+\ell=i}
 p_1^{-1}\nbigc^{k,\ell}_{L^2}(\nbigp_{\ast}E,\theta)
 \otimes_{p_1^{-1}\nbigo_{T^{\lor}}}
 p_{12}^{\ast}\poincare 
\]
We have $\delta^i:\nbigctilde_{L^2}^i\lrarr\nbigctilde_{L^2}^{i+1}$
induced by $\delbar_{E,z}+\theta_w$.
We have a natural inclusion
$\Phi:p_{23\ast}\nbigctilde^{\bullet}_{L^2}
\lrarr
 \underline{\nbiga^{\bullet}}$
of complexes on $T\times\cnum$.
According to the results in \S5.1 of \cite{Mochizuki-wild},
$\Phi$ is a quasi-isomorphism.
We also have the following natural isomorphisms
in $D^b(\nbigo_{T\times\cnum})$:
\[
 p_{23\ast}\nbigctilde_{L^2}
\simeq
 Rp_{23\ast}\Bigl(
 p_{1}^{\ast}\bigl(
 \nbigc_{L^2}^{\bullet}(\nbigp_{\ast}E,\theta)
 \bigr)
\otimes
 p_{12}^{\ast}\poincare
\Bigr)
\stackrel{\Psi}{\simeq}
 Rp_{23\ast}\Bigl(
 p_{1}^{\ast}\bigl(
 \nbigc^{\bullet}(\nbigp_{\ast}E,\theta)
 \bigr)
\otimes
 p_{12}^{\ast}\poincare
\Bigr)
 \simeq
 \Nahm(\nbigp_{\ast}E,\theta)_{|T\times\cnum}
\]
(See Lemma \ref{lem;13.1.26.10}
for $\Psi$.)
Thus, we obtain the desired isomorphism
$E_1\simeq
\Nahm(\nbigp_{\ast}E,\theta)_{|T\times\cnum}$,
by which we shall identity them.

\begin{lem}
\label{lem;13.1.23.101}
To prove Theorem {\rm\ref{thm;13.1.23.20}},
we have only to show 
$\Nahm_a(\nbigp_{\ast}E,\theta)
\subset \nbigp_{a}E_1$
for any $a$.
\end{lem}
\pf
By Proposition \ref{prop;13.1.22.3},
we have
$\deg(\Nahm_{\ast}(\nbigp_{\ast}E,\theta))
=\deg(\nbigp_{\ast}E,\theta)=0$.
By Proposition \ref{prop;13.1.24.3},
we also have
$\deg(\nbigp_{\ast}E_1)=0$.
Hence, 
$\Nahm_a(\nbigp_{\ast}E,\theta)
\subset\nbigp_{a}E_1$
implies 
$\Nahm_a(\nbigp_{\ast}E,\theta)
=\nbigp_{a}E_1$.
\hfill\qed

\vspace{.1in}

To show 
$\Nahm_a(\nbigp_{\ast}E,\theta)
\subset \nbigp_{a}E_1$,
we need an estimate of the upper bound
of the norms of sections of
$\Nahm_a(\nbigp_{\ast}E,\theta)$.
We use an argument of scaling in \cite{Szabo}.
Because we need only the upper bound,
we will not consider more precise estimates
for harmonic representatives or their approximation.

Let $U_{\tau}\subset\proj^1$ be a neighbourhood of $\infty$
with the coordinate $\tau=w^{-1}$.
If $U_{\tau}$ is sufficiently small,
we have the decomposition
$\Nahm_{\ast}(\nbigp_{\ast}E,\theta)
=\bigoplus_{P\in D}
 \Nahm_{\ast}(\nbigp_{\ast}E,\theta)_P$
by the spectrum
on $T\times U_{\tau}$.
We have the refined decomposition
\[
 \Nahm_{\ast}(\nbigp_{\ast}E,\theta)_P
=\bigoplus_{\vecgminio\in\vecIrr(\theta,P)}
 \bigoplus_{\alpha\in\cnum}
 \Nahm_{\ast}(\nbigp_{\ast}E,\theta)_{P,\vecgminio,\alpha},
\]
according to the decomposition of the filtered Higgs bundle
$(\nbigp_{\ast}E,\theta)
=\bigoplus_{\vecgminio\in\vecIrr(\theta,P)}
 \bigoplus_{\alpha\in\cnum}
 (\nbigp_{\ast}E_{P,\vecgminio,\alpha},\theta_{P,\vecgminio,\alpha})$ 
around each $P\in D$.
We have only to show that
$\Nahm_{a}(\nbigp_{\ast}E,\theta)_{P,\vecgminio,\alpha}
\subset
 \nbigp_{a}E_1$.
We shall argue the case $P=\{0\}$ in the following.
The other case can be established similarly.
We omit the subscript $P$.
We take a small neighbourhood
$U_{\zeta}\subset T^{\lor}$ of $\{0\}$.

Let us consider the case 
$(\vecgminio,\alpha)\neq (0,0)$.
Take $\gminia\in\vecgminio$.
For each $c\in\real$,
we have the frame of 
$\Nahm_{c}(\nbigp_{\ast}E,\theta)_{P,\vecgminio,\alpha}$
in Lemma \ref{lem;13.1.23.30}.
We have only to show 
\begin{equation}
\label{eq;13.1.23.31}
\Bigl|
 \bigl[
 \zeta_{\vecgminio}^j
 v_{\vecgminio,i}d\zeta_{\vecgminio}/\zeta_{\vecgminio}
 \bigr]
\Bigr|_{h_1}
=O\bigl(
 |w|^{(b-j-m_{\vecgminio}/2)(p_{\vecgminio}+m_{\vecgminio})^{-1}}
 \bigr)
\end{equation}
Here, $b$ is the parabolic degree of $v_{\vecgminio,i}$.

We give a preliminary.
We have the expression
$\zeta_{\vecgminio}\del_{\zeta_{\vecgminio}}\gminia
+p_{\vecgminio}\alpha 
=\sum_{j=0}^{m_{\vecgminio}}
 \alpha_j \zeta_{\vecgminio}^{-j}
=:G(\zeta_{\vecgminio})$.
We fix a complex number $\gamma$
such that 
$\alpha_{m_{\vecgminio}}
+p_{\vecgminio} \gamma^{p_{\vecgminio}+m_{\vecgminio}}
=0$.
Take a covering
$U_{\eta}\lrarr U_{\tau}$
given by $\eta\longmapsto \eta^{p_{\vecgminio}+m_{\vecgminio}}$.
If $U_{\tau}$ is sufficiently small,
we can take holomorphic functions $u^{(i)}_0(\eta)$
$(i=1,\ldots,p_{\vecgminio+m_{\vecgminio}})$
satisfying the following:
\[
 G\bigl(
 u^{(i)}_{0}(\eta)
 \bigr)
+p_{\vecgminio}u^{(i)}_0(\eta)^{p_{\vecgminio}}
 \eta^{-p_{\vecgminio}-m_{\vecgminio}}=0,
\quad\quad\quad
 \lim_{\eta\to 0}
 u^{(i)}_{0}(\eta)/\eta
=\gamma
 \exp\bigl(2\pi\sqrt{-1} i/(m_{\vecgminio}+p_{\vecgminio})\bigr)
\]
There exist $C_1>0$ and $\epsilon_1>0$
such that 
\[
 \Bigl|
 \eta^{-1} u^{(i)}_0(\eta)
 -\gamma
  \exp\bigl(2\pi\sqrt{-1} i/(m_{\vecgminio}+p_{\vecgminio})\bigr)
 \Bigr|
\leq
 C_1|\eta|^{\epsilon_1}.
\]

\begin{lem}
\label{lem;13.1.22.100}
Let $Z_{\eta}$ denote the support of
$\Cok\bigl(
 \eta^{p_{\vecgminio}+m_{\vecgminio}}
 \zeta_{\gminio}^{m_{\vecgminio}}
 \theta^{\vecgminio}_{\gminia,\alpha}
 +p_{\vecgminio}\zeta_{\vecgminio}^{p_{\vecgminio}+m_{\vecgminio}}
 d\zeta_{\vecgminio}/\zeta_{\vecgminio}\bigr)$
on $U_{\zeta_{\vecgminio}}$.
If $U_{\zeta}$ and $U_{\tau}$ are sufficiently small,
there exists a decomposition
$Z_{\eta}=\coprod_{i=1}^{p_{\vecgminio}+m_{\vecgminio}}
 Z^{(i)}_{\eta}$ such that
the following holds for any $u\in Z^{(i)}_{\eta}$:
\[
 \bigl|
 u_0^{(i)}(\eta)-u
 \bigr|
\leq C|\eta|^{1+m_{\vecgminio}+\epsilon}
\]
Here $C$ and $\epsilon$ are positive constants
which are independent of $\eta$.
\end{lem}
\pf
Take $u_1\in Z_{\eta}$.
There exists a possibly multi-valued holomorphic $1$-form 
$\nu(\zeta_{\vecgminio})\,d\zeta_{\vecgminio}/\zeta_{\vecgminio}$
obtained as the eigenvalue
of $\theta^{\vecgminio}_{\gminia}$,
such that
$\nu(u_1)
 +\eta^{-p_{\vecgminio}-m_{\vecgminio}}
 p_{\vecgminio}u_1^{p_{\vecgminio}}=0$.
Because
$\nu(\zeta_{\vecgminio})-G(\zeta_{\vecgminio})
=O\bigl(\zeta_{\vecgminio}^{\epsilon}\bigr)$,
there exist $C_2>0$ and $\epsilon_2>0$,
independently from $\eta$,
such that 
the following holds for some unique $i$:
\begin{equation}
\label{eq;13.2.27.1}
 \Bigl|
 \eta^{-1}u_1-\gamma
 \exp\bigl(2\pi\sqrt{-1}i/(p_{\vecgminio}+m_{\vecgminio})\bigr)
\Bigr|
\leq C_2|\eta|^{\epsilon_2}.
\end{equation}
We obtain a decomposition of 
$Z_{\eta}=\coprod Z^{(i)}_{\eta}$
by the condition (\ref{eq;13.2.27.1}).

Let $u_1\in Z^{(i)}_{\eta}$.
We set $Q_q(x,y):=\sum_{i+j=q} x^iy^j$.
We have
\begin{multline}
\Bigl(
 \bigl(u_0^{(i)}(\eta)/\eta\bigr)^{-1}
-(u_1/\eta)^{-1}
\Bigr)
 \times \\
\left(
 \sum_{j=1}^{m_{\vecgminio}}
 \alpha_j\eta^{m_{\vecgminio}-j}
 Q_{j-1}\bigl((u_0^{(i)}/\eta)^{-1},(u_1/\eta)^{-1}\bigr)
-\bigl(u^{(i)}_0/\eta\bigr)
 \bigl(u_1/\eta\bigr)
 p_{\vecgminio}Q_{p_{\vecgminio}-1}\bigl(
 u_0^{(i)}/\eta,u_1/\eta
 \bigr)
\right)
 \\
=O\Bigl(
 |u_1/\eta|^{\epsilon}
\,|\eta|^{m_{\vecgminio}+\epsilon}
 \Bigr)
\end{multline}
We obtain 
$\bigl|
 \bigl(u_0^{(i)}(\eta)/\eta\bigr)^{-1}
-(u_1/\eta)^{-1}
\bigr|
=O(|\eta|^{m_{\vecgminio}+\epsilon})$.
Then, we obtain the desired estimate.
\hfill\qed

\vspace{.1in}

Let $\rho$ be an $\real_{\geq 0}$-valued function on $\cnum_{\eta}$
such that $\rho(\eta)=1$ for $|\eta|<1/2$
and $\rho(\eta)=0$ for $|\eta|>1$.
We set $u_0:=u_0^{(0)}$.
We consider the following $C^{\infty}$-sections
of $E^{\vecgminio}_{\gminia,\alpha}
 \otimes\Omega^1_{X^{\vecgminio}}$:
\[
 \mu_1(v_{\vecgminio,i},\xi):=
 \rho\Bigl(
 |\xi|^{1+m_{\vecgminio}/2}\bigl(\zeta_{\vecgminio}-u_0(\xi)\bigr)
 \Bigr)
 \,v_{\vecgminio,i}\,d\zeta_{\vecgminio}/\zeta_{\vecgminio}
\]
\[
 \mu_2(v_{\vecgminio,i},\xi):=
 \bigl(
 \theta^{\vecgminio}_{\gminia}
+\xi^{p_{\vecgminio}+m_{\vecgminio}}
 d\zeta_{\vecgminio}^{p_{\vecgminio}}
 \bigr)^{-1}
 \Bigl(
 \delbar_E\mu_1(v_{\vecgminio,i},\xi)
 \Bigr)
\]
By Lemma \ref{lem;13.1.22.100},
if $|\xi|$ is sufficiently large,
$\rho\Bigl(
 |\xi|^{1+m_{\vecgminio}/2}\bigl(\zeta_{\vecgminio}-u_0(\xi)\bigr)
 \Bigr)$ is constantly $1$ around $Z_{\xi}$.
Hence, the tuple
$\vecmu(v_{\vecgminio,i},\xi)=
 \bigl(
 \mu_1(v_{\vecgminio,i},\xi),
 \mu_2(v_{\vecgminio,i},\xi)
 \bigr)$
gives a representative of
$[v_{\vecgminio,i} d\zeta_{\vecgminio}/\zeta_{\vecgminio}]$.

By an elementary change of variables,
we obtain the following for any $\delta>0$:
\[
 \int \bigl|\mu_1(v_{\vecgminio,i},\xi)\bigr|_h^2
\leq
 \int \rho\bigl(|\xi|^{1+m_{\vecgminio}/2}
 (\zeta_{\vecgminio}-u_0(\xi))\bigr)^2
 |\zeta_{\vecgminio}|^{-2(b+\delta)-2}
 \bigl|d\zeta_{\vecgminio}\,d\zetabar_{\vecgminio}\bigr|
\leq
 C_{1\delta}|\xi|^{2(b+\delta)-m_{\vecgminio}}
\]
Note that 
we have
$|\zeta_{\vecgminio}-u_0(\xi)|\sim
 |\xi|^{-1-m_{\vecgminio}/2}$
for $\zeta_{\vecgminio}$
such that 
$\delbar
 \rho\bigl(|\xi|^{1+m_{\vecgminio}/2}
 (\zeta_{\vecgminio}-u_0(\xi))\bigr)\neq 0$.
Hence, we also have the following:
\begin{multline}
 \int \bigl|\mu_2(v_{\vecgminio,i},\xi)\bigr|_h^2
\leq
 C_{2\delta}\int\bigl|
 \delbar
 \rho\bigl(|\xi|^{1+m_{\vecgminio}/2}
 (\zeta_{\vecgminio}-u_0(\xi))\bigr)
 \bigr|^2
 |\zeta_{\vecgminio}|^{-2(b+\delta)}
 \frac{1}{\Bigl|
 \zeta_{\vecgminio}\del_{\zeta_{\vecgminio}}\gminia
 +p_{\vecgminio}\alpha
 +p_{\vecgminio}\xi^{p_{\vecgminio}+m_{\vecgminio}}
 \zeta_{\vecgminio}^{p_{\vecgminio}}
 \Bigr|^2}
 \\
\leq
 C_{3\delta}|\xi|^{2(b+\delta)-2(m_{\vecgminio}+1)+2(1+m_{\vecgminio}/2)}
=C_{3\delta}|\xi|^{2(b+\delta)-m_{\vecgminio}}
\end{multline}
By the construction of $h_1$,
we have 
$\bigl|
 [v_{\vecgminio,i} d\zeta_{\vecgminio}/\zeta_{\vecgminio}]
 \bigr|_{h_1}^2
\leq 
 \int \Bigl(
 \bigl|\mu_1(v_{\vecgminio,i},\xi)\bigr|_h^2
+\bigl|\mu_2(v_{\vecgminio,i},\xi)\bigr|_h^2
 \Bigr)$.
Hence, we obtain the desired estimate (\ref{eq;13.1.23.31}) for
$[v_{\vecgminio,i} d\zeta_{\vecgminio}/\zeta_{\vecgminio}]$.
We obtain the estimate for 
$[v_{\vecgminio,i} \zeta^jd\zeta_{\vecgminio}/\zeta_{\vecgminio}]$
similarly.

\vspace{.1in}

Let us consider the case
$(\vecgminio,\alpha)=(\{0\},0)$.
The following lemma is easy to see.
\begin{lem}
\label{lem;13.1.23.100}
Let $Z_{w}$ denote the support of
$\Cok\bigl(\theta_{0,0}
 +w d\zeta\bigr)$.
There exist $C>0$ and $\epsilon>0$
such that 
$|u|\leq C|w|^{-1-\epsilon}$
holds for any 
$u\in Z_w$.
\hfill\qed
\end{lem}

For a holomorphic section $s$ of
$\nbigc^1(\nbigp_{\ast} E_{P,\{0\},0}\otimes
 \Omega^{\bullet},\theta)$
(see \S\ref{subsection;13.1.21.1}),
we consider the following $C^{\infty}$-sections
of $E_{P,\{0\},0}\otimes\Omega^1$:
\[
 \mu_1(s,w):=\bigl(\rho(\zeta)-\rho(w\zeta)\bigr)\,s\,d\zeta/\zeta,
\quad\quad\quad
 \mu_2(s,w):=
 \bigl(
 \theta_{P,\{0\},0}+w\,d\zeta
 \bigr)^{-1}(\delbar\mu_1(s))
\]
\[
  \mu'_1(s,w):=\rho(\zeta)\,s\,d\zeta/\zeta,
\quad\quad\quad
 \mu'_2(s,w):=
 \bigl(
 \theta_{P,\{0\},0}+w\,d\zeta
 \bigr)^{-1}(\delbar\mu'_1(s))
\]
By Lemma \ref{lem;13.1.23.100},
$\mu_2$ and $\mu_2'$ are well defined.
The tuples $\vecmu(s,w)=\bigl(\mu_1(s,w),\mu_2(s,w)\bigr)$
and $\vecmu'(s,w)=\bigl(\mu'_1(s,w),\mu'_2(s,w)\bigr)$
naturally induce the same holomorphic section of
$\Nahm(\nbigp_{\ast}E)_{P}$.
If $s$ is a section of $\nbigp_cE_{\{0\},0}$,
then it is elementary to show the following for any $\delta>0$:
\[
 \int\bigl|\mu_i(s,w)\bigr|_h^2
\leq
 C_{\delta}|w|^{2(c+\delta)}
\]
We obtain 
$\bigl|
 \vecmu'(s,w)
 \bigr|_{h_1}
\leq
 C|w|^{c+\delta}$
for any $\delta>0$.
By the construction of
$\Nahm_{\ast}(\nbigp_{\ast}E)_{P,\{0\},0}$,
we obtain
$\Nahm_{\ast}(\nbigp_{\ast}E)_{P,\{0\},0}
\subset
 \nbigp_{\ast}E_1$.
Thus, the proof of Theorem \ref{thm;13.1.23.20}
is finished.
\hfill\qed

\subsubsection{Proof of Theorem \ref{thm;13.1.23.100}}

Let us construct an isomorphism of the Higgs bundles
$(E,\delbar_E,\theta)
\simeq
 \Nahm(\nbigp_{\ast}E_1)_{|T^{\lor}\setminus D}$.
Let us recall the monad construction of $\Nahm(E_1,\nabla_1)$.
Let $\nbiga^{0,i}$ denote the space of
sections $\varphi$ of $E_1\otimes\Omega^{0,i}$
on $T\times\cnum$,
such that 
$\varphi$ and $\delbar_{E_1}\varphi$
are $L^2$ with respect to $h_1$
and the Euclidean metric.
Let $\underline{\nbiga}^{0,i}$ denote 
the sheaf of holomorphic sections of
the product bundle
$\nbiga^{0,i}\times (T^{\lor}\setminus D)$
over $T^{\lor}\setminus D$.
We have the morphism
$\delta^i:\underline{\nbiga}^{0,i}
\lrarr\underline{\nbiga}^{0,i+1}$
induced by 
$\delbar_{E_1}-\zeta\,d\zbar$,
and the sheaf of holomorphic sections of $(E,\delbar_E)$
is isomorphic to
$\Ker\delta^1/\Image\delta^0$.

For $I\subset\{1,2,3\}$,
let $p_I$ denote the projection of
$T^{\lor}\times T\times\cnum$
onto the product of the $i$-th components.
By the construction,
we have a natural morphism
$Rp_{1\ast}(p_{23}^{\ast}\nbigp_{<-1}E
 \otimes p_{12}^{\ast}\poincare)
\lrarr
 \underline{\nbiga}^{0,\bullet}$.
By the results in \S\ref{subsection;13.1.24.1},
it is a quasi-isomorphism.
Hence, we obtain a holomorphic isomorphism
$E\simeq \Nahm(\nbigp_{\ast}E_1)_{|T^{\lor}\setminus D}$,
by which we identify them.
The Higgs fields are equal,
because they are induced by the multiplication of $-w$.

\vspace{.1in}

We give a preliminary.
Let $U\subset\proj^1$ be a small neighbourhood of $\infty$.
On $T\times U$, we have the following decomposition
\begin{equation}
 \label{eq;13.2.14.3}
 \nbigp_{\ast}E_1
=\bigoplus_{P\in\Sp_{\infty}(E_1)}
 \bigoplus_{\vecgminio,\alpha}
\nbigp_{\ast}(E_{1})_{P,\vecgminio,\alpha}.
\end{equation}
Fix a lift of $\Sp_{\infty}(E_1)\subset T^{\lor}$
to $\Sptilde_{\infty}(E_1)\subset\cnum$.
We have the filtered bundles
with an endomorphism $(\nbigp_{\ast}V,g)$ on $U$,
corresponding to $\nbigp_{\ast}E_1$.
It has a decomposition
$(\nbigp_{\ast}V,g)=
 \bigoplus(\nbigp_{\ast}V_{P,\vecgminio,\alpha},
 g_{P,\vecgminio,\alpha})$.

Let $\nbigu\subset\nbigp_{-1}(E_1)$ be the subsheaf
such that
$\nbigu_{|T\times\cnum}=\nbigp_{-1}(E_1)_{|T\times\cnum}$
and 
\begin{equation}
 \label{eq;13.2.13.20}
 \nbigu=\bigoplus_{P}
 \Bigl(
 \nbigp_{-1}(E_1)_{P,\{0\},0}
\oplus
 \bigoplus_{(\vecgminio,\alpha)\neq(\{0\},0)}
 \nbigp_{<-1}(E_1)_{P,\vecgminio,\alpha}
 \Bigr)
\end{equation}
around $T\times\{\infty\}$.
We use the notation in \S\ref{subsection;13.1.21.111}.
\begin{lem}
\label{lem;13.2.14.2}
We have $N(\nbigu)\subset \nbigp_0E$.
\end{lem}
\pf
We give an argument around $0\in T^{\lor}$,
by supposing $0\in D$.
The other case can be proved similarly.
We may suppose the lift of $0\in D$ is $0\in\cnum$.
Let $t$ be a holomorphic section of $N(\nbigu)$
around $0\in T^{\lor}$.
We have to show 
$|t|_{h}=O(|\zeta|^{-\delta})$
for any $\delta>0$.
It is represented by 
a family of $C^{\infty}$-sections 
$\kappa(\zeta)=\kappa^{1}(\zeta)d\zbar+\kappa^{2}(\zeta)d\wbar$ 
of
$\nbigp_{-1}E_1\otimes
 \Omega^{0,1}_{T\times\proj^1}\otimes L^{-1}_{\zeta}$.
According to the decomposition (\ref{eq;13.2.13.20}),
we have
\[
 \kappa^i(\zeta)
=\sum_{P,\vecgminio,\alpha}
 \kappa^i(\zeta)_{P,\vecgminio,\alpha}.
\]
If $P\neq 0$,
we may assume 
$\kappa^i(\zeta)_{P,\vecgminio,\alpha}=0$
on $U$.
(See the proof of Proposition \ref{prop;13.2.13.100}.)
Let $\dvol:=|dz d\zbar dw d\wbar|$.

We take a $C^{\infty}$-metric $h_2$ of $\nbigu$.
Note $h_1=O(h_2|w|^{-2+\delta})$ for any $\delta>0$
on $\nbigp_{-1}(E_1)_{P,\{0\},0}$,
and 
$h_1=O(h_2|w|^{-1-\epsilon})$ for some $\epsilon>0$
on $\nbigp_{<-1}(E_1)_{P,\vecgminio,\alpha}$
for $(\vecgminio,\alpha)\neq (\{0\},0)$.
If $P=0$ and $(\vecgminio,\alpha)\neq (\{0\},0)$,
we have the following finiteness uniformly for $\zeta$:
\[
  \int_{T\times U}
 \bigl|\kappa^i(\zeta)_{0,\vecgminio,\alpha}\bigr|^2_{h_1}
 \dvol
\leq
C_0\int_{T\times U}
 \bigl|\kappa^i(\zeta)_{0,\vecgminio,\alpha}\bigr|_{h_2}^2
 |w|^{-2-\epsilon}
 \dvol
<\infty
\]

We have
$|g_{0,\{0\},0}|_{h_1}\leq C_1|w|^{-1}$ for some $C_1$.
We take a sufficiently small $C_2>0$,
and we put
$H_{\zeta}:=\bigl\{w\,\big|\,|w|^{-1}<C_2|\zeta|\bigr\}$.
We can find a unique family of $C^{\infty}$-sections $\mu(\zeta)$
of $\nbigp_{-1}E\otimes 
 L_{\zeta}^{-1}$ on $H_{\zeta}$
such that
\[
(\delbar_E+\zeta\,d\zbar)\mu(\zeta)
=\Bigl(
 \kappa^1(\zeta)_{0,\{0\},0}d\zbar
+\kappa^2(\zeta)_{0,\{0\},0}d\wbar
 \Bigr)_{|H_{\zeta}}.
\]
There exists $C_3>0$ such that the following holds:
\[
 \int_{T\times\{w\}}
 |\mu(\zeta)|_{h_2}^2\,|dz d\zbar|
\leq
 C_3|\zeta|^{-2}
 \int_{T\times\{w\}}
 |\kappa^1(\zeta)_{0,\{0\},0}|_{h_2}^2\,
 |dz d\zbar|
\]
Let $\chi(w)$ be an $\real_{\geq 0}$-valued $C^{\infty}$-function
such that
$\chi(w)=1$ if $|w|^{-1}\leq C_2/4$
and $\chi(w)=0$ if $|w|^{-1}\geq C_2/2$.
We set
\[
 \kappatilde^{1}(\zeta)
=\kappa^1(\zeta)_{0,\{0\},0}
 -\del_{\zbar}\bigl(\chi(w\zeta)\mu(\zeta)\bigr)
=(1-\chi(w\zeta))\kappa^1(\zeta)_{0,\{0\},0}
\]
\[
  \kappatilde^{2}(\zeta)
=\kappa^2(\zeta)_{0,\{0\},0}
-\del_{\wbar}\bigl(
 \chi(w\zeta)\mu(\zeta)
 \bigr)
=(1-\chi(w\zeta))\kappa^2(\zeta)_{0,\{0\},0}
-(\del_{\wbar}\chi)(w\zeta)
 \cdot \zeta\cdot
 \mu(\zeta)
\]
For any $\delta>0$,
we have the following finiteness,
which is uniform for $\zeta$:
\[
 \int_{T\times U}
 \bigl(
 |\kappatilde^1(\zeta)|^2_{h_2}
+|\kappatilde^2(\zeta)|^2_{h_2}
 \bigr)
 \,|dz d\zbar|\,
 \frac{|dw\,d\wbar|}{|w|^{2+\delta}}
\leq
 C_{1,\delta}
\]
Hence, we have the following for any $\delta>0$:
\[
 \int_{T\times U}
\bigl(
 |\kappatilde_1(\zeta)|_{h_1}^2
+|\kappatilde_2(\zeta)|_{h_1}^2
 \bigr)
 |dz d\zbar dw d\wbar|
\leq
 C_{2,\delta}
 \int_{T\times U}
\bigl(
|\kappatilde^1(\zeta)|_{h_2}^2
+|\kappatilde^2(\zeta)|_{h_2}^2
 \bigr)
|dz d\zbar|\,\frac{|dw\,d\wbar|}{|w|^{2+\delta}}
 |\zeta|^{-2\delta}
\leq C_{3,\delta}|\zeta|^{-2\delta}
\]
Hence, we obtain
$|t(\zeta)|_h\leq 
 C_{4\epsilon}|\zeta|^{-\delta}$ for any $\delta>0$.
Thus,
the proof of Lemma \ref{lem;13.2.14.2} is finished.
\hfill\qed

\vspace{.1in}

Let us prove 
$\Nahm_{\ast}(\nbigp_{\ast}E_1)
=\nbigp_{\ast}E$.
We have the following,
which is similar to Lemma \ref{lem;13.1.23.101}.
\begin{lem}
We have only to show 
$\Nahm_a(\nbigp_{\ast}E_1)
\subset \nbigp_{a}E$
for any $a$.
\hfill\qed
\end{lem}

Around each $P\in D$,
we have the decomposition
\begin{equation}
 \Nahm_{\ast}(\nbigp_{\ast}E_1)
=\bigoplus_{\vecgminio,\alpha}
 \Nahm_{\ast}(\nbigp_{\ast}E_1)_{P,\vecgminio,\alpha},
\end{equation}
according to the decomposition (\ref{eq;13.2.14.3}).
We have only to prove
$\Nahm_a(\nbigp_{\ast}E_1)_{P,\vecgminio,\alpha}
\subset\nbigp_a(E)$.
We shall argue the case $P=0$
in the following.
The other case can be proved similarly.
We shall omit the subscript $P$.
We take a small neighbourhood $U_{\zeta}$ of $P$.

Let us consider the case
$(\vecgminio,\alpha)\neq (0,0)$.
Let $U_{\tau}\subset\proj^1$ be a small neighbourhood of $\infty$
with the coordinate $\tau=w^{-1}$.
Take $\gminia\in\vecgminio$.
For each $c\in\real$,
we have the frame of Lemma \ref{lem;13.1.23.102}.
We have only to show
\begin{equation}
 \label{eq;13.1.23.121}
 \bigl|
 [\tau_{\vecgminio}^jv_{\vecgminio,i}]
 \bigr|_h
=O\Bigl(
 |\zeta|^{
  -(b-j+p_{\vecgminio}-m_{\vecgminio}/2)
  (p_{\vecgminio}-m_{\vecgminio})^{-1}}
 \Bigr).
\end{equation}
Here, $b$ is the parabolic degree of
$v_{\vecgminio,i}$.

We give a preliminary.
We take a ramified covering $U_u\lrarr U_{\zeta}$
given by $\zeta=u^{p_{\vecgminio}-m_{\vecgminio}}$.
We put
$G(\tau_{\vecgminio}):=
 \del_w\gminia(\tau_{\vecgminio})
-\alpha p_{\vecgminio}\tau_{\vecgminio}^{p_{\vecgminio}}
=\sum_{j=0}^{m_{\vecgminio}}\beta_j\tau_{\vecgminio}^{p_{\vecgminio}-j}$.
We take a complex number $\gamma$
such that 
$\beta_{m_{\vecgminio}}
\gamma^{p_{\vecgminio}-m_{\vecgminio}}
-1=0$.
If $U_{\zeta}$ is sufficiently small,
we can take holomorphic functions
$\eta_0^{(i)}(u)$ $(i=1,\ldots,p_{\vecgminio}-m_{\vecgminio})$
on $U_{\zeta}$
satisfying
\[
G\bigl(\eta_0^{(i)}(u)\bigr)
-u^{p_{\vecgminio}-m_{\vecgminio}}=0,
\quad\quad
 \lim_{u\to 0}
 u^{-1}\eta^{(i)}_0(u)
=\gamma\exp\bigl(2\pi\sqrt{-1}i/(p_{\vecgminio}-m_{\vecgminio})\bigr)
\]
There exist $C_1>0$ and $\epsilon_1>0$ such that
$\Bigl|
 u^{-1}\eta_0^{(i)}(u)
-\gamma\exp\bigl(2\pi\sqrt{-1}i
 \big/(p_{\vecgminio}-m_{\vecgminio})\bigr)
\Bigr|
\leq C_1|u|^{\epsilon_1}$.
The following lemma is similar to 
Lemma \ref{lem;13.1.22.100}.
\begin{lem}
\label{lem;13.1.23.110}
Let $Z_{u}$ denote the support of
$\Cok\bigl(
 g_{\gminia,\alpha}-u^{p_{\vecgminio}-m_{\vecgminio}}
 \bigr)$
on $U_{\tau_{\vecgminio}}$.
If $U_{\tau}$  and $U_{\zeta}$ are sufficiently small,
we have a decomposition
$Z_u=\coprod_{i=1}^{p_{\vecgminio}-m_{\vecgminio}}
 Z_u^{(i)}$
and positive constants $C$ and $\epsilon$
such that 
$\bigl|
 \eta^{(i)}_0(u)-\eta_1
 \bigr|
\leq
 C|u|^{1+m_{\vecgminio}+\epsilon}$
for any $\eta_1\in Z^{(i)}_u$.
\hfill\qed
\end{lem}

We set $d:=1+m_{\vecgminio}/2$.
We consider the following sections of
$E^{\vecgminio }_{\gminia}\otimes\Omega^{0,1}$:
\[
 \mu_1(v_{\vecgminio,i},u):=
 \rho\Bigl(|u|^d\bigl(\tau_{\vecgminio}-\eta_0(u)\bigr)\Bigr)
 \,v_{\vecgminio,i}\,d\zbar
\]
\[
 \mu_2(v_{\vecgminio,i},u):=
 \bigl(g_{\gminia,\alpha}-u^{p_{\vecgminio}-m_{\vecgminio}}\bigr)^{-1}
 \Bigl(
 \delbar
 \rho\Bigl(|u|^d\bigl(\tau_{\vecgminio}-\eta_0(u)\bigr)\Bigr)
 \Bigr)\,v_{\vecgminio,i}
\]
The tuple
$\vecmu(v_{\vecgminio,i},u):=\bigl(
 \mu_1(v_{\vecgminio,i},u),
  \mu_2(v_{\vecgminio,i},u)\bigr)$
induces a section of
$\Nahm(\nbigp_{\ast}E_1)_{P,\vecgminio,\alpha}$.
By Lemma \ref{lem;13.1.23.110},
$\rho\bigl(|u|^d\bigl(\tau_{\vecgminio}-\eta_0(u)\bigr)\bigr)$
is constantly $1$ around $Z_u$.
Hence,
$\vecmu(v_{\vecgminio,i},u)$
induces $\bigl[v_{\vecgminio,i}\bigr]$.

By an elementary change of variables,
we obtain the following for any $\delta>0$:
\begin{multline}
 \int |\mu_1(v_{\vecgminio,i},u)|_h^2
\leq
 \int \bigl|
 \rho(|u|^d(\tau_{\vecgminio}-\eta_0(u)))
 \bigr|^2
 |\tau_{\vecgminio}|^{-2(b+\delta)}
 |dz\,d\zbar|
 \,|dw\,d\wbar|
 \\
\leq
 C_{1\delta}|u|^{-2(b+\delta)-2p_{\vecgminio}-2+2d}
=C_{1\delta}|u|^{-2(b+\delta+p_{\vecgminio}-m_{\vecgminio}/2)}
\end{multline}
We also have the following:
\begin{multline}
 \int\bigl|\mu_2(v_{\vecgminio,i},u)\bigr|_h^2
\leq
\int
 \int\bigl|
 \delbar\rho\bigl(|u|^d(\tau_{\vecgminio}-\eta_0(u))\bigr)
 \bigr|^2
 |\tau_{\vecgminio}|^{-2(b+\delta)}
 \frac{1}{\bigl|\del_w\gminia(\zeta_{\vecgminio})
 -\alpha \zeta_{\vecgminio}^{p_{\vecgminio}}
 -u^{p_{\vecgminio}-m_{\vecgminio}}\bigr|^2}
\\
\leq
 C_{2\delta}
 |u|^{-2(b+\delta)-2(p_{\vecgminio}-m_{\vecgminio}-1)-2d}
=
 C_{2\delta}|u|^{-2(b+\delta+p_{\vecgminio}-m_{\vecgminio}/2)}
\end{multline}
Hence, we obtain the estimate (\ref{eq;13.1.23.121}).

\vspace{.1in}

Let us consider the case $(\gminio,\alpha)=(\{0\},0)$.
Note that 
$N(\nbigp_{-1}E_1)_{0,\{0\},0}=N(\nbigu)_{0,\{0\},0}
\subset
 \Nahm_0(\nbigp_{\ast}E_1)$.
Let $\nu\in 
 \Nahm_{1+c}(\nbigp_{\ast}E_1)_{0,\{0\},0}
 \big/N(\nbigp_{-1}E_1)_{0,\{0\},0}$
for $-1<c\leq 0$.
We take $v\in\nbigp_{c}V_{0,\{0\},0}$ which represents $\nu$.
We naturally regard $v$ as a $C^{\infty}$-section of
$\nbigp_c(E_{1})_0$.
Fix a sufficiently small number $b>0$,
and let $\rho$ be a $\real_{\geq 0}$-valued
$C^{\infty}$-function on $\cnum_{\tau}$
such that 
$\rho(\tau)=1$ if $|\tau|\leq b/2$
and $\rho(\tau)=0$ if $|\tau|\geq b$.
We obtain a $C^{\infty}$-section 
$\delbar \bigl(\rho(\tau)vd\zbar\bigr)$
of $\nbigp_{-1}(E_1)_0\otimes\Omega^{0,2}$.
By using
$H^2(T\times\proj^1,\nbigu\otimes L_{-\zeta})=0$
for any $\zeta$,
we can take a holomorphic family of $C^{\infty}$-forms
$\kappa(\zeta)=
 \kappa^1(\zeta)d\zbar+\kappa^2(\zeta)d\wbar$
of $\nbigu\otimes \Omega^{0,1}$
such that 
$\delbar_{E\otimes L_{-\zeta}}\kappa(\zeta)
=\delbar \bigl(\rho(\tau)v d\zbar\bigr)$.
Then, $\rho(\tau)v d\zbar-\kappa(\zeta)$
induces a holomorphic section $\nutilde$
of $\Nahm_{1+c}(\nbigp_{\ast}E_1)$ around $P$
which induces $\nu$
in $\Nahm_{1+c}(\nbigp_{\ast}E_1)/N(\nbigu)$.

We consider the following sections:
\[
\mu_1(v,\zeta):=
\bigl(\rho(\tau)-\rho(\zeta^{-1} \tau)\bigr) v\,d\zbar
\]
\[
\mu_2(v,\zeta):=
\delbar\bigl(\rho(\tau)-\rho(\zeta^{-1}\tau)\bigr) 
 (g_{0,\{0\},0}-\zeta)^{-1}(v)
\]
Then, 
$\mu_1(v,\zeta)+\mu_2(v,\zeta)
-\kappa(\zeta)$
induces the same section $\nutilde$.

We have the following for any $\delta>0$:
\[
 \int \bigl|\mu_1\bigr|_{h_1}^2\,|dw\,d\wbar|
\leq
 C_{\delta}\int_{|\tau|\geq A|\zeta|}
 |\tau|^{-2(c+\delta)-4}
 |d\tau\,d\taubar|
\leq
 C_{\delta}|\zeta|^{-2(c+1+\delta)}
\]
We also have the following:
\[
 \int\bigl|\mu_2\bigr|_{h_1}^2\,|dz\,d\zbar|
\leq
 C_{\delta}\int
 \bigl|\delbar\rho(\zeta^{-1}\tau)\bigr|^2
 |\zeta|^{-2}|\tau|^{-2(c+\delta)}
\leq
 C_{\delta}|\zeta|^{-2(c+1+\delta)}
\]
Because the support of $\delbar(\rho(\tau)v\,d\zbar)$
is compact,
we obtain
$\int|\kappa|_{h_1}^2\dvol=O(|\zeta|^{-\delta})$
for any $\delta>0$,
by the argument in the proof of Lemma \ref{lem;13.2.14.2}.
We obtain
$|\nutilde|_h\leq C_{\delta}|\zeta|^{-(c+1+\delta)}$
for any $\delta>0$.
Thus, we obtain
$\Nahm_{1+c}(\nbigp_{\ast}E_1)_{0,\{0\},0}
\subset
 \nbigp_{1+c}E$,
and the proof of Theorem \ref{thm;13.1.23.100}
is finished.
\hfill\qed

\subsection{Kobayashi-Hitchin correspondence
for $L^2$-instantons}
\subsubsection{Statements}

Let $\nbigp_{\ast}E_1$ be a good filtered bundle 
on $(T\times\proj^1,T\times\{\infty\})$
of degree $0$ satisfying the conditions {\bf(A3)}.
(See \S\ref{subsection;13.2.14.11}
for good filtered bundles.)

\begin{prop}
\label{prop;13.1.24.2}
$\nbigp_{\ast}E_1$ is stable,
if and only if
$\Nahm_{\ast}(\nbigp_{\ast}E_1)$ is 
a stable filtered Higgs bundle.
(See {\rm\S\ref{subsection;13.10.15.200}}
for the stability condition of
$\nbigp_{\ast}E_1$.)
\end{prop}

Before going to the proof,
we give a consequence.

\begin{thm}
\label{thm;13.1.24.1}
Let $\nbigp_{\ast}E_2$ be a stable good filtered bundle
on $(T\times\proj^1,T\times\{\infty\})$
with $\deg(\nbigp_{\ast}E_2)=0$.
We set $E_2:=(\nbigp_aE_2)_{|T\times\cnum}$
which is independent of $a$.
Then, there exists a Hermitian-Einstein metric $h$
of $E_2$ on $T\times\cnum$ such that
(i) its curvature is $L^2$ with respect to
$h$ and the Euclidean metric,
(ii) it is adapted to the filtered bundle $\nbigp_{\ast}E_2$.
(See \S\ref{subsection;13.10.28.3}
for adaptedness.)
Such a metric is unique up to the multiplication
of positive constants.
\end{thm}
\pf
If $\rank E_2=1$,
then $E_2$ is the pull back of a line bundle $L$
of degree $0$ on $T$
by the projection
$T\times\cnum\lrarr T$,
and the parabolic structure is the natural one,
as in Remark {\rm\ref{rem;13.10.20.1}}.
The flat metric of $L$ induces
the Hermitian-Einstein metric of $E_2$
adapted to $\nbigp_{\ast}E_2$.

Suppose $\rank E_2>1$.
By Proposition \ref{prop;13.1.24.2},
$\Nahm(\nbigp_{\ast}E_2)$ is stable.
By Corollary \ref{cor;13.1.29.31},
we have
\[
 \deg\Nahm(\nbigp_{\ast}E_2)=
\deg(\nbigp_{\ast}E_2)=0.
\]
By Corollary \ref{cor;13.1.24.21},
$\Nahm(\nbigp_{\ast}E_1)$ is a good filtered Higgs bundle.
Hence, 
by the Kobayashi-Hitchin correspondence
for wild harmonic bundles on curves
\cite{biquard-boalch},
we obtain an adapted harmonic metric
for $\Nahm(\nbigp_{\ast}E_1)$.
Its Nahm transform induces a Hermitian-Einstein metric
of $E_1$ adapted to the filtered bundle
$\nbigp_{\ast}E_1$,
by Theorem \ref{thm;13.1.23.20}
and Proposition \ref{prop;13.1.22.2}.
\hfill\qed

\vspace{.1in}
Note that the converse is given in 
Proposition \ref{prop;13.1.24.3}

\begin{rem}
This proof of Theorem {\rm\ref{thm;13.1.24.1}}
is based on the idea mentioned in
Remark {\rm 5.13} of {\rm\cite{Biquard-Jardim}}.
\hfill\qed
\end{rem}

\subsubsection{Proof of Proposition \ref{prop;13.1.24.2}}

Let us prove the ``if'' part in Proposition \ref{prop;13.1.24.2}.
Suppose $\Nahm(\nbigp_{\ast}E_1)$ is stable.
By the  Kobayashi-Hitchin correspondence
for wild harmonic bundles on curves
\cite{biquard-boalch} (see also \cite{Mochizuki-wild}
for the case of good filtered flat bundles),
we have an adapted harmonic metric for 
$\Nahm(\nbigp_{\ast}E_1)$.
By Theorem \ref{thm;13.1.23.20},
its Nahm transform gives an adapted Hermitian-Einstein metric
for $\nbigp_{\ast}E_1$.
By Proposition \ref{prop;13.1.24.3},
$\nbigp_{\ast}E_1$ is polystable.
If it is not stable,
the decomposition into the stable components
induces a decomposition of 
$\Nahm(\nbigp_{\ast}E_1)$,
which contradicts with the stability of
$\Nahm(\nbigp_{\ast}E_1)$.
Hence, $\nbigp_{\ast}E_1$ is stable.

Let us show the ``only if'' part in Proposition \ref{prop;13.1.24.2}.
Let $(\nbigp_{\ast}E,\theta):=\Nahm(\nbigp_{\ast}E_1)$.
Let $(\nbigp_{\ast}E',\theta')\subset(\nbigp_{\ast}E,\theta)$
be a strict filtered Higgs subbundle
with $0<\rank E'<\rank E$.
We obtain a subcomplex
$\nbigctilde^{\bullet}(\nbigp_{\ast}E',\theta')
\subset
 \nbigctilde^{\bullet}(\nbigp_{\ast}E,\theta)$
on $T^{\lor}\times T\times\proj^1$.
Let $\nbigytilde^{\bullet}=
 \bigl(\nbigytilde^0\lrarr\nbigytilde^1\bigr)$
be the quotient complex.

\begin{lem}
\label{lem;13.2.27.2}
The induced morphism
$R^1p_{23\ast}\bigl(
 \nbigctilde^{\bullet}(\nbigp_{\ast}E',\theta')
 \bigr)
\lrarr
 R^1p_{23\ast}\bigl(
 \nbigctilde^{\bullet}(\nbigp_{\ast}E,\theta)\bigr)$
is injective.
\end{lem}
\pf
By the construction,
$\nbigytilde^0$ is locally free.
Hence, we obtain that 
$R^0p_{23\ast}\nbigytilde^0$ 
is torsion-free.
Because
\[
 R^0p_{23\ast}\bigl(
 \nbigytilde^{\bullet}
 \bigr)
\lrarr R^0p_{23\ast}\nbigytilde^0
\]
is injective,
we obtain that
$R^0p_{23\ast}\nbigytilde^{\bullet}$
is torsion-free.
On a small neighbourhood $U\subset\proj^1$ of $\infty$,
we have 
$R^i p_{23\ast}\bigl(
 \nbigctilde^{\bullet}(\nbigp_{\ast}E',\theta')
 \bigr)
=R^i p_{23\ast}\bigl(
 \nbigctilde^{\bullet}(\nbigp_{\ast}E,\theta)\bigr)
=0$ unless $i=1$.
It is easy to check that
\[
 R^1 p_{23\ast}\bigl(
 \nbigctilde^{\bullet}(\nbigp_{\ast}E',\theta')
 \bigr)_{|T\times\{\infty\}}
\lrarr
R^1 p_{23\ast}\bigl(
 \nbigctilde^{\bullet}(\nbigp_{\ast}E,\theta)\bigr)
 _{|T\times\{\infty\}} 
\]
is injective.
Hence,
$R^1 p_{23\ast}\bigl(
 \nbigctilde^{\bullet}(\nbigp_{\ast}E',\theta')
 \bigr)_{|T\times U}
\lrarr
R^1 p_{23\ast}\bigl(
 \nbigctilde^{\bullet}(\nbigp_{\ast}E,\theta)\bigr)
 _{|T\times U}$
is injective.
Because 
\[
0\lrarr
R^0p_{23\ast}\nbigytilde^{\bullet}
\lrarr
R^1 p_{23\ast}\bigl(
 \nbigctilde^{\bullet}(\nbigp_{\ast}E',\theta')
 \bigr)
\lrarr
R^1 p_{23\ast}\bigl(
 \nbigctilde^{\bullet}(\nbigp_{\ast}E,\theta)\bigr)
\]
is exact,
we obtain $R^0p_{23\ast}\nbigytilde^{\bullet}=0$,
and 
$R^1 p_{23\ast}\bigl(
 \nbigctilde^{\bullet}(\nbigp_{\ast}E',\theta')
 \bigr)
\lrarr
R^1 p_{23\ast}\bigl(
 \nbigctilde^{\bullet}(\nbigp_{\ast}E,\theta)\bigr)$
is injective.
\hfill\qed

\vspace{.1in}

We define the parabolic structure of
$R^1p_{23\ast}\bigl(\nbigctilde^{\bullet}
 (\nbigp_{\ast}E',\theta')\bigr)$
as in \S\ref{subsection;13.2.14.10}.
The filtered sheaf is denoted by
$\nbigp_{\ast}\nbigv_1$.
We have a naturally defined injective morphism
$\nbigp_{\ast}\nbigv_1\lrarr \nbigp_{\ast}E_1$.
Hence, we have
$\deg(\nbigp_{\ast}\nbigv_1)\leq 0$.
By the argument in \S\ref{subsection;13.1.24.10},
we obtain
\[
 \int_{T\times\proj^1}
 c_1\bigl(\nbigp_{\ast}\nbigv_1\bigr)
 \omega_{T}
-\int_{T\times\proj^1}
 c_1(R^2p_{23\ast}\nbigctilde^{\bullet}(\nbigp_{\ast}E',\theta'))
 \omega_T
=\deg(\nbigp_{\ast}E')
\]
Because 
$R^2p_{23\ast}\nbigctilde^{\bullet}(\nbigp_{\ast}E',\theta')$
is a torsion sheaf,
we obtain 
$\int_{T\times\proj^1}
 c_1(R^2p_{23\ast}\nbigctilde^{\bullet}(\nbigp_{\ast}E',\theta'))
 \omega_T\geq 0$.
Hence,
we obtain
$\deg(\nbigp_{\ast}E')\leq 0$,
i.e.,
$(\nbigp_{\ast}E,\theta)$ is semistable.

We have 
$(\nbigp_{\ast}E',\theta')\subset(\nbigp_{\ast}E,\theta)$
such that 
$(\nbigp_{\ast}E',\theta')$
is stable of degree $0$.
If $(\nbigp_{\ast}E',\theta')$ has no singularity,
it is isomorphic to
$\nbigo_T^{\lor}$
with a Higgs field $\alpha d\zeta$
$(\alpha\in\cnum)$,
and hence
$R^1p_{23\ast}\nbigctilde^{\bullet}(\nbigp_{\ast}E',\theta')$
is a non-zero torsion subsheaf of $E_1$,
which contradicts with Lemma \ref{lem;13.2.27.2}. 
Therefore, $(\nbigp_{\ast}E',\theta')$
has a singularity,
and $\Nahm_{\ast}(\nbigp_{\ast}E',\theta')\neq 0$
is a good filtered subbundle of $\nbigp_{\ast}E_1$.
By the stability of $\nbigp_{\ast}E_1$,
we have 
$\rank\Nahm_{\ast}(\nbigp_{\ast}E',\theta')=\rank E_1$.
Because
$\deg\Nahm_{\ast}(\nbigp_{\ast}E',\theta')
=\deg(\nbigp_{\ast}E)$,
we have 
$\Nahm_{\ast}(\nbigp_{\ast}E',\theta')
=\nbigp_{\ast}E_1$
in codimension one.
Because both of them are filtered bundles,
we have
$\Nahm_{\ast}(\nbigp_{\ast}E',\theta')
=\nbigp_{\ast}E_1$
on $T\times\proj^1$.
Then, we obtain
$(\nbigp_{\ast}E',\theta')
=(\nbigp_{\ast}E,\theta)$
by the involutivity of the algebraic Nahm transforms.
\hfill\qed

\noindent
{\em Address\\
Research Institute for Mathematical Sciences,
Kyoto University,
Kyoto 606-8502, Japan,\\
takuro@kurims.kyoto-u.ac.jp
}


\begin{thebibliography}{99}

\bibitem{Aker-Szabo}
K. Aker, Sz. Szabo,
{\em Algebraic Nahm transform for Parabolic Higgs Bundles on $P^1$},
math.arXiv.0610301

\bibitem{a}
L. V. Ahlfors, {\it An extension of Schwarz's lemma},
Trans. Amer. Math. Soc. {\bf 43} (1938), 359--364.

\bibitem{Andreotti-Vesentini}
A. Andreotti, E. Vesentini,
{\em Carleman estimates for the Laplace-Beltrami equation 
on complex manifolds},
Inst. Hautes \'{E}tudes Sci. Publ. Math. {\bf 25},
(1965), 81--130.

\bibitem{Bartotti-Bruzzo-Ruiperez}
C. Bartocci, 
U. Bruzzo, 
D. Hern\'{a}ndez Ruipérez, 
{\em Fourier-Mukai and Nahm transforms 
in geometry and mathematical physics},
Progress in Mathematics, {\bf 276}. Birkhäuser Boston, 
Inc., Boston, MA, 2009. 

\bibitem{Beilison-Bloch-Deligne-Esnault}
A. Beilinson,
S. Bloch,
P. Deligne,
H. Esnault,
{\em Periods for irregular connections on curves},
preprint.

\bibitem{b}
O. Biquard,
{\it Fibr\'es de Higgs et connexions int\'egrables:
le cas logarithmique (diviseur lisse)},
Ann. Sci. \'Ecole Norm. Sup. {\bf 30} (1997), 41--96.

\bibitem{biquard-boalch}
O. Biquard,
P. Boalch,
{\em Wild non-abelian Hodge theory on curves},
Compos. Math. {\bf 140} (2004),
179--204.

\bibitem{Biquard-Jardim}
O. Biquard,
M. Jardim,
{\em Asymptotic behaviour and the moduli space of 
doubly-periodic instantons}. J. Eur. Math. Soc. {\bf 3}
(2001), 335--375.

\bibitem{Bonsdorff1}
J. Bonsdorff,
{\em A Fourier transformation for Higgs bundles},
J. Reine Angew. Math. {\bf 591} (2006), 21--48. 

\bibitem{Bonsdorff2}
J. Bonsdorff, 
{\em Autodual connection in the Fourier transform of a Higgs bundle.}
Asian J. Math. {\bf 14} (2010), 153--173.

\bibitem{Braam-Baal}
P. Braam,
P. van Baal,
{\em Nahm's transformation for instantons},
Comm. Math. Phys. {\bf 122} (1989), 267--280.

\bibitem{Charbonneau}
B. Charbonneau,
{\em Analytic aspects of periodic instantons},
Thesis (Ph.D.)–Massachusetts Institute of Technology. 2004


\bibitem{Donaldson-Nahm}
S. K. Donaldson, 
{\em Nahm's equations and the classification of monopoles},
Comm. Math. Phys. {\bf 96} (1984), 387--407

\bibitem{Donaldson-Kronheimer}
S. K. Donaldson,
P. B. Kronheimer,
{\em The geometry of four-manifolds},
Oxford Science Publications. 
The Clarendon Press, 
Oxford University Press, New York, 1990.

\bibitem{Bloch-Esnault}
S. Bloch, H. Esnault,
{\em Local Fourier transforms and rigidity for D-modules},
Asian J. Math. {\bf 8} (2004), 587--605. 

\bibitem{Esnault-Viehweg}
H. Esnault, E. Viehweg, 
{\em Logarithmic de Rham complexes and vanishing theorems},
Invent. Math. {\bf 86} (1986), 161--194. 

\bibitem{Fang}
J. Fang, 
{\em Calculation of local Fourier transforms 
for formal connections},
Sci. China Ser. A {\bf 52} (2009),
2195--2206.

\bibitem{Ford-Pawlowski}
C. Ford,
J. M. Pawlowski,
{\em Doubly periodic instantons and their constituents},
Phys. Rev. D (3) {\bf 69} (2004), 065006, 12 pp.

\bibitem{Frejlich-Jardim}
P. Frejlich, M. Jardim, 
{\em Nahm transform for Higgs bundles},
J. Geom. Phys. {\bf 58} (2008), 1221--1230.

\bibitem{Fu}
L. Fu,
{\em Calculation of $\ell$-adic local Fourier transformations},
Manuscripta Math. {\bf 133} (2010), 
409--464.

\bibitem{Fujiki}
A. Fujiki,
{\em An L2 Dolbeault lemma and its applications},
Publ. Res. Inst. Math. Sci. {\bf 28} (1992), 845--884

\bibitem{Garcia-Lopez}
R. Garc\'{\i}a L\'{o}pez, 
{\em Microlocalization and stationary phase},
Asian J. Math. {\bf 8} (2004), 747--768. 

\bibitem{Graham-Squire}
A. Graham-Squire,
{\em Calculation of local formal Fourier transforms},
Arkiv f\"{o}r Matematik, 1--14, \\
10.1007/s11512-011-0156-2

\bibitem{Hitchin-construction-monopole}
N. Hitchin,
{\em Construction of monopoles},
Commun. Math. Phys. {\bf 89}, (1983), 145--190.

\bibitem{Hitchin-self-duality}
N. Hitchin,
{\em The self-duality equations on a Riemann surface},
Proc. London Math. Soc. (3) {\bf 55} (1987), 59--126.
\bibitem{Hitchin-lecture}
N. Hitchin,
{\em Monopoles, minimal surfaces and algebraic curves.}
Presses de l'Université de Montréal, Montreal, QC, 1987. 

\bibitem{Jardim1}
M. Jardim, 
{\em Construction of doubly-periodic instantons,}
Comm. Math. Phys. {\bf 216} (2001), 1--15. 

\bibitem{Jardim2}
M. Jardim, 
{\em Nahm transform and spectral curves 
for doubly-periodic instantons},
Comm. Math. Phys. {\bf 225} (2002), 639--668.

\bibitem{Jardim3}
M. Jardim, 
{\em Classification and existence of doubly-periodic instantons.} 
Q. J. Math. {\bf 53} (2002), 431--442

\bibitem{Jardim4}
M. Jardim,
{\em A survey on Nahm transform},
J. Geom. Phys. {\bf 52} (2004), 313--327

\bibitem{Lang}
S. Lang,
{\em Real Analysis},
Second edition. Addison-Wesley Publishing Company, 
Advanced Book Program, Reading, MA, 1983


\bibitem{Li-Narasimhan}
J. Li,
M. S. Narasimhan,
{\em Hermitian-Einstein metrics on parabolic stable bundles},
Acta Math. Sin. (Engl. Ser.) {\bf 15} (1999), 
93--114.

\bibitem{li2}
J. Li,
{\em Hermitian-Einstein metrics and Chern number inequalities
on parabolic stable bundles over K\"{a}hler manifolds,} 
Comm. Anal. Geom. {\bf 8} (2000), 445--475.

\bibitem{Laumon}
G. Laumon,
{\em Transformation de Fourier, constantes d'équations fonctionnelles 
et conjecture de Weil},
Inst. Hautes \'{E}tudes Sci. Publ. Math.
{\bf 65} (1987), 131--210

\bibitem{Laumon-Fourier-Mukai}
G. Laumon,
{\em Transformation de Fourier generalisee},
arXiv:alg-geom/9603004

\bibitem{Malgrange-book}
B. Malgrange,
{\em
\'{E}quations diff\'{e}rentielles \`{a} coefficients
	polynomiaux},
Progress in Mathematics, {\bf 96},
Birkh\"{a}user Boston, Inc., Boston, MA, 1991. 

\bibitem{Maruyama-Yokogawa}
M. Maruyama,
K. Yokogawa,
{\em Moduli of parabolic stable sheaves},
Math. Ann. {\bf 293} (1992),  77--99. 

\bibitem{mochi2}
T. Mochizuki,
{\em Asymptotic behaviour of tame harmonic bundles
and an application to pure twistor $D$-modules I, II},
Mem. AMS. {\bf 185} (2007)

\bibitem{Mochizuki-Fourier}
T. Mochizuki, 
{\em Note on the Stokes structure of Fourier transform},
Acta Math. Vietnam. {\bf 35} (2010), 107--158. 

\bibitem{Mochizuki-wild}
T. Mochizuki,
{\em Wild harmonic bundles and wild pure twistor $D$-modules},
Ast\'{e}risque {\bf 340},
Soci\'{e}t\'{e} Math\'{e}matique de France, Paris,
2011.

\bibitem{Mukai}
S. Mukai,
{\em Duality between $D(X)$ and $D(\widehat{X})$
 with its application to Picard sheaves},
Nagoya Math. J. {\bf 81} (1981), 153--175. 

\bibitem{Nakajima}
H. Nakajima,
{\em Monopoles and Nahm's equations},
in {\em Einstein metrics and Yang-Mills connections (Sanda, 1990)}, 
Lecture Notes in Pure and Appl. Math., {\bf 145}, Dekker, New York, (1993),
193--211.

\bibitem{Rothstein}
M. Rothstein, 
{\em Sheaves with connection on abelian varieties,}
Duke Math. J. {\bf 84} (1996), 565--598.

\bibitem{sabbah2}
C. Sabbah,
{\em Polarizable twistor $D$-modules},
Ast\'{e}risque, {\bf 300},
Soci\'{e}t\'{e} Math\'{e}matique
de France, Paris, 2005.


\bibitem{Sabbah-local-Fourier}
C. Sabbah, 
{\em An explicit stationary phase formula 
for the local formal Fourier-Laplace transform},
in {\em Singularities I}, 
Contemp. Math., {\bf 474}, Amer. Math. Soc., (2008),
309--330.

\bibitem{Sacks-Uhlenbeck}
J. Sacks, K. Uhlenbeck,
{\em The existence of minimal immersions of 2-spheres},
Ann. of Math. (2) {\bf 113} (1981), 1--24.

\bibitem{Schenk}
H. Schenk,
{\em On a generalised Fourier transform of instantons over flat tori},
Comm. Math. Phys. {\bf 116} (1988), 177--183.

\bibitem{Simpson88}
C. T. Simpson,
{\em Constructing variations of Hodge structure 
using Yang-Mills theory and applications to uniformization},
J. Amer. Math. Soc. {\bf 1} (1988), 867--918.

\bibitem{Simpson90}
C. T. Simpson,
{\em Harmonic bundles on noncompact curves},
J. Amer. Math. Soc. {\bf 3} (1990), 713--770.

\bibitem{s5}
C. T. Simpson,
{\em  Higgs bundles and local systems},
Publ. I.H.E.S., {\bf 75} (1992),  5--95.

\bibitem{s3}
C. T. Simpson,
{\it Mixed twistor structures},
math.AG/9705006.

\bibitem{Siu-extension}
S. T. Siu,
{\em Techniques of extension of analytic objects},
Lecture Notes in Pure and Applied Mathematics, 
Marcel Dekker, Inc., New York, 1974. 

\bibitem{Szabo}
Sz. Szab\'{o}, 
{\em Nahm transform for integrable connections on the Riemann sphere},
M\'{e}m. Soc. Math. Fr. (N.S.) {\bf 110} (2008).

\bibitem{Uhlenbeck1}
K. Uhlenbeck,
{\em Removable singularities in Yang-Mills fields},
Comm. Math. Phys. {\bf 83} (1982), 11--29.

\bibitem{Uhlenbeck2}
K. Uhlenbeck,
{\em Connections with $L^{p}$ bounds on curvature,} 
Comm. Math. Phys. {\bf 83} (1982), 31--42.

\bibitem{Wehrheim}
K. Wehrheim, 
{\em Energy identity for anti-self-dual instantons on 
$\cnum\times\Sigma$},
Math. Res. Lett. {\bf 13} (2006), 
161--166.

\bibitem{Zucker}
S. Zucker,
{\em Hodge theory with degenerating coefficients. $L2$ cohomology 
in the Poincar\'{e} metric},
Ann. of Math. (2) {\bf 109} (1979), 
415--476.

\end{thebibliography}
\end{document}